\documentclass[spanningrule]{cambridge7A} 

\usepackage{latexsym}
\usepackage{amssymb,amsmath}
\usepackage{bbm}
\usepackage{mathrsfs}
\usepackage{graphicx}
\usepackage{color}
\usepackage{multicol}
\usepackage{txfonts}

\usepackage{url}
\usepackage{makeidx}
\makeindex

\usepackage[all,2cell]{xy}
\UseAllTwocells
\xyoption{arc}
\xyoption{rotate}

\usepackage{hyperref}   
\hypersetup{colorlinks, 
breaklinks,             
linkcolor=black,        
urlcolor=black}         

\def\makeRRlabeldot#1{\hss\llap{#1}}    
\renewcommand{\theenumi}{\rm(\alph{enumi})}
\renewcommand{\theequation}{\thechapter.\arabic{equation}}

\newcommand{\demph}[1]{\textbf{#1}}     

\newenvironment{displaytext}
{\begin{list}{}%
{\setlength{\leftmargin}{2em}%
\setlength{\rightmargin}{\leftmargin}}%
\item}%
{\end{list}}

\newenvironment{setprop}{\begin{displaytext}\itshape}{\end{displaytext}}

\newenvironment{slogan}{\begin{displaytext}\itshape}{\end{displaytext}}

\newenvironment{citedsource}{\begin{displaytext}}{\end{displaytext}}

\newcommand{\minihead}[1]{\subsection*{#1}}

\newcommand{\subjectchange}{\bigskip}

\newcommand{\exs}{\minihead{Exercises}}
\newcommand{\exone}{\minihead{Exercise}}

\newcommand{\cat}[1]{\mathscr{#1}}      
\newcommand{\scat}[1]{\mathbf{#1}}      
\newcommand{\fcat}[1]{\mathbf{#1}}      

\newcommand{\ovln}[1]{\overline{#1}}
\newcommand{\lwr}[1]{\mathbf{#1}}       
\newcommand{\rreg}[1]{\underline{#1}}   
\newcommand{\ynt}[1]{\tilde{#1}}        
\newcommand{\yel}[1]{\hat{#1}}          
\newcommand{\wideyel}[1]{\widehat{#1}}  

\newcommand{\epsln}{\varepsilon}        

\newcommand{\nat}{\mathbb{N}}           
\newcommand{\integers}{\mathbb{Z}}
\newcommand{\rationals}{\mathbb{Q}}
\newcommand{\reals}{\mathbb{R}}
\newcommand{\complexes}{\mathbb{C}}

\DeclareMathOperator{\cset}{\cat{C}}    
\DeclareMathOperator{\Cl}{Cl}           
\DeclareMathOperator{\Cone}{Cone}       
\DeclareMathOperator{\Eq}{Eq}           
\DeclareMathOperator{\ev}{ev}           
\DeclareMathOperator{\Hom}{Hom}         
\DeclareMathOperator{\HOM}{\fcat{Hom}}  
\DeclareMathOperator{\im}{im}           
\DeclareMathOperator{\Lan}{Lan}         
\DeclareMathOperator{\ob}{ob}           
\DeclareMathOperator{\oset}{\cat{O}}    
\DeclareMathOperator{\pr}{pr}           
\DeclareMathOperator{\pset}{\cat{P}}    
\DeclareMathOperator{\Sub}{Sub}         

\newcommand{\h}{H}                      
\newcommand{\op}{\mathrm{op}}           
\newcommand{\true}{\texttt{true}}       
\newcommand{\false}{\texttt{false}}     
        
\newcommand{\One}{\mathbf{1}}           
\newcommand{\Two}{\mathbf{2}}           
\newcommand{\Ab}{\fcat{Ab}}             
\newcommand{\Bilin}{\fcat{Bilin}}       
\newcommand{\Cat}{\fcat{Cat}}           
\newcommand{\CAT}{\fcat{CAT}}           
\newcommand{\CptHff}{\fcat{CptHff}}     
\newcommand{\CRing}{\fcat{CRing}}       
\newcommand{\FDVect}{\fcat{FDVect}}     
\newcommand{\Field}{\fcat{Field}}       
\newcommand{\FinSet}{\fcat{FinSet}}     
\newcommand{\Fix}{\fcat{Fix}}           
\newcommand{\Grp}{\fcat{Grp}}           
\newcommand{\Mt}{\fcat{Mat}}            
\newcommand{\Mon}{\fcat{Mon}}           
\newcommand{\Ord}{\fcat{Ord}}           
\newcommand{\Ring}{\fcat{Ring}}         
\newcommand{\Set}{\fcat{Set}}           
\newcommand{\Sym}{\fcat{Sym}}           
\newcommand{\Tp}{\fcat{Top}}            
\newcommand{\Toph}{\fcat{Toph}}         
\newcommand{\Vect}{\fcat{Vect}}         

\newcommand{\such}{\mathrel{|}}         
\newcommand{\Bigsuch}{\mathrel{\big|}}  
\newcommand{\without}{\setminus}        

\newcommand{\emptybk}{\hspace*{.5em}}   
\newcommand{\blank}{(\emptybk)}         
\newcommand{\dashbk}{-}                 
\newcommand{\bl}{\bullet}               

\newcommand{\from}{{\colon}\linebreak[0]}       

\newcommand{\iso}{\cong}                
\newcommand{\eqv}{\simeq}               
\newcommand{\sub}{\subseteq}            
\newcommand{\ladj}{\dashv}              
\newcommand{\divides}{\mathbin{\mid}}   
\newcommand{\of}{\mathbin{\circ}}       

\newcommand{\meet}{\wedge}              
\newcommand{\Meet}{\bigwedge}
\newcommand{\join}{\vee}
\renewcommand{\Join}{\bigvee}

\newcommand{\ftrcat}[2]{[#1,#2]}                
\newcommand{\pshf}[1]{\ftrcat{#1^\op}{\Set}}    
\newcommand{\psh}[1]{\hat{#1}}                  

\newcommand{\elt}[1]{\scat{E}(#1)}              

\newcommand{\colt}[1]{\lim\limits_{\rightarrow #1}}     
\newcommand{\lt}[1]{\lim\limits_{\leftarrow #1}}        

\newcommand{\abel}[1]{#1_{\mathrm{ab}}}                 
\newcommand{\comma}[2]{(#1 \mathbin{\Rightarrow} #2)}   
\newcommand{\qer}[2]{#1/\mathord{#2}}   
\newcommand{\crd}[1]{\left|#1\right|}   

\newcommand{\oppairu}{\rightleftarrows}         
\newcommand{\parpairu}{\rightrightarrows}       
\newcommand{\ot}{\leftarrow}                    
\newcommand{\toby}[1]{\stackrel{#1}{\longrightarrow}}   
\newcommand{\otby}[1]{\stackrel{#1}{\longleftarrow}}    
\newcommand{\longto}{\longrightarrow}           
\newcommand{\toiso}{\toby{\textstyle{}_\sim}}   
\newcommand{\incl}{\hookrightarrow}             
\newcommand{\mapsfrom}{\mathrel{\reflectbox{\ensuremath{\mapsto}}}}
\newcommand{\textiff}{\Longleftrightarrow}      
\newcommand{\textif}{\Longleftarrow}            
\newcommand{\textonlyif}{\Longrightarrow}       

\newcommand{\nent}{\rotatebox{45}{$\Rightarrow$}} 
 
\newcommand{\oppair}[4]{%
\xymatrix{%
#1 \ar@<.5ex>[r]^-{#3} &{#2}\ar@<.5ex>[l]^-{#4}%
}}

\newcommand{\oppairi}[4]{%
\xymatrix@1{%
#1 \ar@<.5ex>[r]^{#3} &{#2}\ar@<.5ex>[l]^{#4}%
}}

\newcommand{\parpair}[4]{%
\xymatrix{%
#1 \ar@<.5ex>[r]^{#3} \ar@<-.5ex>[r]_{#4} &#2%
}}

\newcommand{\parpairi}[4]{%
\xymatrix@1{%
#1 \ar@<.5ex>[r]^{#3} \ar@<-.5ex>[r]_{#4} &#2%
}}

\newcommand{\adjn}[4]{%
\xymatrix{
#1 \ar@{}[d]|\ladj \ar@<1ex>[d]^{#4} \\
#2 \ar@<1ex>[u]^{#3}
}}

\newcommand{\hadjnli}[4]{%
\xymatrix@1{
#1 \ar@<1.1ex>[r]^-{#3} \ar@{}[r]|-\bot &#2 \ar@<1.1ex>[l]^-{#4}}}

\newcommand{\hadjnri}[4]{%
\xymatrix@1{
#1 \ar@<1.1ex>[r]^-{#3} \ar@{}[r]|-\top &#2 \ar@<1.1ex>[l]^-{#4}}}

\newcommand{\searrows}{\makebox[1em]{$\downarrow$%
\hspace*{-.2em}\raisebox{-1ex}{$\rightarrow$}}}

\newcommand{\esarrows}{\makebox[1em]{%
\raisebox{1ex}{$\rightarrow$}\hspace*{-.2em}$\downarrow$}}

\newcommand{\cell}[4]{\put(#1,#2){\makebox(0,0)[#3]{\ensuremath{#4}}}}

\newcommand{\bref}[1]{\ref{#1}}         

\newcommand{\citestyle}[1]{#1}  
\newcommand{\citeGow}{\citestyle{Gowers (2002)}}
\newcommand{\citeKel}{\citestyle{Kelly (1982)}}
\newcommand{\citeLR}{\citestyle{Lawvere and Rosebrugh (2003)}}
\newcommand{\citeRST}{\citestyle{Leinster (2014)}}
\newcommand{\citeCWM}{\citestyle{Mac~Lane (1971)}}

\newcommand{\ntn}[1]{\label{notn:#1}}   
\newcommand{\nref}[1]{\pageref{notn:#1}}
\newcommand{\nuse}[2]{#2, #1\\}         
\usepackage[amsmath,thmmarks]{ntheorem}

\theoremseparator{\hspace{0.5em}}       

\newtheorem{thm}{Theorem}[section]
\newtheorem{propn}[thm]{Proposition}
\newtheorem{lemma}[thm]{Lemma}
\newtheorem{cor}[thm]{Corollary}

\newtheorem{ilemma}{Lemma}[chapter]     
\newtheorem{alemma}{Lemma}[chapter]     

\theorembodyfont{\normalfont}

\newtheorem{defn}[thm]{Definition}
\newtheorem{example}[thm]{Example}
\newtheorem{examples}[thm]{Examples}
\newtheorem{remark}[thm]{Remark}
\newtheorem{remarks}[thm]{Remarks}
\newtheorem{warning}[thm]{Warning}
\newtheorem{notn}[thm]{Notation}
\newtheorem{constn}[thm]{Construction}
\newtheorem{question}[thm]{\hspace{-.5ex}}      

\newtheorem{iexample}[ilemma]{Example}          
\newtheorem{iquestion}[ilemma]{\hspace{-.5ex}}  
\newtheorem{aquestion}[alemma]{\hspace{-.5ex}}  

\theoremstyle{nonumberplain}
\theoremsymbol{\ensuremath{\Box}}
\qedsymbol{\ensuremath{\Box}}

\newtheorem{pf}{Proof}

\newcommand{\theoremtobeproved}{}
\newtheorem{pfoftheorem}{Proof of \theoremtobeproved}
\newenvironment{pfof}[1]
{
\renewcommand{\theoremtobeproved}{#1}
\begin{pfoftheorem}
}
{\end{pfoftheorem}}

\begin{document}

%
%
%

\thispagestyle{empty}

{\centering
\vspace*{20mm}

\scalebox{1.2}{\textbf{\Huge Basic Category Theory}}

\vspace*{12mm}

{\large T\,O\,M \, L\,E\,I\,N\,S\,T\,E\,R}\\
\textit{University \,of \,Edinburgh}

}

\newpage
\thispagestyle{empty}

{\centering\small

First published as \emph{Basic Category Theory}, Cambridge Studies in
Advanced Mathematics, Vol.~143, Cambridge University Press, Cambridge,
2014.\\  
ISBN 978-1-107-04424-1 (hardback).

\bigskip

Information on this title:\\
\href{http://www.cambridge.org/9781107044241}{http://www.cambridge.org/9781107044241}
%

\bigskip

\copyright{} Tom Leinster 2014

\vspace*{20mm}

This arXiv version is published under a Creative Commons
Attribution-NonCommercial-ShareAlike 4.0 International licence\\ 
(CC BY-NC-SA 4.0).

\bigskip

Licence information:\\
\href{https://creativecommons.org/licenses/by-nc-sa/4.0/}{https://creativecommons.org/licenses/by-nc-sa/4.0}

\bigskip

\copyright{} Tom Leinster 2014, 2016

}

%
%
%

\frontmatter

\setcounter{page}{3}

\chapter*{Preface to the arXiv version}

This book was first published by Cambridge University Press in 2014, and is
now being published on the arXiv by mutual agreement.  CUP has consistently
supported the mathematical community by allowing authors to make free
versions of their books available online.  Readers may, in turn, wish to
support CUP by buying the printed version, available at
\href{http://www.cambridge.org/9781107044241}{\url{http://www.cambridge.org/9781107044241}}.

This electronic version is not only free; it is also freely
\emph{editable}.  For instance, if you would like to teach a course using
this book but some of the examples are unsuitable for your class, you can
remove them or add your own.  Similarly, if there is notation that you
dislike, you can easily change it; or if you want to reformat the text for
reading on a particular device, that is easy too.

In legal terms, this text is released under the Creative Commons
Attribution-NonCommercial-ShareAlike 4.0 International licence (CC BY-NC-SA
4.0).  The licence terms are available at the Creative Commons website,
\href{https://creativecommons.org/licenses/by-nc-sa/4.0/}{\url{https://creativecommons.org/licenses/by-nc-sa/4.0}}.
Broadly speaking, the licence allows you to edit and redistribute this text
in any way you like, as long as you include an accurate statement about
authorship and copyright, do not use it for commercial purposes, and
distribute it under this same licence.

In technical terms, all you need to do in order to edit this book is to
download the source files from the \href{https://arxiv.org}{arXiv} and use
\LaTeX{} (or pdflatex) in the usual way.  

This version is identical to the printed version except for the correction
of a small number of minor errors.  Thanks to all those who pointed these
out, including Martin Brandenburg, Miguel Couto, Bradley Hicks, Thomas
Moeller, and Yaokun Wu.

\bigskip

\hfill\textit{Tom Leinster, December 2016}

\newpage


\thispagestyle{empty}

\  
%
%
%

\title{Basic Category Theory}

\author{Tom Leinster\\[1ex]
\emph{University of Edinburgh}}

\bookabstract{This short introduction to category theory is for readers
with relatively little mathematical background.  At its heart is the
concept of a universal property, important throughout mathematics.  After
a chapter giving the basic definitions, the three main chapters present
three ways of expressing universal properties: via adjoint functors,
representable functors, and limits.  A final chapter ties the three
together.

For each new categorical concept, a generous supply of examples is
provided, taken from different parts of mathematics.  At points where the
leap in abstraction is particularly great (such as the Yoneda lemma), the
reader will find careful and extensive explanations.}

\bookkeywords{Category, functor, adjoint, limit, universal property.  MSC
  2010: 18A (primary), 03E (secondary).}

\tableofcontents


\addtocontents{toc}{\vspace{-4\baselineskip}} 

\clearpage 
\thispagestyle{empty}
\ 
%
%
%

\chapter*{Note to the reader}
\label{ch:preface}

This is not a sophisticated text.  In writing it, I have assumed no more
mathematical knowledge than might be acquired from an undergraduate degree
at an ordinary British university, and I have not assumed that you are used
to learning mathematics by reading a book rather than attending lectures.
Furthermore, the list of topics covered is deliberately short, omitting all
but the most fundamental parts of category theory.  A `further reading'
section points to suitable follow-on texts.

There are two things that every reader should know about this book.  One
concerns the examples, and the other is about the exercises.

Each new concept is illustrated with a generous supply of examples, but it
is not necessary to understand them all.  In courses I have taught based on
earlier versions of this text, probably no student has had the background to
understand every example.  All that matters is to understand enough
examples that you can connect the new concepts with mathematics that you
already know.

As for the exercises, I join every other textbook author in exhorting you
to do them; but there is a further important point.  In subjects such as
number theory and combinatorics, some questions are simple to state but
extremely hard to answer.  Basic category theory is not like that.  To
understand the question is very nearly to know the answer.  In most of the
exercises, there is only one possible way to proceed.  So, if you are stuck
on an exercise, a likely remedy is to go back through each term in the
question and make sure that you understand it \emph{in full}.  Take your
time.  Understanding, rather than problem solving, is the main challenge of
learning category theory.

Citations such as \citeCWM\ refer to the sources listed in `Further
reading'.

This book developed out of master's-level courses taught several times at
the University of Glasgow and, before that, at the University of Cambridge.
In turn, the Cambridge version was based on Part~III courses taught for
many years by Martin Hyland and Peter Johnstone.  Although this text is
significantly different from any of their courses, I am conscious that
certain exercises, lines of development and even turns of phrase have
persisted through that long evolution.  I would like to record my
indebtedness to them, as well as my thanks to Fran\c{c}ois Petit, my past
students, the anonymous reviewers, and the staff of Cambridge University
Press.


\mainmatter

%
%
%

\chapter*{Introduction}
\label{ch:intro}

Category theory takes a bird's eye view of mathematics.  From high in the sky,
details become invisible, but we can spot patterns that were impossible to
detect from ground level.  How is the lowest common multiple of two numbers
like the direct sum of two vector spaces?  What do discrete topological
spaces, free groups, and fields of fractions have in common?  We will
discover answers to these and many similar questions, seeing patterns in
mathematics that you may never have seen before.

The most important concept in this book is that of \emph{universal%
\index{universal!property|(}
property}.  The further you go in mathematics, especially pure mathematics,
the more universal properties you will meet.  We will spend most of our
time studying different manifestations of this concept.

Like all branches of mathematics, category theory has its own special
vocabulary, which we will meet as we go along.  But since the idea of
universal property is so important, I will use this introduction to explain
it with no jargon at all, by means of examples.

Our first example of a universal property is very simple.

\begin{iexample} 
\label{eg:univ-terminal-set}
Let $1$%
\ntn{oneset}
denote a set with one%
\index{set!one-element}
element.  (It does not matter what this element is called.)  Then $1$ has
the following property:
\begin{displaytext}
for all sets $X$, there exists a unique map from $X$ to $1$.
\end{displaytext}
(In this context, the words `map', `mapping' and `function' all mean the
same thing.)

Indeed, let $X$ be a set.  There \emph{exists} a map $X \to 1$, because we
can define $f\from X \to 1$ by taking $f(x)$ to be the single element of
$1$ for each $x \in X$.  This is the \emph{unique} map $X \to 1$, because
there is no choice in the matter: any map $X \to 1$ must send each element
of $X$ to the single element of $1$.
\end{iexample}

Phrases of the form `there exists a unique%
\index{uniqueness}
such-and-such satisfying some condition' are common in category theory.
The phrase means that there is one and only one such-and-such satisfying
the condition.  To prove the existence part, we have to show that there is
at least one.  To prove the uniqueness part, we have to show that there is
at most one; in other words, any two such-and-suches satisfying the
condition are equal.

Properties such as this are called `universal' because they state how the
object being described (in this case, the set $1$) relates to the entire
universe in which it lives (in this case, the universe of sets).  The
property begins with the words `\emph{for all} sets $X$', and therefore
says something about the relationship between $1$ and \emph{every} set $X$:
namely, that there is a unique map from $X$ to $1$.

\begin{iexample}        
\label{eg:univ-Z}
This example involves rings,%
\index{ring}
which in this book are always taken to have a multiplicative identity,
called $1$.  Similarly, homomorphisms of rings are understood to preserve
multiplicative identities.

The ring $\integers$%
\index{Z@$\integers$ (integers)!ring@as ring}
has the following property: for all rings $R$, there exists a unique
homomorphism $\integers \to R$.

To prove existence, let $R$ be a ring.  Define a function $\phi \from
\integers \to R$ by
\[
\phi(n)
=
\begin{cases}
\underbrace{1 + \cdots + 1}_n   &\text{ if } n > 0,     \\
0                               &\text{ if } n = 0,     \\
-\phi(-n)                       &\text{ if } n < 0
\end{cases}
\]
($n \in \integers$).  A series of elementary checks confirms that $\phi$ is
a homomorphism.

To prove uniqueness, let $R$ be a ring and let $\psi\from \integers \to R$
be a homomorphism.  We show that $\psi$ is equal to the homomorphism $\phi$
just defined.  Since homomorphisms preserve multiplicative identities,
$\psi(1) = 1$.  Since homomorphisms preserve addition,
\[
\psi(n) 
= 
\psi(\underbrace{1 + \cdots + 1}_n)
=
\underbrace{\psi(1) + \cdots + \psi(1)}_n
=
\underbrace{1 + \cdots + 1}_n
=
\phi(n)
\]
for all $n > 0$.  Since homomorphisms preserve zero, $\psi(0) = 0 = \phi(0)$.
Finally, since homomorphisms preserve negatives, $\psi(n) = -\psi(-n) =
-\phi(-n) = \phi(n)$ whenever $n < 0$.
\end{iexample}

Crucially, there can be essentially only \emph{one} object satisfying a
given universal property.  The word `essentially' means that two objects
satisfying the same universal property need not literally be equal, but
they are always isomorphic.  For example: 

\begin{ilemma}  
\label{lemma:Z-unique}
\index{universal!property!determines object uniquely}
Let $A$ be a ring with the following property: for all rings $R$, there exists
a unique homomorphism $A \to R$.  Then $A \iso \integers$.
\end{ilemma}

\begin{pf}
Let us call a ring with this property `initial'.  We are given that $A$
is initial, and we proved in Example~\ref{eg:univ-Z} that $\integers$ is
initial.  

Since $A$ is initial, there is a unique homomorphism $\phi\from A \to
\integers$.  Since $\integers$ is initial, there is a unique homomorphism
$\phi'\from \integers \to A$.  Now $\phi' \of \phi \from A \to A$ is a
homomorphism, but so too is the identity map $1_A\from A \to A$; hence,
since $A$ is initial, $\phi' \of \phi = 1_A$.  (This follows from the
uniqueness part of initiality, taking `$R$' to be $A$.)  Similarly, $\phi
\of \phi' = 1_\integers$.  So $\phi$ and $\phi'$ are mutually inverse, and
therefore define an isomorphism between $A$ and $\integers$.
\end{pf}

This proof has very little to do with rings.  It really belongs at a higher
level of generality.  To properly understand this, and to convey more fully
the idea of universal property, it will help to consider some more
complex examples.

\begin{iexample} 
\label{eg:univ-basis}
Let $V$ be a vector%
\index{vector space}
space with a basis $(v_s)_{s \in S}$.  (For example, if $V$ is
finite-dimensional then we might take $S = \{1, \ldots, n\}$.)  If $W$ is
another vector space, we can specify a linear map from $V$ to $W$ simply by
saying where the basis elements go.  Thus, for any $W$, there is a natural
one-to-one correspondence between
\begin{displaytext}
linear maps $V \to W$
\end{displaytext}
and
\begin{displaytext}
functions $S \to W$.
\end{displaytext}
This is because any function defined on the basis elements extends
uniquely to a linear map on $V$.

Let us rephrase this last statement.  Define a function $i\from S \to V$ by
$i(s) = v_s$ ($s \in S$).  Then $V$ together with $i$ has the following
universal property:
\[ 
\xymatrix{
S \ar[r]^i \ar[dr]_{\forall \text{ functions } f}  &
V \ar@{.>}[d]^{\exists ! \text{ linear } \bar{f}} \\
&
\forall W.
}
\]
This diagram means that for all vector spaces $W$ and all functions $f\from
S \to W$, there exists a unique linear map $\bar{f}: V \to W$ such that
$\bar{f} \of i = f$.  The symbol $\forall$%
\ntn{forall}
means `for all', and the symbols $\exists !$%
\ntn{ei}
 mean `there exists a unique'.%
\index{uniqueness}

Another way to say `$\bar{f} \of i = f$' is `$\bar{f}(v_s) = f(s)$ for all
$s \in S$'.  So, the diagram asserts that every function $f$ defined on the
basis elements extends uniquely to a linear map $\bar{f}$ defined on the
whole of $V$.  In other words still, the function
\[
\begin{array}{ccc}
\{ \text{linear maps } V \to W \}     &
\ \to \             &
\{ \text{functions } S \to W \}       \\
\bar{f}        &\mapsto        &\bar{f} \of i       
\end{array}
\]
is bijective. 
\end{iexample}

\begin{iexample} 
\label{eg:univ-discrete}
Given a set $S$, we can build a topological space $D(S)$%
\ntn{discrete-space}
by equipping $S$ with the \demph{discrete%
\index{topological space!discrete} 
topology}: all subsets are open.  With this topology, \emph{any} map from
$S$ to a space $X$ is continuous.

Again, let us rephrase this.  Define a function $i\from S \to D(S)$ by
$i(s) = s$ ($s \in S$).  Then $D(S)$ together with $i$ has the following
universal property:
\[
\xymatrix{
S \ar[r]^-i \ar[dr]_{\forall \text{ functions } f}       &
D(S) \ar@{.>}[d]^{\exists! \text{ continuous } \bar{f}}   \\
&
\forall X.
}
\]
In other words, for all topological spaces $X$ and all functions $f\from S \to
X$, there exists a unique continuous map $\bar{f}\from D(S) \to X$ such that
$\bar{f} \of i = f$.  The continuous map $\bar{f}$ is the same thing as the
function $f$, except that we are regarding it as a continuous map between
topological spaces rather than a mere function between sets.

You may feel that this universal property is almost too trivial to mean
anything.  But if we change the definition of $D(S)$~-- say from the
discrete to the indiscrete topology, in which the only open sets are
$\emptyset$ and $S$~-- then the property becomes false.  So this property
really does say something about the discrete topology.  What it says is
that all maps out of a discrete space are continuous.

Indeed, given $S$, the universal property determines $D(S)$ and $i$
uniquely (or rather, uniquely up to isomorphism; but who could want more?).
The proof of this is similar to that of Lemma~\ref{lemma:Z-unique} above
and Lemma~\ref{lemma:tensor-unique} below.
\end{iexample}

\begin{iexample} 
\label{eg:univ-tensor}
Given vector%
\index{vector space}
spaces $U$, $V$ and $W$, a \demph{bilinear%
\index{map!bilinear}
map} $f\from U \times V \to W$ is a function $f$ that is linear in each
variable:
\begin{align*}
f(u, v_1 + \lambda v_2)         &=      
f(u, v_1) + \lambda f(u, v_2),  \\
f(u_1 + \lambda u_2, v)         &=      
f(u_1, v) + \lambda f(u_2, v)   
\end{align*}
for all $u, u_1, u_2 \in U$, $v, v_1, v_2 \in V$, and scalars $\lambda$.  A
good example is the scalar product (dot product), which is a bilinear map
\[
\begin{array}{ccc}
\reals^n \times \reals^n        &\to            &\reals                 \\
(\mathbf{u}, \mathbf{v})        &\mapsto        &\mathbf{u}.\mathbf{v}
\end{array}
\]
of real vector spaces.  The vector product (cross product) $\reals^3 \times
\reals^3 \to \reals^3$ is also bilinear.  

Let $U$ and $V$ be vector spaces.  It is a fact that there is a `universal
bilinear map out of $U \times V$'.  In other words, there exist a certain
vector space $T$%
\index{tensor product|(}
and a certain bilinear map $b\from U \times V \to T$ with the following
universal property:
\begin{equation}        
\label{eq:univ-tensor} 
\begin{array}{c}
\xymatrix{
U \times V \ar[r]^-b \ar[dr]_{\forall \text{ bilinear } f} &
T \ar@{.>}[d]^{\exists! \text{ linear } \bar{f}}     \\
&
\forall W.
}
\end{array}
\end{equation}
Roughly speaking, this property says that bilinear maps out of $U \times V$
correspond one-to-one with linear maps out of $T$.

Even without knowing that such a $T$ and $b$ exist, we can immediately
prove that this universal property determines $T$ and $b$ uniquely up to
isomorphism.  The proof is essentially the same as that of
Lemma~\ref{lemma:Z-unique}, but looks more complicated because of the more
complicated universal property.
\end{iexample}

\begin{ilemma}   
\label{lemma:tensor-unique}
\index{universal!property!determines object uniquely}
Let $U$ and $V$ be vector spaces.  Suppose that $b\from U \times V \to T$
and $b'\from U \times V \to T'$ are both universal bilinear maps out of $U
\times V$.  Then $T \iso T'$.  More precisely, there exists a unique
isomorphism $j\from T \to T'$ such that $j \of b = b'$.
\end{ilemma}

In the proof that follows, it does not actually matter what `bilinear',
`linear' or even `vector space' mean.  The hard part is getting the logic
straight.  That done, you should be able to see that there is really only
one possible proof.  For instance, to use the universality of $b$, we will
have to choose some bilinear map $f$ out of $U\times V$.  There are only
two in sight, $b$ and $b'$, and we use each in the appropriate place.

\begin{pf}
In diagram~\eqref{eq:univ-tensor}, take $\Bigl(U \times V \toby{f} W\Bigr)$
to be $\Bigl(U \times V \toby{b'} T'\Bigr)$.  This gives a linear map
$j\from T \to T'$ satisfying $j \of b = b'$.  Similarly, using the
universality of $b'$, we obtain a linear map $j': T' \to T$ satisfying $j'
\of b' = b$:
\[
\xymatrix{
        &T \ar[d]^-j     \\
U \times V \ar[ur]^-b \ar[r]|-{b'} \ar[dr]_b      &
T' \ar[d]^{j'}  \\
&
T.
}
\]
Now $j' \of j\from T \to T$ is a linear map satisfying $(j' \of j) \of b =
b$; but also, the identity map $1_T\from T \to T$ is linear and satisfies
$1_T \of b = b$.  So by the uniqueness part of the universal property of
$b$, we have $j' \of j = 1_T$.  (Here we took the `$f$'
of~\eqref{eq:univ-tensor} to be $b$.)  Similarly, $j \of j' = 1_{T'}$.  So
$j$ is an isomorphism.  
\end{pf}

In Example~\ref{eg:univ-tensor}, it was stated that given vector spaces $U$
and $V$, there exists a pair $(T, b)$ with the universal property
of~\eqref{eq:univ-tensor}.  We just proved that there is essentially only
one such pair $(T, b)$.  The vector space $T$ is called the \demph{tensor
product} of $U$ and $V$, and is written as $U \otimes V$.%
\ntn{tensor}
Tensor products are very important in algebra.  They reduce the study of
bilinear maps to the study of linear maps, since a bilinear map out of $U
\times V$ is really the same thing as a linear map out of $U \otimes V$.

However, tensor products will not be important in this book.  The real
lesson for us is that it is safe to speak of \emph{the} tensor product, not
just \emph{a} tensor product, and the reason for that is
Lemma~\ref{lemma:tensor-unique}.  This is a general point that applies to
anything satisfying a universal property.

Once you know a universal property of an object, it often does no harm to
forget how it was constructed. For instance, if you look through a pile of
algebra books, you will find several different ways of constructing the
tensor product of two vector spaces.  But once you have proved that the
tensor product satisfies the universal property, you can forget the
construction.%
\index{tensor product|)} 
The universal property tells you all you need to know, because it
determines the object uniquely up to isomorphism.

\begin{iexample} 
\label{eg:univ-kernel}
Let $\theta\from G \to H$ be a homomorphism of groups.%
\index{group}
Associated with $\theta$ is a diagram
\begin{equation}        
\label{eq:eq-groups}
\index{kernel}
\begin{array}{c}
\xymatrix@M+.5ex{
\ker(\theta) \ar@{^{(}->}[r]^-{\iota} &
G \ar@<.5ex>[r]^\theta \ar@<-.5ex>[r]_\epsln      &
H,
}
\end{array}
\end{equation}
where $\iota$ is the inclusion of $\ker(\theta)$ into $G$ and $\epsln$ is
the trivial homomorphism.  `Inclusion'%
\index{inclusion}
means that $\iota(x) = x$ for all $x \in \ker(\theta)$, and `trivial' means
that $\epsln(g) = 1$ for all $g \in G$.  The symbol $\incl$%
\ntn{incl}
is often used for inclusions; it is a combination of a subset symbol
$\subset$ and an arrow.

The map $\iota$ into $G$ satisfies $\theta \of \iota = \epsln \of \iota$,
and is universal as such.  Exercise~\ref{ex:ker-groups} asks you to make
this precise.
\end{iexample}

Here is our final example of a universal property.

\begin{iexample} 
\label{eg:univ-van-Kampen}
Take a topological%
\index{topological space}
space covered by two open subsets: $X = U \cup V$.  The diagram
\[
\xymatrix@M+0.5ex{
U \cap V \ar@{^{(}->}[r]^-{i} \ar@{^{(}->}[d]_-{j}      &
U \ar@{^{(}->}[d]^-{j'} \\
V \ar@{^{(}->}[r]_-{i'} &
X
}
\]
of inclusion maps has a universal property in the world of topological spaces
and continuous maps, as follows:
\begin{equation} 
\label{eq:pushout-vK}
\begin{array}{c}
\xymatrix@M+.5ex{
U \cap V \ar@{^{(}->}[r]^-{i} \ar@{^{(}->}[d]_-{j}      &
U \ar@{^{(}->}[d]^-{j'} \ar@/^/[ddr]^{\forall f} &
\\
V \ar@{^{(}->}[r]_-{i'} \ar@/_/[drr]_{\forall g} &
X \ar@{.>}[dr]|{\exists! h}     &
\\
&
&
\forall Y.
}
\end{array}
\end{equation}
The diagram means that given $Y$, $f$ and $g$ such that $f \of i = g \of
j$, there is exactly one continuous map $h\from X \to Y$ such that $h \of
j' = f$ and $h \of i' = g$.

Under favourable conditions, the induced diagram 
\[
\xymatrix@M+0.5ex{
\pi_1(U \cap V) \ar[r]^-{i_*} \ar[d]_-{j_*}      &
\pi_1(U) \ar[d]^-{j'_*} \\
\pi_1(V) \ar[r]_-{i'_*} &
\pi_1(X)
}
\]
of fundamental%
\index{group!fundamental}
groups has the same property in the world of groups and group
homomorphisms.  This is \emph{van%
\index{van Kampen's theorem}
Kampen's theorem}.  In fact, van Kampen stated his theorem in a much more
complicated way.  Stating it transparently requires some categorical
language, but he was working in the 1930s, before category theory had been
born.
\end{iexample}

You have now seen several examples of universal properties.  As this book
progresses, we will develop different ways of talking about them.  Once we
have set up the basic vocabulary of categories and functors, we will study
\emph{adjoint functors}, then \emph{representable functors}, then
\emph{limits}.  Each of these provides an approach to universal properties,
and each places the idea in a different light.  For instance,
Examples~\ref{eg:univ-basis} and~\ref{eg:univ-discrete} can most readily be
described in terms of adjoint functors, Example~\ref{eg:univ-tensor} via
representable functors, and Examples~\ref{eg:univ-terminal-set},
\ref{eg:univ-Z}, \ref{eg:univ-kernel} and~\ref{eg:univ-van-Kampen} in terms
of limits.%
\index{universal!property|)}

\exs

\begin{iquestion}
Let $S$ be a set.  The \demph{indiscrete}%
\index{topological space!indiscrete}%
\index{indiscrete space} 
topological space $I(S)$%
\ntn{indiscrete-space}
is the space whose set of points is $S$ and whose only open subsets are
$\emptyset$ and $S$ itself.  Imitating Example~\ref{eg:univ-discrete}, find
a universal property satisfied by the space $I(S)$.
\end{iquestion}

\begin{iquestion}        
\label{ex:ker-groups}
Fix a group homomorphism $\theta\from G \to H$.  Find a universal property
satisfied by the pair $(\ker(\theta), \iota)$%
\index{kernel}
of diagram~\eqref{eq:eq-groups}.  (This property can~-- indeed, must~-- make
reference to $\theta$.)
\end{iquestion}

\begin{iquestion}
Verify the universal property shown in diagram~\eqref{eq:pushout-vK}.
\end{iquestion}

\begin{iquestion}       
\label{ex:Zx}
Denote by $\integers[x]$%
\ntn{poly-one}
the polynomial%
\index{ring!polynomial}
ring over $\integers$ in one variable.
\begin{enumerate}[(b)]
\item   
\label{part:Zx-main}
Prove that for all rings $R$ and all $r \in R$, there exists a unique ring
homomorphism $\phi\from \integers[x] \to R$ such that $\phi(x) = r$.

\item 
Let $A$ be a ring and $a \in A$.  Suppose that for all rings $R$ and all
$r \in R$, there exists a unique ring homomorphism $\phi\from A \to R$ such that
$\phi(a) = r$.  Prove that there is a unique isomorphism $\iota\from
\integers[x] \to A$ such that $\iota(x) = a$.  
\end{enumerate}
\end{iquestion}

\begin{iquestion}
Let $X$ and $Y$ be vector spaces.  
\begin{enumerate}[(b)]
\item   
\label{part:vs-prod}
For the purposes of this exercise only, a \emph{cone} is a triple $(V, f_1,
f_2)$ consisting of a vector space $V$, a linear map $f_1\from V \to X$, and
a linear map $f_2\from V \to Y$.  Find a cone $(P, p_1,
p_2)$ with the following property: for all cones $(V, f_1, f_2)$, there exists
a unique linear map $f\from V \to P$ such that $p_1 \of f = f_1$ and $p_2 \of f
= f_2$.

\item 
Prove that there is essentially only one cone with the property stated
in~\bref{part:vs-prod}.  That is, prove that if $(P, p_1, p_2)$ and $(P',
p'_1, p'_2)$ both have this property then there is an isomorphism $i\from P
\to P'$ such that $p'_1 \of i = p_1$ and $p'_2 \of i = p_2$.

\item   
\label{part:vs-coprod} 
For the purposes of this exercise only, a \emph{cocone} is a triple $(V,
f_1, f_2)$ consisting of a vector space $V$, a linear map $f_1\from X \to
V$, and a linear map $f_2\from Y \to V$.  Find a cocone $(Q, q_1, q_2)$
with the following property: for all cocones $(V, f_1, f_2)$, there exists
a unique linear map $f\from Q \to V$ such that $f \of q_1 = f_1$ and $f \of
q_2 = f_2$.

\item 
Prove that there is essentially only one cocone with the property stated
 in~\bref{part:vs-coprod}, in a sense that you should make precise.  
\end{enumerate}
\end{iquestion}


%
%
%

\chapter{Categories, functors and natural transformations}
\label{ch:cfnt}

A category is a system of related objects.  The objects do not live in
isolation: there is some notion of map between objects, binding them
together.

Typical examples of what `object' might mean are `group' and `topological
space', and typical examples of what `map' might mean are `homomorphism'
and `continuous map', respectively.  We will see many examples, and we will
also learn that some categories have a very different flavour from the two
just mentioned.  In fact, the `maps' of category theory need not be
anything like maps in the sense that you are most likely to be familiar
with.

Categories are \emph{themselves} mathematical objects, and with that in
mind, it is unsurprising that there is a good notion of `map between
categories'.  Such maps are called functors.  More surprising, perhaps, is
the existence of a third level: we can talk about maps between
\emph{functors}, which are called natural transformations.  These, then,
are maps between maps between categories.

In fact, it was the desire to formalize the notion of natural
transformation that led to the birth of category theory.  By the early
1940s, researchers in algebraic topology had started to use the phrase
`natural transformation', but only in an informal way.  Two mathematicians,
Samuel Eilenberg%
\index{Eilenberg, Samuel}
and Saunders Mac Lane,%
\index{Mac~Lane, Saunders}
saw that a precise definition was needed.  But before they could define
natural transformation, they had to define functor; and before they could
define functor, they had to define category.  And so the subject was born.

Nowadays, the uses of category theory have spread far beyond algebraic
topology.  Its tentacles extend into most parts of pure mathematics.  They
also reach some parts of applied mathematics; perhaps most notably,
category theory has become a standard tool in certain parts of computer%
\index{computer science}
science.  Applied%
\index{applied mathematics}
mathematics is more than just applied differential equations!

\section{Categories}
\label{sec:cats}

\begin{defn}
A \demph{category}%
\index{category}
$\cat{A}$ consists of:
\begin{itemize}
\item 
a collection $\ob(\cat{A})$%
\ntn{ob}
of \demph{objects};%
\index{object}

\item 
for each $A, B \in \ob(\cat{A})$, a collection $\cat{A}(A, B)$%
\ntn{hom-set-default}
of \demph{maps}%
\index{map}
or \demph{arrows}%
\index{arrow}
or \demph{morphisms}%
\index{morphism}
from $A$ to $B$;

\item 
for each $A, B, C \in \ob(\cat{A})$, a function
\[
\begin{array}{ccc}
\cat{A}(B, C) \times \cat{A}(A, B) &
\to	&
\cat{A}(A, C)	\\
(g, f)	&
\mapsto	&
g \of f,%
\ntn{of}
\end{array}
\]
called \demph{composition};%
\index{composition}

\item 
for each $A \in \ob(\cat{A})$, an element $1_A$%
\ntn{id-map}
of $\cat{A}(A, A)$, called the \demph{identity}%
\index{identity}
on $A$,
\end{itemize}
satisfying the following axioms:
\begin{itemize}
\item 
\demph{associativity}:%
\index{associativity}
for each $f \in \cat{A}(A, B)$, $g \in \cat{A}(B, C)$ and $h \in \cat{A}(C,
D)$, we have $(h \of g) \of f = h \of (g \of f)$;

\item 
\demph{identity%
\index{identity}
laws}: for each $f \in \cat{A}(A, B)$, we have $f \of 1_A = f = 1_B \of f$.
\end{itemize}
\end{defn}

\begin{remarks}  
\label{rmks:defn-cat}
\begin{enumerate}[(b)]
\item   
\label{item:defn-cat-notn}
We often write:
\begin{displaytext}
\begin{tabular}{rcl}
$A \in \cat{A}$			&to mean	&$A \in \ob(\cat{A})$;	\\
$f\from A\to B$ or $A \toby{f} B$&to mean 	&$f \in \cat{A}(A,B)$;%
\ntn{arrow}
\\
$gf$				&to mean	&$g \of f$.%
\ntn{juxt}
\end{tabular}
\end{displaytext}
People also write $\cat{A}(A, B)$ as $\Hom_{\cat{A}}(A, B)$%
\ntn{Hom}
or $\Hom (A, B)$.  The notation `$\Hom$' stands for homomorphism, from one
of the earliest examples of a category.

\item	
\label{rmk:defn-cat:loosely}
The definition of category is set up so that in general, from each string
\[
A_0 \toby{f_1} 
A_1 \toby{f_2}
\ \cdots \ 
\toby{f_n} A_n
\]
of maps in $\cat{A}$, it is possible to construct exactly one%
\index{uniqueness!constructions@of constructions}
map
\[
A_0 \to A_n
\]
(namely, $f_n f_{n - 1} \cdots f_1$).  If we are given extra information
then we may be able to construct other maps $A_0 \to A_n$; for instance, if
we happen to know that $A_{n - 1} = A_n$, then $f_{n - 1} f_{n - 2} \cdots
f_1$ is another such map.  But we are speaking here of the \emph{general}
situation, in the absence of extra information.

For example, a string like this with $n = 4$ gives rise to maps
\[
\xymatrix@=8em{
A_0
\ar@<1ex>[r]^{((f_4 f_3)f_2)f_1}
\ar@<-1ex>[r]_{(f_4(1_{A_3} f_3))((f_2 f_1)1_{A_0})} &
A_4,
}
\]
but the axioms imply that they are equal.  It is safe to omit the brackets
and write both as $f_4 f_3 f_2 f_1$.

Here it is intended that $n \geq 0$.  In the case $n = 0$, the statement is
that for each object $A_0$ of a category, it is possible to construct
exactly one map $A_0 \to A_0$ (namely, the identity $1_{A_0}$).  An
identity map can be thought of as a zero-fold%
\index{identity!zero-fold composite@as zero-fold composite}
composite, in much the same way that the number $1$ can be thought of as
the product of zero numbers.

\item 
We often speak of \demph{commutative%
\index{diagram!commutative}
diagrams}.  For instance, given objects and maps
\[
\xymatrix{
A \ar[rr]^f \ar[d]_h    &               &B \ar[d]^g     \\
C \ar[r]_i              &D \ar[r]_j     &E
}
\]
in a category, we say that the diagram \demph{commutes}%
\index{commutes}
if $gf = jih$.  Generally, a diagram is said to commute if whenever there
are two paths from an object $X$ to an object $Y$, the map from $X$ to $Y$
obtained by composing along one path is equal to the map obtained by
composing along the other.

\item 
The slightly vague word `collection'%
\index{collection}
means \emph{roughly} the same as `set', although if you know about such
things, it is better to interpret it as meaning `class'.%
\index{class}
We come back to this in Chapter~\ref{ch:sets}.

\item 
If $f \in \cat{A}(A, B)$, we call $A$ the \demph{domain}%
\index{domain}
and $B$ the \demph{codomain}%
\index{codomain}
of $f$.  Every map in every category has a definite domain and a definite
codomain.  (If you believe it makes sense to form the intersection of an
arbitrary pair of abstract sets, you should add to the definition of
category the condition that $\cat{A}(A, B) \cap \cat{A}(A', B') =
\emptyset$ unless $A = A'$ and $B = B'$.)
\end{enumerate}
\end{remarks}

\begin{examples}[Categories \:of \:mathematical \:structures]        
\label{egs:cats-of}
\begin{enumerate}[(b)]
\item 
There is a category $\Set$%
\ntn{Set}
described as follows.  Its objects are sets.%
\index{set!category of sets}
Given sets $A$ and $B$, a map from $A$ to $B$ in the category $\Set$ is
exactly what is ordinarily called a map (or mapping, or function) from $A$
to $B$.  Composition in the category is ordinary composition of functions,
and the identity maps are again what you would expect.

In situations such as this, we often do not bother to specify the
composition and identities.  We write `the category of sets and functions',
leaving the reader to guess the rest.  In fact, we usually go further and
call it just `the category of sets'.

\item 
There is a category $\Grp$%
\ntn{Grp}
of groups,%
\index{group!category of groups}
whose objects are groups and whose maps are group homomorphisms.

\item 
Similarly, there is a category $\Ring$%
\ntn{Ring}
of rings%
\index{ring!category of rings}
and ring homomorphisms.

\item 
For each field $k$, there is a category $\Vect_k$%
\ntn{Vect}
of vector%
\index{vector space!category of vector spaces}
spaces over $k$ and linear maps between them.

\item 
There is a category $\Tp$%
\ntn{Top}
of topological%
\index{topological space!category of topological spaces}
spaces and continuous maps.
\end{enumerate}
\end{examples}

This chapter is mostly about the interaction \emph{between} categories,
rather than what goes on \emph{inside} them.  We will, however, need the
following definition.

\begin{defn}    
\label{defn:isomorphism}
A map $f\from A \to B$ in a category $\cat{A}$ is an \demph{isomorphism}%
\index{isomorphism}
if there exists a map $g\from B \to A$ in $\cat{A}$ such that $gf = 1_A$
and $fg = 1_B$.
\end{defn}

In the situation of Definition~\ref{defn:isomorphism}, we call $g$ the
\demph{inverse}%
\index{inverse}
of $f$ and write $g = f^{-1}$.%
\ntn{inverse}
(The word `the' is justified by Exercise~\ref{ex:unique-inverse}.)  If
there exists an isomorphism from $A$ to $B$, we say that $A$ and $B$ are
\demph{isomorphic} and write $A \iso B$.%
\ntn{iso}

\begin{example}
\label{eg:iso-Set}
The isomorphisms in $\Set$%
\index{set!category of sets!isomorphisms in}
are exactly the bijections.  This\linebreak statement is not quite a logical
triviality.  It amounts to the assertion that a function has a two-sided
inverse if and only if it is injective and surjective.
\end{example}

\begin{example}
The isomorphisms in $\Grp$%
\index{group!category of groups!isomorphisms in}
are exactly the isomorphisms of groups.  Again, this is not quite trivial,
at least if you were taught that the definition of group isomorphism is
`bijective homomorphism'.  In order to show that this is equivalent to
being an isomorphism in $\Grp$, you have to prove that the inverse of a
bijective homomorphism is also a homomorphism.

Similarly, the isomorphisms in $\Ring$%
\index{ring!category of rings!isomorphisms in}
are exactly the isomorphisms of rings.
\end{example}

\begin{example}
The isomorphisms in $\Tp$%
\index{topological space!category of topological spaces!isomorphisms in}
are exactly the homeomorphisms.  Note that, in contrast to the situation in
$\Grp$ and $\Ring$, a bijective map in $\Tp$ is not necessarily an
isomorphism.  A classic example is the map
\[
\begin{array}{ccc}
[0, 1)  &\to            &\{z \in \complexes \such \left|z\right| = 1\}     \\
t       &\mapsto        &e^{2\pi i t},
\end{array}
\]
which is a continuous bijection but not a homeomorphism.
\end{example}

The examples of categories mentioned so far are important, but could give a
false impression.  In each of them, the objects of the category are sets
with structure (such as a group structure, a topology, or, in the case
of $\Set$, no structure at all).  The maps are the functions preserving the
structure, in the appropriate sense.  And in each of them, there is a
clear sense of what the elements of a given object are.

However, not all categories are like this.  In general, the objects of a
category are not `sets equipped with extra stuff'.  Thus, in a general
category, it does not make sense to talk about the `elements' of an object.
(At least, it does not make sense in an immediately obvious way; we return
to this in Definition~\ref{defn:gen-elt}.)  Similarly, in a general
category, the maps need not be mappings or functions in the usual sense.
So:

\begin{slogan}
The objects of a category need not be remotely like sets.%
\index{object!need not resemble set}
\end{slogan}
\begin{slogan}
The maps in a category need not be remotely like functions.
\index{map!need not resemble function}
\end{slogan}
The next few examples illustrate these points.  They also show that,
contrary to the impression that might have been given so far, categories
need not be enormous.  Some categories are small, manageable structures in
their own right, as we now see.

\begin{examples}[Categories \,as \,mathematical \,structures]        
\label{egs:cats-as}
\begin{enumerate}[(b)]
\item   
\label{eg:cats-as:graphs}
A category can be specified%
\index{category!drawing of}
by saying directly what its objects, maps, composition and identities are.
For example, there is a category $\emptyset$%
\ntn{empty-cat}
with no objects or maps at all.  There is a category $\One$%
\ntn{terminal-cat}
 with one object
and only the identity map.  It can be drawn like this:
\[
\bullet
\]
(Since every object is required to have an identity map on it, we usually
do not bother to draw the identities.)  There is another category that can
be drawn as
\[
\bullet \to \bullet 
\qquad
\text{or}
\qquad
A \toby{f} B,
\]
with two objects and one non-identity map, from the first object to the
second.  (Composition is defined in the only possible way.)  To reiterate
the points made above, it is not obvious what an `element' of $A$ or $B$
would be, or how one could regard $f$ as a `function' of any sort.

It is easy to make up more complicated examples.  For instance, here are three
more categories:
\[
\begin{array}{c}
\xymatrix{
\bullet \ar@<.5ex>[r] \ar@<-.5ex>[r] &\bullet
}
\end{array}
\qquad
\begin{array}{c}
\xymatrix{
        &B \ar[dr]^g    &       \\
A \ar[ur]^f \ar[rr]_{gf} & &C    
}
\end{array}
\qquad
\begin{array}{c}
\xymatrix{
        &\bullet \ar[dl]_{kj} \ar[r]^f \ar[dr]|{hj=gf} \ar[d]_{j} &
\bullet \ar[d]^g        \\
\bullet &\bullet \ar[l]^k \ar[r]_{h}   &\bullet
}
\end{array}
\]

\item   
\label{eg:cats-as:discrete} 
Some categories contain no maps at all apart from identities (which, as
categories, they are obliged to have).  These are called \demph{discrete}%
\index{category!discrete}
categories.  A discrete category amounts to just a class of objects.  More
poetically, a category is a collection of objects related to one another to
a greater or lesser degree; a discrete category is the extreme case in
which each object is totally isolated from its companions.

\item   
\label{eg:cats-as:groups}
A group is essentially the same thing as a category that has only one%
\index{category!one-object|(}%
\index{group!one-object category@as one-object category}
object and in which all the maps are isomorphisms.

To understand this, first consider a category $\cat{A}$ with just one
object.  It is not important what letter or symbol we use to denote the
object; let us call it $A$.  Then $\cat{A}$ consists of a set (or class)
$\cat{A}(A, A)$, an associative composition function
\[
\of\from  \cat{A}(A, A) \times \cat{A}(A, A) \to \cat{A}(A, A),
\]
and a two-sided unit $1_A \in \cat{A}(A, A)$.  This would make $\cat{A}(A,
A)$ into a group, except that we have not mentioned inverses.  However, to
say that every map in $\cat{A}$ is an isomorphism is exactly to say that
every element of $\cat{A}(A, A)$ has an inverse with respect to $\of$.

If we write $G$ for the group $\cat{A}(A, A)$, then the situation is this:
\begin{displaytext}
\begin{tabular}{l@{\hspace{2em}}l}
\emph{category $\cat{A}$ with single object $A$}       &
\emph{corresponding group $G$}        \\[1ex]
maps in $\cat{A}$       &elements of $G$                \\
$\of$ in $\cat{A}$      &$\cdot$ in $G$                 \\
$1_A$                   &$1 \in G$      \\
\end{tabular}
\end{displaytext}
The category $\cat{A}$ looks something like this:
\[
\SelectTips{cm}{}
\xymatrix{A \ar@(r,u)@<-.5ex>[]_{} \ar@(u,l) \ar@(ld,rd)[]_{}}
\]
The arrows represent different maps $A \to A$, that is, different elements of
the group $G$. 

What the object of $\cat{A}$ is called makes no difference.  It matters
exactly as much as whether we choose $x$ or $y$ or $t$ to denote some
variable in an algebra problem, which is to say, not at all.  Later we will
define `equivalence' of categories, which will enable us to make a precise
statement: the category of groups is equivalent to the category of (small)
one-object categories in which every map is an isomorphism
(Example~\ref{eg:mon-one-obj-eqv}).

The first time one meets the idea that a group is a kind of category, it is
tempting to dismiss it as a coincidence or a trick.  But it is not; there
is real content.

To see this, suppose that your education had been shuffled and that you
already knew about categories before being taught about groups.  In your
first group theory class, the lecturer declares that a group is supposed to
be the system of all symmetries of an object.  A symmetry of an object $X$,
she says, is a way of mapping $X$ to itself in a reversible or invertible
manner.  At this point, you realize that she is talking about a very
special type of category.  In general, a category is a system consisting of
\emph{all} the mappings (not usually just the invertible ones) between
\emph{many} objects (not usually just one).  So a group is just a category
with the special properties that all the maps are invertible and there is
only one object.

\item   
\label{eg:cats-as:monoids}
The inverses played no essential part in the previous example, suggesting that
it is worth thinking about `groups without inverses'.  These are called
monoids.  

Formally, a \demph{monoid}%
\index{monoid}
is a set equipped with an associative binary operation and a two-sided unit
element.  Groups describe the reversible transformations, or symmetries,
that can be applied to an object; monoids describe the
not-necessarily-reversible transformations.  For instance, given any set
$X$, there is a group consisting of all bijections $X \to X$, and there is
a monoid consisting of all functions $X \to X$.  In both cases, the binary
operation is composition and the unit is the identity function on $X$.
Another example of a monoid is the set $\nat = \{0, 1, 2, \ldots\}$%
\ntn{nat}
of natural%
\index{natural numbers}
numbers, with $+$ as the operation and $0$ as the unit.  Alternatively, we
could take the set $\nat$ with $\cdot$ as the operation and $1$ as the
unit.

A category with one%
\index{monoid!one-object category@as one-object category}
object is essentially the same thing as a monoid, by the same argument as
for groups.  This is stated formally in Example~\ref{eg:mon-one-obj-eqv}.%
\index{category!one-object|)}

\item   
\label{eg:cats-as:orders}
A \demph{preorder}%
\index{preorder}
is a reflexive transitive binary relation.  A \demph{preordered set} $(S,
\mathord{\leq})$%
\ntn{leq}
is a set $S$ together with a preorder $\leq$ on it.  Examples: $S = \reals$
and $\leq$ has its usual meaning; $S$ is the set of subsets of $\{1,
\ldots, 10\}$ and $\leq$ is $\sub$ (inclusion); $S = \integers$ and $a \leq
b$ means that $a$ divides $b$.

A preordered set can be regarded as a category $\cat{A}$ in which, for each
$A, B \in \cat{A}$, there is at most one map from $A$ to $B$.  To see this,
consider a category $\cat{A}$ with this property.  It is not important what
letter we use to denote the unique map from an object $A$ to an object $B$;
all we need to record is which pairs $(A, B)$ of objects have the property
that a map $A \to B$ does exist.  Let us write $A \leq B$ to mean that
there exists a map $A \to B$.

Since $\cat{A}$ is a category, and categories have composition, if $A \leq
B \leq C$ then $A \leq C$.  Since categories also have identities, $A \leq
A$ for all $A$.  The associativity and identity axioms are automatic.  So,
$\cat{A}$ amounts to a collection of objects equipped with a transitive
reflexive binary relation, that is, a preorder.  One can think of the
unique map $A \to B$ as the statement or assertion that $A \leq B$.

An \demph{order}%
\index{ordered set}
on a set is a preorder $\leq$ with the property that if $A \leq B$ and $B
\leq A$ then $A = B$.  (Equivalently, if $A \iso B$ in the corresponding
category then $A = B$.)  Ordered sets are also called \demph{partially
ordered sets}%
\index{partially ordered set} 
or \demph{posets}.%
\index{poset}
An example of a preorder that is not%
\index{ordered set!preordered set@vs.\ preordered set}
an order is the divisibility relation $\divides$ on $\integers$: for there
we have $2 \divides {-2}$ and $-2 \divides 2$ but $2 \neq -2$.
\end{enumerate}
\end{examples}

Here are two ways of constructing new categories from old. 

{\sloppy
\begin{constn} 
\label{constn:op-cat}
Every category $\cat{A}$ has an \demph{opposite}%
\index{category!opposite}
or \demph{dual}%
\index{duality}
category $\cat{A}^\op$,%
\ntn{op}
defined by reversing the arrows.  Formally, $\ob(\cat{A}^\op) =
\ob(\cat{A})$ and $\cat{A}^\op(B, A) = \cat{A}(A, B)$ for all objects $A$
and $B$.  Identities in $\cat{A}^\op$ are the same as in $\cat{A}$.
Composition in $\cat{A}^\op$ is the same as in $\cat{A}$, but with the
arguments reversed.  To spell this out: if $A \toby{f} B \toby{g} C$ are
maps in $\cat{A}^\op$ then $A \otby{f} B \otby{g} C$ are maps in $\cat{A}$;
these give rise to a map $A \otby{f \of g} C$ in $\cat{A}$, and the
composite of the original pair of maps is the corresponding map $A \to C$
in $\cat{A}^\op$.

So, arrows $A \to B$ in $\cat{A}$ correspond to arrows $B \to A$ in
$\cat{A}^\op$.  According to the definition above, if $f\from A \to B$ is
an arrow in $\cat{A}$ then the corresponding arrow $B \to A$ in
$\cat{A}^\op$ is also called $f$.  Some people prefer to give it a
different name, such as $f^\op$.
\end{constn}
}

\begin{remark}  
\label{rmk:principle-duality}
The \demph{principle of duality}%
\index{duality!principle of}
is fundamental to category theory.  Informally, it states that every
categorical definition, theorem and proof has a \demph{dual}, obtained by
reversing all the arrows.  Invoking the principle of duality can save work:
given any theorem, reversing the arrows throughout its statement and proof
produces a dual theorem.  Numerous examples of duality appear throughout
this book.
\end{remark}

\begin{constn}  
\label{constn:prod-cat}
Given categories $\cat{A}$ and $\cat{B}$, there is a \demph{product%
\index{category!product of categories}
category} $\cat{A} \times \cat{B}$,%
\ntn{prod-cat}
 in which
\begin{align*}
\ob(\cat{A} \times \cat{B})     &
=      
\ob(\cat{A}) \times \ob(\cat{B}),\\
(\cat{A} \times \cat{B})((A, B), (A', B'))      &
=      
\cat{A}(A, A') \times \cat{B}(B, B').
\end{align*}
Put another way, an object of the product category $\cat{A} \times \cat{B}$
is a pair $(A, B)$ where $A \in \cat{A}$ and $B \in \cat{B}$.  A map $(A,
B) \to (A', B')$ in $\cat{A} \times \cat{B}$ is a pair $(f, g)$ where
$f\from A \to A'$ in $\cat{A}$ and $g\from B \to B'$ in $\cat{B}$.  For the
definitions of composition and identities in $\cat{A} \times \cat{B}$, see
Exercise~\ref{ex:prod-cat}.
\end{constn}

\exs

\begin{question}
Find three examples of categories not mentioned above.
\end{question}

\begin{question}        
\label{ex:unique-inverse}
Show that a map in a category can have at most one inverse.  That is, given
a map $f\from A \to B$, show that there is at most one map $g\from B \to A$
such that $gf = 1_A$ and $fg = 1_B$.  
\end{question}

\begin{question}        
\label{ex:prod-cat}
Let $\cat{A}$ and $\cat{B}$ be categories.
Construction~\ref{constn:prod-cat} defined the product category $\cat{A}
\times \cat{B}$, except that the definitions of composition and identities
in $\cat{A} \times \cat{B}$ were not given.  There is only one sensible way
to define them; write it down.
\end{question}

\begin{question}
There is a category $\Toph$%
\ntn{Toph}
whose objects are topological spaces and whose
maps $X \to Y$ are homotopy%
\index{homotopy}
classes of continuous maps from $X$ to $Y$.  What do you need to know about
homotopy in order to prove that $\Toph$ is a category?  What does it mean,
in purely topological terms, for two objects of $\Toph$ to be isomorphic?
\end{question}

\section{Functors}
\label{sec:ftrs}

One of the lessons of category theory is that whenever we meet a new type
of mathematical object, we should always ask whether there is a sensible
notion of `map' between such objects.  We can ask this about categories
themselves.  The answer is yes, and a map between categories is called a
functor.

\begin{defn}
Let $\cat{A}$ and $\cat{B}$ be categories.  A \demph{functor}%
\index{functor}
$F\from \cat{A} \to \cat{B}$ consists of:
\begin{itemize}
\item 
a function
\[
\ob(\cat{A}) \to \ob(\cat{B}),
\]
written as $A \mapsto F(A)$;

\item 
for each $A, A' \in \cat{A}$, a function
\[
\cat{A}(A, A') \to \cat{B}(F(A), F(A')),
\]
written as $f \mapsto F(f)$,
\end{itemize}
satisfying the following axioms:
\begin{itemize}
\item 
$F(f' \of f) = F(f') \of F(f)$ whenever $A \toby{f} A' \toby{f'} A''$ in
$\cat{A}$;

\item 
$F(1_A) = 1_{F(A)}$ whenever $A \in \cat{A}$.
\end{itemize}
\end{defn}

\begin{remarks} 
\label{rmks:defn-ftr}
\begin{enumerate}[(b)]       
\item 
\label{rmk:defn-ftr:loosely}
The definition of functor is set up so that from each string
\[
A_0 \toby{f_1} \ \cdots\ \toby{f_n} A_n
\]
of maps in $\cat{A}$ (with $n \geq 0$), it is possible to construct exactly
one%
\index{uniqueness!constructions@of constructions}
map
\[
F(A_0) \to F(A_n)
\]
in $\cat{B}$.  For example, given maps
\[
A_0 \toby{f_1} A_1 \toby{f_2} A_2 \toby{f_3} A_3 \toby{f_4} A_4
\]
in $\cat{A}$, we can construct maps
\[
\xymatrix@=10em{
F(A_0) 
\ar@<1ex>[r]^{F(f_4 f_3) F(f_2 f_1)}
\ar@<-1ex>[r]_{F(1_{A_4}) F(f_4) F(f_3 f_2) F(f_1)}     &
F(A_4)
}
\]
in $\cat{B}$, but the axioms imply that they are equal.  

\item   
\label{rmk:defn-ftr:comp}
We are familiar with the idea that structures and the structure-preserving
maps between them form a category (such as $\Grp$, $\Ring$, etc.).  In
particular, this applies to categories and functors: there is a category
$\CAT$%
\index{category!category of categories}%
\ntn{CAT}
whose objects are categories and whose maps are functors.

One part of this statement is that functors can be composed.%
\index{functor!composition of functors}
That is, given functors $\cat{A} \toby{F} \cat{B} \toby{G} \cat{C}$, there
arises a new functor $\cat{A} \toby{G \of F} \cat{C}$,%
\ntn{of-ftr}
defined in the obvious way.  Another is that for every category $\cat{A}$,
there is an identity%
\index{functor!identity}
functor $1_\cat{A}\from \cat{A} \to \cat{A}$.%
\ntn{id-ftr}
\end{enumerate}
\end{remarks}

\begin{examples} 
\label{egs:forgetful-functors}
Perhaps the easiest examples of functors are the so-called \demph{forgetful%
\index{functor!forgetful}
functors}.  (This is an informal term, with no precise definition.)  For
instance:
\begin{enumerate}[(b)]
\item   
\label{eg:forgetful-groups}
There is a functor $U\from \Grp \to \Set$ defined as follows: if $G$ is a
group then $U(G)$ is the underlying%
\index{underlying}
set of $G$ (that is, its set of elements), and if $f\from G \to H$ is a
group homomorphism then $U(f)$ is the function $f$ itself.  So $U$ forgets
the group structure of groups and forgets that group homomorphisms are
homomorphisms.

\item   
\label{eg:forgetful-ring-vs}
Similarly, there is a functor $\Ring \to \Set$ forgetting the ring
structure on rings, and (for any field $k$) there is a functor $\Vect_k
\to \Set$ forgetting the vector space structure on vector spaces.

\item   
\label{eg:forgetful-part}
Forgetful functors do not have to forget \emph{all} the structure.  For
example, let $\Ab$%
\ntn{Ab}
be the category of abelian groups.  There is a functor $\Ring \to \Ab$ that
forgets the multiplicative structure, remembering just the underlying
additive group.  Or, let $\Mon$%
\ntn{Mon}
be the category of monoids.  There is a functor $U\from \Ring \to \Mon$
that forgets the additive structure, remembering just the underlying%
\index{underlying}
multiplicative monoid.  (That is, if $R$ is a ring then $U(R)$ is the set
$R$ made into a monoid via $\cdot$ and $1$.)

\item   
\label{eg:forgetful-ab}
There is an inclusion functor $U\from \Ab \to \Grp$ defined by $U(A) = A$
for any abelian group $A$ and $U(f) = f$ for any homomorphism $f$ of abelian
groups.  It forgets that abelian groups are abelian.
\end{enumerate}

The forgetful functors in examples
\bref{eg:forgetful-groups}--\bref{eg:forgetful-part} forget
\emph{structure} on the objects, but that of
example~\bref{eg:forgetful-ab} forgets a \emph{property}.  Nevertheless,
it turns out to be convenient to use the same word, `forgetful', in both
situations.

Although forgetting is a trivial operation, there are situations in which
it is powerful.  For example, it is a theorem that the order of any finite
field is a prime power.  An important step in the proof is to simply forget
that the field is a field, remembering only that it is a vector space over
its subfield $\{0, 1, 1 + 1, 1 + 1 + 1, \ldots\}$.
\end{examples}

\begin{examples}        
\label{egs:free-functors}
\demph{Free%
\index{functor!free}%
\index{free functor}
functors} are in some sense dual to forgetful functors (as we will see in
the next chapter), although they are less elementary.  Again, `free
functor' is an informal but useful term.

\begin{enumerate}[(b)]
\item 
\label{eg:free-group} 
Given any set $S$, one can build the \demph{free%
\index{group!free}
group} $F(S)$ on $S$.  This is a group containing $S$ as a subset and with
no further properties other than those it is forced to have, in a sense
made precise in Section~\ref{sec:adj-basics}.  Intuitively, the group
$F(S)$ is obtained from the set $S$ by adding just enough new elements that
it becomes a group, but without imposing any equations other than those
forced by the definition of group.

A little more precisely, the elements of $F(S)$ are formal expressions or
\demph{words}%
\index{word}
such as $x^{-4} y x^2 z y^{-3}$ (where $x, y, z \in S$).  Two such words
are seen as equal if one can be obtained from the other by the usual
cancellation rules, so that, for example, $x^3 x y$, $x^4 y$, and $x^2
y^{-1} y x^2 y$ all represent the same element of $F(S)$.  To multiply two
words, just write one followed by the other; for instance, $x^{-4} y x$
times $x z y^{-3}$ is $x^{-4} y x^2 z y^{-3}$.

This construction assigns to each set $S$ a group $F(S)$.  In fact, $F$ is
a functor: any map of sets $f\from S \to S'$ gives rise to a homomorphism
of groups $F(f)\from F(S) \to F(S')$.  For instance, take the map of sets
\[
f\from \{w, x, y, z\} \to \{u, v\}
\]
defined by $f(w) = f(x) = f(y) = u$ and $f(z) = v$.  This gives rise to a
homomorphism 
\[
F(f)\from F(\{w, x, y, z\}) \to F(\{u, v\}),
\]
which maps $x^{-4} y x^2 z y^{-3} \in F(\{w, x, y, z\})$ to 
\[
u^{-4} u u^2 v u^{-3} 
=
u^{-1} v u^{-3}
\in
F(\{u, v\}).
\]

\item 
\label{eg:free-ring} 
Similarly, we can construct the free commutative ring $F(S)$ on a set $S$,
giving a functor $F$ from $\Set$ to the category $\CRing$%
\ntn{CRing}
of commutative rings.  In fact, $F(S)$ is something familiar, namely, the
ring of polynomials%
\index{ring!polynomial}
over $\integers$ in commuting variables $x_s$ ($s \in S$).  (A polynomial
is, after all, just a formal expression built from the variables using
the ring operations $+$, $-$ and $\cdot$.)  For example, if $S$ is a
two-element set then $F(S) \iso \integers[x, y]$.

\item   
\label{eg:free-vs}
We can also construct the free%
\index{vector space!free}
vector space on a set.  Fix a field $k$.  The free functor $F\from \Set \to
\Vect_k$ is defined on objects by taking $F(S)$ to be a vector space with
basis $S$.  Any two such vector spaces are isomorphic; but it is perhaps
not obvious that there is any such vector space at all, so we have to
construct one.  Loosely, $F(S)$ is the set of all formal $k$-linear
combinations of elements of $S$, that is, expressions
\[
\sum_{s \in S} \lambda_s s
\]
where each $\lambda_s$ is a scalar and there are only finitely many values
of $s$ such that $\lambda_s \neq 0$.  (This restriction is imposed because
one can only take \emph{finite} sums in a vector space.)  Elements of
$F(S)$ can be added:
\[
\sum_{s \in S} \lambda_s s + \sum_{s \in S} \mu_s s 
=
\sum_{s \in S} (\lambda_s + \mu_s) s.
\]
There is also a scalar multiplication on $F(S)$:
\[
c \cdot \sum_{s \in S} \lambda_s s
=
\sum_{s \in S} (c \lambda_s) s
\]
($c \in k$).  In this way, $F(S)$ becomes a vector space.

To be completely precise and avoid talking about `expressions', we
can define $F(S)$ to be the set of all functions $\lambda\from S \to k$ such
that $\{ s \in S \such \lambda(s) \neq 0\}$ is finite.  (Think of such a
function $\lambda$ as corresponding to the expression $\sum_{s \in S}
\lambda(s) s$.)  To define addition on $F(S)$, we must define for each
$\lambda, \mu \in F(S)$ a sum $\lambda + \mu \in F(S)$; it is given by
\[
(\lambda + \mu)(s) = \lambda(s) + \mu(s)
\]
($s \in S$). Similarly, the scalar multiplication is given by $(c \cdot
\lambda)(s) = c\cdot \lambda(s)$ ($c \in k$, $\lambda \in F(S)$, $s \in
S$).  
\end{enumerate}

Rings and vector spaces have the special property that it is relatively
easy to write down an explicit formula for the free functor.  The case of
groups is much more typical.  For most types of algebraic structure,
describing the free functor requires as much fussy work as it does for
groups.  We return to this point in Example~\ref{egs:adjns-alg} and
Example~\ref{eg:gaft-free-alg} (where we see how to avoid the fussy work
entirely).
\end{examples}

\begin{examples}[Functors in algebraic topology]
\label{egs:functors:homo-homo}
\index{algebraic topology}   
Historically, some of the first examples of functors arose in algebraic
topology.  There, the strategy is to learn about a space by extracting data
from it in some clever way, assembling that data into an algebraic
structure, then studying the algebraic structure instead of the original
space.  Algebraic topology therefore involves many functors from categories
of spaces to categories of algebras.
\begin{enumerate}[(b)]
\item 
Let $\Tp_*$%
\ntn{Top-star}
 be the category of topological spaces equipped with a
basepoint, together with the continuous basepoint-preserving maps.  There
is a functor $\pi_1\from \Tp_* \to \Grp$%
\ntn{pi-1}
assigning to each space $X$ with basepoint $x$ the fundamental%
\index{group!fundamental}
group $\pi_1(X, x)$ of $X$ at $x$.  (Some texts use the simpler notation
$\pi_1(X)$, ignoring the choice of basepoint.  This is more or less safe if
$X$ is path-connected, but strictly speaking, the basepoint should always
be specified.)

That $\pi_1$ is a functor means that it not only assigns to each
space-with-basepoint $(X, x)$ a group $\pi_1(X, x)$, but also assigns to
each basepoint-pre\-ser\-ving continuous map 
\[
f\from (X, x) \to (Y, y)
\]
a homomorphism 
\[
\pi_1(f)\from \pi_1(X, x) \to \pi_1(Y, y).
\]
Usually $\pi_1(f)$ is written as $f_*$.  The functoriality axioms say that $(g
\of f)_* = g_* \of f_*$ and $(1_{(X, x)})_* = 1_{\pi_1(X, x)}$.  

\item 
For each $n \in \nat$, there is a functor $H_n\from \Tp \to \Ab$%
\ntn{homology}
assigning to a space its $n$th homology%
\index{homology}
group (in any of several possible senses).
\end{enumerate}
\end{examples}

\begin{example}
Any system of polynomial%
\index{polynomial}%
\index{simultaneous equations}
equations such as
\begin{align}
2x^2 + y^2 - 3z^2       &
= 
1
\label{eq:scheme-1}     \\
x^3 + x     &
=
y^2       
\label{eq:scheme-2}
\end{align}
gives rise to a functor $\CRing \to \Set$.  Indeed, for each commutative
ring $A$, let $F(A)$ be the set of triples $(x, y, z) \in A \times A \times
A$ satisfying equations~\eqref{eq:scheme-1} and~\eqref{eq:scheme-2}.
Whenever $f\from A \to B$ is a ring homomorphism and $(x, y, z) \in F(A)$,
we have $(f(x), f(y), f(z)) \in F(B)$; so the map of rings $f\from A \to B$
induces a map of sets $F(f) \from F(A) \to F(B)$.  This defines a functor
$F\from \CRing \to \Set$.

In algebraic%
\index{algebraic geometry}
geometry, a \demph{scheme}%
\index{scheme}
is a functor $\CRing \to \Set$ with certain properties.  (This is not the
most common way of phrasing the definition, but it is equivalent.)  The
functor $F$ above is a simple example.
\end{example}

\begin{example}
\label{eg:ftrs-between-monoids}
Let $G$ and $H$ be monoids (or groups, if you prefer), regarded as
one-object%
\index{monoid!homomorphism of monoids}
categories $\cat{G}$ and $\cat{H}$.  A functor $F\from \cat{G} \to \cat{H}$
must send the unique object of $\cat{G}$ to the unique object of $\cat{H}$,
so it is determined by its effect on maps.  Hence, the functor $F \from
\cat{G} \to \cat{H}$ amounts to a function $F\from G \to H$ such that $F(g'
g) = F(g') F(g)$ for all $g', g \in G$, and $F(1) = 1$.  In other words, a
functor $\cat{G} \to \cat{H}$ is just a homomorphism $G \to H$.
\end{example}

\begin{example}
\label{eg:functor-action}
Let $G$ be a monoid,%
\index{monoid!action of}
regarded as a one-object category $\cat{G}$.  A functor $F\from \cat{G} \to
\Set$ consists of a set $S$ (the value of $F$ at the unique object of
$\cat{G}$) together with, for each $g \in G$, a function $F(g)\from S \to
S$, satisfying the functoriality axioms.  Writing $(F(g))(s) = g \cdot s$,
we see that the functor $F$ amounts to a set $S$ together with a function
\[
\begin{array}{ccc}
G \times S      &\to            &S      \\
(g, s)          &\mapsto        &g \cdot s      
\end{array}
\]
satisfying $(g' g) \cdot s = g' \cdot (g \cdot s)$ and $1 \cdot s = s$ for
all $g, g' \in G$ and $s \in S$.  In other words, a functor $\cat{G} \to
\Set$ is a set equipped with a left action by $G$: a \demph{left $G$-set},%
\index{G-set@$G$-set}
for short.

Similarly, a functor $\cat{G} \to \Vect_k$ is exactly a $k$-linear
representation%
\index{representation!group or monoid@of group or monoid!linear}
of $G$, in the sense of representation theory.  This can reasonably be
taken as the \emph{definition} of representation.
\end{example}

\begin{example}
\label{eg:functor-orders}
When $A$ and $B$ are (pre)ordered sets, a functor between the corresponding
categories is exactly an \demph{order-preserving%
\index{order-preserving}%
\index{map!order-preserving}
map}, that is, a function
$f\from A \to B$ such that $a \leq a' \implies f(a) \leq f(a')$.
Exercise~\ref{ex:functor-orders} asks you to verify this.
\end{example}

Sometimes we meet functor-like operations that reverse the arrows, with a
map $A \to A'$ in $\cat{A}$ giving rise to a map $F(A) \ot F(A')$ in
$\cat{B}$.  Such operations are called contravariant functors.

\begin{defn}    
\label{defn:contravariant}
Let $\cat{A}$ and $\cat{B}$ be categories.  A \demph{contravariant%
\index{functor!contravariant}%
\index{contravariant}
functor} from $\cat{A}$ to $\cat{B}$ is a functor $\cat{A}^\op \to
\cat{B}$.
\end{defn}

To avoid confusion, we write `a contravariant functor from $\cat{A}$ to
$\cat{B}$' rather than `a contravariant functor $\cat{A} \to \cat{B}$'.

Functors $\cat{C} \to \cat{D}$ correspond one-to-one with functors
$\cat{C}^\op \to \cat{D}^\op$, and $(\cat{A}^\op)^\op = \cat{A}$, so a
contravariant functor from $\cat{A}$ to $\cat{B}$ can also be described as
a functor $\cat{A} \to \cat{B}^\op$.  Which description we use is not
enormously important, but in the long run, the convention in
Definition~\ref{defn:contravariant} makes life easier.

An ordinary functor $\cat{A} \to \cat{B}$ is sometimes called a
\demph{covariant%
\index{functor!covariant}%
\index{covariant}
functor} from $\cat{A}$ to $\cat{B}$, for emphasis.

\begin{example}
\label{eg:contra-fn-spaces}
\index{topological space!functions on}%
\index{ring!functions@of functions}
We can tell a lot about a space by examining the functions on it.  The
importance of this principle in twentieth- and twenty-first-century
mathematics can hardly be exaggerated.

For example, given a topological space $X$, let $C(X)$%
\ntn{cts-ring}
be the ring of continuous real-valued functions on $X$.  The ring
operations are defined `pointwise':%
\index{pointwise}
for instance, if $p_1, p_2\from X \to \reals$ are continuous maps then the
map $p_1 + p_2 \from X \to \reals$ is defined by
\[
(p_1 + p_2)(x) = p_1(x) + p_2(x)
\]
($x \in X$).  A continuous map $f\from X \to Y$ induces a ring homomorphism
$C(f)\from C(Y) \to C(X)$, defined at $q \in C(Y)$ by taking $(C(f))(q)$ to be
the composite map
\[
X \toby{f} Y \toby{q} \reals.
\]
Note that $C(f)$ goes in the opposite direction from $f$.  After checking
some axioms (Exercise~\ref{ex:contra-fn-spaces}), we conclude that $C$ is a
contravariant functor from $\Tp$ to $\Ring$.

While this particular example will not play a large part in this text, it
is worth close attention.  It illustrates the important idea of a
structure whose elements are maps (in this case, a ring whose elements are
continuous functions).  The way in which $C$ becomes a functor, via
composition, is also important.  Similar constructions will be crucial in
later chapters.

For certain classes of space, the passage from $X$ to $C(X)$ loses no
information: there is a way of reconstructing the space $X$ from the ring
$C(X)$.  For this and related reasons, it is sometimes said that `algebra
is dual%
\index{duality!algebra--geometry}
to geometry'.
\end{example}

\begin{example}
\label{eg:fns-on-vs}
Let $k$ be a field.  For any two vector spaces $V$ and $W$ over $k$, there
is a vector space
\[
\HOM(V, W) 
=
\{ \text{linear maps } V \to W \}.%
\index{vector space!linear maps@of linear maps}
\ntn{HOM}
\]
The elements of this vector space are themselves maps, and the vector space
operations (addition and scalar multiplication) are defined pointwise, as
in the last example.

Now fix a vector space $W$.  Any linear map $f\from V \to V'$ induces a linear
map 
\[
f^*\from \HOM(V', W) \to \HOM(V, W),%
\ntn{dual-map}
\]
defined at $q\in \HOM(V', W)$ by taking $f^*(q)$ to be the composite map
\[
V \toby{f} V' \toby{q} W.
\]
This defines a functor
\[
\HOM(\dashbk, W)\from \Vect_k^\op \to \Vect_k.
\]
The symbol `$\dashbk$'%
\ntn{dashbk}
is a blank or placeholder, into which arguments can be inserted.  Thus, the
value of $\HOM(\dashbk, W)$ at $V$ is $\HOM(V, W)$.  Sometimes we use a
blank space%
\ntn{empty-blank}
instead of $\dashbk$, as in $\HOM(\hspace*{1em}, W)$.

An important special case is where $W$ is $k$, seen as a one-dimensional
vector space over itself.  The vector space $\HOM(V, k)$ is called the
\demph{dual}%
\index{vector space!functions on}%
\index{vector space!dual}%
\index{duality!vector spaces@for vector spaces}
of $V$, and is written as $V^*$.%
\ntn{dual-vs}
So there is a contravariant functor
\[
\blank^* = \HOM(\dashbk, k)\from \Vect_k^\op \to \Vect_k
\]
sending each vector space to its dual.
\end{example}

\begin{example}
For each $n\in\nat$, there is a functor $H^n\from  \Tp^\op
\to \Ab$%
\ntn{cohomology}
assigning to a space its $n$th cohomology%
\index{cohomology}
group.
\end{example}

\begin{example}
\label{eg:contra-functors:actions}
Let $G$ be a monoid, regarded as a one-object category $\cat{G}$.  A
functor $\cat{G}^\op \to \Set$ is a \emph{right} $G$-set,%
\index{monoid!action of}
for essentially the same reasons as in Example~\ref{eg:functor-action}.

That left actions are covariant functors and right actions are
contravariant functors is a consequence of a basic notational choice: we
write the value of a function $f$ at an element $x$ as $f(x)$, not $(x)f$.
\end{example}

Contravariant functors whose codomain is $\Set$ are important enough to have
their own special name.

\begin{defn}    
\label{defn:presheaf}
Let $\cat{A}$ be a category.  A \demph{presheaf}%
\index{presheaf|(}
on $\cat{A}$ is a functor $\cat{A}^\op \to \Set$.
\end{defn}

The name comes from the following special case.  Let $X$ be a topological
space.  Write $\oset(X)$%
\ntn{oset}
for the poset of open subsets of $X$, ordered by inclusion.  View
$\oset(X)$ as a category, as in
Example~\ref{egs:cats-as}\bref{eg:cats-as:orders}.  Thus, the objects of
$\oset(X)$ are the open subsets of $X$, and for $U, U' \in \oset(X)$, there
is one map $U \to U'$ if $U \sub U'$, and there are none otherwise.  A
\demph{presheaf} on the space $X$ is a presheaf on the category $\oset(X)$.
For example, given any space $X$, there is a presheaf $F$ on $X$ defined by
\[
F(U) = \{ \text{continuous functions } U \to \reals \}
\index{topological space!functions on}
\]
($U \in \oset(X)$) and, whenever $U \sub U'$ are open subsets of $X$, by
taking the map $F(U') \to F(U)$ to be restriction.  Presheaves, and a
certain class of presheaves called sheaves,%
\index{sheaf}
play an important role in modern geometry.%
\index{presheaf|)}

\subjectchange

We know very well that for \emph{functions} between \emph{sets}, it is
sometimes useful to consider special kinds of function such as injections,
surjections and bijections.  We also know that the notions of injection and
subset are related: for instance, whenever $B$ is a subset of $A$, there is
an injection $B \to A$ given by inclusion.  In this section and the next,
we introduce some similar notions for \emph{functors} between
\emph{categories}, beginning with the following definitions.

\begin{defn}
A functor $F\from  \cat{A} \to \cat{B}$ is \demph{faithful}%
\index{functor!faithful}%
\index{faithful}
(respectively, \demph{full})%
\index{functor!full}
if for each $A, A' \in \cat{A}$, the function
\[
\begin{array}{ccc}
\cat{A}(A, A')  &\to            &\cat{B}(F(A), F(A'))   \\
f               &\mapsto        &F(f)
\end{array}
\]
is injective (respectively, surjective).
\end{defn}

\begin{warning}
Note the roles of $A$ and $A'$ in the definition.  Faithfulness does
\emph{not} say that if $f_1$ and $f_2$ are distinct maps in $\cat{A}$ then
$F(f_1) \neq F(f_2)$ (Exercise~\ref{ex:faithful-not-inj}).
\begin{figure}
\[
\cat{A}
\begin{array}{c}
\fbox{ 
\xymatrix{
\ &A\vphantom{F(A)} \ar@{.>}[d]&\ \\ &A'\vphantom{F(A')}
}
}
\end{array}
\qquad%
\toby{F}%
\qquad%
\begin{array}{c}
\fbox{ 
\xymatrix{
\ &F(A) \ar[d]^g&\  \\ & F(A')
}
}
\end{array}
\cat{B}
\]
\caption{Fullness and faithfulness.}
\label{fig:ff}
\end{figure}
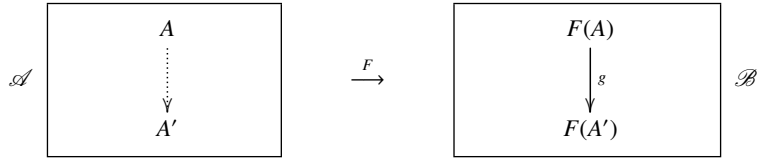
In the situation of Figure~\ref{fig:ff}, $F$ is faithful if for each $A$,
$A'$ and $g$ as shown, there is at most one dotted arrow that $F$ sends to
$g$.  It is full if for each such $A$, $A'$ and $g$, there is at least one
dotted arrow that $F$ sends to $g$.
\end{warning}

\begin{defn}
Let $\cat{A}$ be a category.  A \demph{subcategory}%
\index{subcategory!full}
$\cat{S}$ of $\cat{A}$ consists of a subclass $\ob(\cat{S})$ of
$\ob(\cat{A})$ together with, for each $S, S' \in \ob(\cat{S})$, a subclass
$\cat{S}(S, S')$ of $\cat{A}(S, S')$, such that $\cat{S}$ is closed under
composition and identities.  It is a \demph{full}%
\index{subcategory!full}
subcategory if $\cat{S}(S, S') = \cat{A}(S, S')$ for all $S, S'
\in \ob(\cat{S})$.
\end{defn}

A full subcategory therefore consists of a selection of the objects, with
all of the maps between them.  So, a full subcategory can be specified
simply by saying what its objects are.  For example, $\Ab$ is the full
subcategory of $\Grp$ consisting of the groups that are abelian.

Whenever $\cat{S}$ is a subcategory of a category $\cat{A}$, there is an
inclusion functor $I: \cat{S} \to \cat{A}$ defined by $I(S) = S$ and $I(f)
= f$.  It is automatically faithful, and it is full if and only if
$\cat{S}$ is a full subcategory.

\begin{warning}
The image%
\index{functor!image of}%
\index{image!functor@of functor}
of a functor need not be a subcategory.  For example, consider the functor
\[
\Bigl(
\begin{array}{c}
\xymatrix@1{A \ar[r]^-f&B\quad B' \ar[r]^-g&C}
\end{array}
\Bigr)
\qquad
\toby{F}
\qquad
\left(
\begin{array}{c}
\xymatrix{
        &Y \ar[dr]^q    &       \\
X \ar[ur]^p \ar[rr]_{qp} &&Z
}
\end{array}
\right)
\]
defined by $F(A) = X$, $F(B) = F(B') = Y$, $F(C) = Z$, $F(f) = p$, and $F(g) =
q$.  Then $p$ and $q$ are in the image of $F$, but $qp$ is not.
\end{warning}

\exs

\begin{question}
Find three examples of functors not mentioned above.
\end{question}

\begin{question}        
\label{ex:ftrs-pres-iso}
Show that functors preserve isomorphism.%
\index{isomorphism!preserved by functors}
That is, prove that if $F\from \cat{A} \to \cat{B}$ is a functor and $A, A'
\in \cat{A}$ with $A \iso A'$, then $F(A) \iso F(A')$.
\end{question}

\begin{question}
\label{ex:functor-orders}
Prove the assertion made in Example~\ref{eg:functor-orders}.  In other
words, given ordered sets $A$ and $B$, and denoting by $\cat{A}$ and
$\cat{B}$ the corresponding categories, show that a functor $\cat{A} \to
\cat{B}$ amounts to an order-preserving%
\index{order-preserving}%
\index{map!order-preserving}
map $A \to B$.
\end{question}

\begin{question}
Two categories $\cat{A}$ and $\cat{B}$ are \demph{isomorphic},%
\index{isomorphism!categories@of categories}%
\index{category!isomorphism of categories}
written as $\cat{A} \iso \cat{B}$,%
\ntn{iso-cat}
if they are isomorphic as objects of $\CAT$.  
\begin{enumerate}[(b)]
\item 
Let $G$ be a group, regarded as a one-object%
\index{group!opposite}
category all of whose maps are isomorphisms.  Then its opposite $G^\op$ is
also a one-object category all of whose maps are isomorphisms, and can
therefore be regarded as a group too.  What is $G^\op$, in purely
group-theoretic terms?  Prove that $G$ is isomorphic to $G^\op$.

\item 
Find a monoid%
\index{monoid!opposite}
not isomorphic to its opposite.
\end{enumerate}
\end{question}

\begin{question}
Is there a functor $Z \from \Grp \to \Grp$ with the property that $Z(G)$ is
the centre%
\index{centre}
of $G$ for all groups $G$?
\end{question}

\begin{question}        
\label{ex:ftr-on-product}
Sometimes we meet functors whose domain is a product%
\index{category!product of categories}
$\cat{A} \times \cat{B}$ of categories.  Here you will show that such a
functor can be regarded as an interlocking pair of families of functors,
one defined on $\cat{A}$ and the other defined on $\cat{B}$.  (This is very
like the situation for bilinear and linear maps.)
\begin{enumerate}[(b)]
\item   
\label{part:prod-ftr-compts} 
Let $F\from \cat{A} \times \cat{B} \to \cat{C}$ be a functor.  
Prove that for each $A \in \cat{A}$, there is a functor $F^A\from
\cat{B} \to \cat{C}$ defined on objects $B \in \cat{B}$ by $F^A(B) = F(A,
B)$ and on maps $g$ in $\cat{B}$ by $F^A(g) = F(1_A, g)$.  Prove that for
each $B \in \cat{B}$, there is a functor $F_B \from \cat{A} \to \cat{C}$
defined similarly.  

\item   
\label{part:prod-ftr-condns}
Let $F \from \cat{A} \times \cat{B} \to \cat{C}$ be a functor.  With
notation as in~\bref{part:prod-ftr-compts}, show that
the families of functors $(F^A)_{A \in \cat{A}}$ and $(F_B)_{B \in
  \cat{B}}$ satisfy the following two conditions:
\begin{itemize}
\item 
if $A \in \cat{A}$ and $B \in \cat{B}$ then $F^A(B) = F_B(A)$;

\item 
if $f\from A \to A'$ in $\cat{A}$ and $g\from B \to B'$ in $\cat{B}$ then
$F^{A'}(g) \of F_B(f) = F_{B'}(f) \of F^A(g)$.
\end{itemize}

\item
Now take categories $\cat{A}$, $\cat{B}$ and $\cat{C}$, and take families
of functors $(F^A)_{A \in \cat{A}}$ and $(F_B)_{B \in \cat{B}}$ satisfying
the two conditions in~\bref{part:prod-ftr-condns}.  Prove that there is a
unique functor $F\from \cat{A} \times \cat{B} \to \cat{C}$ satisfying the
equations in~\bref{part:prod-ftr-compts}.  (`There is a unique functor'
means in particular that there \emph{is} a functor, so you have to prove
existence as well as uniqueness.)
\end{enumerate}
\end{question}

\begin{question}
\label{ex:contra-fn-spaces}
Fill in the details of Example~\ref{eg:contra-fn-spaces}, thus constructing
a functor $C\from \Tp^\op \to \Ring$.  
\end{question}

\begin{question}
\label{ex:faithful-not-inj}
Find an example of a functor $F\from \cat{A} \to \cat{B}$ such that $F$ is
faithful%
\index{functor!faithful}%
\index{faithful}
but there exist distinct maps $f_1$ and $f_2$ in $\cat{A}$ with $F(f_1) =
F(f_2)$.
\end{question}

\begin{question}
\begin{enumerate}[(b)]
\item
Of the examples of functors appearing in this section, which are faithful
and which are full?

\item 
Write down one example of a functor that is both full and faithful, one
that is full but not faithful, one that is faithful but not full, and one
that is neither.
\end{enumerate}
\end{question}

\begin{question}
\begin{enumerate}[(b)]
\item
What are the subcategories of an ordered set?  Which are full?

\item
What are the subcategories of a group?  (Careful!)  Which are full?
\end{enumerate}
\end{question}

\section{Natural transformations}
\label{sec:nts}

We now know about categories.  We also know about functors, which
are maps between categories.  Perhaps surprisingly, there is 
a further notion of `map between functors'.  Such maps are called natural
transformations.  This notion only applies when the functors have
the same domain and codomain: 
\[
\parpair{\cat{A}}{\cat{B}}{F}{G}.
\]

To see how this might work, let us consider a special case.  Let $\cat{A}$
be the discrete category
(Example~\ref{egs:cats-as}\bref{eg:cats-as:discrete}) whose objects are
the natural numbers $0, 1, 2$, \ldots.  A functor $F$ from $\cat{A}$ to
another category $\cat{B}$ is simply a sequence $(F_0, F_1, F_2, \ldots)$
of objects of $\cat{B}$.  Let $G$ be another functor from $\cat{A}$ to
$\cat{B}$, consisting of another sequence $(G_0, G_1, G_2, \ldots)$ of
objects of $\cat{B}$.  It would be reasonable to define a `map' from $F$ to
$G$ to be a sequence
\[
\Bigl(
F_0 \toby{\alpha_0} G_0, \ 
F_1 \toby{\alpha_1} G_1, \ 
F_1 \toby{\alpha_2} G_2, 
\ \ldots\ 
\Bigr)
\]
of maps in $\cat{B}$.  The situation can be depicted as follows:
\[
\cat{A} 
\begin{array}{c}
\fbox{
\xymatrix@C=2ex{0 & 1 & 2 & \cdots}
}
\end{array}
\qquad
\qquad
\begin{array}{c}
\fbox{
\xymatrix@C=2ex{
&
F_0 \ar[d]^{\alpha_0} &
F_1 \ar[d]^{\alpha_1} &
F_2 \ar[d]^{\alpha_2} &
{ } \ar@{}[d]^{\textstyle\cdots} &
&
\\
&
G_0 &
G_1 &
G_2 &
}
}
\end{array}
\cat{B} 
\]
(The right-hand diagram should not be understood too literally.  Some of
the objects $F_i$ or $G_i$ might be equal, and there might be much else in
$\cat{B}$ besides what is shown.)

This suggests that in the general case, a natural transformation between
functors $\parpairi{\cat{A}}{\cat{B}}{F}{G}$ should consist of maps
$\alpha_A\from F(A) \to G(A)$, one for each $A \in \cat{A}$.  In the
example above, the category $\cat{A}$ had the special property of not
containing any nontrivial maps.  In general, we demand some kind of
compatibility between the maps in $\cat{A}$ and the maps $\alpha_A$.

\begin{defn}
Let $\cat{A}$ and $\cat{B}$ be categories and let
$\parpairi{\cat{A}}{\cat{B}}{F}{G}$ be functors.  A \demph{natural%
\index{natural transformation}
transformation} $\alpha\from F \to G$ is a family $\Bigl( F(A)
\toby{\alpha_A} G(A) \Bigr)_{A \in \cat{A}}$%
\ntn{nt-comp}
of maps in $\cat{B}$ such that
for every map $A \toby{f} A'$ in $\cat{A}$, the square
\begin{equation}        
\label{eq:nat}
\begin{array}{c}
\xymatrix{
F(A) \ar[r]^{F(f)} \ar[d]_{\alpha_A}    &
F(A') \ar[d]^{\alpha_{A'}}      \\
G(A) \ar[r]_{G(f)}      &
G(A')
} 
\end{array}
\end{equation}
commutes.  The maps $\alpha_A$ are called the \demph{components}%
\index{component!natural transformation@of natural transformation}
of $\alpha$.
\end{defn}

\begin{remarks} 
\label{rmks:defn-nt}
\begin{enumerate}[(b)]
\item   
\label{rmk:defn-nt:loosely}
The definition of natural transformation is set up so that from each map $A
\toby{f} A'$ in $\cat{A}$, it is possible to construct exactly one%
\index{uniqueness!constructions@of constructions}
map $F(A) \to G(A')$ in $\cat{B}$.  When $f = 1_A$, this map is $\alpha_A$.
For a general $f$, it is the diagonal of the square~\eqref{eq:nat}, and
`exactly one' implies that the square commutes.

\item 
We write
\[
\xymatrix{
\cat{A} \rtwocell^F_G{\alpha} &\cat{B}
}%
\ntn{double-arrow}
\]
to mean that $\alpha$ is a natural transformation from $F$ to $G$.
\end{enumerate}
\end{remarks}

\begin{example}
Let $\cat{A}$ be a discrete%
\index{category!discrete!functor out of}
category, and let $F, G\from \cat{A} \to \cat{B}$ be functors.  Then $F$
and $G$ are just families $(F(A))_{A \in \cat{A}}$ and $(G(A))_{A \in
  \cat{A}}$ of objects of $\cat{B}$.  A natural transformation $\alpha\from
F \to G$ is just a family $\Bigl( F(A) \toby{\alpha_A} G(A) \Bigr)_{A \in
  \cat{A}}$ of maps in $\cat{B}$, as claimed above in the case $\ob \cat{A}
= \nat$.  In principle, this family must satisfy the naturality
axiom~\eqref{eq:nat} for every map $f$ in $\cat{A}$; but the only maps in
$\cat{A}$ are the identities, and when $f$ is an identity, this axiom holds
automatically.
\end{example}

\begin{example}
Recall from Examples~\ref{egs:cats-as} that a group (or more generally, a
monoid) $G$ can be regarded as a one-object%
\index{monoid!action of}
category.  Also recall from Example~\ref{eg:functor-action} that a functor
from the category $G$ to $\Set$ is nothing but a left $G$-set.  (Previously
we used $\cat{G}$ to denote the category corresponding to the group $G$;
from now on we use $G$%
\index{monoid!one-object category@as one-object category}
\ntn{group-one-obj}
to denote them both.)  Take two $G$-sets, $S$ and $T$.  Since $S$ and $T$
can be regarded as functors $G \to \Set$, we can ask: what is a natural
transformation
\[
\xymatrix{
G \rtwocell^S_T{\alpha} &\Set,
}
\]
in concrete terms?

Such a natural transformation consists of a single map in $\Set$ (since $G$
has just one object), satisfying some axioms.  Precisely, it is a function
$\alpha\from S \to T$ such that $\alpha(g\cdot s) = g\cdot \alpha(s)$ for
all $s \in S$ and $g \in G$.  (Why?)  In other words, it is just a map of
$G$-sets, sometimes called a \demph{$G$-equivariant}%
\index{equivariant}
map.
\end{example}

\begin{example}
Fix a natural number $n$.  In this example, we will see how `determinant%
\index{determinant}
of an $n \times n$ matrix' can be understood as a natural transformation.

For any commutative ring $R$, the $n \times n$ matrices with entries in $R$
form a monoid $M_n(R)$ under multiplication.  Moreover, any ring
homomorphism $R \to S$ induces a monoid homomorphism $M_n(R) \to M_n(S)$.
This defines a functor $M_n\from \CRing \to \Mon$ from the category of
commutative rings to the category of monoids.

Also, the elements of any ring $R$ form a monoid $U(R)$ under multiplication,
giving another functor $U\from \CRing \to \Mon$.

Now, every $n \times n$ matrix $X$ over a commutative ring $R$
has a determinant ${\det}_R(X)$, which is an element of $R$.
Familiar properties of determinant~--
\[
{\det}_R(XY) = {\det}_R(X) {\det}_R(Y),
\qquad
{\det}_R(I) = 1
\]
-- tell us that for each $R$, the function ${\det}_R\from M_n(R) \to U(R)$
is a monoid homomorphism.  So, we have a family of maps
\[
\Bigl( M_n(R) \toby{\det_R} U(R) \Bigr)_{R \in \CRing},
\]
and it makes sense to ask whether they define a natural transformation
\[
\xymatrix{
\CRing \rrtwocell^{M_n}_U{\hspace{.8em}\det} &&\Mon.
}
\]
Indeed, they do.  That the naturality squares commute (check!)\ reflects
the fact that determinant is defined in the same way for all rings.  We do
not use one definition of determinant for one ring and a different
definition for another ring.  Generally speaking, the naturality
axiom~\eqref{eq:nat} is supposed to capture the idea that the family
$(\alpha_A)_{A \in \cat{A}}$ is defined in a uniform way across all $A \in
\cat{A}$.
\end{example}

\begin{constn}  
\label{constn:comp-nat}
Natural transformations are a kind of map, so we would expect to be able to
compose%
\index{natural transformation!composition of}
them.  We can.  Given natural transformations
\vspace{-2ex}
\[
\xymatrix@C+1.5em{
\cat{A} 
\ruppertwocell<8>^F{\alpha}
\ar[r]|G 
\rlowertwocell<-8>_H{\beta}
&
\cat{B}
},
\]
\vspace{-3ex}\\
there is a composite natural transformation
\[
\xymatrix@C+1.5em{
\cat{A} \rtwocell^F_H{\hspace{.8em}\beta\of\alpha} &\cat{B}
}%
\ntn{of-nt}
\]
defined by $(\beta\of\alpha)_A = \beta_A \of \alpha_A$ for all $A \in
\cat{A}$.  There is also an identity%
\index{natural transformation!identity}
natural transformation
\[
\xymatrix@C+1em{
\cat{A} \rtwocell^F_F{\hspace{.5em}1_F} &\cat{B}
}%
\ntn{id-nt}
\]
on any functor $F$, defined by $(1_F)_A = 1_{F(A)}$.  So for any two
categories $\cat{A}$ and $\cat{B}$, there is a category whose objects are
the functors from $\cat{A}$ to $\cat{B}$ and whose maps are the natural
transformations between them.  This is called the \demph{functor%
\index{functor!category}
category} from $\cat{A}$ to $\cat{B}$, and written as
$\ftrcat{\cat{A}}{\cat{B}}$%
\ntn{ftr-cat-bkts}
 or
$\cat{B}^\cat{A}$.%
\ntn{ftr-cat-power}
\end{constn}

\begin{example}        
Let $2$%
\ntn{two-disc-cat}
be the discrete%
\index{category!discrete!functor out of}
category with two objects.  A functor from $2$ to a category $\cat{B}$ is a
pair of objects of $\cat{B}$, and a natural transformation is a pair of
maps.  The functor category $\ftrcat{2}{\cat{B}}$ is therefore isomorphic
to the product category $\cat{B} \times \cat{B}$
(Construction~\ref{constn:prod-cat}).  This fits well with the alternative
notation $\cat{B}^2$ for the functor category.
\end{example}

\begin{example} 
Let $G$ be a monoid.%
\index{monoid!action of}
Then $\ftrcat{G}{\Set}$%
\ntn{GSet}
is the category of left $G$-sets, and $\ftrcat{G^\op}{\Set}$ is the
category of right $G$-sets (Example~\ref{eg:contra-functors:actions}).
\end{example}

\begin{example}
\label{eg:ftr-cats-orders}
Take ordered%
\index{ordered set}
sets $A$ and $B$, viewed as categories (as in
Example~\ref{egs:cats-as}\bref{eg:cats-as:orders}).  Given order-preserving
maps $\parpairi{A}{B}{f}{g}$, viewed as functors (as in
Example~\ref{eg:functor-orders}), there is at most one natural
transformation
\[
\xymatrix{
A \rtwocell^f_g &B,
}
\]
and there is one if and only if $f(a) \leq g(a)$ for all $a \in A$.  (The
naturality axiom~\eqref{eq:nat} holds automatically, because in an ordered
set, all diagrams commute.)  So $\ftrcat{A}{B}$ is an ordered set too; its
elements are the order-preserving maps from $A$ to $B$, and $f \leq g$ if
and only if $f(a) \leq g(a)$ for all $a \in A$.
\end{example}

Everyday phrases such as `\emph{the}%
\index{uniqueness}
cyclic group of order $6$' and `\emph{the} product of two spaces' reflect
the fact that given two isomorphic objects of a category, we usually
neither know nor care%
\label{p:care}
whether they are actually equal.  This is enormously important.  

In particular, the lesson applies when the category concerned is a functor
category.  In other words, given two functors $F, G\from \cat{A} \to
\cat{B}$, we usually do not care whether they are literally equal.
(Equality would imply that the objects $F(A)$ and $G(A)$ of $\cat{B}$ were
equal for all $A \in \cat{A}$, a level of detail in which we have just
declared ourselves to be uninterested.)  What really matters is whether
they are naturally isomorphic.

\begin{defn}    
\label{defn:nat-iso}
Let $\cat{A}$ and $\cat{B}$ be categories.  A \demph{natural%
\index{isomorphism!natural}
isomorphism} between functors from $\cat{A}$ to $\cat{B}$ is an isomorphism
in $\ftrcat{\cat{A}}{\cat{B}}$.
\end{defn}

An equivalent form of the definition is often useful:

\begin{lemma}   
\label{lemma:nat-iso-compts}
Let $\xymatrix@1{\cat{A}\rtwocell^F_G{\alpha} &\cat{B}}$ be a natural
transformation.  Then $\alpha$ is a natural isomorphism if and only if
$\alpha_A\from F(A) \to G(A)$ is an isomorphism for all $A \in \cat{A}$.
\end{lemma}

\begin{pf}
Exercise~\ref{ex:nat-iso-compts}.
\end{pf}

Of course, we say that functors $F$ and $G$ are \demph{naturally
isomorphic} if there exists a natural isomorphism from $F$ to $G$.  Since
natural isomorphism is just isomorphism in a particular category (namely,
$\ftrcat{\cat{A}}{\cat{B}}$), we already have notation for this: $F \iso
G$.%
\ntn{iso-ftr}

\begin{defn}	
\label{defn:nat-in}
Given functors $\parpair{\cat{A}}{\cat{B}}{F}{G}$, we say that
\[
F(A) \iso G(A) \text{ \demph{naturally in} } A 
\index{naturally}
\]
if $F$ and $G$ are naturally isomorphic.  
\end{defn}

This alternative terminology can be understood as follows.  If $F(A) \iso
G(A)$ naturally in $A$ then certainly $F(A) \iso G(A)$ for each individual
$A$, but more is true: we can choose isomorphisms $\alpha_A\from F(A) \to
G(A)$ in such a way that the naturality axiom~\eqref{eq:nat} is satisfied.

\begin{example}
Let $F, G\from \cat{A} \to \cat{B}$ be functors from a discrete%
\index{category!discrete!functor out of}
category $\cat{A}$ to a category $\cat{B}$.  Then $F \iso G$ if and only if
$F(A) \iso G(A)$ for all $A \in \cat{A}$.

So in \emph{this} case, $F(A) \iso G(A)$ naturally in $A$ if and only if
$F(A) \iso G(A)$ for all $A$.  But this is only true because $\cat{A}$ is
discrete.  In general, it is emphatically false.  There are many examples
of categories and functors $\parpairi{\cat{A}}{\cat{B}}{F}{G}$ such that
$F(A) \iso G(A)$ for all $A \in \cat{A}$, but not \emph{naturally} in $A$.
Exercise~\ref{ex:species} gives an example from combinatorics.
\end{example}

\begin{example}
Let $\FDVect$%
\ntn{FDVect}
be the category of finite-dimensional vector spaces over some field $k$.
The dual%
\index{vector space!dual}%
\index{duality!vector spaces@for vector spaces}
vector space construction defines a contravariant functor from $\FDVect$ to
itself (Example~\ref{eg:fns-on-vs}), and the double dual construction
therefore defines a covariant functor from $\FDVect$ to itself.

Moreover, we have for each $V \in \FDVect$ a canonical isomorphism
$\alpha_V\from V \to V^{**}$.  Given $v \in V$, the element $\alpha_V(v)$
of $V^{**}$ is `evaluation%
\index{evaluation}
at $v$'; that is, $\alpha_V(v)\from V^* \to k$ maps $\phi \in V^*$ to
$\phi(v) \in k$.  That $\alpha_V$ is an isomorphism is a standard result in
the theory of finite-dimensional vector spaces.

This defines a natural transformation
\[
\xymatrix{
\FDVect 
\rtwocell<5>^{1_\FDVect}_{\blank^{**}}{\alpha} &
\FDVect
}
\]
from the identity functor to the double dual functor.  By
Lemma~\ref{lemma:nat-iso-compts}, $\alpha$ is a natural isomorphism.  So
$1_\FDVect \iso \blank^{**}$.  Equivalently, in the language of
Definition~\ref{defn:nat-in}, $V \iso V^{**}$ naturally in $V$.

This is one of those occasions on which category theory makes an intuition
precise.  In some informal sense, evident before you learn anything about
category theory, the isomorphism between a finite-dimensional vector space
and its double dual is `natural' or `canonical': no arbitrary choices are
needed in order to define it.  In contrast, to specify an isomorphism
between $V$ and its single dual $V^*$, we need to make an arbitrary choice
of basis, and the isomorphism really does depend on the basis that we
choose.
\end{example}

In the example on vector spaces, the word \demph{canonical}%
\index{canonical}
was used.  It is an informal word, meaning something like `God-given' or
`defined without making arbitrary choices'.  For example, for any two sets
$A$ and $B$, there is a canonical bijection $A \times B \to B \times A$
defined by $(a, b) \mapsto (b, a)$, and there is a canonical function $A
\times B \to A$ defined by $(a, b) \mapsto a$.  But the function $B \to A$
defined by `choose an element $a_0 \in A$ and send everything to $a_0$' is
not canonical, because the choice of $a_0$ is arbitrary.

\subjectchange

The concept of natural isomorphism leads unavoidably to another central
concept: equivalence of categories.

\index{sameness|(}
Two elements of a set are either equal or not.  Two objects of a category
can be equal, not equal but isomorphic, or not even isomorphic.  As
explained before Definition~\ref{defn:nat-iso}, the notion of equality
between two objects of a category is unreasonably strict; it is usually
isomorphism that we care about.  So:
\begin{itemize}
\item
the right notion of sameness of two elements of a set is
equality;

\item
the right notion of sameness of two objects of a category is
isomorphism.
\end{itemize}
When applied to a functor category $\ftrcat{\cat{A}}{\cat{B}}$, the second
point tells us that:
\begin{itemize}
\item
the right notion of sameness of two functors $\cat{A} \parpairu
\cat{B}$ is natural isomorphism. 
\end{itemize}
But what is the right notion of sameness of two \emph{categories}?
Isomorphism is unreasonably strict, as if $\cat{A} \iso \cat{B}$ then
there are functors
\begin{equation}        
\label{eq:oppair}
\oppair{\cat{A}}{\cat{B}}{F}{G}
\end{equation}
such that 
\begin{equation}
\label{eq:inverse-functors}
G \of F = 1_\cat{A}
\qquad
\text{and}
\qquad
F \of G = 1_\cat{B}, 
\end{equation}
and we have just seen that the notion of equality between functors is too
strict.  The most useful notion of sameness of categories, called
`equivalence', is looser than isomorphism.  To obtain the definition, we
simply replace the unreasonably strict equalities
in~\eqref{eq:inverse-functors} by isomorphisms.  This gives
\[
G \of F \iso 1_\cat{A}
\qquad
\text{and}
\qquad
F \of G \iso 1_\cat{B}.
\index{sameness|)}
\]

\begin{defn}    
\label{defn:eqv}
An \demph{equivalence}%
\index{equivalence of categories}%
\index{category!equivalence of categories}
between categories $\cat{A}$ and $\cat{B}$ consists
of a pair~\eqref{eq:oppair} of functors together with natural isomorphisms
\[
\eta\from 1_\cat{A} \to G \of F,
\qquad
\epsln\from F \of G \to 1_\cat{B}.
\]
If there exists an equivalence between $\cat{A}$ and $\cat{B}$, we say that
$\cat{A}$ and $\cat{B}$ are \demph{equivalent}, and write $\cat{A} \eqv
\cat{B}$.%
\ntn{eqv}
We also say that the functors $F$ and $G$ are \demph{equivalences}.
\end{defn}

The directions of $\eta$ and $\epsln$ are not very important, since they
are isomorphisms anyway.  The reason for this particular choice will become
apparent when we come to discuss adjunctions (Section~\ref{sec:adj-units}).

\begin{warning}
The symbol $\iso$ is used for isomorphism of objects of a category, and in
particular for isomorphism of categories (which are objects of $\CAT$).
The symbol $\eqv$ is used for equivalence of categories.  At least, this is
the convention used in this book and by most category theorists, although
it is far from universal in mathematics at large.
\end{warning}

There is a very useful alternative characterization of those functors that
are equivalences.  First, we need a definition.

\begin{defn} 
A functor $F\from \cat{A} \to \cat{B}$ is \demph{essentially%
\index{essentially surjective on objects}%
\index{functor!essentially surjective on objects}
surjective on objects} if for all $B \in \cat{B}$, there exists $A \in
\cat{A}$ such that $F(A) \iso B$.
\end{defn}

\begin{propn}   
\label{propn:eqv-ffeso}
A functor is an equivalence if and only if it is full, faithful and
essentially surjective on objects.
\end{propn}

\begin{pf}
Exercise~\ref{ex:eqv-ffeso}.
\end{pf}

This result can be compared to the theorem that every bijective group
homomorphism is an isomorphism (that is, its inverse is also a
homomorphism), or that a natural transformation whose components are
isomorphisms is itself an isomorphism (Lemma~\ref{lemma:nat-iso-compts}).
Those two results are useful because they allow us to show that a map is an
isomorphism without directly constructing an inverse.
Proposition~\ref{propn:eqv-ffeso} provides a similar service, enabling us
to prove that a functor $F$ is an equivalence without actually constructing
an `inverse' $G$, or indeed an $\eta$ or an $\epsln$ (in the notation of
Definition~\ref{defn:eqv}).

A corollary of Proposition~\ref{propn:eqv-ffeso} invites us to view full
and faithful%
\index{functor!full and faithful}
functors as, essentially, inclusions of full subcategories:

\begin{cor}     
\label{cor:ff-emb}
Let $F\from \cat{C} \to \cat{D}$ be a full and faithful functor.  Then
$\cat{C}$ is equivalent to the full subcategory $\cat{C}'$ of $\cat{D}$
whose objects are those of the form $F(C)$ for some $C \in \cat{C}$.
\end{cor}

\begin{pf}
The functor $F'\from \cat{C} \to \cat{C}'$ defined by $F'(C) = F(C)$ is
full and faithful (since $F$ is) and essentially surjective on objects (by
definition of $\cat{C}'$).
\end{pf}

This result is true, with the same proof, whether we interpret `of the form
$F(C)$' to mean `equal to $F(C)$' or `isomorphic to $F(C)$'.

\begin{example}	
\label{eg:equivs-skellish}
Let $\cat{A}$ be any category, and let $\cat{B}$ be any full subcategory
containing at least one object from each isomorphism class of $\cat{A}$.
Then the inclusion functor $\cat{B} \incl \cat{A}$ is faithful (like any
inclusion of subcategories), full, and essentially surjective on objects.
Hence $\cat{B} \eqv \cat{A}$.

So if we take a category and remove some (but not all) of the objects in
each isomorphism class, the slimmed-down%
\index{category!slimmed-down}
version is equivalent to the original.  Conversely, if we take a category
and throw in some more objects, each of them isomorphic to one of the
existing objects, it makes no difference: the new, bigger, category is
equivalent to the old one.

For example, let $\FinSet$%
\ntn{FinSet}
be the category of finite%
\index{set!finite}
sets and functions between them.  For each natural number $n$, choose a set
$\lwr{n}$ with $n$ elements, and let $\cat{B}$ be the full subcategory of
$\FinSet$ with objects $\lwr{0}, \lwr{1}$, \ldots.  Then $\cat{B} \eqv
\FinSet$, even though $\cat{B}$ is in some sense much smaller than
$\FinSet$.
\end{example}

\begin{example}
\label{eg:equivs-mon}
In Example~\ref{egs:cats-as}\bref{eg:cats-as:monoids}, we saw that monoids
are essentially the same thing as one-object%
\index{monoid!one-object category@as one-object category}
categories.  With the definition of equivalence in hand, we are nearly
ready to make this statement precise.  We are missing some set-theoretic
language, and we will return to this result once we have that language
(Example~\ref{eg:mon-one-obj-eqv}), but the essential point can be stated
now.

Let $\cat{C}$ be the full subcategory of $\CAT$ whose objects are the
one-object categories.  Let $\Mon$ be the category of monoids.  Then
$\cat{C} \eqv \Mon$.  To see this, first note that given any object $A$ of
any category, the maps $A \to A$ form a monoid under composition (at least,
subject to some set-theoretic restrictions).  There is, therefore, a
canonical functor $F: \cat{C} \to \Mon$ sending a one-object category to
the monoid of maps from the single object to itself.  This functor $F$ is
full and faithful (by Example~\ref{eg:ftrs-between-monoids}) and
essentially surjective on objects.  Hence $F$ is an equivalence.
\end{example}

\begin{example}
An equivalence of the form $\cat{A}^\op \eqv \cat{B}$ is sometimes called a
\demph{duality}%
\index{duality}
between $\cat{A}$ and $\cat{B}$.  One says that $\cat{A}$ is \demph{dual}
to $\cat{B}$.  There are many famous dualities in which $\cat{A}$ is a
category of algebras and $\cat{B}$ is a category of spaces; recall the
slogan `algebra is dual%
\index{duality!algebra--geometry}
to geometry' from Example~\ref{eg:contra-fn-spaces}.

Here are some quite advanced examples, well beyond the scope of this book.
\begin{itemize}
\item 
Stone duality:%
\index{duality!Stone}
the category of Boolean%
\index{Boolean algebra}
algebras is dual to the category of totally disconnected compact Hausdorff
spaces.

\item 
Gelfand--Naimark duality:%
\index{duality!Gelfand--Naimark}
the category of commutative unital $C^*$-algebras%
\index{C-algebra@$C^*$-algebra}
is dual to the category of compact Hausdorff spaces.  ($C^*$-algebras are
certain algebraic structures important in functional analysis.)

\item 
Algebraic geometers%
\index{algebraic geometry}
have several notions of `space', one of which is `affine variety'.%
\index{variety}
Let $k$ be an algebraically closed field.  Then the category of affine
varieties over $k$ is dual to the category of finitely generated
$k$-algebras with no nontrivial nilpotents.

\item 
Pontryagin duality:%
\index{duality!Pontryagin}
the category of locally compact abelian topological%
\index{topological group}%
\index{group!topological}\linebreak
groups is dual to itself.  As the words `topological group' suggest, both
sides of the duality are algebraic \emph{and} geometric.  Pontryagin
duality is an abstraction of the properties of the Fourier%
\index{Fourier analysis}
transform.
\end{itemize}
\end{example}

\begin{example}
It is rarely useful to consider a category of structured objects in
which the maps do not respect that structure.  For instance, let $\cat{A}$
be the category whose objects are groups%
\index{group!non-homomorphisms of groups}
and whose maps are \emph{all} functions between them, not necessarily
homomorphisms.  Let $\Set_{\neq\emptyset}$ be the category of nonempty
sets.  The forgetful functor $U\from \cat{A} \to \Set_{\neq\emptyset}$ is
full and faithful.  It is a (not profound) fact that every nonempty set can
be given at least one group structure, so $U$ is essentially surjective on
objects.  Hence $U$ is an equivalence.  This implies that the category
$\cat{A}$, although defined in terms of groups, is really just the category
of nonempty sets.
\end{example}

\begin{remarks} 
\label{rmks:2-cat-CAT}
Here is a kind of review of the chapter so far.  We have defined:
\begin{itemize}
\item 
categories (Section~\ref{sec:cats});

\item 
functors between categories (Section~\ref{sec:ftrs});

\item 
natural transformations between functors (Section~\ref{sec:nts});

\item 
composition of functors 
\[
\cdot \to \cdot \to \cdot
\]
and the identity functor on any category
(Remark~\ref{rmks:defn-ftr}\bref{rmk:defn-ftr:comp});

\item 
composition of natural transformations%
\index{natural transformation!composition of|(}
\vspace{-3ex}
\[
\xymatrix{
\cdot
\ruppertwocell
\ar[r]
\rlowertwocell
&
\cdot
}
\]
\vspace{-4ex}\\ 
and the identity natural transformation on any functor
(Construction~\ref{constn:comp-nat}).
\end{itemize}
This composition of natural transformations is sometimes called
\demph{vertical%
\index{composition!vertical}
composition}.  There is also \demph{horizontal%
\index{composition!horizontal} 
composition}, which takes natural transformations
\[
\xymatrix@C+.5em{
\cat{A} \rtwocell<4>^F_G{\alpha}   &
\cat{A}' \rtwocell<4>^{F'}_{G'}{\alpha'}   &
\cat{A}''
}
\]
and produces a natural transformation
\[
\xymatrix@C+.5em{
\cat{A} \rtwocell<4>^{F' \of F}_{G' \of G} &\cat{A}'',
}
\]
traditionally written as $\alpha' * \alpha$.%
\ntn{horiz-comp}
The component of $\alpha' * \alpha$ at $A \in \cat{A}$ is defined to be
the diagonal of the naturality square
\[
\xymatrix{
F'(F(A)) \ar[r]^{F'(\alpha_A)} \ar[d]_{\alpha'_{F(A)}}  &
F'(G(A)) \ar[d]^{\alpha'_{G(A)}}        \\
G'(F(A)) \ar[r]_{G'(\alpha_A)}  &
G'(G(A)).
}
\]
In other words, $(\alpha' * \alpha)_A$ can be defined as either 
$\alpha'_{G(A)} \of F'(\alpha_A)$ or $G'(\alpha_A) \of \alpha'_{F(A)}$; it
makes no difference which, since they are equal.  

The special%
\label{p:special-cases}
cases of horizontal composition where either $\alpha$ or $\alpha'$ is an
identity are especially important, and have their own notation.  Thus,
\[
\xymatrix@C+.5em{
\cat{A}
\ar[r]^F        &
\cat{A}'
\rtwocell^{F'}_{G'}{\hspace{.3em}\alpha'}    &
\cat{A}''
}
\qquad
\text{gives rise to}
\qquad
\xymatrix@C+1.5em{
\cat{A} \rtwocell<4>^{F'\of F}_{G'\of F}{\hspace{.8em}\alpha' F}      &\cat{A}''
}%
\ntn{whisker-right}
\]
where $(\alpha' F)_A = \alpha'_{F(A)}$, and 
\[
\xymatrix@C+.5em{
\cat{A}
\rtwocell^F_G{\alpha}   &
\cat{A}'
\ar[r]^{F'}     &
\cat{A}''
}
\qquad
\text{gives rise to}
\qquad
\xymatrix@C+1.5em{
\cat{A} \rtwocell<4>^{F'\of F}_{F' \of G}{\hspace{.8em}F' \alpha} &
\cat{A}''
}%
\ntn{whisker-left}
\]
where $(F'\alpha)_A = F'(\alpha_A)$.  

Vertical and horizontal composition interact well: natural transformations
\[
\xymatrix{
\cat{A} 
\ruppertwocell<8>^F{\alpha}
\ar[r]|G
\rlowertwocell<-8>_H{\beta} &
\cat{A}'
\ruppertwocell<8>^{F'}{\hspace{.2em}\alpha'}
\ar[r]|{G'}
\rlowertwocell<-8>_{H'}{\hspace{.2em}\beta'} &
\cat{A}''
}
\]
obey the \demph{interchange law},%
\index{interchange law}
\[
(\beta' \of \alpha') * (\beta \of \alpha)
=
(\beta' * \beta) \of (\alpha' * \alpha)
\from
F' \of F \to H' \of H.
\]
As usual, a statement on composition is accompanied by a statement on
identities: $1_{F'} * 1_F = 1_{F' \of F}$ too.

All of this enables us to construct, for any categories $\cat{A}$,
$\cat{A}'$ and $\cat{A}''$, a functor
\[
\ftrcat{\cat{A}'}{\cat{A}''}
\times
\ftrcat{\cat{A}}{\cat{A}'}
\to
\ftrcat{\cat{A}}{\cat{A}''},
\index{functor!category}
\]
given on objects by $(F', F) \mapsto F' \of F$ and on maps by $(\alpha',
\alpha) \mapsto \alpha' * \alpha$.  In particular, if $F' \iso G'$ and $F
\iso G$ then $F' \of F \iso G' \of G$, since functors preserve isomorphism
(Exercise~\ref{ex:ftrs-pres-iso}).

(The existence of this functor is similar to the fact that \emph{inside} a
category $\cat{C}$, we have, for any objects $A$, $A'$ and $A''$, a
funct\emph{ion}
\[
\cat{C}(A', A'') \times \cat{C}(A, A') \to \cat{C}(A, A''),
\]
given by $(f', f) \mapsto f' \of f$.)

The diagrams above contain not only objects (0-dimensional) and arrows
$\to$ (1-dimensional), but also double arrows $\Rightarrow$ sweeping out
2-dimensional regions between arrows.  What we are implicitly doing is
called 2-category%
\index{two-category@2-category}
theory.  There is a 2-category%
\index{category!two-category of categories@2-category of categories}
of categories, functors and natural transformations, whose anatomy we have
just been describing.  If we are really serious about categories, we have
to get serious about 2-categories.  And if we are really serious about
2-categories, we have to get serious about 3-categories\ldots%
\index{n-category@$n$-category}
and before we know it, we are studying $\infty$-categories.%
\index{infinity-category@$\infty$-category}
But in this book, we climb no higher than the first rung or two of this
infinite ladder.%
\index{natural transformation!composition of|)}

\end{remarks}

\exs

\begin{question}
Find three examples of natural transformations not mentioned above.
\end{question}

\begin{question}        
\label{ex:nat-iso-compts}
Prove Lemma~\ref{lemma:nat-iso-compts}.
\end{question}

\begin{question}
Let $\cat{A}$ and $\cat{B}$ be categories.  Prove that
$\ftrcat{\cat{A}^\op}{\cat{B}^\op} \iso \ftrcat{\cat{A}}{\cat{B}}^\op$.
\end{question}

\begin{question}
Let $A$ and $B$ be sets, and denote by $B^A$ the set of functions from $A$
to $B$.  Write down:
\begin{enumerate}[(b)]
\item 
a canonical function $A \times B^A \to B$;%
\index{canonical}

\item 
a canonical function $A \to B^{(B^A)}$.
\end{enumerate}
(Although in principle there could be many such canonical functions,
in both these cases there is only one.)
\end{question}

\begin{question}        
\label{ex:nat-iso-on-product}
Here we consider natural transformations between functors whose domain
is a product%
\index{category!product of categories}
category $\cat{A} \times \cat{B}$.  Your task is to show that naturality in
two variables simultaneously is equivalent to naturality in each variable
separately.

Take functors $F, G\from \cat{A} \times \cat{B} \to \cat{C}$.  For each $A
\in \cat{A}$, there are functors $F^A, G^A\from \cat{B} \to \cat{C}$, as in
Exercise~\ref{ex:ftr-on-product}.  Similarly, for each $B\in\cat{B}$,
there are functors $F_B, G_B\from \cat{A} \to \cat{C}$.

Let $\bigl(\alpha_{A, B}\from F(A, B) \to G(A, B)\bigr)_{A \in \cat{A}, B
  \in \cat{B}}$ be a family of maps.  Show that this family is a natural
transformation $F \to G$ if and only if it satisfies the following two
conditions: 
\begin{itemize}
\item 
for each $A \in \cat{A}$, the family $\bigl(\alpha_{A, B}\from F^A(B) \to
G^A(B)\bigr)_{B \in \cat{B}}$ is a natural transformation $F^A \to G^A$;

\item 
for each $B \in \cat{B}$, the family $\bigl(\alpha_{A, B}\from F_B(A) \to
G_B(A)\bigr)_{A \in \cat{A}}$ is a natural transformation $F_B \to G_B$.
\end{itemize}
\end{question}

\begin{question}
Let $G$ be a group.%
\index{group!isomorphism of elements of}
For each $g \in G$, there is a unique homomorphism
$\phi\from \integers \to G$%
\index{Z@$\integers$ (integers)!group@as group}
satisfying $\phi(1) = g$.  Thus, elements of $G$ are essentially the same
thing as homomorphisms $\integers \to G$.  When groups are regarded as
one-object categories, homomorphisms $\integers \to G$ are in turn the same
as functors $\integers \to G$.  Natural isomorphism defines an equivalence
relation on the set of functors $\integers \to G$, and, therefore, an
equivalence relation on $G$ itself.  What is this equivalence relation, in
purely group-theoretic terms?

(First have a guess.  For a general group $G$, what equivalence
relations on $G$ can you think of?)
\end{question}

\begin{question}        
\label{ex:species}
A \demph{permutation}%
\index{permutation}
of a set $X$ is a bijection $X \to X$.  Write $\Sym(X)$ for the set of
permutations of $X$.  A \demph{total%
\index{total order}%
\index{ordered set!totally}
order} on a set $X$ is an order $\leq$ such that for all $x, y \in X$,
either $x \leq y$ or $y \leq x$; so a total order on a finite set amounts
to a way of placing its elements in sequence.  Write $\Ord(X)$ for the set
of total orders on $X$.

Let $\cat{B}$ denote the category of finite sets and bijections.

\begin{enumerate}[(b)]
\item 
Give a definition of $\Sym$ on maps in $\cat{B}$ in such a way that $\Sym$
becomes a functor $\cat{B} \to \Set$.  Do the same for $\Ord$.  Both your
definitions should be canonical (no arbitrary choices).

\item 
Show that there is no natural transformation $\Sym \to \Ord$.
(Hint: consider identity permutations.)

\item 
For an $n$-element set $X$, how many elements do the sets $\Sym(X)$
and\linebreak $\Ord(X)$ have?
\end{enumerate}

Conclude that $\Sym(X) \iso \Ord(X)$ for all $X \in \cat{B}$, but not
\emph{naturally} in $X \in \cat{B}$.  (The moral is that for each finite
set $X$, there are exactly as many permutations of $X$ as there are total
orders on $X$, but there is no natural way of matching them up.)
\end{question}

\begin{question} 
\label{ex:eqv-ffeso}
In this exercise, you will prove Proposition~\ref{propn:eqv-ffeso}.  Let $F
\from \cat{A} \to \cat{B}$ be a functor.
\begin{enumerate}[(b)]
\item 
Suppose that $F$ is an equivalence.  Prove that $F$ is full, faithful and
essentially surjective on objects.  (Hint: prove faithfulness before
fullness.)

\item 
Now suppose instead that $F$ is full, faithful and essentially surjective
on objects.  For each $B \in \cat{B}$, choose an object $G(B)$ of $\cat{A}$
and an isomorphism $\epsln_B\from F(G(B)) \to B$.  Prove that $G$ extends
to a functor in such a way that $(\epsln_B)_{B \in \cat{B}}$ is a natural
isomorphism $FG \to 1_{\cat{B}}$.  Then construct a natural isomorphism
$1_{\cat{A}} \to GF$, thus proving that $F$ is an equivalence.
\end{enumerate}
\end{question}

\begin{question}
This exercise makes precise the idea that linear algebra can equivalently
be done with matrices%
\index{matrix}
or with linear maps.

Fix a field $k$.  Let $\Mt$ be the category whose objects are the natural
numbers and with
\[
\Mt(m, n) 
=
\{ n \times m \text{ matrices over } k \}.
\]
Prove that $\Mt$ is equivalent to $\FDVect$, the category of
finite-dimensional vector%
\index{vector space}
spaces over $k$.  Does your equivalence involve a \emph{canonical} functor
from $\Mt$ to $\FDVect$, or from $\FDVect$ to $\Mt$?

(Part of the exercise is to work out what composition in the category $\Mt$
is supposed to be; there is only one sensible possibility.
Proposition~\ref{propn:eqv-ffeso} makes the exercise easier.)
\end{question}

\begin{question}        
\label{ex:eqv-eq-reln}
Show that equivalence of categories is an equivalence relation.  (Not as
obvious as it looks.)
\end{question}
%
%
%

\chapter{Adjoints}
\label{ch:adj}

The slogan of Saunders Mac Lane's book \emph{Categories for the Working
  Mathematician} is:
\begin{slogan}
Adjoint functors arise everywhere.
\end{slogan}
We will see the truth of this, meeting examples of adjoint functors from
diverse parts of mathematics.  To complement the understanding provided by
examples, we will approach the theory of adjoints from three different
directions, each of which carries its own intuition.  Then we will prove
that the three approaches are equivalent.

Understanding adjointness gives you a valuable addition to your
mathematical toolkit.  Most professional pure mathematicians know what
categories and functors are, but far fewer know about adjoints.  More
should: adjoint functors are both common and easy, and knowing about
adjoints helps you to spot patterns in the mathematical landscape.

\section{Definition and examples}
\label{sec:adj-basics}

Consider a pair of functors in opposite directions, $F\from \cat{A} \to
\cat{B}$ and $G\from \cat{B} \to \cat{A}$.  Roughly speaking, $F$ is said
to be left adjoint to $G$ if, whenever $A \in \cat{A}$ and $B \in \cat{B}$,
maps $F(A) \to B$ are essentially the same thing as maps $A \to G(B)$.

\begin{defn}    
\label{defn:adjn}
Let $\oppairi{\cat{A}}{\cat{B}}{F}{G}$ be categories and functors.  We say
that $F$ is \demph{left adjoint}%
\index{adjunction}
to $G$, and $G$ is \demph{right adjoint} to $F$, and write $F \ladj G$,%
\ntn{ladj}
 if
\begin{equation}        
\label{eq:adjn}
\cat{B}(F(A), B)
\iso
\cat{A}(A, G(B))
\end{equation}
naturally in $A \in \cat{A}$ and $B \in \cat{B}$.  The meaning of
`naturally' is defined below.  An \demph{adjunction} between $F$ and $G$ is
a choice of natural isomorphism~\eqref{eq:adjn}.
\end{defn}

`Naturally in $A \in \cat{A}$ and $B \in \cat{B}$' means that there is a
specified bijection~\eqref{eq:adjn} for each $A \in \cat{A}$ and $B \in
\cat{B}$, and that it satisfies a naturality axiom.  To state it, we need
some notation.  Given objects $A \in \cat{A}$ and $B \in \cat{B}$, the
correspondence~\eqref{eq:adjn} between maps $F(A) \to B$ and $A \to G(B)$
is denoted by a horizontal bar,%
\ntn{adj-bar}
in both directions:
\[
\begin{array}{ccc}
\Bigl(F(A) \toby{g} B\Bigr)             &
\mapsto	&
\Bigl(A \toby{\bar{g}} G(B)\Bigr),	\\
\Bigl(F(A) \toby{\bar{f}} B\Bigr)	&
\mapsfrom	&
\Bigl(A \toby{f} G(B)\Bigr).
\end{array}
\]
So $\bar{\bar{f}} = f$ and $\bar{\bar{g}} = g$.  We call
$\bar{f}$ the \demph{transpose}%
\index{transpose}
of $f$, and similarly for $g$.  The naturality%
\index{adjunction!naturality axiom for}
axiom has two parts:
\begin{equation}        
\label{eq:adj-nat-a}
\ovln{\Bigl(F(A) \toby{g} B \toby{q} B'\Bigr)}
\quad
=
\quad
\Bigl(A \toby{\bar{g}} G(B) \toby{G(q)} G(B')\Bigr)
\end{equation}
(that is, $\ovln{q \of g} = G(q) \of \bar{g}$) for all $g$ and $q$, and 
\begin{equation}        
\label{eq:adj-nat-b}
\ovln{\Bigl(A' \toby{p} A \toby{f} G(B)\Bigr)}
\quad
=
\quad
\Bigl(F(A') \toby{F(p)} F(A) \toby{\bar{f}} B\Bigr)
\end{equation}
for all $p$ and $f$.  It makes no difference whether we put the long bar
over the left or the right of these equations, since bar is self-inverse.

\begin{remarks} 
\label{rmks:adjts}
\begin{enumerate}[(b)]
\item           
\label{rmk:adjts:nat}
The naturality axiom might seem ad hoc, but we will see in
Chapter~\ref{ch:rep} that it simply says that two particular functors are
naturally isomorphic.  In this section, we ignore the naturality axiom
altogether, trusting that it embodies our usual intuitive idea of
naturality: something defined without making any arbitrary choices.

\item 
The naturality axiom implies that from each array of maps
\[
A_0 \to \cdots \to A_n,
\quad
F(A_n) \to B_0,
\quad
B_0 \to \cdots \to B_m,
\]
it is possible to construct exactly one%
\index{uniqueness!constructions@of constructions}
map
\[
A_0 \to G(B_m).
\]
Compare the comments on the definitions of category, functor and natural
transformation (Remarks~\ref{rmks:defn-cat}\bref{rmk:defn-cat:loosely},
\ref{rmks:defn-ftr}\bref{rmk:defn-ftr:loosely},
and~\ref{rmks:defn-nt}\bref{rmk:defn-nt:loosely}).

\item   
\label{rmk:Lie-ass}
Not only do adjoint functors arise everywhere; better, whenever you see a
pair of functors $\cat{A}\oppairu \cat{B}$, there is an excellent chance
that they are adjoint (one way round or the other).

For example, suppose you get talking to a mathematician who tells you that
her work involves Lie%
\index{algebra!associative|(}%
\index{associative algebra|(}%
\index{Lie algebra|(}
algebras and associative algebras.  You try to object
that you don't know what either of those things is, but she carries on
talking anyway, explaining that there's a way of turning any Lie algebra
into an associative algebra, and also a way of turning any associative
algebra into a Lie algebra.  At this point, even without knowing what she's
talking about, you should bet her that one process is adjoint to the other.
This almost always works.

\item   
\label{rmks:adjts:uniqueness}
A given functor $G$ may or may not have a left adjoint, but if it does, it is
unique%
\index{adjunction!uniqueness of adjoints}
up to isomorphism, so we may speak of `\emph{the} left adjoint of $G$'.
The same goes for right adjoints.  We prove this later
(Example~\ref{eg:yon-adjts-unique}).

You might ask `what do we gain from knowing that two functors are adjoint?'
The uniqueness is a crucial part of the answer.  Let us return to the
example of~\bref{rmk:Lie-ass}.  It would take you only a few minutes to
learn what Lie algebras are, what associative algebras are, and what the
standard functor $G$ is that turns an associative algebra into a Lie
algebra.  What about the functor $F$ in the opposite direction?  The
description of $F$ that you will find in most algebra books (under
`universal%
\index{universal!enveloping algebra} 
enveloping algebra') takes much longer to understand.  However, you can
bypass that process completely, just by knowing that $F$ is the left
adjoint of $G$.  Since $G$ can have only \emph{one} left adjoint, this
characterizes $F$ completely.  In a sense, it tells you all you need to
know.%
\index{algebra!associative|)}%
\index{associative algebra|)}%
\index{Lie algebra|)}
\end{enumerate}
\end{remarks}

\begin{examples}[Algebra: free $\ladj$ forgetful]      
\label{egs:adjns-alg}
\index{adjunction!free--forgetful|(}
Forgetful%
\index{functor!forgetful!left adjoint to}
functors between categories of algebraic structures usually have left
adjoints.  For instance:
\begin{enumerate}[(b)]
\item   
\label{eg:adjns-alg:vs}
Let $k$ be a field.  There is an adjunction
\[
\adjn{\Vect_k}{\Set,}{F}{U}
\index{vector space!free}
\]
where $U$ is the forgetful functor of
Example~\ref{egs:forgetful-functors}\bref{eg:forgetful-ring-vs} and $F$ is
the free functor of Example~\ref{egs:free-functors}\bref{eg:free-vs}.
Adjointness says that given a set $S$ and a vector space $V$, a linear map
$F(S) \to V$ is essentially the same thing as a function $S \to U(V)$.

We saw this in Example~\ref{eg:univ-basis}, but let us now check it in
detail.

Fix a set $S$ and a vector space $V$.  Given a linear map $g\from F(S) \to
V$, we may define a map of sets $\bar{g}\from S \to U(V)$ by $\bar{g}(s) =
g(s)$ for all $s \in S$.  This gives a function
\[
\begin{array}{ccc}
\Vect_k(F(S), V)        &\to            &\Set(S, U(V))  \\
g                       &\mapsto        &\bar{g}.
\end{array}
\]
In the other direction, given a map of sets $f\from S \to U(V)$, we may
define a linear map $\bar{f}\from F(S) \to V$ by $\bar{f} \bigl( \sum_{s
  \in S} \lambda_s s \bigr) = \sum_{s \in S} \lambda_s f(s)$ for all formal
linear combinations $\sum \lambda_s s \in F(S)$.  This gives a function
\[
\begin{array}{ccc}
\Set(S, U(V))  &\to            &\Vect_k(F(S), V)       \\
f              &\mapsto        &\bar{f}.
\end{array}
\]
These two functions `bar' are mutually inverse: for any linear map $g\from
F(S) \to V$, we have
\[
\bar{\bar{g}} \Biggl( \sum_{s \in S} \lambda_s s \Biggr)
=
\sum_{s \in S} \lambda_s \bar{g}(s)
=
\sum_{s \in S} \lambda_s g(s)
=
g \Biggl( \sum_{s \in S} \lambda_s s \Biggr)
\]
for all $\sum \lambda_s s \in F(S)$, so $\bar{\bar{g}} = g$, and for any
map of sets $f\from S \to U(V)$, we have
\[
\bar{\bar{f}}(s)
=
\bar{f}(s)
=
f(s)
\]
for all $s \in S$, so $\bar{\bar{f}} = f$.  We therefore have a canonical
bijection between $\Vect_k(F(S), V)$ and $\Set(S, U(V))$ for each $S \in
\Set$ and $V \in \Vect_k$, as required.

Here we have been careful to distinguish between the vector space $V$ and
its underlying set $U(V)$.  Very often, though, in category theory as in
mathematics at large, the symbol for a forgetful functor is omitted.  In
this example, that would mean dropping the $U$ and leaving the reader to
figure out whether each occurrence of $V$ is intended to denote the vector
space itself or its underlying set.  We will soon start using such
notational shortcuts ourselves.

\item   
\label{egs:adjns-alg:gp}
In the same way, there is an adjunction
\[
\adjn{\Grp}{\Set}{F}{U}
\index{group!free}
\]
where $F$ and $U$ are the free and forgetful functors of
Examples~\ref{egs:forgetful-functors}\bref{eg:forgetful-groups}
and~\ref{egs:free-functors}\bref{eg:free-group}.  

The free group functor is tricky to construct explicitly.  In
Chapter~\ref{ch:arl}, we will prove a result (the general adjoint functor
theorem) guaranteeing that $U$ and many functors like it all have left
adjoints.  To some extent, this removes the need to construct $F$
explicitly, as observed in
Remark~\ref{rmks:adjts}\bref{rmks:adjts:uniqueness}.  The point can be
overstated: for a group theorist, the more descriptions of free groups that
are available, the better.  Explicit%
\index{explicit description}
constructions really can be useful.  But it is an important general
principle that forgetful functors of this type always have left adjoints.

\item 
There is an adjunction
\[
\adjn{\Ab}{\Grp}{F}{U}
\]
where $U$ is the inclusion functor of
Example~\ref{egs:forgetful-functors}\bref{eg:forgetful-ab}.  If $G$ is a
group then $F(G)$ is the \demph{abelianization}%
\index{abelianization}%
\index{group!abelianization of}
$\abel{G}$%
\ntn{abel}
of $G$.  This is an abelian quotient group of $G$, with the property that
every map from $G$ to an abelian group factorizes uniquely through
$\abel{G}$:
\[
\xymatrix{
G \ar[r]^-\eta \ar[dr]_{\forall\phi}	&
\abel{G} \ar@{.>}[d]^{\exists!\bar{\phi}}	\\
&
\forall A.
}
\]
Here $\eta$ is the natural map from $G$ to its quotient $\abel{G}$, and $A$
is any abelian group.  (We have adopted the abuse of notation advertised in
example~\bref{eg:adjns-alg:vs}, omitting the symbol $U$ at several places
in this diagram.)  The bijection
\[
\Ab(\abel{G}, A) \iso \Grp(G, U(A))
\]
is given in the left-to-right direction by $\psi \mapsto \psi\of\eta$, and
in the right-to-left direction by $\phi \mapsto \bar{\phi}$.

(To construct $\abel{G}$, let $G'$ be the smallest normal subgroup of $G$
containing $xyx^{-1}y^{-1}$ for all $x, y \in G$, and put $\abel{G} =
G/G'$.  The kernel of any homomorphism from $G$ to an abelian group
contains $G'$, and the universal property follows.)

\item   
\label{eg:adjns-alg:gp-mon}
There are adjunctions
\[
\xymatrix{
\Grp \ar[d]|{U\vphantom{gl'}}	\\
\Mon \ar@<3ex>[u]^F_\ladj \ar@<-3ex>[u]_R^\ladj
}
\index{group!free on monoid}%
\index{monoid!free group on}
\]
between the categories of groups and monoids.  The middle functor $U$ is
inclusion.  The left adjoint $F$ is, again, tricky to describe explicitly.
Informally, $F(M)$ is obtained from $M$ by throwing in an inverse to every
element.  (For example, if $M$ is the additive monoid of natural numbers
then $F(M)$ is the group of integers.)  Again, the general adjoint functor
theorem (Theorem~\ref{thm:gaft}) guarantees the existence of this adjoint.

This example is unusual in that forgetful functors do not usually have
\emph{right} adjoints.  Here, given a monoid $M$, the group $R(M)$ is the
submonoid of $M$ consisting of all the invertible elements.

The category $\Grp$ is both a \demph{reflective}%
\index{subcategory!reflective}%
\index{reflective}
and a \demph{coreflective}%
\index{coreflective}
subcategory of $\Mon$.  This means, by definition, that the inclusion
functor $\Grp \incl \Mon$ has both a left and a right adjoint.  The
previous example tells us that $\Ab$ is a reflective subcategory of $\Grp$.

\item   
\label{egs:adjns-alg:fields}
Let $\Field$%
\ntn{Field}
be the category of fields,%
\index{field}
with ring homomorphisms as the maps.  The forgetful functor $\Field \to \Set$
does \emph{not} have a left adjoint.  (For a proof, see
Example~\ref{eg:no-free-field}.)  The theory of fields is unlike the
theories of groups, rings, and so on, because the operation $x \mapsto
x^{-1}$ is not defined for \emph{all} $x$ (only for $x \neq 0$).
\index{adjunction!free--forgetful|)}
\end{enumerate}
\end{examples}

\begin{remark}  
\label{rmk:alg-thy}
At several points in this book, we make contact with the idea of an
\demph{algebraic%
\index{algebraic theory}
theory}.  You already know several examples: the theory of groups is an
algebraic theory, as are the theory of rings, the theory of vector spaces
over $\reals$, the theory of vector spaces over $\complexes$, the theory of
monoids, and (rather trivially) the theory of sets.  After reading the
description below, you might conclude that the word `theory' is overly
grand, and that `definition' would be more appropriate.  Nevertheless, this
is the established usage.

We will not need to define `algebraic theory' formally, but it will be
important to have the general idea.  Let us begin by considering the theory
of groups.

A group can be defined as a set $X$ equipped with a function $\cdot \from X
\times X \to X$ (multiplication), another function $\blank^{-1}\from X \to
X$ (inverse), and an element $e \in X$ (the identity), satisfying a
familiar list of equations.  More systematically, the three pieces of
structure on $X$ can be seen as maps of sets
\[
\cdot\from X^2 \to X, 
\qquad
\blank^{-1}\from X^1 \to X,
\qquad
e\from X^0 \to X,
\]
where in the last case, $X^0$ is the one-element set $1$ and we are using
the observation that a map $1 \to X$ of sets is essentially the same thing
as an element of $X$.

(You may be more familiar with a definition of group in which only the
multiplication and perhaps the identity are specified as pieces of
\emph{structure}, with the existence of inverses required as a
\emph{property}.  In that approach, the definition is swiftly followed by a
lemma on uniqueness of inverses, guaranteeing that it makes sense to speak
of \emph{the} inverse of an element.  The two approaches are equivalent,
but for many purposes, it is better to frame the definition in the way
described in the previous paragraph.) 

An algebraic theory consists of two things: first, a collection of
operations, each with a specified arity%
\index{arity} 
(number of inputs), and second, a collection of equations.  For example,
the theory of groups has one operation of arity $2$, one of arity $1$, and
one of arity $0$.  An \demph{algebra}%
\index{algebra!algebraic theory@for algebraic theory}
or \demph{model}%
\index{model}
for an algebraic theory consists of a set $X$ together with a specified map
$X^n \to X$ for each operation of arity $n$, such that the equations hold
everywhere.  For example, an algebra for the theory of groups is exactly a
group.

A more subtle example is the theory of vector spaces over $\reals$.  This is
an algebraic theory with, among other things, an infinite number of
operations of arity $1$: for each $\lambda \in \reals$, we have the
operation $\lambda \cdot \dashbk\from X \to X$ of scalar multiplication by
$\lambda$ (for any vector space $X$).  There is nothing special about the
field $\reals$ here; the only point is that it was chosen in advance.  The
theory of vector spaces over $\reals$ is different from the theory of
vector spaces over $\complexes$, because they have different operations of
arity $1$.

In a nutshell, the main property of algebras for an algebraic theory is
that the operations are defined everywhere on the set, and the equations
hold everywhere too.  For example, \emph{every} element of a group has a
specified inverse, and \emph{every} element $x$ satisfies the equation $x
\cdot x^{-1} = 1$.  This is why the theories of groups, rings, and so on,
are algebraic theories, but the theory of fields is not.
\end{remark}

\begin{example}
\label{eg:adjn:spaces}
There are adjunctions
\[
\xymatrix{
\Tp \ar[d]|{U\vphantom{gl'}}	\\
\Set \ar@<3ex>[u]^D_\ladj \ar@<-3ex>[u]_I^\ladj
}
\]
where $U$ sends a space to its set of points, $D$ equips a set with the
discrete%
\index{topological space!discrete}
topology, and $I$ equips a set with the indiscrete%
\index{topological space!indiscrete}%
\index{indiscrete space}
topology.
\end{example}

\begin{figure}
\centering
\setlength{\unitlength}{1em}%
\begin{picture}(10,11.5)(-5,-.8)
\cell{0}{0}{b}{\includegraphics[width=10em]{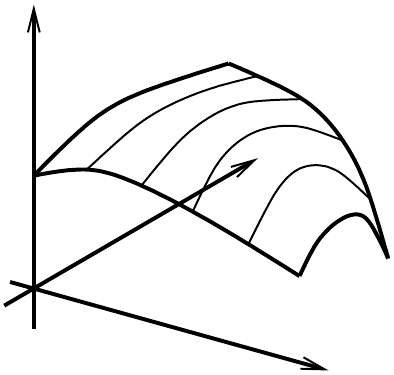}}
\cell{3.8}{-0.2}{c}{A}
\cell{2.2}{5.8}{c}{B}
\cell{-4.1}{10.1}{c}{C}
\end{picture}
\caption{In $\Set$, a map $A \times B \to C$ can be seen as a way of
  assigning to each element of $A$ a map $B \to C$.}
\label{fig:curry}  
\end{figure}

\begin{example}
\label{eg:adjn:cc}
Given sets $A$ and $B$, we can form their (cartesian) product%
\index{set!category of sets!products in}
$A \times B$.  We can also form the set $B^A$%
\index{set!functions@of functions}%
\index{function!set of functions}
of functions from $A$ to $B$.  This is the same as the set $\Set(A, B)$,
but we tend to use the notation $B^A$ when we want to emphasize that it is
an object of the same category as $A$ and $B$.

Now fix a set $B$.  Taking the product with $B$ defines a functor
\[
\begin{array}{cccc}
\dashbk \times B\from	&\Set	&\to	&\Set   \\
			&A	&\mapsto&A \times B.
\end{array}
\]
(Here we are using the blank notation introduced in
Example~\ref{eg:fns-on-vs}.)  There is also a functor
\[
\begin{array}{cccc}
(\dashbk)^B\from	&\Set	&\to	&\Set   \\
			&C	&\mapsto&C^B.
\end{array}
\]
Moreover, there is a canonical bijection
\[
\Set(A \times B, C)
\iso 
\Set(A, C^B)
\]
for any sets $A$ and $C$.  It is defined by simply changing the
punctuation: given a map $g\from A \times B \to C$, define $\bar{g}\from A
\to C^B$ by
\[
(\bar{g}(a))(b) = g(a, b)
\]
($a \in A$, $b \in B$), and in the other direction, given $f\from A \to
C^B$, define $\bar{f}\from A \times B \to C$ by
\[
\bar{f}(a, b) = (f(a))(b)
\]
($a \in A$, $b \in B$).  Figure~\ref{fig:curry} shows an example with $A =
B = C = \reals$.  By slicing up the surface as shown, a map $\reals^2 \to
\reals$ can be seen as a map from $\reals$ to $\{\text{maps } \reals \to
\reals\}$.  

Putting all this together, we obtain an adjunction
\[
\adjn{\Set}{\Set}{\dashbk\times B}{(\dashbk)^B}
\]
for every set $B$.
\end{example}

\begin{defn}    
\label{defn:init-term}
Let $\cat{A}$ be a category.  An object $I \in \cat{A}$ is \demph{initial}%
\index{object!initial}
if for every $A \in \cat{A}$, there is exactly one map $I \to A$.  An
object $T \in \cat{A}$ is \demph{terminal}%
\index{object!terminal}
if for every $A \in \cat{A}$, there is exactly one map $A \to T$.
\end{defn}

For example, the empty set is initial in $\Set$, the trivial group is
initial in $\Grp$, and $\integers$%
\index{Z@$\integers$ (integers)!ring@as ring}
is initial in $\Ring$ (Example~\ref{eg:univ-Z}).  The one-element set is
terminal in $\Set$, the trivial group is terminal (as well as initial) in
$\Grp$, and the trivial (one-element) ring is terminal in $\Ring$.  The
terminal object of $\CAT$ is the category $\One$ containing just one object
and one map (necessarily the identity on that object).

A category need not have an initial object, but if it does have one, it is
unique%
\index{object!initial!uniqueness of}
up to isomorphism.  Indeed, it is unique up to \emph{unique} isomorphism,
as follows.

\begin{lemma}   
\label{lemma:init-unique}
Let $I$ and $I'$ be initial objects of a category.  Then there is a unique
isomorphism $I \to I'$.  In particular, $I \iso I'$.
\end{lemma}

\begin{pf}
Since $I$ is initial, there is a unique map $f\from I \to I'$.  Since $I'$
is initial, there is a unique map $f'\from I' \to I$.  Now $f' \of f$ and
$1_I$ are both maps $I \to I$, and $I$ is initial, so $f' \of f = 1_I$.
Similarly, $f \of f' = 1_{I'}$.  Hence $f$ is an isomorphism, as required.
\end{pf}

\begin{example}
\label{eg:init-term}
\index{object!initial!adjoint@as adjoint}
Initial and terminal objects can be described as adjoints.  Let $\cat{A}$
be a category.  There is precisely one functor $\cat{A} \to \One$.  Also, a
functor $\One \to \cat{A}$ is essentially just an object of $\cat{A}$
(namely, the object to which the unique object of $\One$ is mapped).
Viewing functors $\One \to \cat{A}$ as objects of $\cat{A}$, a left adjoint
to $\cat{A} \to \One$ is exactly an initial object of $\cat{A}$.

Similarly, a right adjoint to the unique functor $\cat{A} \to \One$ is
exactly a terminal object of $\cat{A}$.
\end{example}

\begin{remark}
In the language introduced in Remark~\ref{rmk:principle-duality}, the
concept of terminal object is dual%
\index{duality!principle of}
to the concept of initial object.  (More generally, the concepts of left
and right adjoint are dual to one another.)  Since any two initial objects
of a category are uniquely isomorphic, the principle of duality implies
that the same is true of terminal objects.
\end{remark}

\begin{remark}
Adjunctions can be composed.%
\index{adjunction!composition of adjunctions}
Take adjunctions
\[
\xymatrix{
\cat{A} \ar@<1.1ex>[r]^F_\bot &
\cat{A}' \ar@<1.1ex>[l]^G \ar@<1.1ex>[r]^{F'}_\bot	&
\cat{A}'' \ar@<1.1ex>[l]^{G'}
}
\]
where the $\bot$%
\ntn{bot}
symbol is a rotated $\ladj$ (thus, $F \ladj G$ and $F' \ladj G'$).  Then we
obtain an adjunction
\[
\xymatrix{
\cat{A} \ar@<1.1ex>[r]^{F' \of F}_\bot &
\cat{A}'', \ar@<1.1ex>[l]^{G \of G'}
}
\]
since for $A \in \cat{A}$ and $A'' \in \cat{A}''$,
\[
\cat{A}''\bigl( F'(F(A)), A''\bigr)
\iso 
\cat{A}'\bigl(F(A), G'(A'')\bigr)
\iso 
\cat{A}\bigl(A, G(G'(A''))\bigr)
\]
naturally in $A$ and $A''$.  
\end{remark}

\exs

\begin{question}
Find three examples of adjoint functors not mentioned above.  Do the same
for initial and terminal objects.
\end{question}

\begin{question}
What can be said about adjunctions between discrete categories?
\end{question}

\begin{question}        
\label{ex:adj-nat-in-one}
Show that the naturality%
\index{adjunction!naturality axiom for}
equations~\eqref{eq:adj-nat-a} and~\eqref{eq:adj-nat-b} can equivalently be
replaced by the single equation
\[
\ovln{\Bigl( A' \toby{p} A \toby{f} G(B) \toby{G(q)} G(B') \Bigr)}
\quad
=
\quad
\Bigl(F(A') \toby{F(p)} F(A) \toby{\bar{f}} B \toby{q} B'\Bigr)
\]
for all $p$, $f$ and $q$. 
\end{question}

\begin{question}
Show that left adjoints preserve initial objects: that is, if
$\hadjnli{\cat{A}}{\cat{B}}{F}{G}$ and $I$ is an initial object of
$\cat{A}$, then $F(I)$ is an initial object of $\cat{B}$.  Dually, show
that right adjoints preserve terminal objects.

(In Section~\ref{sec:adj-lim}, we will see this as part of a bigger
picture: right adjoints preserve limits and left adjoints preserve
colimits.)
\end{question}

\begin{question}        
\label{ex:G-set-adjns}
Let $G$ be a group.
\begin{enumerate}[(b)]
\item
\index{group!action of}%
\index{G-set@$G$-set}
What interesting functors are there (in either direction) between $\Set$
and the category $\ftrcat{G}{\Set}$ of left $G$-sets?  Which of those
functors are adjoint to which?

\item 
\index{representation!group or monoid@of group or monoid!linear}
Similarly, what interesting functors are there between $\Vect_k$ and the
category $\ftrcat{G}{\Vect_k}$ of $k$-linear representations of $G$, and
what adjunctions are there between those functors?
\end{enumerate}
\end{question}

\begin{question}
\label{ex:pshf-adjns}
Fix a topological space $X$, and write $\oset(X)$ for the poset of open
subsets of $X$, ordered by inclusion.  Let 
\[
\Delta \from \Set \to \ftrcat{\oset(X)^\op}{\Set}%
\index{functor!diagonal}%
\ntn{Delta}
\]
be the functor assigning to a set $A$ the presheaf%
\index{presheaf}
$\Delta A$ with constant value $A$.  Exhibit a chain of adjoint functors
\[
\Lambda \ladj \Pi \ladj \Delta \ladj \Gamma \ladj \nabla.
\]
\end{question}

\section{Adjunctions via units and counits}
\label{sec:adj-units}

In the previous section, we met the definition of adjunction.  In this
section and the next, we meet two ways of rephrasing the definition.  The
one in this section is most useful for theoretical purposes, while the one
in the next fits well with many examples.

\index{adjunction!naturality axiom for|(}
To start building the theory of adjoint functors, we have to take seriously
the naturality requirement (equations~\eqref{eq:adj-nat-a}
and~\eqref{eq:adj-nat-b}), which has so far been ignored.  Take an
adjunction $\hadjnli{\cat{A}}{\cat{B}}{F}{G}$.  Intuitively, naturality
says that as $A$ varies in $\cat{A}$ and $B$ varies in $\cat{B}$, the
isomorphism between $\cat{B}(F(A), B)$ and $\cat{A}(A, G(B))$ varies in a
way that is compatible with all the structure already in place.  In other
words, it is compatible with composition in the categories $\cat{A}$ and
$\cat{B}$ and the action of the functors $F$ and $G$.

But what does `compatible' mean?  Suppose, for example, that we have maps
\[
F(A) \toby{g} B \toby{q} B'
\]
in $\cat{B}$.  There are two things we can do with this data: either
compose then take the transpose, which produces a map $\ovln{q \of g}\from
A \to G(B')$, or take the transpose of $g$ then compose it with $G(q)$,
which produces a potentially different map $G(q) \of \bar{g}\from A \to
G(B')$.  Compatibility means that they are equal; and that is the first
naturality equation~\eqref{eq:adj-nat-a}.  The second is its dual, and can
be explained in a similar way.

For each $A \in \cat{A}$, we have a map
\[
\Bigl( A \toby{\eta_A} GF(A) \Bigr)
=
\ovln{\Bigl(F(A) \toby{1} F(A)\Bigr)}.%
\ntn{adj-unit}
\]
Dually, for each $B \in \cat{B}$, we have a map
\[
\Bigl( FG(B) \toby{\epsln_B} B \Bigr)
=
\ovln{\Bigl(G(B) \toby{1} G(B)\Bigr)}.%
\ntn{adj-counit}
\]
(We have begun to omit brackets, writing $GF(A)$ instead of $G(F(A))$,
etc.)%
\index{adjunction!naturality axiom for|)}
These define natural transformations
\[
\eta\from 1_\cat{A} \to G \of F,
\qquad
\epsln\from F \of G \to 1_\cat{B},
\]
called the \demph{unit}%
\index{unit and counit}
and \demph{counit} of the adjunction, respectively.

\begin{example}
Take the usual adjunction $\hadjnri{\Vect_k}{\Set}{U}{F}$.%
\index{vector space!free!unit of}
Its unit $\eta\from 1_\Set \to U \of F$ has components
\[
\begin{array}{ccccl}
\eta_S\from         &S      &\to            &UF(S)	&
\!\!\!\!
= 
\bigl\{ 
\text{formal }k\text{-linear sums } \sum_{s \in S} \lambda_s s 
\bigr\}\\
                &s      &\mapsto        &s
\end{array}
\]
($S \in \Set$).  The component of the counit $\epsln$ at a vector space
$V$ is the linear map
\[
\epsln_V\from FU(V) \to V
\]
that sends a \emph{formal} linear sum $\sum_{v \in V} \lambda_v v$ to its
\emph{actual} value in $V$.  

The vector space $FU(V)$ is enormous.  For instance, if $k = \reals$ and
$V$ is the vector space $\reals^2$, then $U(V)$ is the set $\reals^2$ and
$FU(V)$ is a vector space with one basis element for every element of
$\reals^2$; thus, it is uncountably infinite-dimensional.  Then $\epsln_V$
is a map from this infinite-dimensional space to the $2$-dimensional space
$V$.
\end{example}

\begin{lemma}   
\label{lemma:triangle-ids}
Given an adjunction $F \ladj G$ with unit $\eta$ and counit $\epsln$, the
triangles 
\[
\begin{array}{c}
\xymatrix{
F \ar[r]^-{F\eta} \ar[dr]_{1_F}	&
FGF \ar[d]^{\epsln F}	\\
&
F
}
\end{array}
\qquad
\begin{array}{c}
\xymatrix{
G \ar[r]^-{\eta G} \ar[dr]_{1_G}	&
GFG \ar[d]^{G \epsln}	\\
&
G
}
\end{array}
\]
commute.
\end{lemma}

\begin{remark}
These are called the \demph{triangle%
\index{triangle identities}
identities}.  They are commutative diagrams in the functor categories
$\ftrcat{\cat{A}}{\cat{B}}$ and $\ftrcat{\cat{B}}{\cat{A}}$, respectively.
For an explanation of the notation, see Remarks~\ref{rmks:2-cat-CAT}
(particularly the special cases mentioned on
page~\pageref{p:special-cases}).  An equivalent statement is that the
triangles
\begin{equation}        
\label{eq:triangle-ids}
\begin{array}{c}
\xymatrix{
F(A) \ar[r]^-{F(\eta_A)} \ar[dr]_{1_{F(A)}}	&
FGF(A) \ar[d]^{\epsln_{F(A)}}	\\
&
F(A)
}
\end{array}
\qquad
\begin{array}{c}
\xymatrix{
G(B) \ar[r]^-{\eta_{G(B)}} \ar[dr]_{1_{G(B)}}	&
GFG(B) \ar[d]^{G(\epsln_B)}	\\
&
G(B)
}
\end{array}
\end{equation}
commute for all $A \in \cat{A}$ and $B \in \cat{B}$.  
\end{remark}

\begin{pfof}{Lemma~\ref{lemma:triangle-ids}}
We prove that the triangles~\eqref{eq:triangle-ids} commute.  
Let $A \in \cat{A}$.  Since $\ovln{1_{GF(A)}} = \epsln_{F(A)}$,
equation~\eqref{eq:adj-nat-b} gives
\[
\ovln{\Bigl(A \toby{\eta_A} GF(A) \toby{1} GF(A)\Bigr)}
\quad
=
\quad
\Bigl( F(A) \toby{F(\eta_A)} FGF(A) \toby{\epsln_{F(A)}} F(A) \Bigr).
\]
But the left-hand side is $\ovln{\eta_A} = \ovln{\ovln{1_{F(A)}}} =
1_{F(A)}$, proving the first identity.  The second follows by duality.
\end{pfof}

Amazingly, the unit and counit determine the whole adjunction, even though
they appear to know only the transposes \emph{of identities}.  This is the
main content of the following pair of results.

\begin{lemma}   
\label{lemma:unit-determines-adjn}
\index{unit and counit!adjunction in terms of}
Let $\hadjnli{\cat{A}}{\cat{B}}{F}{G}$ be an adjunction, with unit $\eta$
and counit $\epsln$.  Then 
\[
\bar{g} = G(g) \of \eta_A
\]
for any $g\from F(A) \to B$, and 
\[
\bar{f} = \epsln_B \of F(f)
\]
for any $f\from A \to G(B)$.
\end{lemma}

\begin{pf}
For any map $g \from F(A) \to B$, we have
\begin{align*}
\ovln{\Bigl(F(A) \toby{g} B\Bigr)}      &
=
\ovln{\Bigl(F(A) \toby{1} F(A) \toby{g} B\Bigr)}        \\
&
=
\Bigl( A \toby{\eta_A} GF(A) \toby{G(g)} G(B) \Bigr)
\end{align*}
by equation~\eqref{eq:adj-nat-a}, giving the first statement.  The second
follows by duality.
\end{pf}

\begin{thm}   
\label{thm:adj-triangle}
Take categories and functors $\oppairi{\cat{A}}{\cat{B}}{F}{G}$.  There is
a one-to-one correspondence between:
\index{unit and counit!adjunction in terms of}
\begin{enumerate}[(b)]
\item 
adjunctions between $F$ and $G$ (with $F$ on the left and $G$ on the
right);

\item   
\label{item:nat-triangle}
pairs $\Bigl(1_\cat{A} \toby{\eta} GF, \ FG \toby{\epsln} 1_\cat{B}\Bigr)$ of
natural transformations satisfying the triangle identities.
\end{enumerate}
\end{thm}

(Recall that by definition, an adjunction between $F$ and $G$ is a choice
of isomorphism~\eqref{eq:adjn} for each $A$ and $B$, satisfying the
naturality equations~\eqref{eq:adj-nat-a} and~\eqref{eq:adj-nat-b}.)

\begin{pf}
We have shown that every adjunction between $F$ and $G$ gives rise to a
pair $(\eta, \epsln)$ satisfying the triangle identities.  We now have to
show that this process is bijective.  So, take a pair $(\eta, \epsln)$ of
natural transformations satisfying the triangle identities.  We must show
that there is a unique adjunction between $F$ and $G$ with unit $\eta$ and
counit $\epsln$.

Uniqueness follows from Lemma~\ref{lemma:unit-determines-adjn}.  For
existence, take natural transformations $\eta$ and $\epsln$ as
in~\bref{item:nat-triangle}.  For each $A$ and $B$, define functions
\begin{equation}        
\label{eq:adjn-fns}
\cat{B}(F(A), B) 
\oppairu
\cat{A}(A, G(B)),
\end{equation}
both denoted by a bar, as follows.  Given $g \in \cat{B}(F(A), B)$, put
$\bar{g} = G(g) \of \eta_A \in \cat{A}(A, G(B))$.  Similarly, in the
opposite direction, put $\bar{f} = \epsln_B \of F(f)$.

I claim that for each $A$ and $B$, the two functions $g \mapsto \bar{g}$
and $f \mapsto \bar{f}$ are mutually inverse.  Indeed, given a map $g\from
F(A) \to B$ in $\cat{B}$, we have a commutative diagram
\[
\xymatrix{
F(A) \ar[r]^-{F(\eta_A)} \ar[rd]_1	&
FGF(A) \ar[d]^{\epsln_{F(A)}} \ar[r]^-{FG(g)}	&
FG(B) \ar[d]^{\epsln_B}	\\
&
F(A) \ar[r]_-g	&
B.}
\]
The composite map from $F(A)$ to $B$ by one route around the outside of the
diagram is 
\[
\epsln_B \of FG(g) \of F(\eta_A) 
= 
\epsln_B \of F(\bar{g}) 
=
\bar{\bar{g}}, 
\]
and by the other is $g \of 1 = g$, so $\bar{\bar{g}} = g$.  Dually,
$\bar{\bar{f}} = f$ for any map $f\from A \to G(B)$ in $\cat{A}$.  This
proves the claim.

It is straightforward to check the naturality
equations~\eqref{eq:adj-nat-a} and~\eqref{eq:adj-nat-b}.  The
functions~\eqref{eq:adjn-fns} therefore define an adjunction.  Finally, its
unit and counit are $\eta$ and $\epsln$, since the component of the unit at
$A$ is
\[
\ovln{1_{F(A)}} 
= 
G(1_{F(A)}) \of \eta_A 
= 
1 \of \eta_A 
= 
\eta_A,
\]
and dually for the counit.
\end{pf}

\begin{cor}     
\label{cor:adj-triangle}
Take categories and functors $\oppairi{\cat{A}}{\cat{B}}{F}{G}$.  Then $F
\ladj G$ if and only if there exist natural transformations $1 \toby{\eta}
GF$ and $FG \toby{\epsln} 1$ satisfying the triangle identities.  
\qed
\end{cor}

\begin{example} 
\label{eg:poset-adjn}
An adjunction between ordered sets consists of order-pre\-serv\-ing%
\index{ordered set!adjunction between}
maps $\oppairi{A}{B}{f}{g}$ such that
\begin{equation}        
\label{eq:adjn-order}
\forall a \in A, \, \forall b \in B,
\qquad
f(a) \leq b \iff a \leq g(b).
\end{equation}
This is because both sides of the isomorphism~\eqref{eq:adjn} in the
definition of adjunction are sets with at most one element, so they are
isomorphic if and only if they are both empty or both nonempty.  The
naturality requirements~\eqref{eq:adj-nat-a} and~\eqref{eq:adj-nat-b} hold
automatically, since in an ordered set, any two maps with the same domain
and codomain are equal.

Recall from Example~\ref{eg:ftr-cats-orders} that if $\parpair{C}{D}{p}{q}$
are order-preserving maps of ordered sets then there is at most one natural
transformation from $p$ to $q$, and there is one if and only if $p(c) \leq
q(c)$ for all $c \in C$.  The unit of the adjunction above is the statement
that $a \leq gf(a)$ for all $a \in A$, and the counit is the statement that
$fg(b) \leq b$ for all $b \in B$.  The triangle identities say nothing,
since they assert the equality of two maps in an ordered set with the same
domain and codomain.

In the case of ordered sets, Corollary~\ref{cor:adj-triangle} states that
condition~\eqref{eq:adjn-order} is equivalent to:
\[
\forall a \in A, 
\ 
a \leq gf(a)
\qquad \text{ and } \qquad
\forall b \in B,
\ 
fg(b) \leq b.
\]
This equivalence can also be proved directly
(Exercise~\ref{ex:poset-adjn}). 

For instance, let $X$ be a topological space.%
\index{topological space}
Take the set $\cset(X)$ of closed subsets of $X$ and the set $\pset(X)$%
\ntn{power-set}
of all subsets of $X$, both ordered by $\sub$.  There are order-preserving
maps
\[
\oppair{\pset(X)}{\cset(X)}{\Cl}{i}
\index{closure}
\]
where $i$ is the inclusion map and $\Cl$ is closure.  This is an
adjunction, with $\Cl$ left adjoint to $i$, as witnessed by the fact that 
\[
\Cl(A) \sub B 
\iff
A \sub B
\]
for all $A \sub X$ and closed $B \sub X$.  An equivalent statement is that
$A \sub \Cl(A)$ for all $A \sub X$ and $\Cl(B) \sub B$ for all closed $B
\sub X$.  Either way, we see that the topological operation of closure
arises as an adjoint functor.
\end{example}

\begin{remark}  
\label{rmk:eqvs-vs-adjts}
\index{adjunction!equivalence@vs.\ equivalence}%
\index{equivalence of categories!adjunction@vs.\ adjunction}%
\index{category!equivalence of categories!adjunction@vs.\ adjunction}
Theorem~\ref{thm:adj-triangle} states that an adjunction may be regarded as
a quadruple $(F, G, \eta, \epsln)$ of functors and natural transformations
satisfying the triangle identities.  An equivalence $(F, G, \eta, \epsln)$
of categories (as in Definition \ref{defn:eqv}) is not necessarily an
adjunction.  It \emph{is} true that $F$ is left adjoint to $G$
(Exercise~\ref{ex:eqv-is-adjt}), but $\eta$ and $\epsln$ are not
necessarily the unit and counit (because there is no reason why they should
satisfy the triangle identities).
\end{remark}

\begin{remark}
\label{rmk:triangle-string}
There is a way of drawing natural transformations that makes the triangle
identities intuitively plausible.  Suppose, for instance, that we have
categories and functors
\[
\cat{A}   \toby{F_1}
\cat{C}_1 \toby{F_2}
\cat{C}_2 \toby{F_3}
\cat{C}_3 \toby{F_4}
\cat{B},
\qquad
\cat{A}   \toby{G_1}
\cat{D}_1 \toby{G_2}
\cat{B}
\]
and a natural transformation $\alpha\from F_4 F_3 F_2 F_1 \to G_2 G_1$.  We
usually draw $\alpha$ like this:
\[
\begin{xy}
(-40,5)*+{\cat{A}}="A";
(-20,12)*+{\cat{C}_1}="C1";
(0,15)*+{\cat{C}_2}="C2";
(20,12)*+{\cat{C}_3}="C3";
(40,5)*+{\cat{B}}="B";
(0,0)*+{\cat{D}_1}="D1";
(0,7)*+{\Downarrow\,\alpha};
{\ar^{F_1} "A";"C1"};
{\ar^{F_2} "C1";"C2"};
{\ar^{F_3} "C2";"C3"};
{\ar^{F_4} "C3";"B"};
{\ar_{G_1} "A";"D1"};
{\ar_{G_2} "D1";"B"};
\end{xy}
\]
However, we can also draw $\alpha$ as a \demph{string diagram}:%
\index{string diagram}%
\index{diagram!string}
\[
\xymatrix{
F_1 \ar@/_/@{-}[drr] &
\hspace{1em}F_2 \ar@/_/@<1ex>@{-}[dr]&
&
\hspace{-1em}F_3 \ar@/^/@<-1ex>@{-}[dl] &
F_4 \ar@/^/@{-}[dll]\\
&&
*+<1pc>[o][F-]{\alpha}\\
&
G_1 \ar@/^/@{-}[ur]&
&
G_2 \ar@/_/@{-}[ul] 
}
\]
There is nothing special about $4$ and $2$; we could replace them by any
natural numbers $m$ and $n$.  If $m = 0$ then $\cat{A} = \cat{B}$ and the
domain of $\alpha$ is $1_\cat{A}$ (keeping in mind the last paragraph of
Remark~\ref{rmks:defn-cat}\bref{rmk:defn-cat:loosely}).  In that case, the
disk labelled $\alpha$ has no strings coming into the top.  Similarly, if
$n = 0$ then there are no strings coming out of the bottom.

Vertical composition of natural transformations corresponds to joining
string diagrams together vertically, and horizontal composition corresponds
to put\-ting them side by side.  The identity on a functor $F$ is drawn as a
simple string,
\[
\xymatrix{
F \ar@{-}[d] \\
F}
\]

Now let us apply this notation to adjunctions.  The unit and counit are drawn
as 
\[
\begin{array}{c}
\xymatrix@C-2ex{
&
*+<1pc>[o][F-]{\eta}\\
F \ar@/^/@{-}[ur]&
&
G \ar@/_/@{-}[ul] 
}
\end{array}
\qquad
\text{and}
\qquad
\begin{array}{c}
\xymatrix@C-2ex{
G\ar@/_/@{-}[dr]&
&
F\ar@/^/@{-}[dl] \\
&
*+<1pc>[o][F-]{\epsln}
}
\end{array}
\]
The triangle%
\index{triangle identities}
identities now become the topologically plausible equations
\[
\!\begin{array}{c}
\xymatrix@C-3.5ex@R-1.8ex{
&&& F\ar@/^.8pc/@{-}[ddl] \\
& *+<1pc>[o][F-]{\eta} \ar@{-}[dr]|*{G} &&\\
&& *+<1pc>[o][F-]{\epsln} &\\
F\ar@/^.8pc/@{-}[uur] &&&\\
}
\end{array}
=
\begin{array}{c}
\xymatrix@C-3.5ex@R-1.8ex{
F \ar@{-}[dddd]\\
\\
\\
\\
F
}
\end{array}
\qquad\text{and}\qquad
\begin{array}{c}
\xymatrix@C-3.5ex@R-1.8ex{
G\ar@/_.8pc/@{-}[ddr]&&& \\
&& *+<1pc>[o][F-]{\eta} \ar@{-}[dl]|*{F} &\\
& *+<1pc>[o][F-]{\epsln} &&\\
&&& G\ar@/_.8pc/@{-}[uul] \\
}
\end{array}
=
\begin{array}{c}
\xymatrix@C-3.5ex@R-1.8ex{
G \ar@{-}[dddd]\\
\\
\\
\\
G
}
\end{array}
\]
In both equations, the right-hand side is obtained from the left by simply
pulling the string straight.
\end{remark}

\exs

\begin{question}        
\label{ex:poset-adjn}
Let $\oppairi{A}{B}{f}{g}$ be order-preserving maps between ordered sets.
Prove \emph{directly} that the following conditions are equivalent: 
\index{ordered set!adjunction between}
\begin{enumerate}[(b)]
\item 
for all $a \in A$ and $b \in B$, 
\[
f(a) \leq b \iff a \leq g(b);
\]

\item
$a \leq g(f(a))$ for all $a \in A$ and $f(g(b)) \leq b$ for all $b \in B$.
\end{enumerate}
(Both conditions state that $f \ladj g$; see Example~\ref{eg:poset-adjn}.)
\end{question}

\begin{question}
\begin{enumerate}[(b)]
\item   
\label{part:maxl-eqv}
\index{adjunction!fixed points of}%
\index{fixed point}
Let $\hadjnli{\cat{A}}{\cat{B}}{F}{G}$ be an adjunction with unit $\eta$ and
counit $\epsln$.  Write $\Fix(GF)$ for the full subcategory of $\cat{A}$ whose
objects are those $A \in \cat{A}$ such that $\eta_A$ is an isomorphism, and
dually $\Fix(FG) \sub \cat{B}$.  Prove that the adjunction $(F, G, \eta,
\epsln)$ restricts to an equivalence $(F', G', \eta', \epsln')$ between
$\Fix(GF)$ and $\Fix(FG)$.  

\item 
Part~\bref{part:maxl-eqv} shows that every adjunction restricts to an
equivalence between full subcategories in a canonical way.  Take some
examples of adjunctions and work out what this equivalence is.
\end{enumerate}
\end{question}

\begin{question}
\begin{enumerate}[(b)]
\item   
\label{part:refl-defn}
Show that for any adjunction, the right adjoint is full and faithful if and
only if the counit is an isomorphism.  

\item
An adjunction satisfying the equivalent conditions of
part~\bref{part:refl-defn} is called a \demph{reflection}.%
\index{reflection (adjunction)}
(Compare Example~\ref{egs:adjns-alg}\bref{eg:adjns-alg:gp-mon}.)  Of the
examples of adjunctions given in this chapter, which are reflections?
\end{enumerate}
\end{question}

\begin{question}
\begin{enumerate}[(b)]
\item
Let $f\from K \to L$ be a map of sets, and denote by $f^*\from \pset(L) \to
\pset(K)$ the map sending a subset $S$ of $L$ to its inverse%
\index{inverse!image}
image $f^{-1}S \sub K$.  Then $f^*$ is order-preserving with respect to the
inclusion orderings on $\pset(K)$ and $\pset(L)$, and so can be seen as a
functor.  Find left and right adjoints to $f^*$.

\item
Now let $X$ and $Y$ be sets, and write $p\from X \times Y \to X$ for first
projection.  Regard a subset $S$ of $X$ as a predicate%
\index{predicate}
$S(x)$ in one variable $x \in X$, and similarly a subset $R$ of $X \times
Y$ as a predicate $R(x, y)$ in two variables.  What, in terms of
predicates, are the left and right adjoints to $p^*$?  For each of the
adjunctions, interpret the unit and counit as logical implications.  (Hint:
the left adjoint to $p^*$ is often written as $\exists_Y$,%
\index{quantifiers as adjoints}
and the right adjoint as $\forall_Y$.)
\end{enumerate}
\end{question}

\begin{question}
Given a functor $F \from \cat{A} \to \cat{B}$ and a category $\cat{S}$,
there is a functor $F^*\from \ftrcat{\cat{B}}{\cat{S}} \to
\ftrcat{\cat{A}}{\cat{S}}$ defined on objects $Y \in
\ftrcat{\cat{B}}{\cat{S}}$ by $F^*(Y) = Y\of F$ and on maps $\alpha$ by
$F^*(\alpha) = \alpha F$.  Show that any adjunction
$\hadjnli{\cat{A}}{\cat{B}}{F}{G}$ and category $\cat{S}$ give rise to an
adjunction
\[
\hadjnli{\ftrcat{\cat{A}}{\cat{S}}}%
{\ftrcat{\cat{B}}{\cat{S}}}%
{G^*}{F^*}.
\]
(Hint: use Theorem~\ref{thm:adj-triangle}.)
\end{question}

\section{Adjunctions via initial objects}
\label{sec:adj-init}

We now come to the third formulation of adjointness, which is the one you will
probably see most often in everyday mathematics.  

Consider once more the adjunction
\[
\adjn{\Vect_k}{\Set.}{F}{U}
\index{vector space!free!unit of}
\]
Let $S$ be a set.  The universal property of $F(S)$, the vector space whose
basis is $S$, is most commonly stated like this:
\begin{displaytext}
given a vector space $V$, any function $f\from S \to V$ extends uniquely to a
linear map $\bar{f}\from F(S) \to V$.
\end{displaytext}
As remarked in Example~\ref{egs:adjns-alg}\bref{eg:adjns-alg:vs}, forgetful
functors are often forgotten: in this statement, `$f\from S \to V$' should
strictly speaking be `$f\from S \to U(V)$'.  Also, the word `extends' refers
implicitly to the embedding
\[
\begin{array}{cccc}
\eta_S\from     &S      &\to            &UF(S)  \\
                &s      &\mapsto        &s.
\end{array}
\]
So in precise language, the statement reads:
\begin{displaytext}
for any $V \in \Vect_k$ and $f \in \Set(S, U(V))$, there is a unique $\bar{f}
\in \Vect_k(F(S), V)$ such that the diagram
\begin{equation}        
\label{eq:vect-univ}
\begin{array}{c}
\xymatrix{
S \ar[r]^-{\eta_S} \ar[dr]_f     &U(F(S)) \ar[d]^{U(\bar{f})}  \\
&
U(V)
}
\end{array}
\end{equation}
commutes.  
\end{displaytext}
(Compare Example~\ref{eg:univ-basis}.)  In this section, we show that this
statement is equivalent to the statement that $F$ is left adjoint to $U$ with
unit $\eta$.

To do this, we need a definition.

\begin{defn}
Given categories and functors
\[
\xymatrix{
        &\cat{B} \ar[d]^Q       \\
\cat{A} \ar[r]_P        &
\cat{C},
}
\]
the \demph{comma%
\index{comma category}
category} $\comma{P}{Q}$%
\ntn{comma-cat}
(often written as $(P \mathbin{\downarrow} Q)$) is the category defined as
follows:
\begin{itemize}
\item 
objects are triples $(A, h, B)$ with $A \in \cat{A}$, $B \in \cat{B}$, and
$h\from P(A) \to Q(B)$ in $\cat{C}$;

\item 
maps $(A, h, B) \to (A', h', B')$ are pairs $(f\from A \to A',\, g\from B
\to B')$ of maps such that the square
\[
\xymatrix{
P(A)    \ar[r]^{P(f)} \ar[d]_h  &
P(A')   \ar[d]^{h'}     \\
Q(B) \ar[r]_{Q(g)}      &
Q(B')
}
\]
commutes.
\end{itemize}
\end{defn}

\begin{remark}
Given $\cat{A}$, $\cat{B}$, $\cat{C}$, $P$ and $Q$ as above, there
are canonical functors and a canonical natural transformation as shown:
\[
\xymatrix{
\comma{P}{Q} \ar[r] \ar[d] \ar@{}[dr]|-{\nent}  &
\cat{B} \ar[d]^Q        \\
\cat{A} \ar[r]_P        &
\cat{C} 
}
\]
In a suitable 2-categorical sense, $\comma{P}{Q}$ is universal with this
property.
\end{remark}

\begin{example}
Let $\cat{A}$ be a category and $A \in \cat{A}$.  The \demph{slice%
\index{slice category}
category} of $\cat{A}$ over $A$, denoted by $\cat{A}/A$,%
\ntn{slice}
is the category whose objects are maps into $A$ and whose maps are
commutative triangles.  More precisely, an object is a pair $(X, h)$ with
$X \in \cat{A}$ and $h\from X \to A$ in $\cat{A}$, and a map $(X, h) \to
(X', h')$ in $\cat{A}/A$ is a map $f\from X \to X'$ in $\cat{A}$ making the
triangle
\[
\xymatrix{
X \ar[rr]^f \ar[rd]_h   &       &X' \ar[ld]^{h'}     \\
                        &A
}
\]
commute.  

Slice categories are a special case of comma categories.  Recall from
Example~\ref{eg:init-term} that functors $\One \to \cat{A}$ are just
objects of $\cat{A}$.  Now, given an object $A$ of $\cat{A}$, consider the
comma category $\comma{1_\cat{A}}{A}$, as in the diagram
\[
\xymatrix{
        &\One \ar[d]^A       \\
\cat{A} \ar[r]_{1_\cat{A}}      &
\cat{A}.
}
\]
An object of $\comma{1_\cat{A}}{A}$ is in principle a triple $(X, h, B)$
with $X \in \cat{A}$, $B \in \One$, and $h\from X \to A$ in $\cat{A}$; but
$\One$ has only one object, so it is essentially just a pair $(X, h)$.
Hence the comma category $\comma{1_{\cat{A}}}{A}$ has the same objects as
the slice category $\cat{A}/A$.  One can check that it has the same maps
too, so that $\cat{A}/A \iso \comma{1_{\cat{A}}}{A}$.

Dually (reversing all the arrows), there is a \demph{coslice%
\index{coslice category}%
\index{category!coslice}
category} $A/\cat{A} \iso \comma{A}{1_{\cat{A}}}$,%
\ntn{coslice}
 whose objects are the maps out of
$A$.  
\end{example}

\begin{example}
\label{eg:comma-obj-ftr}
Let $G\from \cat{B} \to \cat{A}$ be a functor and let $A \in \cat{A}$.  We
can form the comma category $\comma{A}{G}$,%
\ntn{comma-fix}
 as in the diagram
\[
\xymatrix{
        &\cat{B} \ar[d]^G       \\
\One \ar[r]_-A           &
\cat{A}.
}
\]
Its objects are pairs $(B \in \cat{B},\, f\from A \to G(B))$.  A map $(B,
f) \to (B', f')$ in $\comma{A}{G}$ is a map $q\from B \to B'$ in $\cat{B}$
making the triangle
\[
\xymatrix{
A \ar[r]^-f \ar[rd]_{f'}         &G(B) \ar[d]^{G(q)}     \\
                                &G(B')
}
\]
commute.  

Notice how this diagram resembles the diagram~\eqref{eq:vect-univ} in the
vector space example.  We will use comma categories $\comma{A}{G}$ to
capture the kind of universal property discussed there.

Speaking casually, we say that $f\from A \to G(B)$ is an object of
$\comma{A}{G}$, when what we should really say is that the pair $(B, f)$ is
an object of $\comma{A}{G}$.  There is potential for confusion here, since
there may be different objects $B, B'$ of $\cat{B}$ with $G(B) = G(B')$.
Nevertheless, we will often use this convention.
\end{example}

We now make the connection between comma categories and adjunctions.

\begin{lemma}   
\label{lemma:adj-implies-init}
\index{adjunction!initial objects@via initial objects|(}%
\index{unit and counit!unit as initial object|(}
Take an adjunction $\hadjnli{\cat{A}}{\cat{B}}{F}{G}$ and an object $A \in
\cat{A}$.  Then the unit map $\eta_A\from A \to GF(A)$ is an initial object of
$\comma{A}{G}$. 
\end{lemma}

\begin{pf}
Let $(B,\, f \from A \to G(B))$ be an object of $\comma{A}{G}$.  We have to
show that there is exactly one map from $(F(A), \eta_A)$ to $(B, f)$.

A map $(F(A), \eta_A) \to (B, f)$ in $\comma{A}{G}$ is a map $q\from F(A) \to
B$ in $\cat{B}$ such that
\begin{equation}        
\label{eq:unit-initial}
\begin{array}{c}
\xymatrix{
A \ar[r]^-{\eta_A} \ar[rd]_f     &GF(A) \ar[d]^{G(q)} \\
                                &G(B)
}
\end{array}
\end{equation}
commutes.  But $G(q) \of \eta_A = \bar{q}$ by
Lemma~\ref{lemma:unit-determines-adjn}, so~\eqref{eq:unit-initial} commutes
if and only if $f = \bar{q}$, if and only if $q = \bar{f}$.  Hence
$\bar{f}$ is the unique map $(F(A), \eta_A) \to (B, f)$ in $\comma{A}{G}$.
\end{pf}

We now meet our third and final formulation of adjointness.

\begin{thm}   
\label{thm:adj-comma}
Take categories and functors $\oppairi{\cat{A}}{\cat{B}}{F}{G}$.  There
is a one-to-one correspondence between:
\begin{enumerate}[(b)]
\item 
adjunctions between $F$ and $G$ (with $F$ on the left and $G$ on the
right);
 
\item   
\label{item:init-transf} 
natural transformations $\eta\from 1_\cat{A} \to GF$ such that $\eta_A\from
A \to GF(A)$ is initial in $\comma{A}{G}$ for every $A \in \cat{A}$.
\end{enumerate}
\end{thm}

\begin{pf}
We have just shown that every adjunction between $F$ and $G$ gives rise to
a natural transformation $\eta$ with the property stated
in~\bref{item:init-transf}.  To prove the theorem, we have to show that
every $\eta$ with the property in~\bref{item:init-transf} is the unit of
exactly one adjunction between $F$ and $G$.

By Theorem~\ref{thm:adj-triangle}, an adjunction between $F$ and $G$
amounts to a pair $(\eta, \epsln)$ of natural transformations satisfying
the triangle identities.  So it is enough to prove that for every $\eta$
with the property in~\bref{item:init-transf}, there exists a unique
natural transformation $\epsln\from FG \to 1_\cat{B}$ such that the pair
$(\eta, \epsln)$ satisfies the triangle identities.

Let $\eta\from 1_\cat{A} \to GF$ be a natural transformation with the property
in~\bref{item:init-transf}.

\paragraph*{Uniqueness} 
Suppose that $\epsln, \epsln'\from FG \to 1_\cat{B}$ are natural
transformations such that both $(\eta, \epsln)$ and $(\eta, \epsln')$
satisfy the triangle identities.  One of the triangle identities states
that for all $B \in \cat{B}$, the triangle
\begin{equation}        
\label{eq:comma-triangle}
\begin{array}{c}
\xymatrix@C+1ex{
G(B) \ar[r]^-{\eta_{G(B)}} \ar[rd]_1    &
G(FG(B)) \ar[d]^{G(\epsln_B)}   \\
        &G(B)
}
\end{array}
\end{equation}
commutes.  Thus, $\epsln_B$ is a map 
\[
\Bigl( FG(B), \ G(B) \toby{\eta_{G(B)}} G(FG(B)) \Bigr)
\quad
\longto
\quad 
\Bigl( B, \ G(B) \toby{1} G(B) \Bigr)
\]
in $\comma{G(B)}{G}$.  The same is true of $\epsln'_B$.  But $\eta_{G(B)}$ is
initial, so there is only one such map, so $\epsln_B = \epsln'_B$.  This holds
for all $B$, so $\epsln = \epsln'$.

\paragraph*{Existence}  
For $B \in \cat{B}$, define $\epsln_B\from FG(B) \to B$ to be the unique map
\[
\bigl(FG(B), \eta_{G(B)}\bigr) 
\to 
\bigl(B, 1_{G(B)}\bigr)
\]
in $\comma{G(B)}{G}$.  (So by definition of $\epsln_B$,
triangle~\eqref{eq:comma-triangle} commutes.)  We show that $(\epsln_B)_{B
  \in \cat{B}}$ is a natural transformation $FG \to 1$ such that $\eta$ and
$\epsln$ satisfy the triangle identities.

To prove naturality, take $B \toby{q} B'$ in $\cat{B}$.  We have commutative
diagrams 
\[
\xymatrix@C+1ex{
G(B) \ar[r]^-{\eta_{G(B)}} \ar[rd]_1 \ar[rdd]_{G(q)}     &
GFG(B) \ar[d]^{G(\epsln_B)}     \\
        &G(B) \ar[d]^{G(q)}     \\
        &G(B')
}
\qquad\qquad
\xymatrix@C+1ex{
G(B) \ar[r]^-{\eta_{G(B)}} \ar[d]^{G(q)} 
\ar@/^-5pc/[rdd]_{G(q)}  &
GFG(B) \ar[d]^{GFG(q)}  \\
G(B') \ar[r]^-{\eta_{G(B')}} \ar[rd]_1   &
GFG(B') \ar[d]^{G(\epsln_{B'})} \\
        &
G(B').
}
\]
So $q \of \epsln_B$ and $\epsln_{B'} \of FG(q)$ are both maps $\eta_{G(B)} \to
G(q)$ in $\comma{G(B)}{G}$, and since $\eta_{G(B)}$ is initial, they must be
equal.  This proves naturality of $\epsln$ with respect to $q$.  Hence
$\epsln$ is a natural transformation.

We have already observed that one of the triangle identities,
equation~\eqref{eq:comma-triangle}, holds.  The other states that for $A \in
\cat{A}$,
\[
\xymatrix@C+1ex{
F(A) \ar[r]^-{F(\eta_A)} \ar[rd]_{1_{F(A)}}      &
FGF(A) \ar[d]^{\epsln_{F(A)}}   \\
        &
F(A)
}
\]
commutes.  To prove it, we repeat our previous technique: there are
commutative diagrams
\[
\xymatrix@R+.5em@C+1ex{
A \ar[r]^-{\eta_A} \ar[rdd]_{\eta_A}     &
GF(A) \ar[dd]^{G(1_{F(A)})}     \\
        &       \\
        &GF(A)
}
\qquad\qquad
\xymatrix@C+1ex{
A \ar[r]^-{\eta_A} \ar[d]^{\eta_A}       
\ar@/^-5.7pc/[rdd]_{\eta_A}  &
GF(A) \ar[d]^{GF(\eta_A)}       \\
GF(A) \ar[r]^-{\eta_{GF(A)}} \ar[rd]_1   &
GFGF(A) \ar[d]^{G(\epsln_{F(A)})}       \\
        &
GF(A),
}
\]
so by initiality of $\eta_A$, we have $\epsln_{F(A)} \of F(\eta_A) =
1_{F(A)}$, as required.
\end{pf}

In Section~\ref{sec:adj-lim} we will meet the adjoint functor theorems,
which state conditions under which a functor is guaranteed to have a left
adjoint.  The following corollary is the starting point for their proofs.

\begin{cor}     
\label{cor:pre-AFT}
Let $G\from \cat{B} \to \cat{A}$ be a functor.  Then $G$ has a left adjoint if
and only if for each $A \in \cat{A}$, the category $\comma{A}{G}$ has an
initial object.
\end{cor}

\begin{pf}
Lemma~\ref{lemma:adj-implies-init} proves `only if'.  To prove `if', let us
choose for each $A \in \cat{A}$ an initial object of $\comma{A}{G}$ and
call it $\bigl(F(A), \, \eta_A \from A \to GF(A)\bigr)$.  (Here $F(A)$ and
$\eta_A$ are just the names we choose to use.)  For each map $f\from A \to
A'$ in $\cat{A}$, let $F(f)\from F(A) \to F(A')$ be the unique map such
that
\[
\xymatrix@R=2ex{
A \ar[rr]^-{\eta_A} \ar[rd]_-f    &       &
G(F(A)) \ar[dd]^{G(F(f))}      \\
        &
A'  \ar[rd]_-{\eta_{A'}} &       \\
        &       &G(F(A'))
}
\]
commutes (in other words, the unique map $\eta_A \to \eta_{A'} \of f$ in
$\comma{A}{G}$).  It is easily checked that $F$ is a functor $\cat{A} \to
\cat{B}$, and the diagram tells us that $\eta$ is a natural transformation
$1 \to GF$.  So by Theorem~\ref{thm:adj-comma}, $F$ is left adjoint to $G$.
\end{pf}

This corollary justifies the claim made at the beginning of the section:
that given functors $F$ and $G$, to have an adjunction $F \ladj G$ amounts
to having maps $\eta_A \from A \to GF(A)$ with the universal property
stated there.
\index{adjunction!initial objects@via initial objects|)}%
\index{unit and counit!unit as initial object|)}

\exs

\begin{question}
What can be said about adjunctions between groups (regarded as one-object
categories)? 
\end{question}

\begin{question}
State the dual of Corollary~\ref{cor:pre-AFT}.  How would you prove your
dual statement?
\end{question}

\begin{question}        
\label{ex:eqv-is-adjt}
Let $(F, G, \eta, \epsln)$ be an equivalence of categories, as in
Definition~\ref{defn:eqv}.  Prove that $F$ is left adjoint to $G$ (heeding
the warning in Remark~\ref{rmk:eqvs-vs-adjts}).
\end{question}

\begin{question}	
\label{ex:eta-inj}
\index{unit and counit!injectivity of unit}
Let $\hadjnri{\cat{A}}{\Set}{U}{F}$ be an adjunction.  Suppose that for at
least one $A \in \cat{A}$, the set $U(A)$ has at least two elements.  Prove
that for each set $S$, the unit map $\eta_S\from S \to UF(S)$ is injective.
What does this mean in the case of the usual adjunction between $\Grp$%
\index{group!free}
and $\Set$?
\end{question}

\begin{question}
Given sets $A$ and $B$, a \demph{partial%
\index{partial function}%
\index{function!partial}
function} from $A$ to $B$ is a pair $(S, f)$ consisting of a subset $S \sub
A$ and a function $S \to B$.  (Think of it as like a function from $A$ to
$B$, but undefined at certain elements of $A$.)  Let $\fcat{Par}$ be the
category of sets and partial functions.

Show that $\fcat{Par}$ is equivalent to $\Set_*$, the category of sets
equipped with a distinguished element and functions preserving
distinguished elements.  Show also that $\Set_*$ can be described as a
coslice category in a simple way.
\end{question}

%
%
%

\chapter{Interlude on sets}
\label{ch:sets}

Sets and functions are ubiquitous in mathematics.  You might have the
impression that they are most strongly connected with the pure end of the
subject, but this is an illusion: think of probability density functions in
statistics, data sets in experimental science, planetary motion in
astronomy, or flow in fluid dynamics.

Category theory is often used to shed light on common constructions and
patterns in mathematics.  If we hope to do this in an advanced context, we
must begin by settling the basic notions of set and function.  That is the
purpose of the first section of this chapter.

The definition of category mentions a `collection' of objects and
`collections' of maps.  We will see in the second section that some
collections are too big to be sets, which leads to a distinction between
`small' and `large' collections.  This distinction will be needed later,
most prominently for the adjoint functor theorems (Chapter~\ref{ch:arl}).

The final section takes a historical look at set theory.  It also explains
why the approach to sets taken in this chapter is more relevant to most of
mathematics than the traditional approach is.  None of this section is
logically necessary for anything that follows, but it may provide useful
perspective.

I do not assume that you have encountered axiomatic set theory of any kind.
If you have, it is probably best to put it out of your mind while reading
this chapter, as the approach to set theory that we take is quite different
from the approach that you are most likely to be familiar with.  A brief
comparison of the traditional and categorical approaches can be found at
the very end of the chapter.

\section{Constructions with sets}
\label{sec:Set-properties}

We have made no definition of `set', nor of `function'.  Nevertheless,
guided by our intuition, we can list some properties that we expect the
world of sets and functions to have.  For instance, we can describe some
of the sets that we think ought to exist, and some ways of building new
sets from old.

Intuitively,%
\index{set!intuitive description of}
a set is a bag of points:
\[
\includegraphics[height=5em]{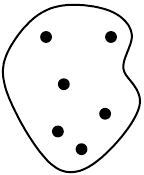}
\]
(There may, of course, be infinitely many.)  These points, or elements, are
not%
\index{set!structurelessness of}
related to one another in any way.  They are not in any order, they do not
come with any algebraic structure (for instance, there is no specified way
of multiplying elements together), and there is no sense of what it means
for one point to be close to another.  In particular examples, we might
have some extra structure in mind; for instance, we often equip the set of
real numbers with an order, a field structure and a metric.  But to view
$\reals$ as a mere \emph{set} is to ignore all that structure, to regard it
as no more than a bunch of featureless points.

Intuitively,%
\index{function!intuitive description of}
a function $A \to B$ is an assignment of a point in bag $B$ to each point
in bag $A$:
\[
\includegraphics[height=5em]{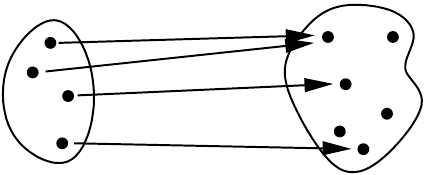}
\]
We can do one function after another: given functions
\[
\includegraphics[height=5em]{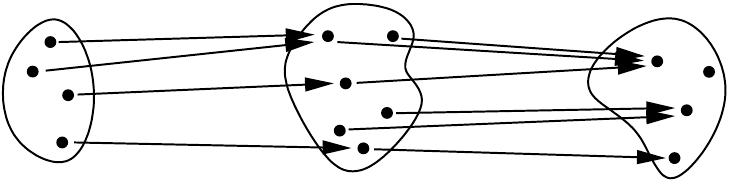}
\]
we obtain a composite function
\[
\includegraphics[height=4.6em]{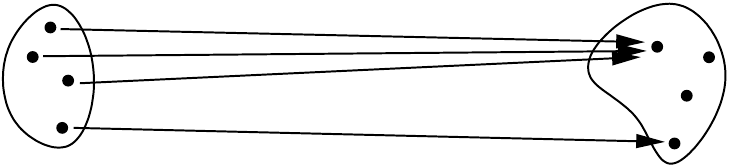}
\]
This composition of functions is associative: $h \of (g \of f) =
(h \of g) \of f$.  There is also an identity function on every set.  Hence:
\begin{setprop}
Sets and functions form a category, denoted by $\Set$.%
\index{set!category of sets}
\end{setprop}

This does not pin things down much: there are many categories, mostly quite
unlike the category of sets.  So, let us list some of the special features
of the category of sets.

\paragraph*{The empty set}  
\index{set!empty}
There is a set $\emptyset$ with no elements.

Suppose that someone hands you a pair of sets, $A$ and $B$, and tells you
to specify a function from $A$ to $B$.  Then your task is to specify for
each element of $A$ an element of $B$.  The larger $A$ is, the longer the
task; the smaller $A$ is, the shorter the task.  In particular, if $A$ is
empty then the task takes no time at all; we have nothing to do.  So there
is a function from $\emptyset$ to $B$ specified by doing nothing.  On the
other hand, there cannot be two different ways to do nothing, so there is
only one function from $\emptyset$ to $B$.  Hence:
\begin{setprop}
$\emptyset$ is an initial object of \hspace{.1em}$\Set$.
\end{setprop}

In case this argument seems unconvincing, here is an alternative.  Suppose
that we have a set $A$ with disjoint subsets $A_1$ and $A_2$ such that $A_1
\cup A_2 = A$.  Then a function from $A$ to $B$ amounts to a function from
$A_1$ to $B$ together with a function from $A_2$ to $B$.  So if all the
sets are finite, we should have the rule
\begin{align*}  
\label{p:num-fns}
\index{function!number of functions}
(\text{number of functions from $A$ to $B$})	&
=	
(\text{number of functions from $A_1$ to $B$})	\\
& \quad\times
(\text{number of functions from $A_2$ to $B$}).
\end{align*}
In particular, we could take $A_1 = A$ and $A_2 = \emptyset$.  This would
force the number of functions from $\emptyset$ to $B$ to be $1$.  So if we
want this rule to hold (and surely we do!), we had better say that there is
exactly one function from $\emptyset$ to $B$.

What about functions \emph{into} $\emptyset$?  There is exactly one function
$\emptyset \to \emptyset$, namely, the identity.  This is a special case of
the initiality of $\emptyset$.  On the other hand, for a set $A$ that is not
empty, there are no functions $A \to \emptyset$, because there is nowhere for
elements of $A$ to go.

\paragraph*{The one-element set}  
\index{set!one-element}
There is a set $1$ with exactly one element.

For any set $A$, there is exactly one function from $A$ to $1$, since every
element of $A$ must be mapped to the unique element of $1$.  That is:
\begin{setprop}
$1$ is a terminal object of \hspace{.1em}$\Set$.
\end{setprop}

A function \emph{from} $1$ \emph{to} a set $B$ is just a choice of an element
of $B$.  In short, the functions $1 \to B$ are the elements of $B$.  Hence:
\begin{slogan}
The concept of element is a special case of the concept of function.%
\index{element!function@as function}
\end{slogan} 

\paragraph*{Products}  Any two sets $A$ and $B$ have a product,%
\index{set!category of sets!products in}
$A \times B$.%
\ntn{prod-set}
  Its elements are the ordered pairs $(a, b)$ with $a \in A$
and $b \in B$.  Ordered pairs are familiar from coordinate geometry.  All
that matters about them is that for $a, a' \in A$ and $b, b' \in B$,
\[
(a, b) = (a', b')
\iff
a = a' \text{ and } b = b'.
\]
More generally, take any set $I$ and any family $(A_i)_{i \in I}$ of sets.
There is a product set $\prod_{i \in I} A_i$,%
\ntn{prod-fam-set}
whose elements are families%
\index{family}
$(a_i)_{i \in I}$ with $a_i \in A_i$ for each $i$.  Just as for ordered
pairs,
\[
(a_i)_{i \in I} = (a'_i)_{i \in I}
\iff
a_i = a'_i \text{ for all } i \in I.
\]

\paragraph*{Sums}  Any two sets $A$ and $B$ have a \demph{sum}%
\index{set!category of sets!sums in}
$A + B$.%
\ntn{sum-set}

Thinking of sets as bags of points, the sum of two sets is obtained by
putting all the points into one big bag:
\[
\begin{array}{c}
\includegraphics[height=5em]{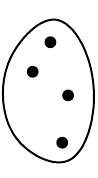}
\end{array}
\, + \,
\begin{array}{c}
\includegraphics[height=5em]{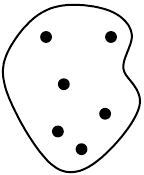}
\end{array}
\ = \
\begin{array}{c}
\includegraphics[height=5em]{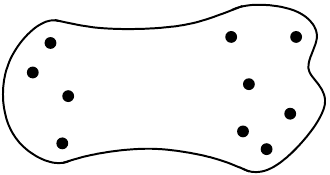}
\end{array}
\]
If $A$ and $B$ are finite sets with $m$ and $n$ elements respectively, then
$A + B$ always has $m + n$ elements.  It makes no difference what the
elements of $A + B$ are called; as usual, we only care what $A + B$ is up
to isomorphism.

There are inclusion functions
\[
A \toby{i} A + B \otby{j} B
\]
such that the union of the images of $i$ and $j$ is all of $A + B$ and the
intersection of the images is empty.

Sum is sometimes called \demph{disjoint union}%
\index{disjoint union}
and written as $\amalg$.%
\ntn{disjt-union}
It is not to be confused with (ordinary) union%
\index{union}
$\cup$.  For a start, we can take the sum of \emph{any} two sets $A$ and
$B$, whereas $A \cup B$ only really makes sense when $A$ and $B$ come as
subsets of some larger set.  (For to say what $A \cup B$ is, we need to
know which elements of $A$ are equal to which elements of $B$.)  And even
if $A$ and $B$ do come as subsets of some larger set, $A + B$ and $A \cup
B$ can be different.  For example, take the subsets $A = \{1, 2, 3\}$ and
$B = \{3, 4\}$ of $\nat$.  Then $A \cup B$ has $4$ elements, but $A + B$
has $3 + 2 = 5$ elements.
  
More generally, any family $(A_i)_{i \in I}$ of sets has a sum $\sum_{i \in
I} A_i$.%
\ntn{sum-fam-set}
If $I$ is finite and each $A_i$ is finite, say with $m_i$
elements, then $\sum_{i \in I} A_i$ has $\sum_{i \in I} m_i$ elements.

\paragraph*{Sets of functions}  
\index{set!functions@of functions}%
\index{function!set of functions}
For any two sets $A$ and $B$, we can form the set $A^B$%
\ntn{exponential-set}
 of functions from $B$ to $A$.  

This is a special case of the product construction: $A^B$ is the product
$\prod_{b \in B} A$ of the constant family $(A)_{b \in B}$.  Indeed, an
element of $\prod_{b \in B} A$ is a family $(a_b)_{b \in B}$ consisting of
one element $a_b \in A$ for each $b \in B$; in other words, it is a function
$B \to A$.

\paragraph*{Digression on arithmetic}%
\label{p:arith}
\index{arithmetic}
We are using notation reminiscent of arithmetic: $A \times B$, $A + B$, and
$A^B$.  There is good reason for this: if $A$ is a finite set with $m$
elements and $B$ a finite set with $n$ elements, then $A \times B$ has $m
\times n$ elements, $A + B$ has $m + n$ elements, and $A^B$ has $m^n$
elements.  Our notation $1$ for a one-element set and the alternative
notation $0$ for the empty set $\emptyset$ also follow this pattern.

All the usual laws of arithmetic have their counterparts for sets:
\begin{align*}
A \times (B + C)        &\iso (A \times B) + (A \times C),      \\
A^{B + C}               &\iso A^B \times A^C,                   \\
(A^B)^C                 &\iso A^{B \times C},
\end{align*}
and so on, where $\iso$ is isomorphism in the category of sets.  (For the
last one, see Example~\ref{eg:adjn:cc}.)  These isomorphisms hold for all
sets, not just finite ones.

\paragraph*{The two-element set}  
\index{set!two-element}
Let $2$%
\ntn{two-set}
be the set $1 + 1$ (a set with two elements!).  For reasons that will soon
become clear, I will write the elements of $2$ as $\true$ and $\false$.

Let $A$ be a set.  Given a subset%
\index{subset}
$S$ of $A$, we obtain a function $\chi_S\from A \to 2$%
\ntn{char-fn}
(the \demph{characteristic%
\index{characteristic function}%
\index{function!characteristic}
function} of $S \sub A$), where
\[
\chi_S(a)
=
\begin{cases}
\true	&\text{if }a \in S,	\\
\false	&\text{if }a \not\in S
\end{cases}
\]
($a \in A$).  Conversely, given a function $f\from A \to 2$, we obtain a
subset
\[
f^{-1}\{\true\} 
= 
\{ a \in A \such f(a) = \true \}
\]
of $A$.  These two processes are mutually inverse; that is, $\chi_S$ is the
unique function $f\from A \to 2$ such that $f^{-1}\{\true\} = S$.  Hence:
\begin{setprop}
Subsets of $A$ correspond one-to-one with functions $A \to 2$.
\end{setprop}
We already know that the functions from $A$ to $2$ form a set, $2^A$.  When we
are thinking of $2^A$ as the set of all subsets of $A$, we call it the
\demph{power%
\index{power!set}
set} of $A$ and write it as $\pset(A)$.%
\ntn{power-set-precise}

\paragraph*{Equalizers}  
\index{equalizer!sets@of sets}%
\index{set!category of sets!equalizers in}
It would be nice if, given a set $A$, we could define a subset $S$ of $A$
by specifying a property that the elements of $S$ are to satisfy:
\[
S
=
\{ a \in A \such \text{some property of } a \text{ holds} \}.
\]
It is hard to give a general definition of `property'.  There is, however, a
special type of property that is easy to handle: equality of two functions.
Precisely, given sets and functions $\parpair{A}{B}{f}{g}$, there is a set
\[
\{ a \in A \such f(a) = g(a) \}.
\]
This set is called the \demph{equalizer} of $f$ and $g$, since it is the
part of $A$ on which the two functions are equal.

\paragraph*{Quotients}
You are probably familiar with quotient groups and quotient rings
(sometimes called factor groups and factor rings) in algebra.  Quotients
also come up everywhere in topology, such as when we glue together opposite
sides of a square to make a cylinder.  But the most basic context for
quotients is that of sets.

Let $A$ be a set and $\sim$ an equivalence%
\index{equivalence relation}
relation on $A$.  There is a set $\qer{A}{\sim}$,%
\ntn{qt-set}
the \demph{quotient of $A$ by $\sim$},%
\index{quotient!set@of set}%
\index{set!quotient of}
whose elements are the equivalence classes.  For example, given a group $G$
and a normal subgroup $N$, define an equivalence relation $\sim$ on $G$ by
$g \sim h \iff gh^{-1} \in N$; then $\qer{G}{\sim} = G/N$.

There is also a canonical map
\[
p\from A \to \qer{A}{\sim},
\]
sending an element of $A$ to its equivalence class.  It is surjective, and has
the property that $p(a) = p(a') \iff a \sim a'$.  In fact, it has a universal
property: any function $f\from A \to B$ such that
\begin{equation}        
\label{eq:er-univ}
\forall a, a' \in A,
\qquad
a \sim a'
\implies
f(a) = f(a')
\end{equation}
factorizes uniquely through $p$, as in the diagram
\[
\xymatrix{
A \ar[r]^-p \ar[rd]_{f}  &\qer{A}{\sim} \ar@{.>}[d]^{\bar{f}}    \\
                        &B.
}
\]
Thus, for any set $B$, the functions $\qer{A}{\sim} \to B$ correspond
one-to-one with the functions $f\from A \to B$
satisfying~\eqref{eq:er-univ}.  This fact is at the heart of the famous
isomorphism theorems of algebra.

\subjectchange

We have now listed the properties of sets and functions that will be most
important for us.  Here are two more.

\paragraph*{Natural numbers}  
\index{natural numbers}
A function with domain $\nat$ is usually called a \demph{sequence}.%
\index{sequence}
A crucial property of $\nat$ is that some sequences can be defined
recursively: given a set $X$, an element $a \in X$, and a function $r\from
X \to X$, there is a unique sequence $(x_n)_{n = 0}^\infty$ of elements of
$X$ such that
\[
x_0 = a,
\qquad
x_{n + 1} = r(x_n) \text{ for all } n \in \nat.
\]
This property refers to two pieces of structure on $\nat$: the element $0$,
and the function $s\from \nat \to \nat$ defined by $s(n) = n + 1$.
Reformulated in terms of functions, and writing $x_n = x(n)$, the property
is this: for any set $X$, element $a \in X$, and function $r\from X \to X$,
there is a unique function $x\from \nat \to X$ such that $x(0) = a$ and
$x\of s = r \of x$.  Exercise~\ref{ex:nno} asks you to show that this is a
\emph{universal} property of $\nat$, $0$ and $s$.

\paragraph*{Choice}  Let $f\from A \to B$ be a map in a category $\cat{A}$.  A
\demph{section}%
\index{section}
(or \demph{right%
\index{inverse!right}
inverse}) of $f$ is a map $i\from B \to A$ in $\cat{A}$ such that $f \of i
= 1_B$.

In the category of sets, any map with a section is certainly surjective.  The
converse statement is called the \demph{axiom%
\index{axiom of choice}
of choice}:
\begin{setprop}
Every surjection has a section.
\end{setprop}
It is called `choice' because specifying a section of $f\from A \to B$ amounts
to choosing, for each $b \in B$, an element of the nonempty set $\{ a \in A
\such f(a) = b\}$.

\subjectchange

\index{foundations|(}%
\index{set!definition of|(}%
The properties listed above are not theorems, since we do not have rigorous
definitions of set and function.  What, then, is their status?

Definitions in mathematics usually depend on previous definitions.  A
vector space is defined as an abelian group with a scalar multiplication.
An abelian group is defined as a group with a certain property.  A group is
defined as a set with certain extra structure.  A set is defined as\ldots\
well, what?

We cannot keep going back indefinitely, otherwise we quite literally would not
know what we were talking about.  We have to start somewhere.  In other words,
there have to be some basic concepts not defined in terms of anything else.
The concept of set is usually taken to be one of the basic ones, which is
why you have probably never read a sentence beginning `Definition: A \demph{set}
is\ldots'.  We will treat function as a basic concept, too.

But now there seems to be a problem.  If these basic concepts are not
defined in terms of anything else, how are we to know what they really are?
How are we going to reason in the watertight, logical way upon which
mathematics depends?  We cannot simply trust our intuitions, since your
intuitive idea of set might be slightly different from mine, and if it came
to a dispute about how sets behave, we would have no way of deciding who
was right.

The problem is solved as follows.  Instead of \emph{defining} a set to be a
such-and-such and a function to be a such-and-such else, we list some
\emph{properties} that we assume sets and functions to have.  In other words,
we never attempt to say what sets and functions \emph{are}; we just say what
you can \emph{do} with them.

In his excellent book \emph{Mathematics: A Very Short Introduction},
Timothy \citeGow\ considers the question: `What is the black%
\index{black king}
king in chess?'%
\index{chess}  
He swiftly points out that this question is rather peculiar.  It is not
important that the black king \emph{is} a small piece of wood, painted a
certain colour and carved into a certain shape.  We could equally well use
a scrap of paper with `BK' written on it.  What matters is what the black
king \emph{does}: it can move in certain ways but not others, according to
the rules of chess.

Similarly, we might not be able to say directly what a set or function
`is', but we agree that they are to satisfy all the properties on the list.
So the list of properties acts as an agreement on how to use the words
`set' and `function', just as the rules of chess act as an agreement on how
to use the chess pieces.

What we are doing is often referred to as \emph{foundations}.  In this
metaphor, the foundation consists of the basic concepts (set and function),
which are not built on anything else, but are assumed to satisfy a stated
list of properties.  On top of the foundations are built some basic
definitions and theorems.  On top of those are built further definitions
and theorems, and so on, towering upwards.

The properties above are stated informally, but they can be formalized
using some categorical language.  (See \citeLR\ or \citeRST.)  In the
formal version, we begin by saying that sets and functions form a category,
$\Set$.  We then list some properties of this category.  For example, the
category is required to have an initial and a terminal object, and the
properties described informally under the headings `Products' and
`Equalizers' are made formal by the statement that $\Set$ `has limits' (a
phrase defined in Chapter~\ref{ch:lims}).

While we were making the list, we were guided by our intuition about sets.
But once it is made, our intuition plays no further official role: any
disputes about the nature of sets are settled by consulting the list of
properties.

(A subtlety arises.  Whatever list of properties one writes down, there
might be some questions that cannot be settled.  In other words, there
might be multiple inequivalent categories satisfying all the properties
listed.  This gets us into the realm of advanced logic: G\"odel
incompleteness, the continuum hypothesis, and so on, all beyond the scope
of this book.)

Now let us look again at the section on the empty%
\index{set!empty}
set.  You might have felt that I was on shaky ground when trying to
convince you that $\emptyset$ is initial.  But the point is that I do not
need to convince you that this is a \emph{true statement}; I only need to
convince you that it is a \emph{convenient assumption}.  Compare the rule
for numbers that $x^0 = 1$.  One can reasonably argue that $0$ copies of
$x$ multiplied together ought to be $1$, but really the best justification
for this rule is convenience: it makes other rules such as $x^{m + n} = x^m
\cdot x^n$ true without exception.  Indeed, it does not even make sense to
ask whether it is `true' that $\emptyset$ is initial until we have written
down our assumptions about how sets and functions behave.  For until then,
what could `true' mean?  There is no physical world of sets against which
to test such statements.

We can make whatever assumptions about sets we like, but some lead to more
interesting mathematics than others.  If, for instance, you want to assume
that there are \emph{no} functions from $\emptyset$ to any other set, you
can, but the tower of mathematics built on that foundation will look
different from what you are used to, and probably not in a good way.  For
example, the `number of functions' rule (page~\pageref{p:num-fns}) will
fail, and there will be further unpleasant surprises higher up the tower.%
\index{foundations|)}%
\index{set!definition of|)}%

\exs

\begin{question}        
\label{ex:diagonal-Set}
The \demph{diagonal%
\index{functor!diagonal}
functor} $\Delta\from \Set \to \Set \times \Set$%
\ntn{diag-set}
is defined by $\Delta(A) = (A, A)$ for all sets $A$.  Exhibit left and
right adjoints to $\Delta$.
\end{question}

\begin{question}        
\label{ex:nno}
In the paragraph headed `Natural numbers', it was observed that the set
$\nat$, together with the element $0$ and the function $s\from \nat \to
\nat$, has a certain property.  This property can be understood as stating
that the triple $(\nat, 0, s)$ is the initial object of a certain category
$\cat{C}$.  Find $\cat{C}$.
\end{question}

\section{Small and large categories}
\label{sec:small-large}

We have now made some assumptions about the nature of sets.  One
consequence of those assumptions is that in many of the categories we have
met, the collection of all objects is too large to form a set.  In fact,
even the collection of \emph{isomorphism classes} of objects is often too
large to form a set.  In this section, I will explain what these statements
mean, and prove them.

This section is not of central importance.  As this book proceeds, I will
say as little as possible about the distinction between sets and
collections too large to be sets.  Nevertheless, the distinction begins to
matter in parts of category theory lying just within the scope of this book
(the adjoint functor theorems), as well as beyond.

Given sets $A$ and $B$, write $\crd{A} \leq \crd{B}$%
\index{set!size of|(}%
\ntn{card-leq}
(or $\crd{B} \geq \crd{A}$) if there exists an injection $A \to B$.  We
give no meaning to the expression `$\crd{A}$' or `$\crd{B}$' in isolation.
(It would perhaps be more logical to write $A \leq B$ rather than $\crd{A}
\leq \crd{B}$, but the notation is well-established.)  In the case of
finite sets, it just means that the number of elements of $A$ is less than
or equal to the number of elements of $B$.

Since identity maps are injective, $\crd{A} \leq \crd{A}$ for all sets
$A$, and since the composite of two injections is an injection,
\[
\crd{A} \leq \crd{B} \leq \crd{C} \implies \crd{A} \leq \crd{C}.
\]
Also, if $A \iso B$ then $\crd{A} \leq \crd{B} \leq \crd{A}$.  Less obvious
is the converse:

\begin{thm}[Cantor--Bernstein]	
\label{thm:cantor-bernstein}
\index{Cantor, Georg!Cantor--Bernstein theorem}
Let $A$ and $B$ be sets.  If $\crd{A} \leq \crd{B} \leq \crd{A}$ then $A
\iso B$. 
\end{thm}

\begin{pf}
Exercise~\ref{ex:cantor-bernstein}.  
\end{pf}

These observations tell us that $\leq$ is a preorder
(Example~\ref{egs:cats-as}\bref{eg:cats-as:orders}) on the collection of
all sets.  It is not a genuine order, since $\crd{A} \leq \crd{B} \leq
\crd{A}$ only implies that $A \iso B$, not $A = B$.  We write $\crd{A} =
\crd{B}$,%
\ntn{card-eq}
and say that $A$ and $B$ \demph{have the same cardinality},%
\index{cardinality}
if $A \iso B$, or equivalently if $\crd{A} \leq \crd{B} \leq \crd{A}$.

As long as we do not confuse equality with isomorphism, the sign $\leq$
behaves as we might imagine.  For example, write $\crd{A} < \crd{B}$ if
$\crd{A} \leq \crd{B}$ and $\crd{A} \neq \crd{B}$.  Then
\begin{equation}	
\label{eq:bernstein-impl}
\crd{A} \leq \crd{B} < \crd{C} \implies \crd{A} < \crd{C}
\end{equation}
for sets $A$, $B$ and $C$.  Indeed, we have already established that
$\crd{A} \leq \crd{C}$, and the strict inequality follows from
Theorem~\ref{thm:cantor-bernstein}.

Here is another fundamental result of set theory.

\begin{thm}[Cantor]  
\label{thm:cantor}
\index{Cantor, Georg!Cantor's theorem}
Let $A$ be a set.  Then $\crd{A} < \crd{\pset(A)}$.  
\end{thm}
Recall that $\pset(A)$ is the power set of $A$.  The lemma is easy for finite
sets, since if $A$ has $n$ elements then $\pset(A)$ has $2^n$ elements, and
$n < 2^n$.

\begin{pf}
Exercise~\ref{ex:cantor-diagonal}.
\end{pf}

\begin{cor}     
\label{cor:no-biggest-set}
For every set $A$, there is a set $B$ such that $\crd{A} < \crd{B}$.
\qed
\end{cor}

In other words, there is no biggest set. 

We now justify the claim made at the beginning of this section: that for
many familiar categories, the collection of isomorphism classes of objects
is too large to form a set.  We begin by doing this for the category $\Set$
itself.  

As a clue to why the collection of isomorphism classes of sets might be too
large to form a set, consider the following statement: the collection of
isomorphism classes of \emph{finite} sets is too large to form a
\emph{finite} set.  This is because there is one isomorphism class of
finite sets for each natural number, but there are infinitely many natural
numbers.

\begin{propn}	
\label{propn:iso-classes-of-sets}
Let $I$ be a set, and let $(A_i)_{i \in I}$ be a family of sets.  Then there
exists a set not isomorphic to any of the sets $A_i$.
\end{propn}

\begin{pf}
Put
\[
A = \pset\Biggl(\sum_{i \in I} A_i\Biggr),
\]
the power set of the sum of the sets $A_i$.  For each $j \in I$, we have
the inclusion function $A_j \to \sum_{i \in I} A_i$, so by
Theorem~\ref{thm:cantor},
\[
\crd{A_j} \leq \biggl| \sum_{i \in I} A_i \biggr| < \crd{A}.
\]
Hence $\crd{A_j} < \crd{A}$ by~\eqref{eq:bernstein-impl}, and in particular,
$A_j \not\iso A$.
\end{pf}
\index{set!size of|)}%

We use the word \demph{class}%
\index{class}
informally to mean any collection of mathematical objects.  All sets are
classes, but some classes (such as the class of all sets) are too big to be
sets.  A class will be called \demph{small}%
\index{small}
if it is a set, and \demph{large}%
\index{large}
otherwise.  For example, Proposition~\ref{propn:iso-classes-of-sets} states
that the class of isomorphism classes of sets is large.  The crucial point
is:
\begin{slogan}
Any \emph{individual} set is small, but the \emph{class} of sets is
large.
\end{slogan}
This is even true if we pretend that isomorphic sets are equal.

Although the `definition' of class is not precise, it will do for our
purposes.  We make a naive distinction between small and large collections,
and implicitly use some intuitively plausible principles (for example, that
any subcollection of a small collection is small).

A category $\cat{A}$ is \demph{small}%
\index{small}%
\index{category!small}
if the class or collection of all maps in $\cat{A}$ is small, and
\demph{large}%
\index{large}%
\index{category!large}
otherwise.  If $\cat{A}$ is small then the class of objects of $\cat{A}$ is
small too, since objects correspond one-to-one with identity maps.

A category $\cat{A}$ is \demph{locally small}%
\index{locally small}%
\index{category!locally small}
if for each $A, B \in \cat{A}$, the class $\cat{A}(A, B)$ is small.  (So,
small implies locally small.)  Many authors take local smallness to be part
of the definition of category.  The class $\cat{A}(A, B)$ is often called
the \demph{hom-set}%
\index{hom-set}
from $A$ to $B$, although strictly speaking, we should only call it this
when $\cat{A}$ is locally small.

\begin{example}
\index{set!category of sets!locally small@is locally small}
$\Set$ is locally small, because for any two sets $A$ and $B$, the
functions from $A$ to $B$ form a set.  This was one of the properties of
sets stated in Section~\ref{sec:Set-properties}.
\end{example}

\begin{example}
\index{vector space!category of vector spaces!locally small@is locally small}%
\index{group!category of groups!locally small@is locally small}%
\index{ring!category of rings!locally small@is locally small}%
\index{topological space!category of topological spaces!locally small@is
  locally small}%
$\Vect_k$, $\Grp$, $\Ab$, $\Ring$ and $\Tp$ are all locally small.  For
example, given rings $A$ and $B$, a homomorphism from $A$ to $B$ is a
function from $A$ to $B$ with certain properties, and the collection of all
functions from $A$ to $B$ is small, so the collection of homomorphisms from
$A$ to $B$ is certainly small.
\end{example}

A category is small if and only if it is locally small and its class of
objects is small.  Again, it may help to consider a similar fact about
finiteness: a category $\cat{A}$ is finite (that is, the class of all maps
in $\cat{A}$ is finite) if and only if it is locally finite (that is, each
class $\cat{A}(A, B)$ is finite) and its class of objects is finite.

\begin{example}	
\label{eg:FinSet-skel}
Consider the category $\cat{B}$ defined in the last paragraph of 
Example~\ref{eg:equivs-skellish}.  Its objects correspond to the natural
numbers, which form a set, so the class of objects of $\cat{B}$ is small.
Each hom-set $\cat{B}(\lwr{m}, \lwr{n})$ is a set (indeed, a finite set),
so $\cat{B}$ is locally small.  Hence $\cat{B}$ is small.
\end{example}

A category is \demph{essentially small}%
\index{essentially small}%
\index{category!essentially small}
if it is equivalent to some small category.  For example, the category of
finite%
\index{set!finite}
sets is essentially small, since by Example~\ref{eg:equivs-skellish}, it is
equivalent to the small category $\cat{B}$ just mentioned.

If two categories $\cat{A}$ and $\cat{B}$ are equivalent, the class of
isomorphism classes of objects of $\cat{A}$ is in bijection with that of
$\cat{B}$.  In a small category, the class of objects is small, so the
class of isomorphism classes of objects is certainly small.  Hence in an
essentially small category, the class of isomorphism classes of objects is
small.  From this we deduce:

\begin{propn}	
\label{propn:Set-large}
\index{set!category of sets!essentially small@is not essentially small}
$\Set$ is not essentially small.
\end{propn}

\begin{pf}
Proposition~\ref{propn:iso-classes-of-sets} states that the class of
isomorphism classes of sets is large.  The result follows.
\end{pf}

By adapting this argument, we can show that many of our standard examples
of categories are not essentially small.  The strategy is to prove that
there are at least as many objects of our category as there are sets.

\begin{example}
For any field $k$, the category $\Vect_k$ of vector%
\index{vector space!category of vector spaces!essentially small@is not
  essentially small} 
spaces over $k$ is not essentially small.  As in the proof of
Proposition~\ref{propn:Set-large}, it is enough to prove that the class of
isomorphism classes of vector spaces is large.  In other words, it is
enough to prove that for any set $I$ and family $(V_i)_{i \in I}$ of vector
spaces, there exists a vector space not isomorphic to any of the spaces
$V_i$.

To show this, write $\hadjnri{\Vect_k}{\Set}{U}{F}$ for the free
and forgetful functors.  As in the proof of
Proposition~\ref{propn:iso-classes-of-sets}, the set
\[
S = \pset \Biggl( \sum_{i \in I} U(V_i) \Biggr) 
\]
has the property that $\crd{U(V_i)} < \crd{S}$ for all $i \in I$.  The free
vector space $F(S)$ on $S$ contains a copy of $S$ as a basis, so $\crd{S}
\leq \crd{UF(S)}$.  Hence $\crd{U(V_i)} < \crd{UF(S)}$ for all $i$, and so
$F(S) \not\iso V_i$ for all $i$, as required.
\end{example}

Similarly, none of the categories $\Grp$,%
\index{group!category of groups!essentially small@is not essentially small}
$\Ab$, $\Ring$%
\index{ring!category of rings!essentially small@is not essentially small}
and $\Tp$%
\index{topological space!category of topological spaces!essentially small@is not
  essentially small}  
is essentially small (Exercise~\ref{ex:not-ess-small}).

Recall that the category of \emph{all} categories and functors is written
as $\CAT$.  

\begin{defn}    
\label{defn:Cat}
We denote by $\Cat$%
\index{category!category of categories}%
\ntn{Cat}
the category of small categories and functors between them.
\end{defn}

\begin{example} 
\label{eg:mon-one-obj-eqv}
Monoids are by definition \emph{sets} equipped with certain structure, so
the one-object%
\index{monoid!one-object category@as one-object category}
categories that they correspond to are small.  Let $\cat{M}$ be the full
subcategory of $\Cat$ consisting of the one-object categories.  Then there
is an equivalence of categories $\Mon \eqv \cat{M}$.  This is proved by the
argument in Example~\ref{eg:equivs-mon}, noting that because each object of
$\cat{M}$ is a \emph{small} one-object category, the collection of maps
from the single object to itself really is a set.
\end{example}

\exs

\begin{question}        
\label{ex:cantor-bernstein}
\begin{enumerate}[(b)]
\item   
\label{part:cb-fix}
Let $A$ be a set.  Let $\theta\from \pset(A) \to \pset(A)$ be a map that is
order-preserving with respect to inclusion.  A \demph{fixed%
\index{fixed point}
point} of $\theta$ is an element $S \in \pset(A)$ such that $\theta(S) =
S$.  By considering
\[
S = \bigcup_{R \in \pset(A) \colon \theta(R) \supseteq R} R,
\]
prove that $\theta$ has at least one fixed point.

\item
Take sets and functions $\oppairi{A}{B}{f}{g}$.  Using~\bref{part:cb-fix},
show that there is some subset $S$ of $A$ such that $g(B \without fS) = A
\without S$.

\item
Deduce the Cantor--Bernstein theorem (Theorem~\ref{thm:cantor-bernstein}). 
\end{enumerate}
\end{question}

\begin{question}        
\label{ex:cantor-diagonal}
\begin{enumerate}[(b)]
\item 
Let $A$ be a set and $f\from A \to \pset(A)$ a function.  By considering
\[
\{ a \in A \such a \not\in f(a) \},
\]
prove that $f$ is not surjective.

\item 
Deduce Cantor's theorem (Theorem~\ref{thm:cantor}): $\crd{A} <
\crd{\pset(A)}$ for all sets $A$.
\end{enumerate}
\end{question}

\begin{question}        
\label{ex:not-ess-small}
\begin{enumerate}[(b)]
\item   
\label{part:big-objects}
Let $\cat{A}$ be a category.  Suppose there exists a functor $U\from
\cat{A} \to \Set$ such that $U$ has a left adjoint and for at least one $A
\in \cat{A}$, the set $U(A)$ has at least two elements.  Prove that for any
set $I$ and any family $(A_i)_{i \in I}$ of objects of $\cat{A}$, there is
some object of $\cat{A}$ not isomorphic to $A_i$ for any $i \in I$.  (Hint:
use Exercise~\ref{ex:eta-inj}.)

\item 
Let $\cat{A}$ be a category satisfying the assumption
of~\bref{part:big-objects}.  Prove that $\cat{A}$ is not essentially
small.

\item 
Deduce that none of the categories $\Set$, $\Vect_k$, $\Grp$, $\Ab$,
$\Ring$, and $\Tp$ is essentially small.
\end{enumerate}
\end{question}

\begin{question}
Which of the following categories are small?  Which are locally small? 
\begin{enumerate}[(b)]
\item 
$\Mon$, the category of monoids;

\item 
$\integers$, the group of integers, viewed as a one-object category;

\item 
$\integers$, the ordered set of integers;

\item 
$\Cat$, the category of small categories;

\item 
the multiplicative monoid of cardinals.
\end{enumerate}
\end{question}

\begin{question}        
\label{ex:cdoi}
\index{category!category of categories!adjunctions with $\Set$}%
\index{category!discrete}
Let $O\from \Cat \to \Set$ be the functor sending a small category to its
set%
\index{object!set of category@-set of category}
of objects.  Exhibit a chain of adjoints $C \ladj D \ladj O \ladj I$.
\end{question}

\section{Historical remarks}
\label{sec:set-history}

\index{set!history|(}
The set theory that we began to develop in Section~\ref{sec:Set-properties}
is rather different from what many mathematicians think of as set theory.
Here I will explain what the socially dominant version of set theory is,
why, despite its dominance, it is the object of widespread suspicion, and
why the kind of set theory outlined here is a more accurate reflection of
how mathematicians use sets in practice.

\paragraph*{Cantor's set theory} The creation of set theory is generally
credited to the German mathematician Georg Cantor,%
\index{Cantor, Georg}
in the late nineteenth century.  Previously, sets had seldom been regarded
as entities worthy of study in their own right; but Cantor, originally
motivated by a problem in Fourier%
\index{Fourier analysis}
analysis, developed an extensive theory.  Among many other things, he
showed that there are different sizes of infinity, proving, for instance,
that there is no bijection between $\nat$ and $\reals$.

Cantor's theory met all the resistance that typically greets a really new
idea.  His work was criticized as nonsensical, as meaningless, as far too
abstract; then later, as all very well but of no use to the mainstream of
mathematics.  Kronecker,%
\index{Kronecker, Leopold}
an important mathematician of the day, called him a charlatan and a
corrupter of youth.  But nowadays, the basics of Cantor's work are on
nearly every undergraduate mathematics syllabus.

Times change.  In the modern style of mathematics, almost every definition,
when unravelled sufficiently, depends on the notion of set.  But
pre-Cantor, this was not so.  It is interesting to try to understand the
outlook of mathematicians of the time, who had successfully developed
sophisticated subjects such as complex analysis and Galois theory without
depending on this notion that we now regard as fundamental.

Before continuing with the history, we need to discuss another fundamental
concept.

\paragraph*{Types}  
\index{type|(}
Suppose someone asks you `is $\sqrt{2} = \pi$?'  Your answer is, of course,
`no'.  Now suppose someone asks you `is $\sqrt{2} = \log$?'  You might
frown and wonder if you had heard right, and perhaps your answer would
again be `no'; but it would be a different kind of `no'.  After all,
$\sqrt{2}$ is a number, whereas $\log$ is a function, so it is
inconceivable that they could be equal.  A better answer would be `your
question makes no sense'.

This illustrates the idea of \demph{types}.  The square root of $2$ is a
real number, $\rationals$ is a field, $S_3$ is a group, $\log$ is a
function from $(0, \infty)$ to $\reals$, and $\frac{d}{dx}$ is an operation
that takes as input one function from $\reals$ to $\reals$ and produces as
output another such function.  One says that the type of $\sqrt{2}$ is
`real number', the type of $\rationals$ is `field', and so on.  We all have
an inbuilt sense of type, and it would not usually occur to us to ask
whether two things of different type were equal.

You may have met this idea before if you have programmed computers.%
\index{computer science}
Many programming languages require you to declare the type of a variable
before you first use it.  For example, you might declare that $x$ is to be
a variable of type `real number', $n$ a variable of type `integer', $M$ a
variable of type `$3 \times 3$ matrix of lists of binary digits', and so
on.

The distinction between different types of object has always been
instinctively understood.  At the beginning of the twentieth century, however,
events took a strange turn.

\paragraph*{Membership-based set theory} 
\index{set!axiomatization of sets|(}
Those who came after Cantor sought to compile a definitive list of
assumptions to be made about sets: an \emph{axiomatization} of set theory.
The list they arrived at, in the early years of the twentieth century, is known
as ZFC%
\index{ZFC (Zermelo--Fraenkel with choice)|(}
(Zermelo--Fraenkel with Choice).  It soon became the standard, and
it is the only kind of axiomatic set theory that most present-day
mathematicians know.

The axiomatization of Zermelo et al.\ was in some ways similar to the one
that we were working towards in the first section of this chapter.  But
there is at least one crucial difference: whereas we took sets and
\emph{functions} as our basic concepts, they took sets and
\emph{membership}.

At first sight, this difference might seem mild.  But when the
membership-based approach is used as a foundation on which to build the
rest of mathematics, several bizarre features become apparent:
\begin{itemize}
\item 
In the Zermelo approach, \emph{everything} is a set.  For instance, a
function is defined as a set with certain properties.  Many other things
that you would not think of as being sets are, nevertheless, treated as
sets: the number $\sqrt{2}$ is a set, the function $\log$ is a set, the
operator $\frac{d}{dx}$ is a set, and so on.

You might wonder how this is possible.  Perhaps it is useful to compare
data storage in a computer,%
\index{computer science}
where files of all different types (text, sound, images, and so on) are
ultimately encoded as sequences of $0$s and $1$s.  To give an example, in
the membership-based set theory presented in most books, the number $4$ is
encoded as the set
\[
\{ \emptyset, \{ \emptyset \}, \{ \emptyset, \{ \emptyset \} \}, 
\{ \emptyset, \{ \emptyset \}, \{ \emptyset, \{ \emptyset \} \} \} \}.
\]

\item 
The virtue of this approach is its simplicity: \emph{everything} is a set!
But the price to be paid is very high: we lose the fundamental notion of
type, precisely because everything is regarded as being of type `set'.  

\item 
In the Zermelo approach, the elements of sets are always sets too.  This
is in conflict%
\index{set!conflicting meaning in ZFC}
with ordinary mathematics.  For instance, in ordinary mathematics, $\reals$
is certainly a set, but real numbers themselves are not regarded as sets.
(After all, what is an element of $\pi$?)

\item
In this approach, membership is a global relation, meaning that for
\emph{any} two sets $A$ and $B$, it makes sense to ask whether $A \in B$.
Since this approach views everything as a set, it makes sense to ask such
apparently nonsensical questions as `is $\rationals \in \sqrt{2}$?'

Further still, the axioms of ZFC imply that we can form the intersection $A
\cap B$ of \emph{any} sets $A$ and $B$.  (Its elements are those sets $C$
for which $C \in A$ and $C \in B$.)  This makes possible further
nonsensical questions such as `does the cyclic group of order $10$ have
nonempty intersection with $\integers$?'

The answers to these nonsensical questions depend on the fine detail of how
mathematical objects (numbers, functions, groups, etc.)\ are encoded as
sets.  Even devotees of the membership-based approach agree that this
encoding is a matter of convention, just like a word processor's encoding
of a document as a string of $0$s and $1$s.  So the answers to these
questions are meaningless.
\end{itemize}

\paragraph*{Set theory today}
It should now be apparent why many modern-day mathematicians are suspicious
of set theory.  However often they are told that it is `the foundation%
\index{foundations}
of mathematics', they feel that much of it is irrelevant to their concerns.

To some extent, this is justified.  But it is also a symptom of the
historical dominance of membership-based set theory: most mathematicians do
not realize that there is any other kind.  This is a shame.  Taking sets
and functions (rather than sets and membership) as the basic concepts leads
to a theory containing all of the meaningful results of Cantor and others,
but with none of the aspects that seem so remote from the rest of
mathematics.  In particular, the function-based approach respects the
fundamental notion of type.%
\index{type|)}

The function-based approach is, of course, categorical, and its advantages
are related to more general points about how mathematics looks through
categorical eyes.  Objects are understood through their place in the
ambient category.  We get inside an object by probing%
\index{object!probing of}
it with maps to or from other objects.  For example, an element of a set
$A$ is a map $1 \to A$, and a subset of $A$ is a map $A \to 2$.  Probing of
this kind is the main theme of the next chapter.

\paragraph*{Footnote for those familiar with ZFC}  
People brought up on traditional axiomatic set theory often have the
following concern when they come across categorical set theory for the
first time.  The objects and maps of a category form a collection of some
kind, perhaps a set, so the notion of category appears to depend on some
prior set-like notion.  How, then, can sets be axiomatized categorically?
Is that not circular?

It is not, because sets can be axiomatized categorically without mentioning
categories once.  To see how, let us first recall the shape of the ZFC
axiomatization of sets.  Informally, it looks like this:
\begin{itemize}
\item 
there are some things called sets;

\item 
there is a binary relation on sets, called membership ($\in$);

\item 
some axioms hold.
\end{itemize}
A categorical axiomatization of sets looks, informally, like this:
\begin{itemize}
\item 
there are some things called sets;

\item 
for each set $A$ and set $B$, there are some things called functions from
$A$ to $B$;

\item 
to each function $f$ from $A$ to $B$ and function $g$ from $B$ to $C$,
there is assigned a function $g\of f$ from $A$ to $C$;

\item 
some axioms hold.
\end{itemize}
Making precise such phrases as `some things' requires delicacy, as will be
familiar to anyone who has done a logic course.  But the difficulties are
no worse for categorical axiomatizations of sets than for membership-based
axiomatizations such as ZFC.

One popular choice of categorical axioms for set theory can be summarized
informally as follows.

\begin{tabular}{@{}r@{\ \ }l}
\ \\
1. &Composition of functions is associative and has identities.\\
2. &There is a terminal set.\\
3. &There is a set with no elements.\\
4. &A function is determined by its effect on elements.\\
5. &Given sets $A$ and $B$, one can form their product $A \times B$.\\
6. &Given sets $A$ and $B$, one can form the set of functions from $A$ to
$B$.\\ 
7. &Given $f\from A \to B$ and $b \in B$, one can form the inverse image
$f^{-1}\{b\}$.\\ 
8. &The subsets of a set $A$ correspond to the functions from $A$ to $\{0,
1\}$.\\ 
9. &The natural numbers form a set.\\
10. &Every surjection has a section.\\
\ \\
\end{tabular}

This informal summary uses terms such as `element' and `inverse image',
which can be defined in terms of the basic concepts of set, function and
composition.  For instance, an element of a set $A$ is defined as a map from
the terminal set to $A$.  

It is certainly \emph{convenient} to express these axioms in terms of
categories.  For example, the first axiom says that sets and functions form
a category, and all ten together can be expressed in categorical jargon as
`sets and functions form a well-pointed topos%
\index{topos}%
\index{set!category of sets!topos@as topos}
with natural numbers object and choice'.  But in order to state the axioms,
it is not \emph{necessary} to appeal to any general notion of category.
They can be expressed directly in terms of sets and functions.  For
details, see \citeLR\ or \citeRST.%
\index{set!history|)}%
\index{set!axiomatization of sets|)}%
\index{ZFC (Zermelo--Fraenkel with choice)|)}%

\exone

\begin{question}
Choose a mathematician at random.  Ask them whether they can accurately
state any axiomatization of sets (without looking it up).  If not, ask them
what operating principles they actually use when handling sets in their
day-to-day work. 
\end{question}
%
%
%

\chapter{Representables}
\label{ch:rep}

A category is a world of objects, all looking at one another.  Each sees
the world from a different viewpoint.  

Consider, for instance, the category of topological spaces, and let us ask
how it looks when viewed from the one-point space $1$.  A map from $1$ to a
space $X$ is essentially the same thing as a point of $X$, so we might say
that $1$ `sees%
\index{functor!seeing@`seeing'}
points'.  Similarly, a map from $\reals$ to a space $X$
could reasonably be called a curve in $X$, and in this sense, $\reals$ sees
curves.

Now consider the category of groups.  A map from the infinite
cyclic group $\integers$%
\index{Z@$\integers$ (integers)!group@as group}
to a group $G$ amounts to an element of $G$.  (For given $g \in G$, there
is a unique homomorphism $\phi\from \integers \to G$ such that $\phi(1) =
g$.)  So, $\integers$ sees elements.  Similarly, if $p$ is a prime number
then the cyclic group $\integers/p\integers$ sees elements of order $1$ or
$p$.

Any ring homomorphism between fields is injective, so in the category of
fields,%
\index{field}
a map $K \to L$ is a way of realizing $L$ as an extension of $K$.  Hence
each field $K$ sees the extensions of itself.  If $K$ and $L$ are fields of
different characteristic then there are no homomorphisms between $K$ and
$L$, so the category of fields is the union of disjoint subcategories
$\Field_0$, $\Field_2$, $\Field_3$, $\Field_5$, \ldots\ consisting of the
fields of characteristics $0, 2, 3, 5$, \ldots.  Each field is blind to the
fields of different characteristic.

In the ordered set $(\reals, \mathord{\leq})$, the object $0$ sees whether
a number is nonnegative.  In other words, if $x$ is nonnegative then
there is one map $0 \to x$, and if not, there are none.

We can also ask the dual question: fixing an object of a category, what are
the maps \emph{into} it?  Let $S$ be the two-element set, for instance.
For an arbitrary set $X$, the maps from $X$ to $S$ correspond to the
subsets of $X$ (as we saw in Section~\ref{sec:Set-properties}).  Now give
$S$ the topology in which one of the singleton subsets is open but the
other is not.  For any topological space $X$, the continuous maps from $X$
into $S$ correspond to the \emph{open} subsets of $X$.

This chapter explores the theme of how each object sees and is seen by the
category in which it lives.  We are naturally led to the notion of
representable functor, which (after adjunctions) provides our second
approach to the idea of universal property.

\section{Definitions and examples}
\label{sec:rep-defns}

Fix an object $A$ of a category $\cat{A}$.  We will consider the totality
of maps out of $A$.  To each $B \in \cat{A}$, there is assigned the set
(or class) $\cat{A}(A, B)$ of maps from $A$ to $B$.  The content of the
following definition is that this assignation is functorial in $B$: any map
$B \to B'$ induces a function $\cat{A}(A, B) \to \cat{A}(A, B')$.

\begin{defn}  
\label{defn:co-rep}
Let $\cat{A}$ be a locally small category and $A \in \cat{A}$.  We define a
functor
\[
\h^A = \cat{A}(A, \dashbk)\from \cat{A} \to \Set%
\ntn{hom-out}
\]
as follows:
\begin{itemize}
\item 
for objects $B \in \cat{A}$, put $\h^A(B) = \cat{A}(A, B)$;

\item 
for maps $B \toby{g} B'$ in $\cat{A}$, define
\[
\h^A(g) = \cat{A}(A, g)\from 
\cat{A}(A, B) \to \cat{A}(A, B')
\]
by 
\[
p 
\mapsto
g \of p
\]
for all $p\from A \to B$.
\end{itemize}
\end{defn}

\begin{remarks}
\begin{enumerate}[(b)]
\item 
Recall that `locally%
\index{locally small}%
\index{category!locally small}
small' means that each class $\cat{A}(A, B)$ is in fact a set.  This
hypothesis is clearly necessary in order for the definition to make sense.

\item 
Sometimes $\h^A(g)$ is written as $g \of \dashbk$%
\ntn{of-blank}
 or $g_*$.%
\ntn{lower-star}
All three forms, as well as $\cat{A}(A, g)$, are in use.
\end{enumerate}
\end{remarks}

\begin{defn}
Let $\cat{A}$ be a locally small category.  A functor $X\from \cat{A} \to
\Set$ is \demph{representable}%
\index{functor!representable}
if $X \iso \h^A$ for some $A \in \cat{A}$.
A \demph{representation}%
\index{representation!functor@of functor}
of $X$ is a choice of an object $A \in \cat{A}$ and an isomorphism between
$\h^A$ and $X$.
\end{defn}
Representable functors are sometimes just called `representables'.  Only
set-valued%
\index{functor!set-valued}%
\index{set!valued functor@-valued functor}
functors (that is, functors with codomain $\Set$) can be representable.

\begin{example}        
\label{eg:co-reps-id-Set}
Consider $\h^1\from \Set \to \Set$, where $1$ is the one-element set.
Since a map from $1$ to a set $B$ amounts to an element of $B$, we have
\[
\h^1(B) \iso B
\]
for each $B \in \Set$.  It is easily verified that this isomorphism is
natural in $B$, so $\h^1$ is isomorphic to the identity functor $1_\Set$.
Hence $1_\Set$ is representable.
\end{example}

\begin{example}
\label{eg:co-reps-seeing}
All of the `seeing'%
\index{functor!seeing@`seeing'}
functors in the introduction to this chapter are representable.  The
forgetful%
\index{functor!forgetful!representable@is representable}
functor $\Tp \to \Set$ is isomorphic to $\h^1 = \Tp(1, \dashbk)$, and the
forgetful functor $\Grp \to \Set$ is isomorphic to $\Grp(\integers,
\dashbk)$.  For each prime $p$, there is a functor $U_p\from \Grp \to \Set$
defined on objects by
\[
U_p(G) = 
\{\text{elements of }G \text{ of order } 1 \text{ or } p \},
\index{group!order of element of}
\]
and as claimed above, $U_p \iso \Grp(\integers/p\integers, \dashbk)$
(Exercise~\ref{ex:cyclic-rep}).  Hence $U_p$ is representable.
\end{example}

\begin{example}
There is a functor $\ob\from \Cat \to \Set$ sending a small category to its
set%
\index{object!set of category@-set of category}
of objects.  (The category $\Cat$ was introduced in
Definition~\ref{defn:Cat}.) It is representable.  Indeed, consider the
terminal category $\One$ (with one object and only the identity map).  A
functor from $\One$ to a category $\cat{B}$ simply picks out an object of
$\cat{B}$.  Thus,
\[
\h^\One(\cat{B}) \iso \ob\cat{B}.
\]
Again, it is easily verified that this isomorphism is natural in $\cat{B}$;
hence $\ob \iso \Cat(\One, \dashbk)$.  It can be shown similarly that the
functor $\Cat \to \Set$ sending a small category to its set of maps is
representable (Exercise~\ref{ex:arrows-rep}).  
\end{example}

\begin{example}
Let $M$ be a monoid, regarded as a one-object%
\index{monoid!action of}
category.  Recall from Example~\ref{eg:functor-action} that a set-valued
functor on $M$ is just an $M$-set.  Since the category $M$ has only one
object, there is only one representable functor on it (up to isomorphism).
As an $M$-set, the unique representable is the so-called \demph{left
regular%
\index{representation!group or monoid@of group or monoid!regular}
representation} of $M$, that is, the underlying set of $M$ acted on by
multiplication on the left.
\end{example}

\begin{example}
Let $\Toph_*$%
\ntn{Toph-star}
be the category whose objects are topological spa\-ces equipped with a
basepoint and whose arrows are homotopy%
\index{homotopy}
classes of basepoint-preserving continuous maps.  Let $S^1 \in \Toph_*$%
\ntn{circle}
be the circle.  Then for any object $X \in \Toph_*$, the maps $S^1 \to X$
in $\Toph_*$ are the elements of the fundamental%
\index{group!fundamental}
group $\pi_1(X)$.  Formally, this says that the composite functor
\[
\Toph_* \toby{\pi_1} \Grp \toby{U} \Set
\]
is isomorphic to $\Toph_*(S^1, \dashbk)$.  In particular, it is
representable.
\end{example}

\begin{example}
\label{eg:co-reps-tensor}
Fix a field $k$ and vector spaces $U$ and $V$ over $k$.  There is a
functor 
\[
\Bilin(U, V; \dashbk)\from \Vect_k \to \Set%
\index{map!bilinear}
\ntn{Bilin}
\]
whose value $\Bilin(U, V; W)$ at $W \in \Vect_k$ is the set of bilinear
maps $U \times V \to W$.  It can be shown that this functor is
representable; in other words, there is a space $T$ with the property that
\[
\Bilin(U, V; W) \iso \Vect_k(T, W)
\]
naturally in $W$.  This $T$ is the tensor%
\index{tensor product}
product $U \otimes V$, which we met just after the proof of
Lemma~\ref{lemma:tensor-unique}.
\end{example}

Adjunctions give rise to representable%
\index{functor!representable!adjoints@and adjoints}
functors in the following way.

\begin{lemma}   
\label{lemma:adj-to-rep}
Let $\hadjnli{\cat{A}}{\cat{B}}{F}{G}$ be locally small categories, and let
$A \in \cat{A}$.  Then the functor
\[
\cat{A}(A, G(\dashbk))\from \cat{B} \to \Set
\]
(that is, the composite $\cat{B} \toby{G} \cat{A} \toby{\h^A} \Set$) is
representable.  
\end{lemma}

\begin{pf}
We have
\[
\cat{A}(A, G(B)) \iso \cat{B}(F(A), B)
\]
for each $B \in \cat{B}$.  If we can show that this isomorphism is natural
in $B$, then we will have proved that $\cat{A}(A, G(\dashbk))$ is isomorphic
to $\h^{F(A)}$ and is therefore representable.  So, let $B \toby{q} B'$ be
a map in $\cat{B}$.  We must show that the square
\[
\xymatrix{
\cat{A}(A, G(B)) \ar[r] \ar[d]_{G(q) \of \dashbk}       &
\cat{B}(F(A), B) \ar[d]^{q \of \dashbk}        \\
\cat{A}(A, G(B')) \ar[r]        &
\cat{B}(F(A), B')
}
\]
commutes, where the horizontal arrows are the bijections provided by the
adjunction.  For $f \from A \to G(B)$, we have
\[
\xymatrix@C+3em{
f \ar@{|->}[r] \ar@{|->}[d]    &
\bar{f} \ar@{|->}[d]    \\
G(q) \of f \ar@{|->}[r]    &
*!<0mm,-1.4ex>+\txt{$q \of \bar{f}$\\
$\ovln{G(q) \of f}$,}
}
\]
so we must prove that $q \of \bar{f} = \ovln{G(q) \of f}$.  This follows
immediately from the naturality condition~\eqref{eq:adj-nat-a} in the
definition of adjunction (with $g = \bar{f}$). 
\end{pf}

You would not expect a randomly-chosen functor into $\Set$ to be
rep\-re\-sen\-table.  In some sense, rather few functors are.  However,
forgetful%
\index{functor!forgetful!representable@is representable|(}
functors do tend to be representable:
\begin{propn}   
\label{propn:ladj-rep}
\hspace*{-2.5pt}Any set-valued functor with a left adjoint is representable.
\end{propn}

\begin{pf}
Let $G\from \cat{A} \to \Set$ be a functor with a left adjoint $F$.  Write
$1$ for the one-point set.  Then
\[
G(A) \iso \Set(1, G(A))
\]
naturally in $A \in \cat{A}$
(by Example~\ref{eg:co-reps-id-Set}), that is, $G \iso
\Set(1, G(\dashbk))$.  So by Lemma~\ref{lemma:adj-to-rep}, $G$ is
representable; indeed, $G \iso \h^{F(1)}$.  
\end{pf}

\begin{example}        
Several of the examples of representables mentioned above arise as in
Proposition~\ref{propn:ladj-rep}.  For instance, $U \from \Tp \to \Set$
has a left adjoint $D$%
\index{topological space!discrete}
(Example~\ref{eg:adjn:spaces}), and $D(1) \iso 1$, so we recover the result
that $U \iso \h^1$.  Similarly, Exercise~\ref{ex:cdoi} asked you to
construct a left adjoint $D$%
\index{category!discrete}
to the objects functor $\ob\from \Cat \to \Set$.  This functor $D$
satisfies $D(1) \iso \One$, proving again that $\ob \iso \h^\One$.
\end{example}

\begin{example}
The forgetful functor $U\from \Vect_k \to \Set$ is representable,\linebreak
since it has a left adjoint.%
\index{vector space!free}
Indeed, if $F$ denotes the left adjoint then $F(1)$ is the $1$-dimensional
vector space $k$, so $U \iso \h^k$.  This is also easy to see directly: a
map from $k$ to a vector space $V$ is uniquely determined by the image of
$1$, which can be any element of $V$; hence $\Vect_k(k, V) \iso U(V)$
naturally in $V$.
\end{example}

\begin{example}
\label{eg:ladj-rep-ring}
Examples~\ref{egs:adjns-alg} began with the declaration that forgetful%
\index{functor!forgetful!left adjoint to}
functors between categories of algebraic structures usually have left
adjoints.  Take the category $\CRing$ of commutative rings%
\index{ring!free}
and the forgetful functor $U\from \CRing \to \Set$.  This general principle
suggests that $U$ has a left adjoint, and Proposition~\ref{propn:ladj-rep}
then tells us that $U$ is representable.

Let us see how this works explicitly.  Given a set $S$, let $\integers[S]$
be the ring of polynomials%
\index{ring!polynomial}
over $\integers$ in commuting variables $x_s$ ($s \in S$).  (This was
called $F(S)$ in Example~\ref{egs:free-functors}\bref{eg:free-ring}.)  Then
$S \mapsto \integers[S]$ defines a functor $\Set \to \CRing$, and this is
left adjoint to $U$.  Hence $U \iso \h^{\integers[x]}$.  Again, this
can be verified directly: for any ring $R$, the maps $\integers[x] \to R$
correspond one-to-one with the elements of $R$ (Exercises~\ref{ex:Zx}
and~\ref{ex:free-ring-one-gen}).
\index{functor!forgetful!representable@is representable|)}
\end{example}

We have defined, for each object $A$ of our category $\cat{A}$, a functor
$\h^A \in \ftrcat{\cat{A}}{\Set}$.  This describes how $A$ sees the world.
As $A$ varies, the view varies.  On the other hand, it is always the same
world being seen, so the different views from different objects are somehow
related.  (Compare aerial%
\index{aerial photography}
photos taken from a moving aeroplane, which agree well enough on their
overlaps that they can be patched together to make one big picture.)  So
the family $\bigl(\h^A\bigr)_{A \in \cat{A}}$ of `views' has some
consistency to it.  What this means is that whenever there is a map between
objects $A$ and $A'$, there is also a map between $\h^A$ and $\h^{A'}$.

Precisely, a map $A' \toby{f} A$ induces a natural transformation
\[
\xymatrix@C+1em{
\cat{A} \rtwocell<4>^{\h^A}_{\h^{A'}}{\hspace{.5em}\h^f} &\Set,
}%
\ntn{hom-out-map}
\]
whose $B$-component (for $B \in \cat{A}$) is the function  
\[
\begin{array}{ccc}
\h^A(B) = \cat{A}(A, B)  &\to            &\h^{A'}(B) = \cat{A}(A', B)     \\
p                       &\mapsto        &p \of f.
\end{array}
\]
Again, $\h^f$
goes by a variety of other names: $\cat{A}(f, \dashbk)$,%
\ntn{hom-out-map-blank}
$f^*$,%
\ntn{upper-star}
and $\dashbk \of f$.%
\ntn{blank-of}

Note the reversal of direction!  Each functor $\h^A$ is covariant, but they
come together to form a \emph{contravariant} functor, as in the following
definition. 

\begin{defn}
Let $\cat{A}$ be a locally small category.  The functor
\[
\h^\bl\from \cat{A}^\op \to \ftrcat{\cat{A}}{\Set}%
\ntn{hom-out-blank}
\]
is defined on objects $A$ by $\h^\bl(A) = \h^A$ and on maps $f$ by $\h^\bl(f)
= \h^f$. 
\end{defn}

The symbol $\bl$ is another type of blank, like $\dashbk$.  

All of the definitions presented so far in this chapter can be dualized.
At the formal level, this is trivial: reverse all the arrows, so that every
$\cat{A}$ becomes an $\cat{A}^\op$ and vice versa.  But in our usual
examples, the flavour is different.  We are no longer asking what objects
\emph{see}, but how they are \emph{seen}.

Let us first dualize Definition~\ref{defn:co-rep}.

\begin{defn}
Let $\cat{A}$ be a locally small category and $A \in \cat{A}$.  We define a
functor
\[
\h_A = \cat{A}(\dashbk, A)\from \cat{A}^\op \to \Set%
\ntn{hom-in}
\]
as follows:
\begin{itemize}
\item 
for objects $B \in \cat{A}$, put $\h_A(B) = \cat{A}(B, A)$;

\item 
for maps $B' \toby{g} B$ in $\cat{A}$, define
\[
\h_A(g) = \cat{A}(g, A) = g^* = \dashbk \of g\from 
\cat{A}(B, A) \to \cat{A}(B', A)
\]
by 
\[
p 
\mapsto
p \of g
\]
for all $p\from B \to A$.
\end{itemize}
\end{defn}

If you know about dual vector spaces, this construction will seem familiar.
In particular, you will not be surprised that a map $B' \to B$ induces a
map in the opposite direction, $\h_A(B) \to \h_A(B')$.

We now define representability for \emph{contravariant} set-valued
functors.  Stri\-ctly speaking, this is unnecessary, as a contravariant
functor on $\cat{A}$ is a covariant functor on $\cat{A}^\op$, and we
already know what it means for a covariant set-valued functor to be
representable.  But it is useful to have a direct definition.

\begin{defn}
Let $\cat{A}$ be a locally small category.  A functor $X\from \cat{A}^\op \to
\Set$ is \demph{representable}%
\index{functor!representable}
if $X \iso \h_A$ for some $A \in \cat{A}$.  A \demph{representation}%
\index{representation!functor@of functor}
of $X$ is a choice of an object $A \in \cat{A}$ and an isomorphism between
$\h_A$ and $X$.
\end{defn}

\begin{example}        
There is a functor
\[
\pset\from \Set^\op \to \Set%
\ntn{power-set-ftr}
\]
sending each set $B$ to its power%
\index{power!set}
set $\pset(B)$, and defined on maps $g\from B' \to B$ by $(\pset(g))(U) =
g^{-1}U$ for all $U \in \pset(B)$.  (Here $g^{-1}U$ denotes the inverse%
\index{inverse!image}
image or preimage of $U$ under $g$, defined by $g^{-1} U = \{ x' \in B'
\such g(x') \in U \}$.)  As we saw in Section~\ref{sec:Set-properties}, a
subset amounts to a map into the two-point%
\index{set!two-element}
set $2$.  Precisely put, $\pset \iso \h_2$.
\end{example}

\begin{example}
\label{eg:contra-rep-sier}
Similarly, there is a functor
\[
\oset\from\Tp^\op \to \Set%
\ntn{oset-ftr}
\index{topological space!open subset of}%
\index{open subset}%
\index{set!open}
\]
defined on objects $B$ by taking $\oset(B)$ to be the set of open subsets
of $B$.  If $S$ denotes the two-point%
\index{topological space!two-point}
topological space in which exactly one of the two singleton subsets is
open, then continuous maps from a space $B$ into $S$ correspond naturally
to open subsets of $B$ (Exercise~\ref{ex:sierpinski-space}).  Hence $\oset
\iso \h_S$, and $\oset$ is representable.
\end{example}

\begin{example}
In Example~\ref{eg:contra-fn-spaces}, we defined a functor $C\from \Tp^\op
\to \Ring$, assigning to each space the ring%
\index{ring!functions@of functions}
of continuous real-valued functions%
\index{topological space!functions on}
on it.  The composite functor
\[
\Tp^\op \toby{C} \Ring \toby{U} \Set
\]
is representable, since by definition, $U(C(X)) = \Tp(X, \reals)$ for
topological spaces $X$.
\end{example}

Previously, we assembled the covariant representables $\bigl( \h^A
\bigr)_{A \in \cat{A}}$ into one big functor $\h^\bl$.  We now do the same
for the contravariant representables $\bigl( \h_A \bigr)_{A \in \cat{A}}$.
Any map $A \toby{f} A'$ in $\cat{A}$ induces a natural transformation
\[
\xymatrix@C+1em{
\cat{A}^\op \rtwocell<4>^{\h_A}_{\h_{A'}}{\hspace{.5em}\h_f} &\Set%
\ntn{hom-in-map}
}
\]
(also called $\cat{A}(\dashbk, f)$,%
\ntn{hom-in-map-blank}
$f_*$%
\ntn{lower-star-bis}
or $f \of \dashbk$),%
\ntn{of-blank-bis}
whose component at an object $B \in \cat{A}$ is
\[
\begin{array}{ccc}
\h_A(B) = \cat{A}(B, A)  &\to               &
\h_{A'}(B) = \cat{A}(B, A')     \\
p                       &\mapsto            &
f \of p.
\end{array}
\]

\begin{defn}    
\label{defn:yon-emb}
Let $\cat{A}$ be a locally small category.  The \demph{Yoneda%
\index{Yoneda embedding}
embedding} of $\cat{A}$ is the functor
\[
\h_\bl\from \cat{A} \to \ftrcat{\cat{A}^\op}{\Set}%
\ntn{hom-in-blank}
\]
defined on objects $A$ by $\h_\bl(A) = \h_A$ and on maps $f$ by $\h_\bl(f) =
\h_f$. 
\end{defn}

Here is a summary of the definitions so far.
\begin{center}
\begin{tabular}{ll}
\ \\[-2ex]
For each $A \in \cat{A}$, we have a functor     &
$\cat{A} \toby{\h^A} \Set$. \\
Putting them all together gives a functor       &
$\cat{A}^\op \toby{\h^\bl} \ftrcat{\cat{A}}{\Set}$.
\vspace*{2ex}\\
For each $A \in \cat{A}$, we have a functor     &
$\cat{A}^\op \toby{\h_A} \Set$. \\
Putting them all together gives a functor       &
$\cat{A} \toby{\h_\bl} \ftrcat{\cat{A}^\op}{\Set}$.\\[-2ex]
\ 
\end{tabular}%
\end{center}
The second pair of functors is the dual of the first.  Both involve
contravariance;%
\index{functor!contravariant}%
\index{contravariant}
it cannot be avoided.

In the theory of representable functors, it does not make much difference
whether we work with the first or the second pair.  Any theorem that we
prove about one dualizes to give a theorem about the other.  We choose to
work with the second pair, the $\h_A$s and $\h_\bl$.  In a sense to be
explained, $\h_\bl$ `embeds' $\cat{A}$ into $\ftrcat{\cat{A}^\op}{\Set}$.
This can be useful, because the category $\ftrcat{\cat{A}^\op}{\Set}$ has
some good properties that $\cat{A}$ might not have.

Exercise~\ref{ex:yoneda-conservative} asks you to prove that $\h_\bl$ is
injective on isomorphism classes of objects.  It is strongly recommended
that you do it before reading on, as it encapsulates the key ideas of the
rest of this chapter.

There is one more functor to define.  It unifies the first and second pairs
of functors shown above.

\begin{defn}
Let $\cat{A}$ be a locally small category.  The functor
\[
\Hom_\cat{A}\from \cat{A}^\op \times \cat{A} \to \Set%
\index{hom-set}
\ntn{Hom-functor}
\]
is defined by
\[
\xymatrix{
(A, B) \ar@<1ex>[d]^g    &
\mapsto \ar@{}[d]|*+{\mapsto}       &
\cat{A}(A, B) \ar[d]^{g \of \dashbk \of f}      \\
(A', B') \ar@<1ex>[u]^f  &
\mapsto &
\cat{A}(A', B').
}
\]
In other words, $\Hom_\cat{A}(A, B) = \cat{A}(A, B)$ and $(\Hom_\cat{A}(f,
g))(p) = g \of p \of f$, whenever $A' \toby{f} A \toby{p} B \toby{g} B'$.
\end{defn}

\begin{remarks}         
\label{rmks:global-hom}
\begin{enumerate}[(b)]
\item 
The existence of the functor $\Hom_\cat{A}$ is something like the fact that for
a metric%
\index{metric space}
space $(X, d)$, the metric is itself a continuous map $d\from X \times X
\to \reals$.  (If we take two points and move each one slightly, the
distance between them changes only slightly.)

\item 
In terms of Exercise~\ref{ex:ftr-on-product}, $\Hom_{\cat{A}}$ is the
functor $\cat{A}^\op \times \cat{A} \to \Set$ corresponding to the families
of functors $\bigl( \h^A \bigr)_{A \in \cat{A}}$ and $\bigl( \h_B \bigr)_{B
  \in \cat{A}}$.

\item   
\label{rmks:global-hom:cc}
In Example~\ref{eg:adjn:cc}, we saw that for any set $B$, there is an
adjunction $(\dashbk \times B) \ladj (\dashbk)^B$ of functors $\Set \to
\Set$.  Similarly, for any category $\cat{B}$, there is an adjunction
$(\dashbk \times \cat{B}) \ladj \ftrcat{\cat{B}}{\dashbk}$ of functors
$\CAT \to \CAT$; in other words, there is a canonical bijection
\[
\CAT(\cat{A} \times \cat{B}, \cat{C})
\iso
\CAT(\cat{A}, \ftrcat{\cat{B}}{\cat{C}})
\]
for $\cat{A}, \cat{B}, \cat{C} \in \CAT$.  Under this bijection, the functors
\[
\Hom_\cat{A}\from \cat{A}^\op \times \cat{A} \to \Set,
\qquad
\h^\bl\from \cat{A}^\op \to \ftrcat{\cat{A}}{\Set}
\]
correspond to one another.  Thus, $\Hom_\cat{A}$ carries the same
information as $\h^\bl$ (or $\h_\bl$), presented slightly differently.
\end{enumerate}
\end{remarks}

\begin{remark}  
\label{rmk:adj-nat}
We can now explain the naturality%
\index{adjunction!naturality axiom for}
in the definition of adjunction (Definition~\ref{defn:adjn}).  Take
categories and functors $\oppairi{\cat{A}}{\cat{B}}{F}{G}$.  They give rise
to functors
\[
\xymatrix{
\cat{A}^\op \times \cat{B} \ar[r]^-{1\times G} \ar[d]_{F^\op \times 1} &
\cat{A}^\op \times \cat{A} \ar[d]^{\Hom_\cat{A}}        \\
\cat{B}^\op \times \cat{B} \ar[r]_-{\Hom_\cat{B}} &
\Set.
}
\]
The composite functor $\searrows$ sends $(A, B)$ to $\cat{B}(F(A), B)$; it
can be written as $\cat{B}(F(\dashbk), \dashbk)$.  The composite
$\esarrows$ sends $(A, B)$ to $\cat{A}(A, G(B))$.
Exercise~\ref{ex:adj-nat} asks you to show that these two functors
\[
\cat{B}(F(\dashbk), \dashbk), 
\ 
\cat{A}(\dashbk, G(\dashbk))
\from
\cat{A}^\op \times \cat{B} \to \Set
\]
are naturally isomorphic if and only if $F$ and $G$ are adjoint.  This
justifies the claim in Remark~\ref{rmks:adjts}\bref{rmk:adjts:nat}: the
naturality requirements~\eqref{eq:adj-nat-a} and~\eqref{eq:adj-nat-b} in
the definition of adjunction simply assert that two particular functors
are naturally isomorphic.
\end{remark}

Objects of an arbitrary category do not have elements in any obvious sense.
However, \emph{sets} certainly have elements, and we have observed that an
element of a set $A$ is the same thing as a map $1 \to A$.  This inspires
the following definition.

\begin{defn}    
\label{defn:gen-elt}
Let $A$ be an object of a category.  A \demph{generalized%
\index{element!generalized}
element} of $A$ is a map with codomain $A$.  A map $S \to A$ is a
generalized element of $A$ of \demph{shape%
\index{shape!generalized element@of generalized element}
$S$}.
\end{defn}

`Generalized element' is nothing more than a synonym of `map', but
sometimes it is useful to think of maps as generalized elements.

For example, when $A$ is a set, a generalized element of $A$ of shape $1$
is an ordinary element of $A$, and a generalized element of $A$ of shape
$\nat$ is a sequence%
\index{sequence}
in $A$.  In the category of topological spaces, the generalized elements of
shape $1$ (the one-point space) are the points, and the generalized
elements of shape $S^1$ (the circle) are, by definition, loops.%
\index{loop|(}
As this suggests, in categories of geometric objects, we might equally well
say `figures of shape $S$'.

In algebra,%
\index{algebra}
we are often interested in solutions to equations such as $x^2 + y^2 = 1$.
Perhaps we begin by being particularly interested in solutions in
$\rationals$, but then realize that in order to study rational solutions,
it will be helpful to study solutions in other rings first.  (This is often
a fruitful strategy.)  Given a ring $A$, a pair $(a, b) \in A \times A$
satisfying $a^2 + b^2 = 1$ amounts to a homomorphism of rings
\[
\integers[x, y]/(x^2 + y^2 - 1) \to A.
\index{algebraic geometry}
\]
Thus, the solutions to our equation (in any ring) can be seen as the
generalized elements of shape $ \integers[x, y]/(x^2 + y^2 - 1)$.

For an object $S$ of a category $\cat{A}$, the functor
\[
\h^S\from \cat{A} \to \Set
\]
sends an object to its set of generalized elements of shape $S$.  The
functoriality tells us that any map $A \to B$ in $\cat{A}$ transforms
$S$-elements of $A$ into $S$-elements of $B$.  For example, taking $\cat{A}
= \Tp$ and $S = S^1$, any continuous map $A \to B$ transforms loops in $A$
into loops in $B$.%
\index{loop|)}

\exs

\begin{question}
Find three examples of representable functors not mentioned above.
\end{question}

\begin{question}        
\label{ex:yoneda-conservative}
Let $\cat{A}$ be a locally small category, and let $A, A' \in \cat{A}$ with
$\h_A \iso \h_{A'}$.  Prove directly that $A \iso A'$.
\end{question}

\begin{question}        
\label{ex:cyclic-rep}
Let $p$ be a prime number.  Show that the functor $U_p\from \Grp \to \Set$
defined in Example~\ref{eg:co-reps-seeing} is isomorphic to
$\Grp(\integers/p\integers, \dashbk)$.  (To check that there is an
isomorphism of functors~-- that is, a \emph{natural} isomorphism~-- you
will first need to define $U_p$ on maps.  There is only one sensible way to
do this.)
\end{question}

\begin{question}        
\label{ex:free-ring-one-gen}
Using the result of Exercise~\ref{ex:Zx}\bref{part:Zx-main}, prove that the
forgetful functor $\CRing \to \Set$ is isomorphic to $\CRing(\integers[x],
\dashbk)$, as in Example~\ref{eg:ladj-rep-ring}.
\end{question}

\begin{question}        
\label{ex:sierpinski-space}
The \demph{Sierpi\'nski%
\index{Sierpi\'nski space}
space} is the two-point topological space $S$ in
which one of the singleton subsets is open but the other is not.  Prove
that for any topological space $X$, there is a canonical bijection between
the open subsets of $X$ and the continuous maps $X \to S$.  Use this to
show that the functor $\oset\from \Tp^\op \to \Set$ of
Example~\ref{eg:contra-rep-sier} is represented by $S$.
\end{question}

\begin{question}        
\label{ex:arrows-rep}
Let $M \from \Cat \to \Set$ be the functor that sends a small category
$\cat{A}$ to the set of all maps in $\cat{A}$.  Prove that $M$ is
representable.
\end{question}

\begin{question}        
\label{ex:adj-nat}
Take locally small categories $\cat{A}$ and $\cat{B}$, and functors
$\oppairi{\cat{A}}{\cat{B}}{F}{G}$.  Show that $F$ is left adjoint to $G$
if and only if the two functors 
\[
\cat{B}(F(\dashbk), \dashbk), 
\ 
\cat{A}(\dashbk, G(\dashbk))
\from
\cat{A}^\op \times \cat{B} \to \Set
\]
of Remark~\ref{rmk:adj-nat} are naturally isomorphic.  (Hint: this
  is made easier by using either Exercise~\ref{ex:nat-iso-on-product} or
  Exercise~\ref{ex:adj-nat-in-one}.)
\end{question}

\section{The Yoneda lemma}

What do representables see?

Recall from Definition~\ref{defn:presheaf} that functors $\cat{A}^\op \to
\Set$ are sometimes called `presheaves' on $\cat{A}$.  So for each $A \in
\cat{A}$ we have a representable presheaf $\h_A$, and we are asking how the
rest of the presheaf category $\pshf{\cat{A}}$ looks from the viewpoint of
$\h_A$.  In other words, if $X$ is another presheaf, what are the maps
$\h_A \to X$?

Newcomers to category theory commonly find that the material presented in
this section is where they first get stuck.  Typically, the core of the
difficulty is in understanding the question just asked.  Let us ask it
again. 

We start by fixing a locally small category $\cat{A}$.  We then take an
object $A \in \cat{A}$ and a functor $X\from \cat{A}^\op \to \Set$.  The
object $A$ gives rise to another functor $\h_A = \cat{A}(\dashbk, A)\from
\cat{A}^\op \to \Set$.  The question is: what are the maps $\h_A \to X$?
Since $\h_A$ and $X$ are both objects of the presheaf category
$\pshf{\cat{A}}$, the `maps' concerned are maps in $\pshf{\cat{A}}$.  So,
we are asking what natural transformations
\begin{equation}        
\label{eq:yoneda-transf}
\begin{array}{c}
\xymatrix{
\cat{A}^\op \rtwocell^{\h_A}_{X} &\Set
}
\end{array}
\end{equation}
there are.  The set of such natural transformations is called
\[
\pshf{\cat{A}}(\h_A, X).
\]
(This is a special case of the notation $\cat{B}(B, B')$ for the set of
maps $B \to B'$ in a category $\cat{B}$.  Here, $\cat{B} = \pshf{\cat{A}}$,
$B = \h_A$, and $B' = X$.)  We want to know what this set is.

There is an informal principle%
\index{uniqueness!constructions@of constructions|(}
of general category theory that allows us to
guess the answer.  Look back at
Remarks~\ref{rmks:defn-cat}\bref{rmk:defn-cat:loosely},
\ref{rmks:defn-ftr}\bref{rmk:defn-ftr:loosely}
and~\ref{rmks:defn-nt}\bref{rmk:defn-nt:loosely} on the definitions of
category, functor and natural transformation.  Each remark is of the form
`from input of one type, it is possible to construct exactly one output of
another type'.  For example, in
Remark~\ref{rmks:defn-cat}\bref{rmk:defn-cat:loosely}, the input is a
sequence of maps $A_0 \toby{f_1} \cdots \toby{f_n} A_n$, the output is a
map $A_0 \to A_n$, and the statement is that no matter what we do with the
input data $f_1, \ldots, f_n$, there is only one map $A_0 \to A_n$ that we
can construct.

Let us apply this principle to our question.  We have just seen how, given
as input an object $A \in \cat{A}$ and a presheaf $X$ on $\cat{A}$, we can
construct a set, namely, $\pshf{\cat{A}}(\h_A, X)$.  Are there any other
ways to construct a set from the same input data $(A, X)$?  Yes: simply
take the set $X(A)$!  The informal principle suggests that these two sets
are the same:
\begin{equation}        
\label{eq:yoneda-pre}
\pshf{\cat{A}}(\h_A, X) \iso X(A)
\end{equation}
for all $A \in \cat{A}$ and $X \in \pshf{\cat{A}}$.  This turns out to be
true; and that is the Yoneda lemma.%
\index{uniqueness!constructions@of constructions|)}

Informally, then, the Yoneda lemma says that for any $A \in \cat{A}$ and
presheaf $X$ on $\cat{A}$:
\begin{slogan}
A natural transformation $\h_A \to X$ is an element of $X(A)$.
\end{slogan}
Here is the formal statement.  The proof follows shortly.

\begin{thm}[Yoneda]   
\label{thm:yoneda}
\index{Yoneda lemma}
Let $\cat{A}$ be a locally small category.  Then
\begin{equation}        
\label{eq:yoneda}
\pshf{\cat{A}}(\h_A, X)
\iso
X(A)
\end{equation}
naturally in $A \in \cat{A}$ and $X \in \pshf{\cat{A}}$.  
\end{thm}

This is exactly what was stated in~\eqref{eq:yoneda-pre}, except that the
word `naturally' has appeared.  Recall from Definition~\ref{defn:nat-in}
that for functors $F, G\from \cat{C} \to \cat{D}$, the phrase `$F(C) \iso
G(C)$ naturally in $C$' means that there is a natural isomorphism $F \iso
G$.  So the use of this phrase in the Yoneda lemma suggests that each side
of~\eqref{eq:yoneda} is functorial in both $A$ and $X$.  This means, for
instance, that a map $X \to X'$ must induce a map
\[
\pshf{\cat{A}}(\h_A, X) 
\to
\pshf{\cat{A}}(\h_A, X'),
\]
and that not only does the isomorphism~\eqref{eq:yoneda} hold for
\emph{every} $A$ and $X$, but also, the isomorphisms can be chosen in a way
that is compatible with these induced maps.  Precisely, the Yoneda lemma
states that the composite functor
\[
\xymatrix@R=0em{
\cat{A}^\op \times \pshf{\cat{A}}  \ar[r]^-{\h_\bl^\op \times 1} &
{\pshf{\cat{A}}}^\op \times \pshf{\cat{A}} 
\ar[r]^-{\Hom_{\pshf{\cat{A}}}} &
\Set    \\
(A, X) \ar@{}[r]|(.43)*+{\longmapsto} &
(\h_A, X) \ar@{}[r]|(.55)*+{\longmapsto} &
\pshf{\cat{A}}(\h_A, X)
}
\]
is naturally isomorphic to the evaluation%
\index{evaluation}
 functor
\[
\begin{array}{ccc}
\cat{A}^\op \times \pshf{\cat{A}}       &\to            &\Set   \\
(A, X)                                  &\mapsto        &X(A).
\end{array}
\]

If the Yoneda lemma were false then the world would look much more complex.
For take a presheaf $X\from \cat{A}^\op \to \Set$, and define a new
presheaf $X'$ by
\[
X' = \pshf{\cat{A}}(\h_\bl, X)
\from 
\cat{A}^\op \to \Set,
\]
that is, $X'(A) = \pshf{\cat{A}}(\h_A, X)$ for all $A \in \cat{A}$.  Yoneda
tells us that $X'(A) \iso X(A)$ naturally in $A$; in other words, $X' \iso
X$.  If Yoneda were false then starting from a single presheaf $X$, we
could build an infinite sequence $X, X', X'', \ldots$ of new presheaves,
potentially all different.  But in reality, the situation is very simple:
they are all the same.

The proof of the Yoneda lemma is the longest proof so far.  Nevertheless,
there is essentially only one way to proceed at each stage.  If you suspect
that you are one of those newcomers to category theory for whom the Yoneda
lemma presents the first serious challenge, an excellent exercise is to
work out the proof before reading it.  No ingenuity is required, only an
understanding of all the terms in the statement.

\begin{pfof}{the Yoneda lemma}
We have to define, for each $A$ and $X$, a bijection between the sets
$\pshf{\cat{A}}(\h_A, X)$ and $X(A)$.  We then have to show that our bijection
is natural in $A$ and $X$.
\pagebreak 

First, fix $A \in \cat{A}$ and $X \in \pshf{\cat{A}}$.  We define functions
\begin{equation}        
\label{eq:yoneda-fns}
\oppair{\pshf{\cat{A}}(\h_A, X)}{X(A)}{\yel{\blank}}{\ynt{\blank}}
\end{equation}
and show that they are mutually inverse.  So we have to do four things: define
the function $\yel{\blank}$, define the function $\ynt{\blank}$, show that
$\yel{\ynt{\blank}}$ is the identity, and show that $\ynt{\yel{\blank}}$ is
the identity.
\begin{itemize}
\item 
Given $\alpha\from \h_A \to X$, define $\yel{\alpha} \in X(A)$ by
$\yel{\alpha} = \alpha_A(1_A)$.  (How else could we possibly define it?)

\item 
Let $x \in X(A)$.  We have to define a natural transformation $\ynt{x}\from
\h_A \to X$.  That is, we have to define for each $B \in \cat{A}$ a
function
\[
\ynt{x}_B\from \h_A(B) = \cat{A}(B, A) \to X(B)
\]
and show that the family $\ynt{x} = (\ynt{x}_B)_{B \in \cat{A}}$ satisfies
naturality.  

Given $B \in \cat{A}$ and $f \in \cat{A}(B, A)$, define
\[
\ynt{x}_B(f) = (X(f))(x) \in X(B).
\]
(How else could we possibly define it?)  This makes sense, since $X(f)$ is a
map $X(A) \to X(B)$.  To prove naturality, we must show that for any map
$B' \toby{g} B$ in $\cat{A}$, the square
\[
\xymatrix@C+4em@R+1ex{
\cat{A}(B, A) \ar[r]^{\h_A(g)\: =\: \dashbk\of g} \ar[d]_{\ynt{x}_B}        &
\cat{A}(B', A) \ar[d]^{\ynt{x}_{B'}} \\
X(B) \ar[r]_{X(g)}        &
X(B')
}
\]
commutes.  To reduce clutter, let us write $X(g)$ as $Xg$, and so on.  Now
for all $f \in \cat{A}(B, A)$, we have
\[
\xymatrix@C+3em{
f \ar@{|->}[r] \ar@{|->}[d]    &
f \of g \ar@{|->}[d]    \\
(Xf)(x) \ar@{|->}[r]    &
*!<0mm,-1.4ex>+\txt{$(X(f\of g))(x)$\\
$(Xg)((Xf)(x))$,}
}
\]
and $X(f \of g) = (Xg) \of (Xf)$ by functoriality, so the square does commute.

\item 
Given $x \in X(A)$, we have to show that $\yel{\ynt{x}} = x$, and indeed,
\[
\yel{\ynt{x}}
=
\ynt{x}_A(1_A)
=
(X 1_A)(x)
=
1_{X(A)}(x)
=
x.
\]

\item 
Given $\alpha\from \h_A \to X$, we have to show that $\ynt{\yel{\alpha}} =
\alpha$.  Two natural transformations are equal if and only if all their
components are equal; so, we have to show that
$\Bigl(\ynt{\yel{\alpha}}\Bigr)_B = \alpha_B$ for all $B \in \cat{A}$.
Each side of this equation is a function from $\h_A(B) = \cat{A}(B, A)$ to
$X(B)$, and two functions are equal if and only if they take equal values
at every element of the domain; so, we have to show that
\[
\Bigl(\ynt{\yel{\alpha}}\Bigr)_B (f)
= 
\alpha_B(f)
\]
for all $B \in \cat{A}$ and $f\from B \to A$ in $\cat{A}$.  The left-hand side
is by definition
\[
\Bigl(\ynt{\yel{\alpha}}\Bigr)_B (f)
=
(Xf)(\yel{\alpha})
=
(Xf)(\alpha_A(1_A)),
\]
so it remains to prove that
\begin{equation}        
\label{eq:recovery}
(Xf)(\alpha_A(1_A)) = \alpha_B(f).  
\end{equation}
By naturality of $\alpha$ (the only tool at our disposal), the square
\[
\xymatrix@C+4em@R+1ex{
\cat{A}(A, A) \ar[r]^{\h_A(f) \:=\: \dashbk \of f} 
\ar[d]_{\alpha_A}    &
\cat{A}(B, A) \ar[d]^{\alpha_B} \\
X(A) \ar[r]_{Xf}        &
X(B)
}
\]
commutes, which when taken at $1_A \in \cat{A}(A, A)$ gives
equation~\eqref{eq:recovery}.   
\end{itemize}

(The proof is not over yet, but it is worth pausing to consider the
significance of the fact that $\ynt{\yel{\alpha}} = \alpha$.  Since
$\yel{\alpha}$ is the value of $\alpha$ at $1_A$, this implies:
\begin{slogan}
A natural transformation $\h_A \to X$ is determined by its
value at $1_A$.
\end{slogan}
Just \emph{how} a natural transformation $\h_A \to X$ is determined by its
value at $1_A$ is described in equation~\eqref{eq:recovery}.)

This establishes the bijection~\eqref{eq:yoneda-fns} for each $A \in
\cat{A}$ and $X \in \ftrcat{\cat{A}^\op}{\Set}$.  We now show that the
bijection is natural in $A$ and $X$.  

We employ two mildly labour-saving devices.  First, in principle we have
to prove naturality of both $\yel{\blank}$ and $\ynt{\blank}$, but by
Lemma~\ref{lemma:nat-iso-compts}, it is enough to prove naturality of just
one of them.  We prove naturality of $\yel{\blank}$.  Second, by
Exercise~\ref{ex:nat-iso-on-product}, $\yel{\blank}$ is natural in the pair
$(A, X)$ if and only if it is natural in $A$ for each fixed $X$ and natural
in $X$ for each fixed $A$.  So, it remains to check these two types of
naturality.  

Naturality in $A$ states that for each $X \in \pshf{\cat{A}}$ and
$B \toby{f} A$ in $\cat{A}$, the square
\[
\xymatrix@C+2em{
\pshf{\cat{A}}(\h_A, X) \ar[r]^{\dashbk \of \h_f} \ar[d]_{\yel{\blank}} &
\pshf{\cat{A}}(\h_B, X) \ar[d]^{\yel{\blank}}   \\
X(A) \ar[r]_{Xf}        &
X(B)
}
\]
commutes.  For $\alpha\from \h_A \to X$, we have
\[
\xymatrix@C+4em{
\alpha \ar@{|->}[r] \ar@{|->}[d]    &
\alpha \of \h_f \ar@{|->}[d]    \\
\alpha_A(1_A) \ar@{|->}[r]    &
*!<0mm,-1.4ex>+\txt{$\bigl(\alpha \of \h_f\bigr)_B(1_B)$\\
$(Xf)(\alpha_A(1_A))$,}
}
\]
so we have to show that $\bigl(\alpha \of \h_f\bigr)_B(1_B)=
(Xf)(\alpha_A(1_A))$.  Indeed,
\begin{align*}
\bigl( \alpha \of \h_f \bigr)_B(1_B)    &
=
\alpha_B ((\h_f)_B (1_B)) \\
&
=
\alpha_B (f \of 1_B)    
=
\alpha_B(f) \\
&
=
(Xf)(\alpha_A(1_A)),
\end{align*}
where the first step is by definition of composition in $\pshf{\cat{A}}$,
the second is by definition of $\h_f$, and the last is by
equation~\eqref{eq:recovery}.

Naturality in $X$ states that for each $A \in \cat{A}$ and map
\[
\xymatrix{
\cat{A}^\op \rtwocell^{X}_{X'}{\theta} &\Set
}
\]
in $\pshf{\cat{A}}$, the square
\[
\xymatrix@C+2em{
\pshf{\cat{A}}(\h_A, X)  \ar[r]^{\theta \of \dashbk} \ar[d]_{\yel{\blank}} &
\pshf{\cat{A}}(\h_A, X') \ar[d]^{\yel{\blank}}  \\
X(A) \ar[r]_{\theta_A}  &
X'(A)
}
\]
commutes.  For $\alpha\from \h_A \to X$, we have
\[
\xymatrix@C+4em{
\alpha \ar@{|->}[r] \ar@{|->}[d]    &
\theta \of \alpha \ar@{|->}[d]    \\
\alpha_A(1_A) \ar@{|->}[r]    &
*!<0mm,-1.4ex>+\txt{$(\theta \of \alpha)_A(1_A)$\\
$\theta_A(\alpha_A(1_A))$,}
}
\]
and $(\theta\of\alpha)_A = \theta_A \of \alpha_A$ by definition of composition
in $\pshf{\cat{A}}$, so the square does commute.  This completes the proof.
\end{pfof}

\exs

\begin{question}
State the dual of the Yoneda lemma.
\end{question}

\begin{question}
One way to understand the Yoneda lemma is to examine some special
cases.  Here we consider one-object%
\index{monoid!Yoneda lemma for monoids}%
\index{Yoneda lemma!monoids@for monoids}
categories.

Let $M$ be a monoid.  The underlying set of $M$ can be given a right
$M$-action by multiplication: $x \cdot m = xm$ for all $x, m \in M$.  This
$M$-set is called the \demph{right regular%
\index{representation!group or monoid@of group or monoid!regular} 
representation} of $M$.  Let us write it as $\rreg{M}$.

\begin{enumerate}[(b)]
\item 
When $M$ is regarded as a one-object category, functors $M^\op \to \Set$
correspond to right $M$-sets (Example~\ref{eg:contra-functors:actions}).
Show that the $M$-set corresponding to the unique representable functor
$M^\op \to \Set$ is the right regular representation.

\item 
Now let $X$ be any right $M$-set.  Show that for each $x \in X$, there is a
unique map $\alpha\from \rreg{M} \to X$ of right $M$-sets such that
$\alpha(1) = x$.  Deduce that there is a bijection between $ \{ \text{maps
} \rreg{M} \to X \text{ of right }M\text{-sets} \} $ and $X$.

\item 
Deduce the Yoneda lemma for one-object categories.  
\end{enumerate}
\end{question}

\section{Consequences of the Yoneda lemma}

The Yoneda lemma is fundamental in category theory.  Here we look at three
important consequences.

\begin{notn}
An arrow decorated with a $\sim$, as in $A \toiso B$,%
\ntn{iso-arrow}
 denotes an isomorphism.
\end{notn}

\minihead{A representation is a universal element}
\index{representation!functor@of functor!universal element@as universal
  element|(} 

\begin{cor}     
\label{cor:rep-univ}
Let $\cat{A}$ be a locally small category and $X\from \cat{A}^\op \to \Set$.
Then a representation of $X$ consists of an object $A \in \cat{A}$ together
with an element $u \in X(A)$ such that:
\begin{equation}        
\label{eq:univ-elt}
\parbox{.8\textwidth}{%
for each $B \in \cat{A}$ and $x \in X(B)$, there is a unique map $\bar{x}\from
B \to A$ such that $(X\bar{x})(u) = x$.%
}
\end{equation}
\end{cor}
To clarify the statement, first recall that by definition, a representation
of $X$ is an object $A \in \cat{A}$ together with a natural isomorphism
$\alpha\from \h_A \toiso X$.  Corollary~\ref{cor:rep-univ} states that such
pairs $(A, \alpha)$ are in natural bijection with pairs $(A, u)$ satisfying
condition~\eqref{eq:univ-elt}.

Pairs $(B, x)$ with $B \in \cat{A}$ and $x \in X(B)$ are sometimes called
\demph{elements}%
\index{element!presheaf@of presheaf}%
\index{presheaf!element of}
of the presheaf $X$.  (Indeed, the Yoneda lemma tells us that $x$ amounts
to a generalized element of $X$ of shape $\h_B$.)  An element $u$
satisfying condition~\eqref{eq:univ-elt} is sometimes called a
\demph{universal}%
\index{element!universal}%
\index{universal!element}
element of $X$.  So, Corollary~\ref{cor:rep-univ} says that a
representation of a presheaf $X$ amounts to a universal element of $X$.

\begin{pf}
By the Yoneda lemma, we have only to show that for $A \in \cat{A}$ and $u
\in X(A)$, the natural transformation $\ynt{u}\from \h_A \to X$ is an
isomorphism if and only if~\eqref{eq:univ-elt} holds.  (Here we are using
the notation introduced in the proof of the Yoneda lemma.)  Now, $\ynt{u}$
is an isomorphism if and only if for all $B \in \cat{A}$, the function
\[
\ynt{u}_B\from \h_A(B) = \cat{A}(B, A) \to X(B)
\]
is a bijection, if and only if for all $B \in \cat{A}$ and $x \in X(B)$, there
is a unique $\bar{x} \in \cat{A}(B, A)$ such that $\ynt{u}_B(\bar{x}) = x$.
But $\ynt{u}_B(\bar{x}) = (X\bar{x})(u)$, so this is exactly
condition~\eqref{eq:univ-elt}.  
\end{pf}

Our examples will use the dual form, for covariant set-valued functors:

\begin{cor}     
\label{cor:rep-univ-dual}
Let $\cat{A}$ be a locally small category and $X\from \cat{A} \to \Set$.  Then
a representation of $X$ consists of an object $A \in \cat{A}$ together with an
element $u \in X(A)$ such that:
\begin{equation}        
\label{eq:univ-elt-dual}
\parbox{.8\textwidth}{%
for each $B \in \cat{A}$ and $x \in X(B)$, there is a unique map $\bar{x}\from
A \to B$ such that $(X\bar{x})(u) = x$.%
}
\end{equation}
\end{cor}

\begin{pf}
Follows immediately by duality.
\end{pf}

\begin{example}
Fix a set $S$ and consider the functor
\[
\begin{array}{cccc}
X = \Set(S, U(\dashbk))\from            &
\Vect_k &\to            &\Set           \\
                                        &
V       &\mapsto        &\Set(S, U(V)).
\end{array}
\]
Here are two familiar (and true!)\ statements about $X$:
\begin{enumerate}[(b)]
\item   
\label{item:vs-rep-direct}
there exist a vector space $F(S)$ and an isomorphism%
\index{unit and counit!unit as initial object}%
\index{adjunction!initial objects@via initial objects}
\begin{equation}        
\label{eq:vect-adj}
\Vect_k(F(S), V) \iso \Set(S, U(V))
\end{equation}
natural in $V \in \Vect_k$
(Example~\ref{egs:adjns-alg}\bref{eg:adjns-alg:vs});

\item   
\label{item:vs-rep-univ}
there exist a vector space $F(S)$ and a function $u\from S \to U(F(S))$%
\index{vector space!free!unit of}
such that:
\begin{displaytext}
for each vector space $V$ and function $f\from S \to U(V)$, there is a unique
linear map $\bar{f}\from F(S) \to V$ such that
\[
\xymatrix{
S \ar[r]^-u \ar[dr]_f    &U(F(S)) \ar[d]^{U(\bar{f})}    \\
        &
U(V)
}
\]
commutes
\end{displaytext}
(as in the introduction to Section~\ref{sec:adj-init}, where $u$ was called by
its usual name, $\eta_S$).
\end{enumerate}

Each of these two statements says that $X$ is representable.
Statement~\bref{item:vs-rep-direct} says that there is an isomorphism
$X(V) \iso \Vect(F(S), V)$ natural in $V$, that is, an isomorphism $X \iso
\h^{F(S)}$.  So $X$ is representable, by definition of representability.
Statement~\bref{item:vs-rep-univ} says that $u \in X(F(S))$ satisfies
condition~\eqref{eq:univ-elt-dual}.  So $X$ is representable, by
Corollary~\ref{cor:rep-univ-dual}.

You will have noticed that the first way of saying that $X$ is representable
is substantially shorter than the second.  Indeed, it is clear that if the
situation of~\bref{item:vs-rep-univ} holds then there is an isomorphism
\[
\Vect_k(F(S), V) \toiso \Set(S, U(V))
\]
natural in $V$, defined by $g \mapsto U(g) \of u$.  But it looks at first
as if~\bref{item:vs-rep-univ} says rather more
than~\bref{item:vs-rep-direct}, since it states that the two functors are
not only naturally isomorphic, but naturally isomorphic in a rather special
way.  Corollary~\ref{cor:rep-univ-dual} tells us that this is an illusion:
all natural isomorphisms~\eqref{eq:vect-adj} arise in this way.  It is the
word `natural'%
\index{adjunction!naturality axiom for}
in~\bref{item:vs-rep-direct} that hides the explicit detail.
\end{example}

\begin{example}
The same can be said for any other adjunction
$\hadjnli{\cat{A}}{\cat{B}}{F}{G}$.  Fix $A \in \cat{A}$ and put
\[
X = \cat{A}(A, G(\dashbk))\from
\cat{B} \to \Set.
\]
Then $X$ is representable, and this can be expressed in either of the
following ways:
\begin{enumerate}[(b)]
\item 
$\cat{A}(A, G(B)) \iso \cat{B}(F(A), B)$ naturally in $B$; in other words,
$X \iso \h^{F(A)}$ (as in Lemma~\ref{lemma:adj-to-rep});

\item 
the unit map $\eta_A\from A \to G(F(A))$ is an initial object of the comma
category $\comma{A}{G}$; that is, $\eta_A \in X(F(A))$ satisfies
condition~\eqref{eq:univ-elt-dual}. 
\end{enumerate}
This observation can be developed into an alternative proof of
Theorem~\ref{thm:adj-comma}, the reformulation of adjointness%
\index{adjunction!initial objects@via initial objects}
in terms of initial objects.
\end{example}

\begin{example} 
\label{eg:gp-fgt-rep}
For any group%
\index{group}
$G$ and element $x \in G$, there is a unique homomorphism $\phi\from
\integers \to G$%
\index{Z@$\integers$ (integers)!group@as group}
 such that $\phi(1) = x$.  This means that $1 \in
U(\integers)$ is a universal element of the forgetful functor $U\from \Grp
\to \Set$; in other words, condition~\eqref{eq:univ-elt-dual} holds when
$\cat{A} = \Grp$, $X = U$, $A = \integers$ and $u = 1$.  So $1 \in
U(\integers)$ gives a representation $\h^\integers \toiso U$ of $U$.

On the other hand, the same is true with $-1$ in place of $1$.  The
isomorphisms $\h^\integers \toiso U$ coming from $1$ and $-1$ are not
equal, because Corollary~\ref{cor:rep-univ-dual} provides a \emph{one-to-one}
correspondence between universal elements and representations.
\end{example}
\index{representation!functor@of functor!universal element@as universal
  element|)} 

\minihead{The Yoneda embedding}
\index{Yoneda embedding|(}

Here is a second corollary of the Yoneda lemma.

\begin{cor}     
\label{cor:yoneda-ff}
For any locally small category $\cat{A}$, the Yoneda embedding
\[
\h_\bl\from \cat{A} \to \pshf{\cat{A}}
\]
is full and faithful.
\end{cor}

Informally, this says that for $A, A' \in \cat{A}$, a map $\h_A \to \h_{A'}$ of
presheaves is the same thing as a map $A \to A'$ in $\cat{A}$.

\begin{pf}
We have to show that for each $A, A' \in \cat{A}$, the function
\begin{equation}        
\label{eq:yoneda-emb-1}
\begin{array}{ccc}
\cat{A}(A, A')  &\to            &\pshf{\cat{A}}(\h_A, \h_{A'})    \\
f               &\mapsto        &\h_f
\end{array}
\end{equation}
is bijective.  By the Yoneda lemma (taking `$X$' to be $\h_{A'}$), the function
\begin{equation}        
\label{eq:yoneda-emb-2}
\ynt{\blank}\from
\h_{A'}(A) \to \pshf{\cat{A}}(\h_A, \h_{A'})
\end{equation}
is bijective, so it is enough to prove that the
functions~\eqref{eq:yoneda-emb-1} and~\eqref{eq:yoneda-emb-2} are equal.
Thus, given $f\from A \to A'$, we have to prove that $\ynt{f} = \h_f$, or
equivalently, $\wideyel{\h_f} = f$.  And indeed,
\[
\wideyel{\h_f}
=
(\h_f)_A(1_A)
=
f \of 1_A 
=
f,
\]
as required.
\end{pf}

In mathematics at large, the word `embedding'%
\index{embedding}
is used (sometimes informally) to mean a map $A \to B$ that makes $A$
isomorphic to its image in $B$.  For example, an injection of sets $i\from
A \to B$ might be called an embedding, because it provides a bijection
between $A$ and the subset $iA$ of $B$.  Similarly, a map $i\from A \to B$
of topological spaces might be called an embedding if it is a homeomorphism
to its image, so that $A \iso iA$.  Corollary~\ref{cor:ff-emb} tells us
that in category theory, a full and faithful functor $\cat{A} \to \cat{B}$
can reasonably be called an embedding, as it makes $\cat{A}$ equivalent to
a full subcategory of $\cat{B}$.

\begin{figure}
\centering
\setlength{\unitlength}{.6ex}
\begin{picture}(46,30)(-15,-15)
\cell{0}{0}{c}{\includegraphics[height=30\unitlength]{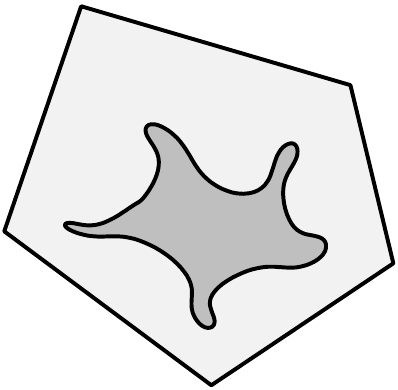}}
\cell{15}{1}{l}{\pshf{\cat{A}}}
\cell{1}{-2.5}{c}{\cat{A}}
\end{picture}
\caption{A category $\cat{A}$ embedded into its presheaf category.}
\label{fig:yoneda-embedding}
\end{figure}

In the case at hand, the Yoneda embedding $\h_\bl\from \cat{A} \to
\pshf{\cat{A}}$ embeds $\cat{A}$ into its own presheaf category
(Figure~\ref{fig:yoneda-embedding}). 
So, $\cat{A}$ is equivalent to the full subcategory of $\pshf{\cat{A}}$
whose objects are the representables.

In general, full%
\index{subcategory!full}
subcategories are the easiest subcategories to handle.  For instance, given
objects $A$ and $A'$ of a full subcategory, we can speak unambiguously of
the `maps' from $A$ to $A'$; it makes no difference whether this is
understood to mean maps in the subcategory or maps in the whole category.
Similarly, we can speak unambiguously of isomorphism of objects of the
subcategory, as in the following lemma.

\begin{lemma}   
\label{lemma:ff-refl-isos}
\index{functor!full and faithful}%
\index{isomorphism!full and faithful functors@and full and faithful functors}
Let $J\from \cat{A} \to \cat{B}$ be a full and faithful functor and $A, A' \in
\cat{A}$.  Then: 
\begin{enumerate}[(b)]
\item 
a map $f$ in $\cat{A}$ is an isomorphism if and only if the map $J(f)$ in
$\cat{B}$ is an isomorphism;

\item 
\label{lemma:ff-refl-isos:medium} 
for any isomorphism $g\from J(A) \to J(A')$ in $\cat{B}$, there is a unique
isomorphism $f\from A \to A'$ in $\cat{A}$ such that $J(f) = g$;

\item   
\label{lemma:ff-refl-isos:weakest} 
the objects $A$ and $A'$ of $\cat{A}$ are isomorphic if and only if the
objects $J(A)$ and $J(A')$ of $\cat{B}$ are isomorphic.
\end{enumerate}
\end{lemma}

\begin{pf}
Exercise~\ref{ex:ff-refl-isos}.
\end{pf}

\begin{example}
In Example~\ref{eg:gp-fgt-rep}, we considered the representations of the
forgetful functor $U\from \Grp \to \Set$,%
\index{group}
and found two different isomorphisms\linebreak $\h^\integers \toiso U$.%
\index{Z@$\integers$ (integers)!group@as group}
Did we find all of them?

Since $\h^\integers \iso U$, there are as many isomorphisms $\h^\integers
\toiso U$ as there are isomorphisms $\h^\integers \toiso \h^\integers$.  By
Corollary~\ref{cor:yoneda-ff} and
Lemma~\ref{lemma:ff-refl-isos}\bref{lemma:ff-refl-isos:medium}, there are
as many of \emph{these} as there are group isomorphisms $\integers \toiso
\integers$.  There are precisely two such (corresponding to the two
generators $\pm 1$ of $\integers$), so we did indeed find all the
isomorphisms $\h^\integers \toiso U$.  Differently put, there are exactly
two universal elements of $U(\integers)$.
\end{example}

In Section~\ref{sec:lim-pshf}, we will see that every presheaf can be built
from representables, in very roughly the same way that every positive
integer can be built from primes.%
\index{Yoneda embedding|)}

\minihead{Isomorphism of representables}
\index{functor!representable!isomorphism of representables|(}

In Exercise~\ref{ex:yoneda-conservative}, you were asked to prove directly
that if $\h_A \iso \h_{A'}$ then $A \iso A'$.  The proof contains all the
main ideas in the proof of the Yoneda lemma.  The result itself can also be
deduced from the Yoneda lemma, as follows.

\begin{cor}     
\label{cor:reps-unique}
Let $\cat{A}$ be a locally small category and $A, A' \in \cat{A}$.  Then
\[
\h_A \iso \h_{A'}
\iff
A \iso A'
\iff
\h^A \iso \h^{A'}.
\]
\end{cor}

\begin{pf}
By duality, it is enough to prove the first `$\textiff$'.  This follows from 
Corollary~\ref{cor:yoneda-ff} and
Lemma~\ref{lemma:ff-refl-isos}\bref{lemma:ff-refl-isos:weakest}.  
\end{pf}

\begin{figure}
\centering
\setlength{\unitlength}{.6ex}
\begin{picture}(39.5,39.5)(-21,-21)
\cell{-1}{-1}{c}{\includegraphics[width=40\unitlength]{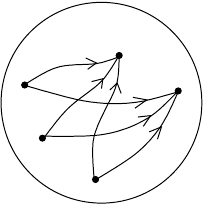}}
\cell{-16}{-1}{c}{B_1}
\cell{-10}{-11}{c}{B_2}
\cell{0.5}{-18}{c}{B_3}
\cell{3}{11}{c}{A}
\cell{14.5}{4}{c}{A'}
\cell{16}{16}{c}{\cat{A}}
\end{picture}%
\caption{If $\cat{A}(B, A) \iso \cat{A}(B, A')$ naturally in $B$, then $A
  \iso A'$.}
\label{fig:view}
\end{figure}

Since functors always preserve isomorphism
(Exercise~\ref{ex:ftrs-pres-iso}), the force of this statement is that
\[
\h_A \iso \h_{A'} \implies A \iso A'.  
\]
In other words, if $\cat{A}(B, A) \iso \cat{A}(B, A')$ naturally in $B$,
then $A \iso A'$.  Thinking of $\cat{A}(B, A)$ as `$A$ viewed from $B$',
the corollary tells us that two objects are the same if and only if they
look the same from all viewpoints (Figure~\ref{fig:view}).  (If it looks
like a duck,%
\index{duck}
walks like a duck, and quacks like a duck, then it probably is a duck.)

\begin{example} 
Consider Corollary~\ref{cor:reps-unique} in the case $\cat{A} = \Grp$.
Take two groups $A$ and $A'$, and suppose someone tells us that $A$ and
$A'$ `look the same from $B$' (meaning that $\h_A(B) \iso \h_{A'}(B)$) for
all groups $B$.  Then, for instance:
\begin{itemize}
\item 
$\h_A(1) \iso \h_{A'}(1)$, where $1$ is the trivial group.  But $\h_A(1) =
\Grp(1, A)$ is a one-element set, as is $\h_{A'}(1)$, no matter what $A$
and $A'$ are.  So this tells us nothing at all.

\item 
$\h_A(\integers) \iso \h_{A'}(\integers)$.  We know that $\h_A(\integers)$
is the underlying set of $A$, and similarly for $A'$.  So $A$ and $A'$
have isomorphic underlying sets.  But for all we know so far, they might
have entirely different group structures.

\item 
$\h_A(\integers/p\integers) \iso \h_{A'}(\integers/p\integers)$ for every
prime $p$, so by Example~\ref{eg:co-reps-seeing}, $A$ and $A'$ have the
same number of elements of each prime order.%
\index{group!order of element of}
\end{itemize}
Each of these isomorphisms gives only partial information about the
similarity of $A$ and $A'$.  But if we know that $\h_A(B) \iso \h_{A'}(B)$
for all groups $B$, and \emph{naturally} in $B$, then $A \iso A'$.
\end{example}

\begin{example}
The category of sets is very unusual in this respect.  For any set $A$, we have
\[
A \iso \Set(1, A) = \h_A(1),
\]
so $\h_A(1) \iso \h_{A'}(1)$ implies $A \iso A'$.  In other words, two
objects of $\Set$ are the same if they look the same from the point of view
of the one-element set.  This is a familiar feature of sets: the only thing
that matters about a set is its elements!

For a general category, Corollary~\ref{cor:reps-unique} tells us that two
objects are the same if they have the same generalized%
\index{element!generalized}
elements of all shapes.  But the category of sets has a special property:
if I choose an object and tell you only what its generalized elements of
shape $1$ are, then you can deduce exactly what my object must be.
\end{example}

\begin{example} 
\label{eg:yon-adjts-unique}
Let $G\from \cat{B} \to \cat{A}$ be a functor, and suppose that both $F$
and $F'$ are left adjoint to $G$.  Then for each $A \in \cat{A}$, we have
\[
\cat{B}(F(A), B)
\iso
\cat{A}(A, G(B))
\iso 
\cat{B}(F'(A), B)
\]
naturally in $B \in \cat{B}$, so $\h^{F(A)} \iso \h^{F'(A)}$, so $F(A) \iso
F'(A)$ by Corollary~\ref{cor:reps-unique}.  In fact, this isomorphism is
natural in $A$, so that $F \iso F'$.  This shows that left adjoints are
unique, as claimed in Remark~\ref{rmks:adjts}\bref{rmks:adjts:uniqueness}.
Dually, right adjoints are unique.  See also Exercise~\ref{ex:adjts-unique}.
\end{example}

\begin{example}
Corollary~\ref{cor:reps-unique} implies that if a set-valued functor is
isomorphic to both $\h^A$ and $\h^{A'}$ then $A \iso A'$.  So the functor
\emph{determines} the representing object, if one exists.  For instance, take
the functor 
\[
\Bilin(U, V; \dashbk)\from \Vect_k \to \Set
\index{map!bilinear}
\]
of Example~\ref{eg:co-reps-tensor}.  Corollary~\ref{cor:reps-unique}
implies that up to isomorphism, there is \emph{at most one} vector space
$T$ such that
\[
\Bilin(U, V; W) \iso \Vect_k(T, W)
\]
naturally in $W$.  It can be shown that there does, in fact, exist such a
vector space $T$.  Since all such spaces $T$ are isomorphic, it is
legitimate to refer to any of them as \emph{the}%
\index{uniqueness}
tensor%
\index{tensor product}
product of $U$ and $V$.
\end{example}
\index{functor!representable!isomorphism of representables|)}

\exs

\begin{question}        
\label{ex:ff-refl-isos}
Prove Lemma~\ref{lemma:ff-refl-isos}.
\end{question}

\begin{question}
Let $\cat{A}$ be a locally small category.  Prove each of the following
statements directly (without using the Yoneda lemma).
\begin{enumerate}[(b)]
\item 
$\h_\bl\from \cat{A} \to \pshf{\cat{A}}$ is faithful.

\item 
$\h_\bl$ is full.

\item
Given $A \in \cat{A}$ and a presheaf $X$ on $\cat{A}$, if $X(A)$ has an
element $u$ that is universal in the sense of Corollary~\ref{cor:rep-univ},
then $X \iso \h_A$.
\end{enumerate}
\end{question}

\begin{question}
Interpret the theory of Chapter~\ref{ch:rep} in the case where the category
$\cat{A}$ is discrete.  For example, what do presheaves look like, and
which ones are representable?  What does the Yoneda lemma tell us?  Does
its proof become any shorter?  What about the corollaries of the Yoneda
lemma?
\end{question}

\begin{question}        
\label{ex:adjts-unique}
Let $\cat{B}$ be a category and $J\from \cat{C} \to \cat{D}$ a functor.
There is an induced functor
\[
J \of \dashbk\from \ftrcat{\cat{B}}{\cat{C}} \to \ftrcat{\cat{B}}{\cat{D}}
\]
defined by composition with $J$.  
\begin{enumerate}[(b)]
\item 
Show that if $J$ is full and faithful then so is $J\of\dashbk$.

\item 
Deduce that if $J$ is full and faithful and $G, G'\from \cat{B} \to
\cat{C}$ with $J \of G \iso J \of G'$ then $G \iso G'$.

\item 
Now deduce that right adjoints are unique:%
\index{adjunction!uniqueness of adjoints}
if $F\from \cat{A} \to \cat{B}$ and $G, G'\from \cat{B} \to \cat{A}$ with
$F \ladj G$ and $F \ladj G'$ then $G \iso G'$.  (Hint: the Yoneda embedding
is full and faithful.)
\end{enumerate}
\end{question}

%
%
%

\chapter{Limits}
\label{ch:lims}

Limits, and the dual concept, colimits, provide our third approach to the
idea of universal property.  

Adjointness is about the relationships \emph{between} categories.
Representability is a property of \emph{set-valued} functors.  Limits are
about what goes on \emph{inside} a category.

The concept of limit unifies many familiar constructions in mathematics.
Whenever you meet a method for taking some objects and maps in a category
and constructing a new object out of them, there is a good chance that you
are looking at either a limit or a colimit.  For instance, in group theory,
we can take a homomorphism between two groups and form its kernel, which
is a new group.  This construction is an example of a limit in the category
of groups.  Or, we might take two natural numbers and form their lowest
common multiple.  This is an example of a colimit in the poset of natural
numbers, ordered by divisibility.

\section{Limits: definition and examples}
\label{sec:lims-basics}

The definition of limit is very general.  We build up to it by first
examining some particularly useful types of limit: products, equalizers,
and pullbacks.

\minihead{Products}

Let $X$ and $Y$ be sets.  The familiar cartesian product%
\index{set!category of sets!products in}
$X \times Y$ is characterized by the property that an element of $X \times
Y$ is an element of $X$ together with an element of $Y$.  Since elements
are just maps from $1$, this says that a map $1 \to X \times Y$ amounts to
a map $1 \to X$ together with a map $1 \to Y$.

A little thought reveals that the same is true when $1$ is replaced
throughout by any set $A$ whatsoever.  (In other words, a generalized
element of $X \times Y$ of shape $A$ amounts to a generalized element of
$X$ of shape $A$ together with a generalized element of $Y$ of shape $A$.)
The bijection between
\[
\text{maps } A \to X \times Y
\]
and
\[
\text{pairs of maps } (A \to X,\ A \to Y)
\]
is given by composing with the projection maps
\[
\begin{array}{ccccc}
X       &\otby{p_1}     &X \times Y     &\toby{p_2}     &Y      \\
x       &\mapsfrom      &(x, y)         &\mapsto        &y.
\end{array}
\]
This suggests the following definition.

\begin{defn}    
\label{defn:bin-prod}
Let $\cat{A}$ be a category and $X, Y \in \cat{A}$.  A \demph{product}%
\index{product}
of $X$ and $Y$ consists of an object $P$ and maps
\[
\xymatrix{
        &P \ar[ld]_{p_1} \ar[rd]^{p_2}  &       \\
X       &                               &Y
}
\]
with the property that for all objects and maps
\begin{equation}        
\label{eq:bin-prod-cone}
\begin{array}{c}
\xymatrix{
        &A \ar[ld]_{f_1} \ar[rd]^{f_2}  &       \\
X       &                               &Y
}
\end{array}
\end{equation}
in $\cat{A}$, there exists a unique map $\bar{f}\from A \to P$ such that 
\begin{equation}        
\label{eq:bin-prod-lim}
\begin{array}{c}
\xymatrix{
        &A \ar[ldd]_{f_1} \ar@{.>}[d]|{\bar{f}\vphantom{\bar{\bar{f}}}} 
\ar[rdd]^{f_2}&       \\
        &P \ar[ld]^{p_1} \ar[rd]_{p_2}                  &       \\
X       &                                               &Y
}
\end{array}
\end{equation}
commutes.  The maps $p_1$ and $p_2$ are called the \demph{projections}.%
\index{projection}
\end{defn}

\begin{remarks} 
\label{rmks:prod}
\begin{enumerate}[(b)]
\item   
\label{rmks:prod:exist}
Products do not always exist.  For example, if $\cat{A}$ is the discrete
two-object category
\[
\fbox{$X\bullet$ \hspace*{2em} $\bullet Y$}
\]
then $X$ and $Y$ do not have a product.  But when objects $X$ and $Y$ of a
category do have a product, it is unique%
\index{product!uniqueness of} 
up to isomorphism.  (This can be proved directly, much as in
Lemma~\ref{lemma:init-unique}.  It also follows from
Corollary~\ref{cor:lims-unique}.)  This justifies talking about \emph{the}
product of $X$ and $Y$.

\item 
Strictly speaking, the product consists of the object $P$ \emph{together
  with} the projections $p_1$ and $p_2$.  But informally,%
\index{product!informal usage}
we often refer to $P$ alone as the product of $X$ and $Y$.  We write $P$ as
$X \times Y$.%
\ntn{prod-gen}
\end{enumerate}
\end{remarks}

\begin{example}
\label{eg:sets}
Any two sets $X$ and $Y$ have a product%
\index{set!category of sets!products in}
in $\Set$.  It is the usual cartesian product $X \times Y$, equipped with
the usual projection maps $p_1$ and $p_2$.

Let us check that this really is a product in the sense of
Definition~\ref{defn:bin-prod}.  Take sets and functions as in
diagram~\eqref{eq:bin-prod-cone}.  Define $\bar{f}\from A \to X \times Y$
by $\bar{f}(a) = (f_1(a), f_2(a))$.  Then $p_i \of \bar{f} = f_i$ for $i =
1, 2$; that is, diagram~\eqref{eq:bin-prod-lim} commutes with $P = X \times
Y$.  Moreover, this is the \emph{only} map making
diagram~\eqref{eq:bin-prod-lim} commute.  For suppose that $\hat{f}\from A
\to X \times Y$, in place of $\bar{f}$, also makes~\eqref{eq:bin-prod-lim}
commute.  Let $a \in A$, and write $\hat{f}(a)$ as $(x, y)$.  Then
\[
f_1(a) = p_1(\hat{f}(a)) = p_1(x, y) = x,
\]
and similarly, $f_2(a) = y$.  Hence $\hat{f}(a) = (f_1(a), f_2(a)) =
\bar{f}(a)$ for all $a \in A$, giving $\hat{f} = \bar{f}$, as required.
\end{example}

In general, in any category, the map $\bar{f}$ of
diagram~\eqref{eq:bin-prod-lim} is usually written as $(f_1, f_2)$.

\begin{example}
\label{eg:prod-spaces}
In the category of topological spaces, any two objects $X$ and $Y$ have a
product.%
\index{topological space!category of topological spaces!products in}
It is the set $X \times Y$ equipped with the product topology and the
standard projection maps.  The product topology is deliberately designed so
that a function
\[
\begin{array}{ccc}
A       &\to            &X \times Y     \\
t       &\mapsto        &(x(t), y(t))
\end{array}
\]
is continuous if and only if it is continuous in each coordinate (that is to
say, both functions
\[
t \mapsto x(t),
\qquad
t \mapsto y(t)
\]
are continuous).  This holds for any space $A$, but the idea is perhaps at
its most intuitively appealing when $A = \reals$ and we think of $t$ as a
time parameter.

A closely related statement is that the product topology is the smallest
topology on $X \times Y$ for which the projections are continuous.  Here
`smallest' means that for any other topology $\mathcal{T}$ on $X \times Y$
such that $p_1$ and $p_2$ are continuous, every subset of $X \times Y$ open
in the product topology is also open in $\mathcal{T}$.  Thus, to define the
product topology, we declare just enough sets to be open that the
projections are continuous.
\end{example}

\begin{example}
\label{eg:prod-vs}
Now let $X$ and $Y$ be vector spaces.  We can form their direct sum,%
\index{vector space!direct sum of vector spaces}%
\index{vector space!category of vector spaces!products in}
$X \oplus Y$,%
\ntn{direct-sum}
whose elements can be written as either $(x, y)$ or $x + y$ (with $x \in
X$ and $y \in Y$), according to taste.  There are linear projection maps
\[
\xymatrix{
        &X \oplus Y \ar[ld]_{p_1} \ar[rd]^{p_2}     &       \\
X       &                                               &Y
}
\qquad
\xymatrix{
        &(x, y) \ar@{|->}[ld] \ar@{|->}[rd]     &       \\
x       &                                       &y.
}
\]
It can be shown that $X \oplus Y$, together with $p_1$ and $p_2$, is the
product of $X$ and $Y$ in the category of vector spaces
(Exercise~\ref{ex:prod-vs}).
\end{example}

\begin{examples}[Elements of ordered sets]     
\label{eg:prod-order}
\index{ordered set!product in|(}
\begin{enumerate}[(b)]
\item 
Let $x, y \in \reals$.  Their minimum $\min\{x, y\}$ satisfies
\[
\min\{x, y\} \leq x,
\qquad
\min\{x, y\} \leq y
\index{minimum}
\]
and has the further property that whenever $a \in \reals$ with
\[
a \leq x,
\qquad 
a \leq y,
\]
we have $a \leq \min\{x, y\}$.  This means exactly that when the poset
$(\reals, \mathord{\leq})$ is viewed as a category, the product of $x, y
\in \reals$ is $\min\{x, y\}$.  The definition of product simplifies when
interpreted in a poset, since all diagrams commute.

\item 
Fix a set $S$.  Let $X, Y \in \pset(S)$.%
\index{power!set}
Then $X \cap Y$%
\index{intersection}
satisfies
\[
X \cap Y \sub X,
\qquad
X \cap Y \sub Y
\]
and has the further property that whenever $A \in \pset(S)$ with
\[
A \sub X,
\qquad
A \sub Y,
\]
we have $A \sub X \cap Y$.  This means that $X \cap Y$ is the product of
$X$ and $Y$ in the poset $(\pset(S), \mathord{\sub})$ regarded as a
category.

\item 
Let $x, y \in \nat$.  Their greatest%
\index{greatest common divisor}
common divisor $\gcd(x, y)$ satisfies
\[
\gcd(x, y) \divides x,
\qquad
\gcd(x, y) \divides y
\]
(it's a common divisor!)\ and has the further property that whenever $a \in
\nat$ with
\[
a \divides x,
\qquad
a \divides y,
\]
we have $a \divides \gcd(x, y)$.  This means that $\gcd(x, y)$ is the
product of $x$ and $y$ in the poset $(\nat, \mathord{\mid})$ regarded as a
category.
\end{enumerate}

Generally, let $(A, \mathord{\leq})$ be a poset and $x, y \in A$.  A
\demph{lower%
\index{lower bound}
bound} for $x$ and $y$ is an element $a \in A$ such that $a
\leq x$ and $a \leq y$.  A \demph{greatest%
\index{greatest lower bound}
lower bound} or \demph{meet}%
\index{meet}
of $x$ and $y$ is a lower bound $z$ for $x$ and $y$ with the further
property that whenever $a$ is a lower bound for $x$ and $y$, we have $a
\leq z$.

When a poset is regarded as a category, meets are exactly products.  They do
not always exist, but when they do, they are unique.  The meet of $x$ and $y$
is usually written as $x \meet y$%
\ntn{meet}
rather than $x \times y$.  Thus, in the three examples above,
\[
x \meet y = \min\{x, y\},
\qquad
X \meet Y = X \cap Y,
\qquad
x \meet y = \gcd(x, y),
\]
the second example being the origin of the notation.%
\index{ordered set!product in|)}
\end{examples}

We have been discussing products $X \times Y$ of \emph{two} objects,
so-called \demph{binary%
\index{product!binary}
products}.  But there is no reason to stick to two.  We can just as well
talk about products $X \times Y \times Z$ of three objects, or of
infinitely many objects.  The definition changes in the most obvious way:

\begin{defn}    
\label{defn:gen-prod}
Let $\cat{A}$ be a category, $I$ a set, and $(X_i)_{i \in I}$ a family of
objects of $\cat{A}$.  A \demph{product}%
\index{product}
of $(X_i)_{i \in I}$ consists of an object $P$ and a family of maps
\[
\Bigl(P \toby{p_i} X_i\Bigr)_{i \in I}
\]
with the property that for all objects $A$ and families of maps
\begin{equation}        
\label{eq:gen-prod-cone}
\Bigl(A \toby{f_i} X_i\Bigr)_{i \in I}
\end{equation}
there exists a unique map $\bar{f}\from A \to P$ such that $p_i \of \bar{f} =
f_i$ for all $i \in I$.
\end{defn}

Remarks~\ref{rmks:prod} apply equally to this definition.  When the product
$P$ exists, we write $P$ as $\prod_{i \in I} X_i$%
\ntn{prod-fam-gen}
and the map $\bar{f}$ as $(f_i)_{i \in I}$.%
\ntn{map-to-prod}
We call the maps $f_i$ the \demph{components}%
\index{component!map into product@of map into product}
of the map $(f_i)_{i \in I}$.  Taking $I$ to be a two-element set, we
recover the special case of binary products.

\begin{example}        
In ordered sets, the extension from binary to arbitrary products works in
the obvious way: given an ordered set $(A, \mathord{\leq})$, a
\demph{lower%
\index{lower bound}
bound} for a family $(x_i)_{i \in I}$ of elements is an element $a \in A$
such that $a \leq x_i$ for all $i$, and a \demph{greatest%
\index{greatest lower bound}
lower bound} or \demph{meet}%
\index{meet}
of the family is a lower bound greater than any other, written as $\Meet_{i
\in I} x_i$.%
\ntn{Meet}
These are the products in $(A, \mathord{\leq})$.

For example, in $\reals$ with its usual ordering, the meet of a family
$(x_i)_{i \in I}$ is $\inf\{x_i \such i \in I\}$%
\index{infimum}
(and one exists if and only if the other does).
\end{example}

\begin{example}
\label{eg:arb-prods-terminal}
What happens to the definition of product when the indexing set $I$ is
empty?%
\index{product!empty}
Let $\cat{A}$ be a category.  In general, an $I$-indexed family
$(X_i)_{i \in I}$ of objects of $\cat{A}$ is a function $I \to
\ob(\cat{A})$.  When $I$ is empty, there is exactly one such function.  In
other words, there is exactly one family $(X_i)_{i \in \emptyset}$, the
\demph{empty%
\index{empty family}%
\index{family!empty}
family}.  Similarly, when $I$ is empty, there is exactly one
family~\eqref{eq:gen-prod-cone} for any given object $A$.

A product of the empty family therefore consists of an object $P$ of $\cat{A}$
such that for each object $A$ of $\cat{A}$, there exists a unique map
$\bar{f}\from A \to P$.  (The condition `$p_i \of \bar{f} = f_i$ for all $i
\in I$' holds trivially.)  In other words, a product of the empty family is
exactly a terminal%
\index{object!terminal}
object.

We have been writing $1$%
\ntn{terminal}
for terminal objects, which was justified by the fact that in categories
such as $\Set$, $\Tp$, $\Ring$ and $\Grp$, the terminal object has one%
\index{set!one-element}
element.  But we have just seen that the terminal object is the product of
no things, which in the context of elementary arithmetic%
\index{arithmetic}
is the number $1$.  This is a second, related, reason for the notation.
\end{example}

\begin{example}
Take an object $X$ of a category $\cat{A}$, and a set $I$.  There is a
constant family $(X)_{i \in I}$.  Its product $\prod_{i \in I} X$, if it
exists, is written as $X^I$%
\ntn{power-of-obj}
and called a \demph{power}%
\index{power}
of $X$.  

We met powers in $\Set$ in Section~\ref{sec:Set-properties}.  When $X$ is a
set, $X^I$ is the set of functions from $I$ to $X$, also written as
$\Set(I, X)$.
\end{example}

\minihead{Equalizers}

To define our second type of limit, we need a preliminary piece of
terminology: a \demph{fork}%
\index{fork}
in a category consists of objects and maps
\begin{equation}        
\label{eq:fork}
\xymatrix{
A \ar[r]^f   &
X \ar@<.5ex>[r]^s \ar@<-.5ex>[r]_t    &
Y
}
\end{equation}
such that $sf = tf$.

\begin{defn}    
\label{defn:equalizer}
Let $\cat{A}$ be a category and let $\parpairi{X}{Y}{s}{t}$ be objects and
maps in $\cat{A}$.  An \demph{equalizer}%
\index{equalizer}
of $s$ and $t$ is an object $E$ together with a map $E \toby{i} X$ such
that
\[
\xymatrix{
E \ar[r]^i   &
X \ar@<.5ex>[r]^s \ar@<-.5ex>[r]_t    &
Y
}
\]
is a fork, and with the property that for any fork~\eqref{eq:fork},
there exists a unique map $\bar{f}\from A \to E$ such that
\begin{equation}        
\label{eq:equalizer-whole}
\begin{array}{c}
\xymatrix{
A \ar@{.>}[d]_{\bar{f}} \ar[rd]^f    &       \\
E \ar[r]_i                      &X
}
\end{array}
\end{equation}
commutes.
\end{defn}

Remarks~\ref{rmks:prod} on products apply to equalizers too.

\begin{example}
\label{eg:equalizers-Set}
We have already met equalizers%
\index{set!category of sets!equalizers in}%
\index{equalizer!sets@of sets}
in $\Set$ (Section~\ref{sec:Set-properties}).  They really are equalizers
in the sense of Definition~\ref{defn:equalizer}.  Indeed, take sets and
functions $\parpairi{X}{Y}{s}{t}$\!, write
\[
E = \{ x \in X \such s(x) = t(x) \},
\]
and write $i\from E \to X$ for the inclusion.  Then $s i = t i$, so we have
a fork, and one can check that it is universal among all forks on $s$ and
$t$.

An equalizer describes the set of solutions of a single equation, but by
combining equalizers with products, we can also describe the solution-set
of any system of simultaneous%
\index{simultaneous equations}
equations.  Take a set $\Lambda$ and a family 
\[
\biggl( 
\parpair{X}{Y_\lambda}{s_\lambda}{t_\lambda} 
\biggr)_{\lambda \in \Lambda} 
\]
of pairs of maps in $\Set$.  Then the solution-set
\[
\{ 
x \in X 
\such
s_\lambda(x) = t_\lambda(x) \text{ for all } \lambda \in \Lambda 
\} 
\]
is the equalizer of the functions
\[
\xymatrix@C+1em{
X
\ar@<.5ex>[r]^-{(s_\lambda)_{\lambda \in \Lambda}}
\ar@<-.5ex>[r]_-{(t_\lambda)_{\lambda \in \Lambda}}        &
\displaystyle \prod_{\lambda \in \Lambda} Y_\lambda
}
\]
(using the notation introduced after Definition~\ref{defn:gen-prod}).  To
see this, observe that for $x \in X$,
\begin{align*}
(s_\lambda)_{\lambda \in \Lambda}(x) =
(t_\lambda)_{\lambda \in \Lambda}(x) &
\iff
\bigl(s_\lambda(x)\bigr)_{\lambda \in \Lambda} = 
\bigl(t_\lambda(x)\bigr)_{\lambda \in \Lambda}  \\
&
\iff
s_\lambda(x) = t_\lambda(x) \text{ for all } \lambda \in \Lambda,
\end{align*}
as required.
\end{example}

\begin{example}
Take continuous maps $\parpairi{X}{Y}{s}{t}$ between topological\linebreak
spaces.  We can form their equalizer%
\index{topological space!category of topological spaces!equalizers in}
$E$ in the category of sets, with inclusion map $i\from E \to X$, say.
Since $E$ is a subset of the space $X$, it acquires the subspace%
\index{topological space!subspace of}
topology from $X$, and $i$ is then continuous.  This space $E$, together
with $i$, is the equalizer of $s$ and $t$.

Showing this amounts to showing that for any fork~\eqref{eq:fork} in
$\Tp$, the induced function $\bar{f}$ is continuous.  This follows from
the definition of the subspace topology, which is the smallest topology
such that the inclusion map is continuous.  Compare the remarks on products
in Example~\ref{eg:prod-spaces}.
\end{example}

\begin{example}
Let $\theta\from G \to H$ be a homomorphism of groups.  As in
Example~\ref{eg:univ-kernel}, the homomorphism $\theta$ gives rise to a fork
\[
\xymatrix{
\ker\theta\      \ar@{^{(}->}[r]^-\iota   &
G \ar@<.5ex>[r]^\theta \ar@<-.5ex>[r]_\epsln    &
H
}
\index{kernel}
\]
where $\iota$ is the inclusion and $\epsln$ is the trivial homomorphism.
This is an equalizer%
\index{group!category of groups!equalizers in}
in $\Grp$.  Showing this amounts to showing that the map that we have been
calling $\bar{f}$ is a homomorphism, which is left to the reader.

Thus, kernels are a special case of equalizers.
\end{example}

\begin{example} 
\label{eg:eq-vect}
Let $\parpairi{V}{W}{s}{t}$ be linear maps between vector spaces.\linebreak  
There is a linear map $t - s \from V \to W$, and the equalizer%
\index{vector space!category of vector spaces!equalizers in}
of $s$ and $t$ in the category of vector spaces is the space $\ker(t - s)$%
\index{kernel}
together with the inclusion map $\ker(t - s) \incl V$.
\end{example}

\minihead{Pullbacks}

We explore one more type of limit before formulating the general
definition. 

\begin{defn}    
\label{defn:pb}
Let $\cat{A}$ be a category, and take objects and maps
\begin{equation}        
\label{eq:pb-corner}
\begin{array}{c}
\xymatrix{
                &Y \ar[d]^t     \\
X \ar[r]_s      &Z
}
\end{array}
\end{equation}
in $\cat{A}$.  A \demph{pullback}%
\index{pullback}
of this diagram is an object $P \in \cat{A}$ together with maps $p_1\from P
\to X$ and $p_2\from P \to Y$ such that
\begin{equation}        
\label{eq:pb}
\begin{array}{c}
\xymatrix{
P \ar[r]^{p_2} \ar[d]_{p_1}     &
Y \ar[d]^t      \\
X \ar[r]_s      &
Z
}
\end{array}
\end{equation}
commutes, and with the property that for any commutative square
\begin{equation}        
\label{eq:pb-cone}
\begin{array}{c}
\xymatrix{
A \ar[r]^{f_2} \ar[d]_{f_1}     &
Y \ar[d]^t      \\
X \ar[r]_s      &
Z
}
\end{array}
\end{equation}
in $\cat{A}$, there is a unique map $\bar{f}\from A \to P$ such that
\begin{equation}        
\label{eq:pb-whole}
\begin{array}{c}
\xymatrix{
A \ar@/^/[rrd]^{f_2} \ar@{.>}[rd]|{\bar{f}} \ar@/_/[rdd]_{f_1}&       
        &       \\
        &
P \ar[r]^{p_2} \ar[d]_{p_1}     &
Y \ar[d]^t      \\
        &
X \ar[r]_s      &Z
}
\end{array}
\end{equation}
commutes.  (For~\eqref{eq:pb-whole} to commute means only that $p_1 \bar{f}
= f_1$ and $p_2 \bar{f} = f_2$, since the commutativity of the square is
already given.)
\end{defn}

Again, Remarks~\ref{rmks:prod} apply.  

We call~\eqref{eq:pb} a \demph{pullback%
\index{pullback!square}
square}.  Another name for pullback is \demph{fibred%
\index{fibred product}
product}.  This name is partially explained by the following fact: when $Z$
is a terminal object (and $s$ and $t$ are the only maps they can possibly
be), a pullback of the diagram~\eqref{eq:pb-corner} is simply a product%
\index{product!pullback@as pullback}
of $X$ and $Y$.

\begin{examples}[Pullbacks in $\Set$]
\label{egs:pb-sets}
The pullback of a diagram~\eqref{eq:pb-corner} in $\Set$ is
\[
P = \{ (x, y) \in X \times Y \such s(x) = t(y) \}
\]
with projections $p_1$ and $p_2$ given by $p_1(x, y) = x$ and $p_2(x, y) =
y$. 

Although you might not be familiar with general pullbacks in $\Set$, there are
at least two instances that you are likely to have met.
\begin{enumerate}[(b)]
\item   
\label{eg:pb-sets-inv}
A basic construction with sets and functions is the formation of inverse%
\index{inverse!image!pullback@as pullback}
images.  They are an instance of pullbacks.  Indeed, given a function
$f\from X \to Y$ and a subset $Y' \sub Y$, we obtain a new set, the inverse
image
\[
f^{-1}Y' 
=
\{ x \in X \such f(x) \in Y' \} \sub X,
\]
and a new function,
\[
\begin{array}{cccc}
f'\from         &f^{-1} Y'      &\to            &Y'     \\
                &x              &\mapsto        &f(x).
\end{array}
\]
We also have the inclusion functions $j\from Y' \incl Y$ and $i\from f^{-1}Y'
\incl X$.  Putting everything together gives a commutative square
\begin{equation}        
\label{eq:pb-inv}
\begin{array}{c}
\xymatrix@M+.5ex{
f^{-1}Y' \ar[r]^{f'} \ar@{^{(}->}[d]_i &
Y' \ar@{^{(}->}[d]^j   \\
X \ar[r]_f      &
Y.
}
\end{array}
\end{equation}
The data we started with was the lower-right part of this square ($X$, $Y$,
$Y'$, $f$ and $j$), and from it we constructed the rest of the square
($f^{-1} Y'$, $f'$ and $i$).

The square~\eqref{eq:pb-inv} is a pullback.  Let us verify this in detail.
Take any commutative square
\[
\xymatrix@M+.5ex{
A \ar[r]^h \ar[d]_g     &
Y' \ar@{^{(}->}[d]^j   \\
X \ar[r]_f      &
Y.
}
\]
We must show that there is a unique map $k\from A \to f^{-1}Y'$ such that
\[
\xymatrix@M+.5ex{
A \ar@/^/[rrd]^h \ar@{.>}[rd]|k \ar@/_/[rdd]_g  &       &       \\
                                        &
f^{-1}Y' \ar[r]_-{f'} \ar@{^{(}->}[d]^-i       &
Y' \ar@{^{(}->}[d]^j   \\
        &
X \ar[r]_f      &Y
}
\]
commutes.  For uniqueness, let $k$ be a map making the diagram commute.
Then for all $a \in A$, we have $i(k(a)) = g(a)$, that is, $k(a) = g(a)$,
and this determines $k$ uniquely.  For existence, first note that for all
$a \in A$ we have $f(g(a)) = j(h(a)) \in Y'$, so $g(a) \in f^{-1}Y'$.
Hence we may define $k\from A \to f^{-1} Y'$ by $k(a) = g(a)$ for all $a
\in A$.  Then for all $a \in A$, we have $i(k(a)) = k(a) = g(a)$ and
\[
f'(k(a)) = f(k(a)) = f(g(a)) = j(h(a)) = h(a).
\]
Hence $i \of k = g$ and $f' \of k = h$, as required.

\item 
Intersection%
\index{intersection!pullback@as pullback}
of subsets provides another example of pullbacks.  Indeed, let $X$ and $Y$
be subsets of a set $Z$.  Then
\[
\xymatrix@M+.5ex{
X \cap Y \ar@{^{(}->}[r] \ar@{^{(}->}[d]        &
Y \ar@{^{(}->}[d]    \\
X \ar@{^{(}->}[r]       &
Z
}
\]
is a pullback square, where all the arrows are inclusions of subsets.  

In fact, this is a special case of~\bref{eg:pb-sets-inv}, since $X \cap Y$ is
the inverse image of $Y \sub Z$ under the inclusion map $X \incl Z$.
\end{enumerate}
\end{examples}

In the situation of Example~\ref{egs:pb-sets}\bref{eg:pb-sets-inv}, where we
have a map $f\from X \to Y$ and a subset $Y'$ of $Y$, people sometimes say
that $f^{-1}Y'$ is obtained by `pulling $Y'$ back' along $f$: hence the name. 

\minihead{The definition of limit}

We have now looked at three constructions: products, equalizers and
pullbacks.  They clearly have something in common.  Each starts with some
objects and (in the case of equalizers and pullbacks) some maps between
them.  In each, we aim to construct a new object together with some maps
from it to the original objects, with a universal property.

Let us analyse this more closely.  What is the starting data in each
construction?  For (binary) products, it is a pair of objects
\begin{equation}        
\label{eq:prod-data}
X \hspace*{3em} Y.
\end{equation}
For equalizers, it is a diagram 
\begin{equation}        
\label{eq:eq-data}
\parpair{X}{Y.}{s}{t}
\end{equation}
For pullbacks, it is a diagram
\begin{equation}        
\label{eq:pb-data}
\begin{array}{c}
\xymatrix{
                &Y \ar[d]^t     \\
X \ar[r]_s      &Z.
}
\end{array}
\end{equation}

In Definition~\ref{defn:gen-elt}, we met the notion of generalized%
\index{element!generalized}
element, and we saw there that the `figures' in a geometric object can often
be described by maps into it.  For instance, a curve in a topological space
$A$ can be thought of as a map $\reals \to A$.  Similarly, an object of a
category $\cat{A}$ amounts to a functor $D\from \One \to \cat{A}$; think of
$\One = \fbox{$\bullet$}$ as an unlabelled object and $D$ as labelling it
with the name of an object of $\cat{A}$.  And similarly again, a map in a
category $\cat{A}$ is a functor $\Two \to \cat{A}$, where $\Two =
\fbox{$\bullet \to \bullet$}$.%
\ntn{Two}
(Here $\Two$ is the category with two objects, say $0$ and $1$, with one
map $0 \to 1$, and with no other maps except for identities.)  Finally, if
we take $\scat{I}$ to be one of the categories
\begin{equation}        
\label{eq:lim-shapes}
\scat{T} = 
\fbox{$\bullet \hspace*{2.5em} \bullet$}\ ,
\quad
\scat{E} = 
\fbox{$\parpair{\bullet}{\bullet}{}{}$}
\quad
\text{or}
\quad
\scat{P} =
\fbox{%
$\begin{array}{c}
\xymatrix{
                &\bullet \ar[d] \\
\bullet \ar[r]  &\bullet
}
\end{array}$}
\end{equation}
then a functor $\scat{I} \to \cat{A}$ consists of data~\eqref{eq:prod-data},
\eqref{eq:eq-data} or~\eqref{eq:pb-data} in $\cat{A}$, respectively. 

We have just begun to use the convention that one typeface ($\scat{A}$,
$\scat{B}$,%
\ntn{small-cat-face}
$\scat{C}$, \ldots) denotes small%
\index{small}%
\index{category!small}
categories, and another ($\cat{A}$, $\cat{B}$,%
\ntn{arb-cat-face}
 $\cat{C}$, \ldots) denotes arbitrary categories.  Although not strictly
necessary, this convention is helpful, since small categories and arbitrary
categories often play different roles in the theory.

\begin{defn}    
\label{defn:diagram}
Let $\cat{A}$ be a category and $\scat{I}$ a small category.  A functor
$\scat{I} \to \cat{A}$ is called a \demph{diagram}%
\index{diagram}
in $\cat{A}$ of \demph{shape}%
\index{shape!diagram@of diagram}
$\scat{I}$.
\end{defn}

So~\eqref{eq:prod-data}, \eqref{eq:eq-data} and~\eqref{eq:pb-data} are
diagrams of shape $\scat{T}$, $\scat{E}$ and $\scat{P}$.  

\begin{figure}
\centering
\setlength{\unitlength}{1em}
\begin{picture}(18.5,12)(1,-1.5)
\cell{0}{0}{bl}{%
\begin{array}{c}
\xymatrix@C=6em{
        &       &D(I)   \\
A \ar[rru]^{f_I} \ar@{.>}[r]|{\,\exists! \bar{f}\,} \ar[rrd]_{f_J}  &
L \ar[ru]_{p_I} \ar[rd]^{p_J}   &       \\
        &       &D(J)
}
\end{array}
}
\put(16,4.5){\oval(4.5,12)}
\cell{19}{-1}{c}{D}
\end{picture}
\caption{The definition of limit.}  
\label{fig:defn-limit}
\end{figure}
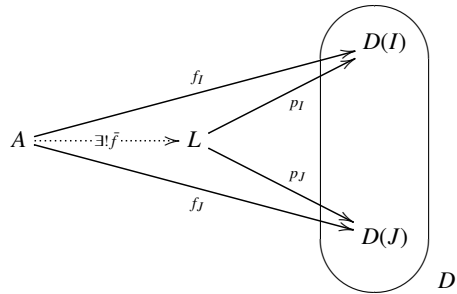

We already have the definitions of product of a diagram of shape
$\scat{T}$, equalizer of a diagram of shape $\scat{E}$, and pullback of a
diagram of shape $\scat{P}$.  We now unify them in the definition of limit
(Figure~\ref{fig:defn-limit}).

\begin{defn}    
\label{defn:lim}
Let $\cat{A}$ be a category, $\scat{I}$ a small category, and $D\from \scat{I}
\to \cat{A}$ a diagram in $\cat{A}$.  
\begin{enumerate}[(b)]
\item
A \demph{cone}%
\index{cone}
on $D$ is an object $A \in \cat{A}$ (the \demph{vertex}%
\index{vertex}
of the cone) together with a family
\begin{equation}        
\label{eq:gen-cone}
\Bigl(
A \toby{f_I} D(I)
\Bigr)_{I \in \scat{I}}
\end{equation}
of maps in $\cat{A}$ such that for all maps $I \toby{u} J$ in $\scat{I}$, the
triangle
\[
\xymatrix@R=1ex{
                                &D(I) \ar[dd]^{Du}      \\
A \ar[ru]^{f_I} \ar[rd]_{f_J}   &                       \\
                                &D(J)
}
\]
commutes.  (Here and later, we abbreviate $D(u)$ as $Du$.)  

\item
A \demph{limit}%
\index{limit}
of $D$ is a cone $\Bigl(L \toby{p_I} D(I)\Bigr)_{I \in \scat{I}}$ with the
property that for any cone~\eqref{eq:gen-cone} on $D$, there exists a
unique map $\bar{f}\from A \to L$%
\ntn{lim-bar}
such that $p_I \of \bar{f} = f_I$ for all $I \in \scat{I}$.  The maps $p_I$
are called the \demph{projections}%
\index{projection}
of the limit.
\end{enumerate}
\end{defn}

\begin{remarks} 
\label{rmks:defn-lim}
\begin{enumerate}[(b)]
\item	
\label{rmk:defn-lim-univ} 
Loosely, the universal property says that for any $A \in \cat{A}$, maps $A
\to L$ correspond one-to-one with cones on $D$ with vertex $A$.  (Any map
$g\from A \to L$ gives rise to a cone $\Bigl( A \toby{p_I g} D(I) \Bigr)_{I
  \in \scat{I}}$, and the definition of limit is that for each $A$, this
process is bijective.)  In Section~\ref{sec:lra}, we will use this thought
to rephrase the definition of limit in terms of representability.  From
this it will follow that limits are unique up to canonical isomorphism,
when they exist (Corollary~\ref{cor:lims-unique}).  Alternatively,
uniqueness can be proved by the usual kind of direct argument, as in
Lemma~\ref{lemma:init-unique}.

\item 
If $\Bigl( L \toby{p_I} D(I) \Bigr)_{I \in \scat{I}}$ is a limit of $D$, we
sometimes abuse language slightly by referring to $L$ (rather than the
whole cone) as the limit%
\index{limit!informal usage}
of $D$.  For emphasis, we sometimes call $\Bigl( L \toby{p_I} D(I)
\Bigr)_{I \in \scat{I}}$ a \demph{limit%
\index{limit!cone}%
\index{cone!limit}
cone}.  We write $L = \lt{\scat{I}} D$.%
\ntn{lim}
Remark~\bref{rmk:defn-lim-univ} can then be stated as:
\begin{slogan}
A map into $\lt{\scat{I}} D$ is a cone on $D$.
\end{slogan}

\item   
\label{rmks:defn-lim:small}
By assuming from the outset that the shape category $\scat{I}$ is small, we
are restricting ourselves to what are officially called \demph{small%
\index{limit!small}%
\index{small}
limits}.  We will seldom be interested in any other kind.
\end{enumerate}
\end{remarks}

\begin{examples}[Limit shapes]        
\label{egs:lims}
Let $\cat{A}$ be any category.  Recall the categories $\scat{T}$,
$\scat{E}$ and $\scat{P}$ of~\eqref{eq:lim-shapes}.
\begin{enumerate}[(b)]
\item
A diagram $D$ of shape $\scat{T}$ in $\cat{A}$ is a pair $(X, Y)$ of
objects of $\cat{A}$.  A cone on $D$ is an object $A$ together with maps
$f_1\from A \to X$ and $f_2\from A \to Y$ (as in
Definition~\ref{defn:bin-prod}), and a limit of $D$ is a product of $X$ and
$Y$.

More generally, let $I$ be a set and write $\scat{I}$ for the discrete
category on $I$.  A functor $D\from \scat{I} \to \cat{A}$ is an $I$-indexed
family $(X_i)_{i \in I}$ of objects of $\cat{A}$, and a limit of $D$ is
exactly a product of the family $(X_i)_{i \in I}$.

In particular, a limit of the unique functor $\emptyset \to \cat{A}$ is a
terminal object of $\cat{A}$, where $\emptyset$ denotes the empty category.

\item	
\label{eg:lim-eq}
A diagram $D$ of shape $\scat{E}$ in $\cat{A}$ is a parallel pair
$\parpairi{X}{Y}{s}{t}$ of maps in $\cat{A}$.  A cone on $D$ consists of
objects and maps
\[
\xymatrix{
        &A \ar[ld]_f \ar[rd]^g  &       \\
X \ar@<.5ex>[rr]^s \ar@<-.5ex>[rr]_t    &       &Y      
}
\]
such that $s \of f = g$ and $t \of f = g$.  But since $g$ is determined by
$f$, it is equivalent to say that a cone on $D$ consists of an object $A$ and
a map $f\from A \to X$ such that 
\[
\xymatrix{
A \ar[r]^f   &
X \ar@<.5ex>[r]^s \ar@<-.5ex>[r]_t    &
Y
}
\]
is a fork.  A limit of $D$ is a universal fork on $s$ and $t$, that is, an
equalizer of $s$ and $t$.

\item
A diagram $D$ of shape $\scat{P}$ in $\cat{A}$ consists of objects and maps
\[
\xymatrix{      
                &Y \ar[d]^t     \\
X \ar[r]_s      &Z
}
\]
in $\cat{A}$.  Performing a simplification similar to that
in~\bref{eg:lim-eq}, we see that a cone on $D$ is a commutative
square~\eqref{eq:pb-cone}.  A limit of $D$ is a pullback.

\item   
\label{eg:lim-seq}
Let $\scat{I} = (\nat, \mathord{\leq})^\op$.  A diagram $D\from \scat{I}
\to \cat{A}$ consists of objects and maps
\[
\cdots \toby{s_3} X_2 \toby{s_2} X_1 \toby{s_1} X_0.
\]
For example, suppose that we have a set $X_0$ and a chain of
subsets
\[
\cdots \sub X_2 \sub X_1 \sub X_0.
\]
The inclusion maps form a diagram in $\Set$ of the type above, and its limit
is $\bigcap_{i \in \nat} X_i$.%
\index{intersection}
In this and similar contexts, limits are sometimes referred to as
\demph{inverse%
\index{inverse!limit}%
\index{limit!inverse}
limits}, although many category theorists regard this usage as
old-fashioned.
\end{enumerate}
\end{examples}

In general, the limit of a diagram $D$ is the terminal object in the
category of cones on $D$, and is therefore an extremal example of a cone on
$D$.  The word `limit' can be understood as meaning `on the boundary',
rather than indicating a limiting process of the type encountered in
analysis.  Nevertheless, the two ideas make contact in
Example~\ref{egs:lims}\bref{eg:lim-seq}.

We have said little so far about which limits exist, except to observe in
Remark~\ref{rmks:prod}\bref{rmks:prod:exist} that they do not exist
always.  We now show that in many familiar categories, all limits do exist;
indeed, we can construct them explicitly.

\begin{example} 
\label{eg:lims-Set}
Let $D\from \scat{I} \to \Set$ and, as a kind of thought%
\index{thought experiment}
experiment, let us ask ourselves what $\lt{\scat{I}} D$%
\index{set!category of sets!limits in}
would have to be if it existed.  (We do not know yet that it does.)  We
would have
\begin{align}
\lt{\scat{I}} D        &
\iso   \Set\biggl(1, \lt{\scat{I}} D\biggr)      
\nonumber       \\
                        &
\iso   
\{ \text{cones on } D \text{ with vertex } 1 \}             
\nonumber       \\
                        &\iso         
\Bigl\{ 
(x_I)_{I \in \scat{I}} 
\Bigsuch
x_I \in D(I) \text{ for all } I \in \scat{I} \text{ and }
(Du)(x_I) = x_J 
\nonumber \\
        &
\phantom{\iso \Bigl\{ (x_I)_{I \in \scat{I}} \Bigsuch {}}
\text{ for all } I \toby{u} J \text{ in }
\scat{I} 
\,\Bigr\},
\label{eq:Set-lim}
\end{align}
where the second isomorphism is by
Remark~\ref{rmks:defn-lim}\bref{rmk:defn-lim-univ} and the third is by
definition of cone.  In fact,~\eqref{eq:Set-lim} really \emph{is} the limit
of $D$ in $\Set$, with projections $p_J\from \lt{\scat{I}} D \to D(J)$
given by $p_J\bigl( (x_I)_{I \in \scat{I}} \bigr) = x_J$
(Exercise~\ref{ex:lims-Set}).  So in $\Set$, all limits exist.
\end{example}

\begin{example} 
\label{eg:lims-alg}
The same formula gives limits in categories of algebras such as $\Grp$,%
\index{group!category of groups!limits in}
$\Ring$,%
\index{ring!category of rings!limits in}
$\Vect_k$,%
\index{vector space!category of vector spaces!limits in}
\ldots.  Of course, we also have to say what the
group/ring/\ldots\ structure on the set~\eqref{eq:Set-lim} is, but this
works in the most straightforward way imaginable.  For instance, in
$\Vect_k$, if $(x_I)_{I \in \scat{I}}, (y_I)_{I \in \scat{I}} \in
\lt{\scat{I}} D$ then
\[
(x_I)_{I \in \scat{I}} + (y_I)_{I \in \scat{I}}
=
(x_I + y_I)_{I \in \scat{I}}.
\]
\end{example}

\begin{example}
The same formula also gives limits in $\Tp$.%
\index{topological space!category of topological spaces!limits in}
The topology on the set~\eqref{eq:Set-lim} is the smallest for which the
projection maps are continuous.
\end{example}

\begin{defn}    
\label{defn:has-lims}
\begin{enumerate}[(b)]
\item 
Let $\scat{I}$ be a small category.  A category $\cat{A}$ \demph{has%
\index{limit!has limits}
limits of shape $\scat{I}$} if for every diagram $D$ of shape $\scat{I}$ in
$\cat{A}$, a limit of $D$ exists.

\item 
A category \demph{has all limits} (or properly, \demph{has small limits})
if it has limits of shape $\scat{I}$ for all small categories $\scat{I}$.
\end{enumerate}
\end{defn}

Thus, $\Set$, $\Tp$, $\Grp$, $\Ring$, $\Vect_k$, \ldots\ all have all
limits. 

Similar terminology can be applied to special classes of limits (for
instance, `has pullbacks').  The class of finite limits is particularly
important.  By definition, a category is \demph{finite}%
\index{category!finite}
if it contains only finitely many maps (in which case it also contains only
finitely many objects).  A \demph{finite%
\index{limit!finite}
limit} is a limit of shape $\scat{I}$ for some finite category $\scat{I}$.
For instance, binary products, terminal objects, equalizers and pullbacks
are all finite limits.

The next result tells us that all limits can be built up from limits of
just a few familiar, basic types.

\begin{propn}	
\label{propn:lim-prod-eq}
\index{limit!products and equalizers@from products and equalizers}
Let $\cat{A}$ be a category.  
\begin{enumerate}[(b)]
\item	
\label{propn:lim-prod-eq:all}
If $\cat{A}$ has all products and equalizers then $\cat{A}$ has all limits.

\item   
\label{propn:lim-prod-eq:fin}
If $\cat{A}$ has binary products, a terminal object and equalizers then
$\cat{A}$ has finite limits.
\end{enumerate}
\end{propn}

To understand the idea, consider formula~\eqref{eq:Set-lim} for limits in
$\Set$.  There, the limit of a diagram $D$ is described as the subset of
the product $\prod_{I \in \scat{I}} D(I)$ consisting of those elements for
which certain equations hold.  We saw in Example~\ref{eg:equalizers-Set}
that the set of solutions to any system of simultaneous%
\index{simultaneous equations}
equations can be described via products and equalizers.  Thus, we can
describe any limit in $\Set$ in terms of products and equalizers.  And in
fact, this same description is valid in any category.

We now examine this idea more closely, in preparation for the proof
(Exercise~\ref{ex:lim-prod-eq}).  First-time readers may wish to skip the
next two paragraphs, resuming at Example~\ref{eg:lim-prod-eq-cpthff}.

Equation~\eqref{eq:Set-lim} states that in $\Set$, the limit of a diagram
$D\from \scat{I} \to \Set$ consists of the elements $(x_I)_{I \in \scat{I}}
\in \prod_{I \in \scat{I}} D(I)$ such that
\[
(Du)(x_J) = x_K 
\]
in $D(K)$ for each map $J \toby{u} K$ in $\scat{I}$.  For each such map
$u$, define maps
\[
\xymatrix{
\displaystyle \prod_{I \in \scat{I}} D(I)
\ar@<.5ex>[r]^-{s_u} \ar@<-.5ex>[r]_-{t_u}      &
D(K)
}
\]
by
\[
s_u\bigl( (x_I)_{I \in \scat{I}} \bigr) 
=
(Du)(x_J),
\qquad
t_u\bigl( (x_I)_{I \in \scat{I}} \bigr) 
=
x_K.
\]
Then $\lt{\scat{I}} D$ is the set of families $x = (x_I)_{I \in \scat{I}}$
satisfying the equation $s_u(x) = t_u(x)$ for each map $u$ in $\scat{I}$.
It follows from Example~\ref{eg:equalizers-Set} that $\lt{\scat{I}} D$ is
the equalizer of
\[
\xymatrix{
\displaystyle 
\prod_{I \in \scat{I}} D(I)
\ar@<1.5ex>[r]^-s \ar@<.5ex>[r]_-t      &
\displaystyle
\prod_{J \toby{u} K \text{ in } \scat{I}} D(K)
}
\]
where $s$ and $t$ are the maps with components $s_u$ and $t_u$,
respectively. 

We have now described any limit in $\Set$ in terms of products and
equalizers.  Although our argument took place entirely in $\Set$, it
suggests how we might proceed in an arbitrary category.  With this in mind,
the proof of Proposition~\ref{propn:lim-prod-eq} is routine, and is left as
Exercise~\ref{ex:lim-prod-eq}.

\begin{example}        
\label{eg:lim-prod-eq-cpthff}
Let $\CptHff$%
\index{topological space!compact Hausdorff}%
\ntn{CptHff}
denote the category of compact Hausdorff spa\-ces and continuous maps.  It
is a classic exercise in topology to show that given continuous maps $s$
and $t$ from a topological space $X$ to a Hausdorff space $Y$, the subset
$\{ x \in X \such s(x) = t(x) \}$ of $X$ is closed.  From this it follows
that $\CptHff$ has equalizers.  Also, Tychonoff's theorem states that any
product (in $\Tp$) of compact spaces is compact, and it is easy to show
that any product (in $\Tp$) of Hausdorff spaces is Hausdorff.  From this it
follows that $\CptHff$ has all products.  Hence by
Proposition~\ref{propn:lim-prod-eq}\bref{propn:lim-prod-eq:all}, $\CptHff$
has all limits.
\end{example}

\begin{example}
Recall from Example~\ref{eg:eq-vect} that kernels provide equalizers in
$\Vect_k$.  By
Proposition~\ref{propn:lim-prod-eq}\bref{propn:lim-prod-eq:fin}, finite
limits%
\index{vector space!category of vector spaces!limits in}%
\index{group!abelian!finite limit of}
in $\Vect_k$ can always be expressed in terms of $\oplus$ (binary direct
sum), $\{0\}$, and kernels.  The same is true in $\Ab$.
\end{example}

\minihead{Monics}

For functions between sets, injectivity is an important concept.  For maps
in an arbitrary category, injectivity does not make sense, but there is a
concept that plays a similar role.

\begin{defn}
Let $\cat{A}$ be a category.  A map $X \toby{f} Y$ in $\cat{A}$ is
\demph{monic}%
\index{monic}
(or a \demph{monomorphism})%
\index{monomorphism}
if for all objects $A$ and maps $\parpairi{A}{X}{x}{x'}$,
\[
f \of x = f \of x'
\implies 
x = x'.
\]
\end{defn}

This can be rephrased suggestively in terms of generalized%
\index{element!generalized}
elements: $f$ is monic if for all generalized elements $x$ and $x'$ of $X$
(of the same shape), $fx = fx' \implies x = x'$.  Being monic is,
therefore, the generalized-element analogue of injectivity.%
\index{function!injective}%
\index{injection}

\begin{example}	
In $\Set$, a map is monic%
\index{set!category of sets!monics in}
if and only if it is injective.  Indeed, if $f$ is injective then certainly
$f$ is monic, and for the converse, take $A = 1$.
\end{example}

\begin{example}
\label{eg:monics-alg} 
In categories of algebras such as $\Grp$,%
\index{group!category of groups!monics in}
$\Vect_k$,%
\index{vector space!category of vector spaces!monics in}
$\Ring$,%
\index{ring!category of rings!monics in}
etc., it is also true that the monic maps are exactly the injections.
Again, it is easy to show that injections are monic.  For the converse,
take $A = F(1)$ where $F$ is the free functor
(Examples~\ref{egs:adjns-alg}).
\end{example}

Why is the definition of monic in a chapter on limits?  Because of this:

\begin{lemma}   
\label{lemma:monic-pb}
A map $X \toby{f} Y$ is monic if and only if the square
\[
\xymatrix{
X \ar[r]^1 \ar[d]_1     &X \ar[d]^f     \\
X \ar[r]_f              &Y
}
\]
is a pullback.
\end{lemma}

\begin{pf}
Exercise~\ref{ex:monic-pb}.
\end{pf}

The significance of this lemma is that whenever we prove a result about
limits, a result about monics will follow.  For example, we will soon show
that the forgetful functors from $\Grp$, $\Vect_k$, etc., to $\Set$
preserve limits (in a sense to be defined), from which it will follow
immediately that they also preserve monics.  This in turn gives an
alternative proof that monics in these categories are injective.

\exs

\begin{question}        
\label{ex:prod-vs}
Verify that in the category of vector spaces, the product of two vector
spaces is their direct sum (Example~\ref{eg:prod-vs}).
\end{question}

\begin{question}
Take objects and maps $\xymatrix@1{E \ar[r]^i &X \ar@<.5ex>[r]^{f}
\ar@<-.5ex>[r]_g &Y}$ in some category.  If this is an equalizer, is the
square 
\[
\xymatrix{
E \ar[r]^i \ar[d]_i     &
X \ar[d]^g      \\
X \ar[r]_f      &
Y
}
\]
necessarily a pullback?%
\index{equalizer!pullback@vs.\ pullback}%
\index{pullback!equalizer@vs.\ equalizer}
What about the converse?  Give proofs or counterexamples.
\end{question}

\begin{question}        
\label{ex:pb-pasting}
Take a commutative diagram 
\[
\xymatrix{
\cdot \ar[r] \ar[d]     &
\cdot \ar[r] \ar[d]     &
\cdot \ar[d]    \\
\cdot \ar[r]    &
\cdot \ar[r]    &
\cdot
}
\index{pullback!pasting of pullbacks}
\]
in some category.  Suppose that the right-hand square is a pullback.  Show
that the left-hand square is a pullback if and only if the outer rectangle is
a pullback.  
\end{question}

\begin{question}        
\label{ex:jointly-monic}
Let $D \from \scat{I} \to \cat{A}$ be a diagram and $\Bigl( L \toby{p_I}
D(I) \Bigr)_{I \in \scat{I}}$ a limit cone on $D$.
\begin{enumerate}[(b)]
\item   
\label{part:j-m-main}
Prove that whenever $\parpairi{A}{L}{h}{h'}$ are maps such that $p_I \of h
= p_I \of h'$ for all $I \in \scat{I}$, then $h = h'$.

\item
What does the result of~\bref{part:j-m-main} mean when $\scat{I}$ is the
two-object discrete category, $\cat{A} = \Set$, and $A = 1$?  Answer
without using any category-theoretic terminology.
\end{enumerate}
\end{question}

\begin{question}        
\label{ex:lims-Set}
Show that the set~\eqref{eq:Set-lim} in Example~\ref{eg:lims-Set} really
is a limit of $D$.  
\end{question}

\begin{question}        
\label{ex:lim-prod-eq}
In this exercise, you will prove Proposition~\ref{propn:lim-prod-eq},
following the plan described after the statement of that proposition.
\begin{enumerate}[(b)]
\item 
Let $\cat{A}$ be a category with all products and equalizers.  Let $D \from
\scat{I} \to \cat{A}$ be a diagram in $\cat{A}$.  Define maps 
\[
\xymatrix{
\displaystyle
\prod_{I \in \scat{I}} D(I)
\ar@<1.5ex>[r]^-s \ar@<.5ex>[r]_-t      &
\displaystyle
\prod_{J \toby{u} K \text{ in } \scat{I}} D(K)
}
\]
as follows: given $J \toby{u} K$ in $\scat{I}$, the $u$-component of
$s$ is the composite
\[
\prod_{I \in \scat{I}} D(I) \toby{\pr_J} D(J) \toby{Du} D(K)
\]
(where $\pr$ denotes a product projection), and the $u$-component of $t$ is
$\pr_K$.  Let $L \toby{p} \prod_{I \in \scat{I}} D(I)$ be the equalizer of
$s$ and $t$, and write $p_I$ for the $I$-component of $p$.  Show that
$\Bigl( L \toby{p_I} D(I) \Bigr)_{I \in \scat{I}}$ is a limit cone on $D$,
thus proving
Proposition~\ref{propn:lim-prod-eq}\bref{propn:lim-prod-eq:all}.

\item   
\label{part:fin-lim-prod-eq}
Adapt the argument to prove
Proposition~\ref{propn:lim-prod-eq}\bref{propn:lim-prod-eq:fin}.
\end{enumerate}
\end{question}

\begin{question}
Prove that a category with pullbacks%
\index{limit!pullbacks and terminal object@from pullbacks and terminal object}
and a terminal object has all finite limits.
\end{question}

\begin{question}        
\label{ex:subobjects}
Let $\cat{A}$ be a category and $A \in \cat{A}$.  A \demph{subobject}%
\index{subobject}
of $A$ is an isomorphism class of monics into $A$.  More precisely, let
$\fcat{Monic}(A)$ be the full subcategory of $\cat{A}/A$ whose objects are
the monics; then a subobject of $A$ is an isomorphism class of objects of
$\fcat{Monic}(A)$.
\begin{enumerate}[(b)]
\item   
\label{part:sub-Set}
Let $X \toby{m} A$ and $X' \toby{m'} A$ be monics in $\Set$.  Show that
$m$ and $m'$ are isomorphic in $\fcat{Monic}(A)$ if and only if they have
the same image.  Deduce that the subobjects of $A$ are in canonical
one-to-one correspondence with the subsets%
\index{subset}
of $A$.

\item 
Part~\bref{part:sub-Set} says that in $\Set$, subobjects are subsets.
What are subobjects in $\Grp$, $\Ring$ and $\Vect_k$?  

\item
What are subobjects in $\Tp$?  (Careful!)
\end{enumerate}
\end{question}

\begin{question}        
\label{ex:monic-pb}
Prove Lemma~\ref{lemma:monic-pb}.
\end{question}

\begin{question}        
\label{ex:pb-monic}
Let
\[
\xymatrix{
X' \ar[r]^{f'} \ar[d]_{m'}      &
X \ar[d]^m      \\
A' \ar[r]_f     &
A
}
\index{pullback!monic@of monic}%
\index{monic!pullback of}
\]
be a pullback square in some category.  Show that if $m$ is monic then so
is $m'$.  (We already know this in the category of sets, by
Example~\ref{egs:pb-sets}\bref{eg:pb-sets-inv}.)
\end{question}

\section{Colimits: definition and examples}

We have seen that examples of limits occur throughout mathematics.  It
therefore makes sense to examine the dual concept, colimit, and ask whether
it is similarly ubiquitous.

By dualizing, we can write down the definition of colimit immediately.  We
then specialize to sums, coequalizers and pushouts, the duals of products,
equalizers and pullbacks.

There are two common conventions for naming dual%
\index{duality!terminology for}
concepts: sometimes we add or subtract the prefix `co' (as in
limit/colimit), and sometimes we use `left' and `right' (as for adjoints).
There are also some irregular names, such as terminal/initial object and
pullback/pushout.

\begin{defn}
Let $\cat{A}$ be a category and $\scat{I}$ a small category.  Let $D\from
\scat{I} \to \cat{A}$ be a diagram in $\cat{A}$, and write $D^\op$ for the
corresponding functor $\scat{I}^\op \to \cat{A}^\op$.  A \demph{cocone}%
\index{cocone}
on $D$ is a cone on $D^\op$, and a \demph{colimit}%
\index{colimit}
of $D$ is a limit of $D^\op$.
\end{defn}

Explicitly, a cocone on $D$ is an object $A \in \cat{A}$ (the
\demph{vertex}%
\index{vertex}
of the cocone) together with a family
\begin{equation}        
\label{eq:gen-cocone}
\Bigl(
D(I) \toby{f_I} A
\Bigr)_{I \in \scat{I}}
\end{equation}
of maps in $\cat{A}$ such that for all maps $I \toby{u} J$ in $\scat{I}$, the
diagram 
\[
\xymatrix@R=1ex{
D(I) \ar[dd]_{Du} \ar[rd]^{f_I} &       \\
                                &A      \\
D(J) \ar[ru]_{f_J}
}
\]
commutes.  A colimit of $D$ is a cocone 
\[
\Bigl(D(I) \toby{p_I} C\Bigr)_{I \in \scat{I}}
\]
with the property that for any cocone~\eqref{eq:gen-cocone} on
$D$, there is a unique map $\bar{f}\from C \to A$%
\ntn{colim-bar}
such that $\bar{f} \of p_I = f_I$ for all $I \in \scat{I}$.  The associated
picture is the mirror image of Figure~\ref{fig:defn-limit}.

Of course, Remarks~\ref{rmks:defn-lim} apply equally here.  We write (the
vertex of) the colimit as $\colt{\scat{I}} D$,%
\ntn{colim}
and call the maps $p_I$ \demph{coprojections}.%
\index{coprojection}

\minihead{Sums}

\begin{defn}
A \demph{sum}%
\index{sum}
or \demph{coproduct}%
\index{coproduct}
is a colimit over a discrete category.  (That is, it is a colimit of shape
$\scat{I}$ for some discrete category $\scat{I}$.)
\end{defn}
Let $(X_i)_{i \in I}$ be a family of objects of a category.  Their sum (if
it exists) is written as $\sum_{i \in I} X_i$%
\ntn{sum-fam-gen}
or $\coprod_{i \in I} X_i$.%
\ntn{disjt-union-fam-gen}
When
$I$ is a finite set $\{1, \ldots, n\}$, we write $\sum_{i \in I} X_i$ as $X_1
+ \cdots + X_n$,%
\ntn{sum-gen}
or as $0$%
\ntn{initial}
if $n = 0$.  

\begin{example}        
\label{eg:sums-init}
By the dual of Example~\ref{eg:arb-prods-terminal}, a sum of the empty%
\index{sum!empty}%
\index{family!empty}%
\index{empty family}
family is exactly an initial%
\index{object!initial}
object.
\end{example}

\begin{example}
Sums in $\Set$%
\index{set!category of sets!sums in}
were described in Section~\ref{sec:Set-properties}.  Let us look in detail
at the universal property, in the case of binary sums.  Take two sets,
$X_1$ and $X_2$.  Form their sum, $X_1 + X_2$, and consider the inclusions
\[
\xymatrix{
X_1 \ar[r]^-{p_1} &X_1 + X_2     &X_2. \ar[l]_-{p_2}
}
\]
This is a colimit cocone.  To prove this, we have to prove the
following universal property: for any diagram 
\[
\xymatrix{
X_1 \ar[r]^{f_1} &A     & X_2 \ar[l]_{f_2}
}
\]
of sets and functions, there is a unique function $\bar{f}\from X_1 + X_2 \to
A$ making
\[
\xymatrix@C+2ex{
X_1 \ar[rrd]^{f_1} \ar[rd]_{p_1}        &       &       \\
        &
X_1 + X_2 \ar@{.>}[r]|{\bar{f}}   &A   \\
X_2 \ar[ru]^{p_2} \ar[rru]_{f_2}
}
\]
commute.  Now, we noted in Section~\ref{sec:Set-properties} that $p_1$ and
$p_2$ are injections whose images partition $X_1 + X_2$.  This means that
every element $x$ of $X_1 + X_2$ is \emph{either} equal to $p_1(x_1)$ for
some $x_1 \in X_1$ (and this $x_1$ is then unique), \emph{or} equal to
$p_2(x_2)$ for some $x_2 \in X_2$ (and this $x_2$ is then unique), but not
both.  So we may define $\bar{f}(x)$ to be equal to $f_1(x_1)$ in the first
case and $f_2(x_2)$ in the second.  This defines a function $\bar{f}$
making the diagram commute, and it is clearly the unique function that does
so. 
\end{example}

\begin{example}
Let $X_1$ and $X_2$ be vector%
\index{vector space!category of vector spaces!sums in}
spaces.  There are linear maps
\begin{equation}        
\label{eq:vs-coprs}
\begin{array}{c}
\xymatrix{
X_1 \ar[r]^-{i_1} &X_1 \oplus X_2        &X_2 \ar[l]_-{i_2}
}
\end{array}
\end{equation}
defined by $i_1(x_1) = (x_1, 0)$ and $i_2(x_2) = (0, x_2)$, and it can be
checked that~\eqref{eq:vs-coprs} is a colimit cocone in $\Vect_k$.  Hence
binary direct%
\index{vector space!direct sum of vector spaces}
sums are sums in the categorical sense.  This is remarkable,
since we saw in Example~\ref{eg:prod-vs} that $X_1 \oplus X_2$ is also the
\emph{product} of $X_1$ and $X_2$!  Contrast this with the category of sets
(or almost any other category), where sums and products are very different.
\end{example}

\begin{example}
Let $(A, \mathord{\leq})$ be an ordered set.%
\index{ordered set!sum in}
\demph{Upper%
\index{upper bound}
bounds} and \demph{least%
\index{least upper bound}
upper bounds} (or \demph{joins})%
\index{join}
in $A$ are defined by dualizing the definitions in
Example~\ref{eg:prod-order}, and, dually, they are sums in the corresponding
category.  The join of a family $(x_i)_{i \in I}$ is written as $\Join_{\! i
\in I} x_i$.%
\ntn{Join}
In the binary case (where $I$ has two elements), the join of $x_1$ and
$x_2$ is written as $x_1 \join x_2$.%
\ntn{join}
A join of the empty family (where $I = \emptyset$) is an initial
object of the category $A$, as in Example~\ref{eg:sums-init}.
Equivalently, it is a \demph{least%
\index{least element}
element} of $A$: an element $0 \in A$ such that $0 \leq a$ for all $a \in
A$.

For instance, in $(\reals, \mathord{\leq})$, join is supremum%
\index{supremum}
and there is no least element.  In a power%
\index{power!set}
set $(\pset(S), \sub)$, join is union%
\index{union}
and the least element is $\emptyset$.  In $(\nat, \divides)$, join is
lowest%
\index{lowest common multiple}
common multiple and the least element is $1$ (since $1$ divides
everything).  So in this order on the natural numbers, $1$ is least; but
also, everything divides $0$, so $0$ is greatest!
\end{example}

\minihead{Coequalizers}

We continue to write $\scat{E}$ for the category $\fbox{$\bullet \parpairu
\bullet$}\,$. 

\begin{defn}
A \demph{coequalizer}%
\index{coequalizer}
is a colimit of shape $\scat{E}$.
\end{defn}

In other words, given a diagram $\parpairi{X}{Y}{s}{t}$, a coequalizer of
$s$ and $t$ is a map $Y \toby{p} C$ satisfying $p \of s = p \of t$ and
universal with this property.

We will see that coequalizers are something like quotients.  But first, we
need some background material on equivalence relations.

\begin{remarks}  
\label{rmks:gen-er}
A binary relation%
\index{relation}
$R$ on a set $A$ can be viewed as a subset $R \sub A \times A$.  Think of
$(a, a') \in R$ as meaning `$a$ and $a'$ are related'.  We can speak of one
relation $S$ on $A$ `containing' another such relation, $R$.  This means
that $R \sub S$: whenever $a$ and $a'$ are $R$-related, they are also
$S$-related.

We will need to use the fact that for any binary relation $R$ on a set $A$,
there is a smallest equivalence relation $\sim$ containing $R$.  This is
called the equivalence relation \demph{generated}%
\index{equivalence relation!generated by relation}%
\index{generated equivalence relation}
by $R$.  `Smallest' means that any equivalence relation containing $R$ also
contains $\sim$.

We can construct $\sim$ as the intersection of all equivalence relations on
$A$ containing $R$, since the intersection of any family of equivalence
relations is again an equivalence relation.  There is also an explicit
construction.  The rough idea is as follows: writing $x \rightarrow y$ to
mean $(x, y) \in R$, we should have $a \sim a'$ if and only if there is a
zigzag such as
\[
a \rightarrow b \leftarrow c \leftarrow d \rightarrow e \leftarrow a'
\]
between $a$ and $a'$.  To make this precise, we first define a relation $S$
on $A$ by
\[
S = 
\{
(a, a') \in A \times A
\such
(a, a') \in R \text{ or } (a', a) \in R
\}
\]
(which enlarges $R$ to a symmetric relation), then define $\sim$ by
declaring that $a \sim a'$ if and only if there exist $n \geq 0$ and $a_0,
\ldots, a_n \in A$ such that
\[
a = a_0, 
\ 
(a_0, a_1) \in S,
\
(a_1, a_2) \in S,
\ 
\ldots,
\ 
(a_{n - 1}, a_n) \in S,
\ 
a_n = a'
\]
(which forces reflexivity and transitivity, while preserving the symmetry). 

Next, recall some facts about equivalence relations from
Section~\ref{sec:Set-properties}.  Given any equivalence relation $\sim$ on
a set $A$, we can construct the set $\qer{A}{\sim}$%
\index{quotient!set@of set}%
\index{set!quotient of}
of equivalence classes and the quotient map $p\from A \to \qer{A}{\sim}$.
This quotient map $p$ is surjective and has the property that $p(a) = p(a')
\iff a \sim a'$, for $a, a' \in A$.  We saw that for any set $B$, the maps
$\qer{A}{\sim} \to B$ correspond one-to-one (via composition with $p$) with
the maps $f\from A \to B$ such that
\begin{equation}        
\label{eq:eq-reln-univ}
\forall a, a' \in A,
\qquad
a \sim a'
\implies
f(a) = f(a').
\end{equation}

Finally, let us consider this universal property in the case where $\sim$
is the equivalence relation generated by some relation $R$.
Condition~\eqref{eq:eq-reln-univ} is then equivalent to:
\begin{equation}        
\label{eq:reln-univ}
\forall a, a' \in A,
\qquad
(a, a') \in R
\implies
f(a) = f(a').
\end{equation}
(Proof: define an equivalence relation $\approx$ on $A$ by $a \approx a'
\iff f(a) = f(a')$.  Condition~\eqref{eq:eq-reln-univ} says that
$\mathord{\sim} \sub \mathord{\approx}$, and condition~\eqref{eq:reln-univ}
that $R \sub \mathord{\approx}$.  But $\sim$ is the smallest equivalence
relation containing $R$, so these statements are equivalent.)  In
conclusion, for any set $B$, the maps $\qer{A}{\sim} \to B$ correspond
one-to-one with the maps $f\from A \to B$ satisfying~\eqref{eq:reln-univ}.
\end{remarks}

\begin{example}
Take sets and functions $\parpairi{X}{Y}{s}{t}$.  To find the
coequalizer%
\index{set!category of sets!coequalizers in}
of $s$ and $t$, we must construct in some canonical way a set $C$ and a
function $p\from Y \to C$ such that $p(s(x)) = p(t(x))$ for all $x \in X$.
So, let $\sim$ be the equivalence relation on $Y$ generated by $s(x) \sim
t(x)$ for all $x \in X$.  (In other words, $\sim$ is generated by the
relation
\[
R = \{(s(x), t(x)) \such x \in X\}
\]
on $Y$.)  Take the quotient map $p\from Y \to \qer{Y}{\sim}$.  By the
correspondence described in Remarks~\ref{rmks:gen-er}, this is indeed the
coequalizer of $s$ and $t$.
\end{example}

\begin{example}
For each pair of homomorphisms $\parpairi{A}{B}{s}{t}$ in $\Ab$,%
\index{group!abelian!coequalizer of}
there is a homomorphism $t - s \from A \to B$, which gives rise to a
subgroup $\im(t - s)$%
\index{image!homomorphism@of homomorphism}
of $B$.  The coequalizer of $s$ and $t$ is the canonical homomorphism $B
\to B/\!\im(t - s)$.  (Compare Example~\ref{eg:eq-vect}.)
\end{example}

\minihead{Pushouts}

\begin{defn}
A \demph{pushout}%
\index{pushout}
is a colimit of shape
\[
\scat{P}^\op
=
\fbox{%
$\begin{array}{c}
\xymatrix{
\bullet \ar[r] \ar[d]   &\bullet        \\
\bullet
}
\end{array}$}\ .
\]
\end{defn}

In other words, the pushout of a diagram
\begin{equation}        
\label{eq:pushout-corner}
\begin{array}{c}
\xymatrix{
X \ar[r]^s \ar[d]_t     &Y      \\
Z
}
\end{array}
\end{equation}
is (if it exists) a commutative square
\[
\xymatrix{
X \ar[r]^s \ar[d]_t     &Y \ar[d]       \\
Z \ar[r]                &\cdot
}
\]
that is universal as such.  In other words still, a pushout in a category
$\cat{A}$ is a pullback in $\cat{A}^\op$. 

\begin{example}
Take a diagram~\eqref{eq:pushout-corner} in $\Set$.  Its pushout%
\index{set!category of sets!pushouts in}
$P$ is $\qer{(Y + Z)}{\sim}$, where $\sim$ is the equivalence relation on
$Y + Z$ generated by $s(x) \sim t(x)$ for all $x \in X$.  The coprojection
$Y \to P$ sends $y \in Y$ to its equivalence class in $P$, and similarly
for the coprojection $Z \to P$.

For example, let $Y$ and $Z$ be subsets of some set $A$.  Then
\[
\xymatrix@M+.5ex{
Y \cap Z \ar@{^{(}->}[r] \ar@{^{(}->}[d] &
Y \ar@{^{(}->}[d] \\
Z \ar@{^{(}->}[r]  &Y \cup Z
}
\index{union!pushout@as pushout}
\]
is a pushout square in $\Set$.  (It is also a pullback%
\index{intersection!pullback@as pullback}
square!  This coincidence is a special property of the category of sets.)
You can check this by verifying the universal property or by using the
formula just stated.  In this case, the formula takes the two sets $Y$ and
$Z$, places them side by side (giving $Y + Z$), then glues the subset $Y
\cap Z$ of $Y$ to the subset $Y \cap Z$ of $Z$ (giving $\qer{(Y + Z)}{\sim}
= Y \cup Z$).
\end{example}

\begin{example}
If $\cat{A}$ is a category with an initial object $0$, and if $Y, Z \in
\cat{A}$, then a pushout of the unique diagram
\[
\xymatrix{
0 \ar[r] \ar[d] &Y      \\
Z
}
\index{sum!pushout@as pushout}
\]
is exactly a sum of $Y$ and $Z$.
\end{example}

\begin{example}
The van Kampen%
\index{van Kampen's theorem}
theorem (Example~\ref{eg:univ-van-Kampen}) says that given a pushout square
in $\Tp$ satisfying certain further hypotheses, the square in $\Grp$
obtained by taking fundamental%
\index{group!fundamental}
groups throughout is also a pushout.
\end{example}

Here is one more shape of colimit, dual to that in
Example~\ref{egs:lims}\bref{eg:lim-seq}.

\begin{example}
A diagram $D\from (\nat, \mathord{\leq}) \to \cat{A}$ consists of objects
and maps
\[
X_0 \toby{s_1} X_1 \toby{s_2} X_2 \toby{s_3} \cdots
\]
in $\cat{A}$.  Colimits of such diagrams are traditionally called
\demph{direct%
\index{direct limit}%
\index{limit!direct}
limits}.  Although the old terms `inverse limit'
(Example~\ref{egs:lims}\bref{eg:lim-seq}) and `direct limit' are made
redundant by the general categorical terms `limit' and `colimit'
respectively, it is worth being aware of them.
\end{example}

With all these examples in mind, we now write down a general formula for
colimits in $\Set$.

\begin{example}
The colimit%
\index{set!category of sets!colimits in}
of a diagram $D\from \scat{I} \to \Set$ is given by
\[
\colt{\scat{I}} D 
= 
\Biggl(
\sum_{I \in \scat{I}} D(I)
\Biggr)
\Bigg/\text{$\sim$}
\]
where $\sim$ is the equivalence relation on $\sum D(I)$ generated by
\[
x \sim (Du)(x)
\]
for all $I \toby{u} J$ in $\scat{I}$ and $x \in D(I)$.  To see this, note
that for any set $A$, the maps
\[
\biggl(\sum D(I)\biggr) \bigg/\text{$\sim$} \: \to A
\]
correspond bijectively with the maps $f\from \sum D(I) \to A$ such that
\[
f(x) = f\bigl((Du)(x)\bigr)
\]
for all $u$ and $x$ (by Remarks~\ref{rmks:gen-er}).  These in turn
correspond to families of maps $\Bigl(D(I) \toby{f_I} A\Bigr)_{I \in
  \scat{I}}$ such that $f_I(x) = f_J\bigl((Du)(x)\bigr)$ for all $u$ and
$x$; but these are exactly the cocones on $D$ with vertex $A$.
\end{example}

\begin{figure}
\centering
\setlength{\unitlength}{1mm}
\begin{picture}(40,48)(10,4)
\cell{30}{10}{b}{\includegraphics[width=39\unitlength]{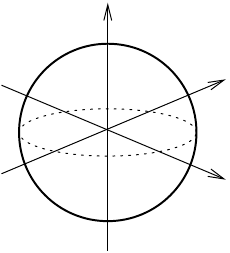}}
\cell{46}{20}{c}{x}
\cell{48}{36}{c}{y}
\cell{31}{51}{c}{z}
\cell{30}{4}{b}{\text{(a)}}
\end{picture}
\hspace*{15mm}
\setlength{\unitlength}{1mm}
\begin{picture}(34,48)(13,4)
\cell{30}{14.5}{b}{\includegraphics[width=33\unitlength]{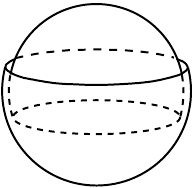}}
\cell{17}{43}{c}{D}
\cell{17}{17}{c}{D}
\cell{30}{4}{b}{\text{(b)}}
\end{picture}
\caption{Sphere as~(a) a limit, and~(b) a colimit.}
\label{fig:lim-colim-spheres}
\end{figure}

There is a kind of duality%
\index{duality}%
\index{limit!colimit@vs.\ colimit}
between the formulas for limits in $\Set$ (Example~\ref{eg:lims-Set}) and
colimits in $\Set$.  Whereas the limit is constructed as a \emph{subset} of
a \emph{product}, the colimit is a \emph{quotient}%
\index{quotient}
of a \emph{sum}.

Figure~\ref{fig:lim-colim-spheres} is intended to convey the difference in
flavour between limits and colimits, in a particular topological context.
In elementary texts, surfaces%
\index{surface|(}
are almost always seen as subsets of Euclidean space $\reals^3$, with the
sphere%
\index{sphere|(}
$S^2$ typically defined as
\[
\bigl\{ (x, y, z) \in \reals^3 \such x^2 + y^2 + z^2 = 1 \bigr\}.
\]
This is a \emph{subspace} of the \emph{product} space $\reals^3 = \reals
\times \reals \times \reals$, which suggests that it is a limit.  Indeed,
the sphere is the equalizer
\[
\xymatrix@M+.5ex{
S^2 \ar@{^{(}->}[r] &\reals^3 \ar@<.5ex>[r]^s \ar@<-.5ex>[r]_t&\reals
}
\]
where the maps $s, t \from \reals^3 \to \reals$ are given by
\[
s(x, y, z) = x^2 + y^2 + z^2,
\qquad
t(x, y, z) = 1.
\]
(An \emph{equa}tion is captured by an \emph{equa}lizer.)%
\index{equalizer}

In more advanced mathematics, however, this point of view is used less
often.  A surface can instead be thought of as the gluing-together of lots
of little patches, each isomorphic to the open unit disk $D$.  For example,
we could in principle construct an entire bicycle%
\index{bicycle inner tube}
inner tube by gluing together a large number of puncture-repair patches.
Figure~\ref{fig:lim-colim-spheres}(b) shows the simpler example of a sphere
made up of two disks glued together.  This realizes the sphere as a
\emph{quotient} (gluing) of the \emph{sum} (disjoint union) of the two
copies of $D$, suggesting that we have constructed the sphere as a colimit.
Indeed, the sphere is the coequalizer
\[
\xymatrix@M+.5ex{
S^1 \times (0, 1) \ar@<.7ex>@{^{(}->}[r] \ar@<-.7ex>@{^{(}->}[r]      &
D + D \ar[r] &S^2
}
\]
where $S^1$ is the circle, the cylinder $S^1 \times (0, 1)$ is the
intersection of the two copies of $D$ (the central belt of
Figure~\ref{fig:lim-colim-spheres}(b)), and the two maps into $D + D$ are
the inclusions of the cylinder into the first and second copies of $D$.%
\index{sphere|)}%
\index{surface|)}

One disadvantage of the limit point of view is that it makes an arbitrary
choice of coordinate system.  It is generally best to think of spaces as
free-standing objects, existing independently of any particular embedding into
Euclidean space.  

One disadvantage of the colimit point of view is that it makes an arbitrary
choice of decomposition.  For example, we could decompose the sphere into
three patches rather than two, or use a different two patches from those
shown.

The colimit point of view has the upper hand in modern geometry.  (If you
are familiar with the definition of manifold,%
\index{manifold}
you will recognize that an atlas is essentially a way of viewing a manifold
as a colimit of Euclidean balls.)  One reason for this is that we are often
concerned with maps \emph{out} of spaces $X$, such as maps $X \to \reals$.
Maps \emph{out} of a colimit are easy; it is in the very definition of
colimit that we know what the maps out of it are.

\minihead{Epics}

\begin{defn}
Let $\cat{A}$ be a category.  A map $X \toby{f} Y$ in $\cat{A}$ is
\demph{epic}%
\index{epic}
(or an \demph{epimorphism})%
\index{epimorphism}
if for all objects $Z$ and maps $\parpairi{Y}{Z}{g}{g'}$,
\[
g \of f = g' \of f
\implies 
g = g'.
\]
\end{defn}

This is the formal dual of the definition of monic.  (In other words, an
epic in $\cat{A}$ is a monic in $\cat{A}^\op$.)  It is in some sense the
categorical version of surjectivity.  But whereas the definition of monic
closely resembles the definition of injective, the definition of epic does
not look much like the definition of surjective.%
\index{function!surjective}%
\index{surjection}
The following examples confirm that in categories where surjectivity makes
sense, it is only sometimes equivalent to being epic.

\begin{example}        
In $\Set$, a map is epic%
\index{set!category of sets!epics in}
if and only if it is surjective.  If $f$ is surjective then certainly $f$
is epic.  To see the converse, take $Z$ to be a two-element set $\{ \true,
\false \}$, take $g$ to be the characteristic function of the image of $f$
(as defined in Section~\ref{sec:Set-properties}), and take $g'$ to be the
function with constant value $\true$.

Any isomorphism in any category is both monic and epic.  In $\Set$, the
converse also holds, since any injective surjective function is invertible
(Example~\ref{eg:iso-Set}). 
\end{example}

\begin{example}
\label{eg:epic-alg}
In categories of algebras, any surjective map is certainly epic.  In some
such categories, including $\Ab$, $\Vect_k$%
\index{vector space!category of vector spaces!epics in}
and $\Grp$,%
\index{group!category of groups!epics in}
the converse also holds.  (The proof is straightforward for $\Ab$ and
$\Vect_k$, but much harder for $\Grp$.)  However, there are other
categories of algebras where it fails.  For instance, in $\Ring$,%
\index{ring!category of rings!epics in}
the inclusion $\integers \incl \rationals$ is epic but not surjective
(Exercise~\ref{ex:epic-Ring}).  This is also an example of a map that is
monic and epic but not an isomorphism.
\end{example}

\begin{example}
In the category of Hausdorff topological spaces%
\index{topological space!category of topological spaces!epics in}%
\index{topological space!Hausdorff}
and continuous maps, any map with dense image is epic.
\end{example}

Of course, there is a dual of Lemma~\ref{lemma:monic-pb}, saying that a map is
epic if and only if a certain square is a pushout.

\exs

\begin{question}
Let $\parpair{X}{Y}{s}{t}$ be maps in some category.  Prove that $s = t$ if
and only if the equalizer of $s$ and $t$ exists and is an isomorphism, if
and only if the coequalizer of $s$ and $t$ exists and is an isomorphism.
\end{question}

\begin{question}
\begin{enumerate}[(b)]
\item 
Let $X$ be a set and $f\from X \to X$ a map.  Describe the coequalizer of
$\parpairi{X}{X}{f}{1}$ in $\Set$ as explicitly as possible. 

\item
Do the same in $\Tp$ rather than $\Set$.  When $X$ is the circle $S^1$,
find an $f$ such that the coequalizer is an uncountable space with the
indiscrete topology.
\end{enumerate}
\end{question}

\begin{question}        
\label{ex:epic-Ring}
\begin{enumerate}[(b)]
\item 
Prove that in the category of monoids,%
\index{monoid!epics between monoids}
the inclusion $(\nat, +, 0) \incl (\integers, +, 0)$ is epic, even though
it is not surjective.

\item 
Prove that in the category of rings, the inclusion $\integers \incl
\rationals$ is epic, even though it is not surjective.
\end{enumerate}
\end{question}

\begin{question}        
\label{ex:qt-objects}
(Compare Exercise~\ref{ex:subobjects}.)  Let $\cat{A}$ be a category and $A
\in \cat{A}$.  Define a \demph{quotient%
\index{quotient}
object} of $A$ to be an isomorphism class of epics out of $A$.  That is,
let $\fcat{Epic}(A)$ be the full subcategory of $A/\cat{A}$ whose objects
are the epics; then a quotient object of $A$ is an isomorphism class of
objects of $\fcat{Epic}(A)$.
\begin{enumerate}[(b)]
\item  
Let $A \toby{e} X$ and $A \toby{e'} X'$ be epics in $\Set$.  Show that $e$
and $e'$ are isomorphic in $\fcat{Epic}(A)$ if and only if they induce the
same equivalence%
\index{equivalence relation}
relation on $A$.  Deduce that the quotient objects of $A$ are in canonical
one-to-one correspondence with the equivalence relations on $A$.

\item 
Assuming the (nontrivial) fact that the epics in $\Grp$ are the
surjections, show that the quotient objects of a group correspond
one-to-one with its normal%
\index{group!normal subgroup of}
subgroups.
\end{enumerate}
(The name `quotient object' is not standard, and indeed there is no
standard name for it.  Arguably, `quotient object' would be more suitable
for an isomorphism class of \emph{regular} epics, as defined in the
following exercises.)
\end{question}

\begin{question}        
\label{ex:reg-split-monic}
A map $m\from A \to B$ is \demph{regular%
\index{monic!regular}
monic} if there exist an object $C$ and maps $B \parpairu C$ of which $m$
is an equalizer.  A map $m\from A \to B$ is \demph{split%
\index{monic!split}
monic} if there exists a map $e\from B \to A$ such that $em = 1_A$.
\begin{enumerate}[(b)]
\item 
Show that split monic $\implies$ regular monic $\implies$ monic.

\item 
In $\Ab$, show that all monics are regular but not all monics are split.
(Hint for the first part: equalizers in $\Ab$ are calculated as in
Example~\ref{eg:eq-vect}.)

\item 
In $\Tp$, describe the regular monics, and find a monic that is not
regular.  
\end{enumerate}
\end{question}

\begin{question}
Dualizing the definitions in Exercise~\ref{ex:reg-split-monic} gives
definitions of \demph{regular}%
\index{epic!regular}
and \demph{split%
\index{epic!split}
epic}.
\begin{enumerate}[(b)]
\item 
We saw in Example~\ref{eg:epic-alg} that a map may be monic and epic but
not an isomorphism.  Prove that in any category, a map is an isomorphism if
and only if it is both monic and \emph{regular} epic.

\item 
Using the assumption that our category of sets satisfies the axiom of
choice (Section~\ref{sec:Set-properties}), show that
\[
\text{epic} \iff \text{regular epic} \iff \text{split epic}
\]
in $\Set$.

\item 
Let us say that a category $\cat{A}$ satisfies the \demph{axiom%
\index{axiom of choice}
of choice} if all epics in $\cat{A}$ are split.  Prove that neither $\Tp$
nor $\Grp$ satisfies the axiom of choice.
\end{enumerate}
\end{question}

\begin{question}
The result of Exercise~\ref{ex:pb-monic} can be phrased as `the class of
monics%
\index{monic!pullback of}%
\index{pullback!monic@of monic}
is stable under pullback'.  It is also a fact that the composite%
\index{monic!composition of monics}
of two monics is always monic; we say that the class of monics is `closed
under composition'.

Consider the following six classes of map:
\begin{displaytext}
monics, regular monics, split monics,
epics, regular epics, split epics.
\end{displaytext}
Determine whether each class is stable under pullback or closed under
composition.
\end{question}

\section{Interactions between functors and limits}
\label{sec:lims-ftrs}

We saw in Example~\ref{eg:lims-alg} that limits in categories such as
$\Grp$, $\Ring$ and $\Vect_k$ can be computed by first taking the limit in
the category of sets, then equipping the result with a suitable algebraic
structure.  On the other hand, colimits in these categories are unlike
colimits in $\Set$.  For example, the underlying set of the initial object
of $\Grp$ (which has one element) is not the initial object of $\Set$
(which has no elements), and the underlying set of the direct sum $X \oplus
Y$ of two vector spaces is not the sum of the underlying sets of $X$ and
$Y$.  So, these forgetful functors interact well with limits and badly with
colimits.

In this section, we develop terminology that will enable us to express
these thoughts precisely.

\begin{defn}
\begin{enumerate}[(b)]
\item   
\label{defn:pres-lims-shape}
Let $\scat{I}$ be a small category.  A functor $F\from \cat{A} \to \cat{B}$
\demph{preserves%
\index{limit!preservation of}
limits of shape $\scat{I}$} if for all diagrams $D\from \scat{I} \to
\cat{A}$ and all cones $\Bigl(A \toby{p_I} D(I)\Bigr)_{I \in \scat{I}}$ on
$D$,
\begin{align*}
        &
\Bigl(A \toby{p_I} D(I)\Bigr)_{I \in \scat{I}}
\text{ is a limit cone on } D \text{ in }\cat{A}\\
\implies        &
\Bigl(F(A) \toby{Fp_I} FD(I)\Bigr)_{I \in \scat{I}}
\text{ is a limit cone on } F \of D \text{ in }\cat{B}. 
\end{align*}

\item   
\label{defn:pres-lims-all}
A functor $F\from \cat{A} \to \cat{B}$ \demph{preserves limits} if it
preserves limits of shape $\scat{I}$ for all small categories $\scat{I}$.

\item 
\demph{Reflection}%
\index{limit!reflection of}%
\index{reflection of limits}
of limits is defined as in~\bref{defn:pres-lims-shape}, but with $\textif$
in place of $\textonlyif$.
\end{enumerate}
\end{defn}
Of course, the same terminology applies to colimits.

Here is a different way to state the definition of preservation.  A functor
$F\from \cat{A} \to \cat{B}$ preserves limits if and only if it has the
following property: whenever $D\from \scat{I} \to \cat{A}$ is a diagram
that has a limit, the composite $F \of D \from \scat{I} \to \cat{B}$ also
has a limit, and the canonical map
\[
F \biggl( \lt{\scat{I}} D \biggr)
\to
\lt{\scat{I}}(F \of D)
\]
is an isomorphism.  Here the `canonical map' has $I$-component 
\[
F \biggl( \lt{\scat{I}} D \biggr)
\toby{F(p_I)}
F(D(I)),
\]
where $p_I$ is the $I$th projection of the limit cone on $D$.

In particular, if $F$ preserves limits then
\begin{equation}
\label{eq:lim-pres-rough}
F \biggl( \lt{\scat{I}} D \biggr)
\iso
\lt{\scat{I}}(F \of D)
\end{equation}
whenever $D$ is a diagram with a limit.  Preservation of limits says more
than \eqref{eq:lim-pres-rough} does: the left- and right-hand sides are
required to be not just isomorphic, but isomorphic \emph{in a particular
  way}.  Nevertheless, we will sometimes omit this check, acting as if
preservation means only that~\eqref{eq:lim-pres-rough} holds.

\begin{example}
The forgetful functor $U\from \Tp \to \Set$%
\index{topological space!category of topological spaces!limits in}%
\index{topological space!category of topological spaces!colimits in}
preserves both limits and colimits.  (As we will see, this follows from the
fact that $U$ has adjoints on both sides.)  It does not reflect all limits
or all colimits.  For instance, choose any non-discrete spaces $X$ and $Y$,
and let $Z$ be the set $U(X) \times U(Y)$ equipped with the discrete
topology.  (All that matters here is that the topology on $Z$ is strictly
larger than the product topology.)  Then we have a cone
\begin{equation}        
\label{eq:discrete-pjn}
X \ot Z \to Y
\end{equation}
in $\Tp$ whose image in $\Set$ is the product cone
\[
U(X) \ot U(X) \times U(Y) \to U(Y).
\]
But~\eqref{eq:discrete-pjn} is not a product cone in $\Tp$, since the
discrete topology on $U(X) \times U(Y)$ is not the product topology.
\end{example}

\begin{example}
In the first paragraph of this section, we observed that the forgetful
functor $\Grp \to \Set$%
\index{group!category of groups!colimits in}
does not preserve initial objects and that the forgetful functor $\Vect_k
\to \Set$%
\index{vector space!category of vector spaces!colimits in}
does not preserve binary sums.  Forgetful functors out of categories of
algebras very seldom preserve all colimits.
\end{example}

\begin{figure}
\centering
\setlength{\unitlength}{1em}
\begin{picture}(12.8,14)(-5.6,-7)
\cell{0}{0}{c}{\includegraphics[height=14\unitlength]{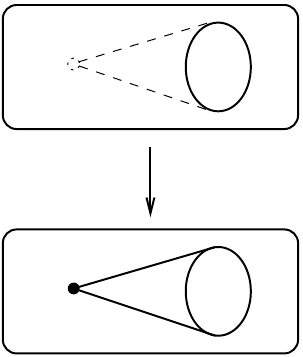}}
\cell{6.7}{4.4}{c}{\cat{A}}
\cell{2.7}{4.4}{c}{D}
\cell{-0.2}{4.3}{c}{p}
\cell{-3.8}{4.4}{c}{A}
\cell{1}{0}{c}{F}
\cell{6.7}{-4.4}{c}{\cat{B}}
\cell{2.7}{-4.4}{c}{F \!\of\! D}
\cell{-0.2}{-4.5}{c}{q}
\cell{-3.8}{-4.4}{c}{B}
\end{picture}%
\caption{Creation of limits.}
\label{fig:creation}
\end{figure}

\begin{example} 
\label{eg:gp-creation}
\index{group!category of groups!limits in|(}%
\index{vector space!category of vector spaces!limits in|(}%
\index{ring!category of rings!limits in|(}
We also saw that (in the examples mentioned) forgetful functors on
categories of algebras do preserve limits.  In fact, something stronger is
true.  Let us examine the case of binary products in $\Grp$, although all
of the following can be said for any limits in any of the categories
$\Grp$, $\Ab$, $\Vect_k$, $\Ring$, etc.

Take groups $X_1$ and $X_2$.  We can form the product set $U(X_1) \times
U(X_2)$, which comes equipped with projections
\[
U(X_1) \otby{p_1} U(X_1) \times U(X_2) \toby{p_2} U(X_2).
\]
I claim that there is exactly one group structure on the set $U(X_1) \times
U(X_2)$ with the property that $p_1$ and $p_2$ are homomorphisms.  To prove
uniqueness, suppose that we have a group structure on $U(X_1) \times U(X_2)$
with this property.  Take elements $(x_1, x_2)$ and $(x'_1, x'_2)$ of
$U(X_1) \times U(X_2)$ and write $(x_1, x_2) \cdot (x'_1, x'_2) = (y_1,
y_2)$.  Since $p_1$ is a homomorphism,
\[
y_1
=
p_1(y_1, y_2)
=
p_1((x_1, x_2) \cdot (x'_1, x'_2))
=
p_1(x_1, x_2) \cdot p_1(x'_1, x'_2)
=
x_1 \cdot x'_1,
\]
and similarly $y_2 = x_2 \cdot x'_2$.  Hence 
\[
(x_1, x_2) \cdot (x'_1, x'_2) = (x_1 x'_1, x_2 x'_2).  
\]
A similar argument shows that $(x_1, x_2)^{-1} = (x_1^{-1}, x_2^{-1})$ and
that the identity element $1$ of the group is $(1, 1)$.  Now, for
existence, define $\cdot$, $\blank^{-1}$ and $1$ by the formulas just
given; it can then be checked that the group axioms are satisfied and that
$p_1$ and $p_2$ are group homomorphisms.  This proves the claim.

Write $L$ for the set $U(X_1) \times U(X_2)$ equipped with this group
structure.  Then we have a cone
\[
X_1 \otby{p_1} L \toby{p_2} X_2
\]
in $\Grp$.  It is easy to check that this is, in fact, a \emph{product}
cone in $\Grp$.

We can summarize this in language that is not tied to group theory.  Given
objects $X_1$ and $X_2$ of $\Grp$,
\begin{itemize}
\item 
for any product cone on $(U(X_1), U(X_2))$ in $\Set$, there is a unique
cone on $(X_1, X_2)$ in $\Grp$ whose image under $U$ is the cone we started
with; 

\item 
this cone on $(X_1, X_2)$ is a product cone.
\end{itemize}
\end{example}
This suggests the following definition (Figure~\ref{fig:creation}).

\begin{defn}    
\label{defn:creates}
A functor $F\from \cat{A} \to \cat{B}$ \demph{creates%
\index{limit!creation of|(}%
\index{creation of limits|(}
limits (of shape $\scat{I}$)} if whenever $D\from \scat{I} \to \cat{A}$ is
a diagram in $\cat{A}$,
\begin{itemize}
\item 
for any limit cone $\Bigl(B \toby{q_I} FD(I)\Bigr)_{I \in \scat{I}}$ on
the diagram $F \of D$, there is a unique cone $\Bigl(A \toby{p_I}
D(I)\Bigr)_{I \in \scat{I}}$ on $D$ such that $F(A) = B$ and $F(p_I) = q_I$
for all $I \in \scat{I}$;

\item 
this cone $\Bigl(A \toby{p_I} D(I)\Bigr)_{I \in \scat{I}}$ is a limit cone
on $D$.
\end{itemize}
\end{defn}
The forgetful functors from $\Grp$, $\Ring$, \ldots\ to $\Set$ all create
limits (Exercise~\ref{ex:fgt-create}).  The word \emph{creates} is
explained by the following result.

\begin{lemma}   
\label{lemma:creates-preserves}
Let $F \from \cat{A} \to \cat{B}$ be a functor and $\scat{I}$ a small
category.  Suppose that $\cat{B}$ has, and $F$ creates, limits of shape
$\scat{I}$.  Then $\cat{A}$ has, and $F$ preserves, limits of shape
$\scat{I}$. 
\end{lemma}

\begin{pf}
Exercise~\ref{ex:creates-preserves}.
\end{pf}

Since $\Set$ has all limits, it follows that all our categories of algebras
have all limits, and that the forgetful functors preserve them.

\begin{remark}
There is something suspicious about Definition~\ref{defn:creates}.  It
refers to \emph{equality} of objects of a category, a relation that, as we
saw on page~\pageref{p:care}, is usually too strict to be appropriate.  It
is almost always better to replace equality by isomorphism.  If we replace
equality by isomorphism throughout the definition of `creates limits', we
obtain a more healthy and inclusive notion.  In the notation of
Definition~\ref{defn:creates}, we ask that if $F \of D$ has a limit then
there exists a cone on $D$ whose image under $F$ is a limit cone, and that
every such cone is itself a limit cone.  

In fact, what we are calling creation of limits should really be called
\emph{strict} creation of limits, with `creation of limits' reserved for
the more inclusive notion.  That is how `creates' is used in most of the
literature.  I have chosen to use the strict version here because it is
slightly simpler to state, and because the examples at hand all satisfy the
stricter condition.
\index{limit!creation of|)}%
\index{creation of limits|)}%
\index{group!category of groups!limits in|)}%
\index{vector space!category of vector spaces!limits in|)}%
\index{ring!category of rings!limits in|)}
\end{remark}

\exs

\begin{question}        
\label{ex:prod-Set-functorial}
Taking the limit is a process that receives as its input a diagram in a
category $\cat{A}$, and produces as its output a new object of $\cat{A}$.
Later, we will see that this process is functorial%
\index{limit!functoriality of}
(Proposition~\ref{propn:lim-const-adjn}).  Here you are asked to prove this
in the case of binary products.%
\index{product!functoriality of}

Let $\cat{A}$ be a category with binary products.  Suppose that we have chosen
for each pair $(X, Y)$ of objects a product cone
\[
\xymatrix@1{
X       &X \times Y \ar[l]_{p^{X, Y}_1} \ar[r]^-{p^{X, Y}_2}     &Y.
}
\]
Construct a functor $\cat{A} \times \cat{A} \to \cat{A}$ given on objects
by $(X, Y) \mapsto X \times Y$.
\end{question}

\begin{question}
Let $\cat{A}$ be a category with binary products.  Prove directly that 
\[
\cat{A}(A, X \times Y)
\iso 
\cat{A}(A, X) \times \cat{A}(A, Y)
\]
naturally in $A, X, Y \in \cat{A}$.  (This presupposes that we have chosen
for each $X$ and $Y$ a product cone on $(X, Y)$.  By
Exercise~\ref{ex:prod-Set-functorial}, the assignment $(X, Y) \mapsto X
\times Y$ is then functorial, which it must be in order for `naturally' to
make sense.)
\end{question}

\begin{question}
Prove that if a functor creates limits then it also reflects them. 
\end{question}

\begin{question}        
\label{ex:fgt-create}
It was shown in Example~\ref{eg:gp-creation} that the forgetful functor
$U\from \Grp \to \Set$%
\index{group!category of groups!limits in}
creates binary products.
\begin{enumerate}[(b)]
\item
Using the formula for limits in $\Set$ (Example~\ref{eg:lims-Set}), prove
that, in fact, $U$ creates arbitrary limits.

\item
Satisfy yourself that the same is true if $\Grp$ is replaced by any other
category of algebras such as $\Ring$,%
\index{ring!category of rings!limits in}
$\Ab$ or $\Vect_k$.%
\index{vector space!category of vector spaces!limits in}
\end{enumerate}
\end{question}

\begin{question}        
\label{ex:creates-preserves}
Prove Lemma~\ref{lemma:creates-preserves}.
\end{question}

\begin{question}
\begin{enumerate}[(b)]
\item 
An object $P$ of a category $\cat{B}$ is \demph{projective}%
\index{projective object}%
\index{object!projective}
if $\cat{B}(P, \dashbk) \from \cat{B} \to \Set$ preserves epics.  (This
means that if $f$ is epic then so is $\cat{B}(P, f)$.)  Let
$\hadjnli{\Set}{\cat{B}}{F}{G}$ be an adjunction in which $G$ preserves
epics.  Prove that $F(S)$ is projective for all sets $S$.

\item 
Find a non-projective object of $\Ab$.

\item 
An object $I$ of a category $\cat{B}$ is \demph{injective}%
\index{injective object}%
\index{object!injective}
if it is projective in $\cat{B}^\op$, or equivalently if $\cat{B}(\dashbk,
I)\from \cat{B}^\op \to \Set$ preserves epics.  Show that all objects of
$\Vect_k$ are injective, and find a non-injective object of $\Ab$.
\end{enumerate}
\end{question}

%
%
%

\chapter{Adjoints, representables and limits}
\label{ch:arl}

We have approached the idea of universal property from three different
angles, producing three different formalisms: adjointness,
representability, and limits.  In this final chapter, we work out the
connections between them.

In principle, anything that can be described in one of the three formalisms
can also be described in the others.  The situation is similar to that of
cartesian and polar coordinates: anything that can be done in polar
coordinates can in principle be done in cartesian coordinates, and vice
versa, but some things are more gracefully done in one system than the
other.

In comparing the three approaches, we will discover many of the fundamental
results of category theory.  Here are some highlights.
\begin{itemize}
\item 
Limits and colimits in functor categories work in the simplest possible
way.

\item 
The embedding of a category $\scat{A}$ into its presheaf category
$\pshf{\scat{A}}$ preserves limits (but not colimits).

\item 
The representables are the prime numbers of presheaves: every presheaf can
be expressed canonically as a colimit of representables.

\item 
A functor with a left adjoint preserves limits.  Under suitable hypotheses,
the converse holds too.

\item 
Categories of presheaves $\pshf{\scat{A}}$ behave very much like the
category of sets, the beginning of an incredible story that brings together
the subjects of logic and geometry.
\end{itemize}

\section{Limits in terms of representables and adjoints}
\label{sec:lra}

There is more than one way to present the definition of limit.  In
Chapter~\ref{ch:lims}, we used an explicit form of the definition that is
particularly convenient for examples.  But we will soon be developing the
\emph{theory} of limits and colimits, and for that, a rephrased form of the
definition is useful.  In fact, we rephrase it in two different ways: once
in terms of representability, and once in terms of adjoints.

We begin by showing that cones are simply natural transformations of a
special kind.  To do this, we need some notation.  Given categories
$\scat{I}$ and $\cat{A}$ and an object $A \in \cat{A}$, there is a functor
$\Delta A\from \scat{I} \to \cat{A}$ with constant value $A$ on objects and
$1_A$ on maps.  This defines, for each $\scat{I}$ and $\cat{A}$, the
\demph{diagonal%
\index{functor!diagonal}
functor}
\[
\Delta\from \cat{A} \to \ftrcat{\scat{I}}{\cat{A}}.%
\ntn{diag-gen}
\]
The name can be understood by considering the case in which $\scat{I}$ is the
discrete category with two objects; then $\ftrcat{\scat{I}}{\cat{A}} = \cat{A}
\times \cat{A}$ and $\Delta(A) = (A, A)$.

Now, given a diagram $D\from \scat{I} \to \cat{A}$ and an object $A \in
\cat{A}$, a cone on $D$ with vertex $A$ is simply%
\index{cone!natural transformation@as natural transformation}
a natural transformation
\[
\xymatrix{
\scat{I} \rtwocell^{\Delta A}_{D} &\cat{A}.
}
\]
Writing $\Cone(A, D)$%
\ntn{Cone}
for the set of cones on $D$ with vertex $A$, we therefore have
\begin{equation}        
\label{eq:cones-as-transfs}
\Cone(A, D)
=
\ftrcat{\scat{I}}{\cat{A}} (\Delta A, D).
\end{equation}
Thus, $\Cone(A, D)$ is functorial in $A$ (contravariantly) and $D$
(covariantly).

Here is our first rephrasing of the definition of limit.

\begin{propn}   
\label{propn:lim-rep}
\index{limit!representation of cone functor@as representation of cone
  functor} 
Let $\scat{I}$ be a small category, $\cat{A}$ a category, and $D\from \scat{I}
\to \cat{A}$ a diagram.  Then there is a one-to-one correspondence between
limit cones on $D$ and representations of the functor
\[
\Cone(\dashbk, D) \from \cat{A}^\op \to \Set,
\]
with the representing objects of \hspace{.05em}$\Cone(\dashbk, D)$ being the
limit objects (that is, the vertices of the limit cones) of $D$.
\end{propn}

Briefly put: a limit of $D$ is a representation of
$\ftrcat{\scat{I}}{\cat{A}}(\Delta\dashbk, D)$. 

\begin{pf}
By Corollary~\ref{cor:rep-univ}, a representation of $\Cone(\dashbk, D)$
consists of a cone on $D$ with a certain universal property.  This is exactly
the universal property in the definition of limit cone.
\end{pf}

The proposition formalizes the thought that cones on a diagram $D$
correspond one-to-one with maps into $\lt{\scat{I}} D$.  It implies that if
$D$ has a limit then
\begin{equation}        
\label{eq:cone-is-map-to-lim}
\Cone(A, D) \iso \cat{A}\biggl(A, \lt{\scat{I}} D\biggr)
\end{equation}
naturally in $A$.  The correspondence is given from left to right by 
\[
(f_I)_{I \in \scat{I}} \mapsto \bar{f}
\]
(in the notation of Definition~\ref{defn:lim}), and from right to left by 
\[
(p_I \of g)_{I \in \scat{I}} \mapsfrom g
\]  
where $p_I\from \lt{\scat{I}} D \to D(I)$ are the projections.

From Proposition~\ref{propn:lim-rep} and Corollary~\ref{cor:reps-unique} we
deduce: 
\begin{cor}     
\label{cor:lims-unique}
\index{limit!uniqueness of}
Limits are unique up to isomorphism.
\qed
\end{cor}

The characterization~\eqref{eq:cones-as-transfs} of cones suggests that we
might consider varying the diagram $D$ as well as the vertex $A$.  We are
naturally led to ask questions such as: given a map $D \to D'$ between
diagrams, is there an induced map between the limits of $D$ and $D'$?  The
answer is yes (Figure~\ref{fig:induced-maps}):

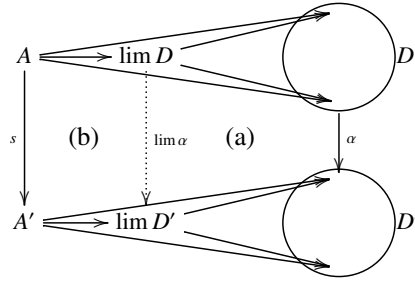
\begin{figure}
\[
\xymatrix@R=2ex{
        &       &\,      \\
A \ar[r] \ar[dddd]_s 
\ar[rru] \ar[rrd]       
\ar@{}[rdddd]|-{\text{\normalsize(b)}}&
\lim D \ar@{.>}[dddd]^{\lim \alpha}       
\ar[ru] \ar[rd]
\ar@{}[rdddd]|-{\text{\normalsize(a)}}&
\save*+<10ex>[o][F]{\,}\restore \hspace*{5.5em}D 
        \\
        &       & \ar[dd]^\alpha        \\
        &       &        \\
        &       &\,      \\
A' \ar[r]   
\ar[rru] \ar[rrd]       &
\lim D' \ar[ru] \ar[rd]         &       
\save*+<10ex>[o][F]{\,}\restore \hspace*{5.8em}D'
\\
        &       &\,       
}
\]
\caption{Illustration of Lemma~\ref{lemma:lim-functorial}.}
\label{fig:induced-maps}
\end{figure}

\begin{lemma}   
\label{lemma:lim-functorial}
Let $\scat{I}$ be a small category and $\xymatrix@1{\scat{I}
  \rtwocell^D_{D'}{\alpha} &\cat{A}}$ a natural transformation.  Let
\[
\biggl(\lt{\scat{I}} D \toby{p_I} D(I)\biggr)_{I \in \scat{I}}
\qquad
\text{and}
\qquad
\biggl(\lt{\scat{I}} D' \toby{p'_I} D'(I)\biggr)_{I \in \scat{I}}
\]
be limit cones.  Then:
\begin{enumerate}[(b)]
\item   
\label{lemma:lim-functorial:induced}
there is a unique map%
\index{limit!map between limits}
$\lt{\scat{I}} \alpha \from \lt{\scat{I}} D \to \lt{\scat{I}} D'$ such that
for all $I \in \scat{I}$, the square
\[
\xymatrix{
\lt{\scat{I}} D \ar[r]^-{p_I} \ar[d]_{\lt{\scat{I}} \alpha}     &
D(I) \ar[d]^{\alpha_I}    \\
\lt{\scat{I}} D' \ar[r]_-{p'_I}         &
D'(I)
}
\]
commutes;

\item   
\label{lemma:lim-functorial:commutes}
given cones $\Bigl(A \toby{f_I} D(I)\Bigr)_{I \in \scat{I}}$ and $\Bigl(A'
\toby{f'_I} D'(I)\Bigr)_{I \in \scat{I}}$ and a map $s\from A \to A'$
such that
\[
\xymatrix{
A \ar[r]^-{f_I} \ar[d]_s        &
D(I) \ar[d]^{\alpha_I}  \\
A' \ar[r]_-{f'_I}       &
D'(I)
}
\]
commutes for all $I \in \scat{I}$, the square
\[
\xymatrix{
A \ar[r]^-{\bar{f}} \ar[d]_s        &
\lt{\scat{I}} D \ar[d]^{\lt{\scat{I}} \alpha}  \\
A' \ar[r]_-{\ovln{f'}}       &
\lt{\scat{I}} D'
}
\]
also commutes.
\end{enumerate}
\end{lemma}

\begin{pf}
Part~\bref{lemma:lim-functorial:induced} follows immediately from the fact
that $\biggl(\lt{\scat{I}} D \toby{\alpha_I p_I} D'(I)\biggr)_{I \in
  \scat{I}}$ is a cone on $D'$.  To
prove~\bref{lemma:lim-functorial:commutes}, note that for each $I \in
\scat{I}$, we have
\[
p'_I \of \biggl(\lt{\scat{I}} \alpha\biggr) \of \bar{f}
=
\alpha_I \of p_I \of \bar{f}
=
\alpha_I \of f_I
=
f'_I \of s
=
p'_I \of \ovln{f'} \of s.
\]
So by Exercise~\ref{ex:jointly-monic}\bref{part:j-m-main},
$\biggl(\lt{\scat{I}} \alpha \biggr) \of \bar{f} = \ovln{f'} \of s$.
\end{pf}

We can now give the second rephrasing of the definition of limit.  It only
applies when the category has \emph{all} limits of the shape concerned.

\begin{propn}   
\label{propn:lim-const-adjn}
\index{limit!adjoint@as adjoint}
Let $\scat{I}$ be a small category and $\cat{A}$ a category with all limits of
shape $\scat{I}$.  Then $\lt{\scat{I}}$ defines a functor
$\ftrcat{\scat{I}}{\cat{A}} \to \cat{A}$, and this functor is right adjoint
to the diagonal functor.  
\end{propn}

\begin{pf}
Choose for each $D \in \ftrcat{\scat{I}}{\cat{A}}$ a limit cone on $D$, and
call its vertex $\lt{\scat{I}} D$.  For each map $\alpha\from D \to D'$ in
$\ftrcat{\scat{I}}{\cat{A}}$, we have a canonical map
$\lt{\scat{I}} \alpha \from \lt{\scat{I}} D \to \lt{\scat{I}} D'$,
defined as in
Lemma~\ref{lemma:lim-functorial}\bref{lemma:lim-functorial:induced}.  This
makes $\lt{\scat{I}}$ into a functor.  Proposition~\ref{propn:lim-rep}
implies that
\[
\ftrcat{\scat{I}}{\cat{A}}(\Delta A, D)
=
\Cone(A, D)
\iso
\cat{A}\biggl(A, \lt{\scat{I}} D\biggr)
\]
naturally in $A$, and taking $s = 1_A$ in
Lemma~\ref{lemma:lim-functorial}\bref{lemma:lim-functorial:commutes} tells
us that the isomorphism is also natural in $D$.
\end{pf}

To define the functor $\lt{\scat{I}}$, we had to \emph{choose}%
\index{limit!uniqueness of}
for each $D$ a limit cone on $D$.  This is a non-canonical choice.
Nevertheless, different choices only affect the functor $\lt{\scat{I}}$ up
to natural isomorphism, by uniqueness of adjoints.

\exs

\begin{question}
Interpret all the theory of this section in the special case where
$\scat{I}$ is the discrete category with two objects.
\end{question}

\begin{question}
What is the content of Proposition~\ref{propn:lim-const-adjn} when
$\scat{I}$ is a group and $\cat{A} = \Set$?  What about the dual of
Proposition~\ref{propn:lim-const-adjn}? 
\end{question}

\section{Limits and colimits of presheaves}
\label{sec:lim-pshf}

What do limits and colimits look like in functor categories
$\ftrcat{\cat{A}}{\cat{B}}$?  In particular, what do they look like in
presheaf categories $\ftrcat{\cat{A}^\op}{\Set}$?  More particularly still,
what about limits and colimits of representables?  Are they, too,
representable?

We will answer all these questions.  In order to do so, we first prove that
representables preserve limits.

\minihead{Representables preserve limits}
\index{functor!representable!preserves limits|(}

Let us begin by recalling that, by definition of product, a map $A \to X
\times Y$ amounts to a pair of maps $(A \to X,\, A \to Y)$.  Here $A$, $X$ and
$Y$ are objects of a category $\cat{A}$ with binary products.  There is,
therefore, a bijection
\begin{equation}
\label{eq:prod-rep}
\cat{A}(A, X \times Y)
\iso 
\cat{A}(A, X) \times \cat{A}(A, Y)
\index{product!map into}
\end{equation}
natural in $A, X, Y \in \cat{A}$.  

Is this a special feature of products, or does some analogous statement
hold for every kind of limit?  Let us try equalizers.  Suppose that
$\cat{A}$ has equalizers, and write
$\Eq\biggl(\parpairi{X}{Y}{s}{t}\biggr)$ for the equalizer of maps $s$ and
$t$.  By definition of equalizer, maps
\begin{equation}
\label{eq:map-to-eq}
A \: \to \: \Eq\biggl(\parpairi{X}{Y}{s}{t}\biggr)
\end{equation}
correspond one-to-one with maps $f\from A \to X$ such that $s \of f = t \of
f$.  Now recall that $s$ induces a map 
\[
s_* = \cat{A}(A, s)\from \cat{A}(A, X) \to \cat{A}(A, Y),
\]
and similarly for $t$.  In this notation, what we have just said is that
maps~\eqref{eq:map-to-eq} correspond one-to-one with elements $f \in
\cat{A}(A, X)$ such that
\[
\bigl(\cat{A}(A, s)\bigr)(f) = \bigl(\cat{A}(A, t)\bigr)(f).
\]
By the explicit formula for equalizers in $\Set$
(Example~\ref{eg:equalizers-Set}), such an $f$ is exactly an element of the
equalizer of $\cat{A}(A, s)$ and $\cat{A}(A, t)$.  So, we have a canonical
bijection
\begin{equation}
\label{eq:equalizer-rep}
\cat{A}
\biggl( 
A,\, \Eq\Bigl(\!\parpair{X}{Y}{s}{t}\!\Bigr)
\biggr)
\iso 
\Eq\biggl(\!
\parpair{\cat{A}(A, X)}{\cat{A}(A, Y)}{\cat{A}(A, s)}{\cat{A}(A, t)}
\!\biggr).
\index{equalizer!map into}
\end{equation}
This looks something like our isomorphism~\eqref{eq:prod-rep} for products.

The isomorphisms~\eqref{eq:prod-rep} and~\eqref{eq:equalizer-rep} suggest
that, more generally, we might have
\begin{equation}        
\label{eq:rep-pres-lims}
\cat{A}\biggl(A, \lt{\scat{I}} D\biggr)
\iso
\lt{\scat{I}} \cat{A}(A, D)
\end{equation}
naturally in $A \in \cat{A}$ and $D \in \ftrcat{\scat{I}}{\cat{A}}$, whenever 
$\cat{A}$ is a category with limits of shape $\scat{I}$.  Here $\cat{A}(A,
D)$%
\ntn{AAD}
is the functor
\[
\begin{array}{cccc}
\cat{A}(A, D)\from      &\scat{I}       &\to        &\Set   \\
                        &I              &\mapsto    &\cat{A}(A, D(I)).
\end{array}
\]
This functor could also be written as $\cat{A}(A, D(\dashbk))$, and is the
composite 
\[
\xymatrix@1@C+1em{
\scat{I} \ar[r]^D       &
\cat{A} \ar[r]^-{\cat{A}(A, \dashbk)}    &
\Set.
}
\]
The conjectured isomorphism~\eqref{eq:rep-pres-lims} states, essentially,
that representables preserve limits.  We now set about proving this.  

\begin{lemma}   
\label{lemma:cone-rep}
Let $\scat{I}$ be a small category, $\cat{A}$ a locally small category, $D\from
\scat{I} \to \cat{A}$ a diagram, and $A \in \cat{A}$.  Then
\[
\Cone(A, D) 
\iso 
\lt{\scat{I}} \cat{A}(A, D)
\index{cone!set of cones as limit}
\]
naturally in $A$ and $D$.
\end{lemma}

\begin{pf}
Like all functors from a small category into $\Set$, the functor
$\cat{A}(A, D)$ does have a limit, given by the explicit
formula~\eqref{eq:Set-lim}.  According to this formula, $\lt{\scat{I}}
\cat{A}(A, D)$ is the set of all families $(f_I)_{I \in \scat{I}}$ such
that $f_I \in \cat{A}(A, D(I))$ for all $I \in \scat{I}$ and
\begin{equation}        
\label{eq:cone-from-lim}
(\cat{A}(A, Du))(f_I) 
= 
f_J
\end{equation}
for all $I \toby{u} J$ in $\scat{I}$.  But equation~\eqref{eq:cone-from-lim}
just says that $(Du) \of f_I = f_J$, so an element of
$\lt{\scat{I}} \cat{A}(A, D)$ is nothing but a cone on $D$ with vertex $A$.  
\end{pf}

\begin{propn}[Representables preserve limits] 
\label{propn:reps-cts}
\hspace*{-5pt}%
Let $\cat{A}$ be a locally small category and $A \in \cat{A}$.  Then
$\cat{A}(A, \dashbk) \from \cat{A} \to \Set$ preserves limits.
\end{propn}

\begin{pf}
Let $\scat{I}$ be a small category and let $D\from \scat{I} \to \cat{A}$ be a
diagram that has a limit.  Then
\[
\cat{A}\biggl(A, \lt{\scat{I}} D\biggr)
\iso
\Cone(A, D)
\iso
\lt{\scat{I}} \cat{A}(A, D)
\]
naturally in $A$.  Here the first isomorphism is
Proposition~\ref{propn:lim-rep} (or more particularly, the
isomorphism~\eqref{eq:cone-is-map-to-lim} that follows it), and the second is
Lemma~\ref{lemma:cone-rep}.
\end{pf}

\begin{remark}  
\label{rmk:rep-pres}
Proposition~\ref{propn:reps-cts} tells us that
\begin{equation}        
\label{eq:lim-out}
\cat{A} \biggl(A, \lt{\scat{I}} D\biggr)
\iso
\lt{\scat{I}} \cat{A}(A, D).
\index{limit!map into}
\end{equation}
To dualize Proposition~\ref{propn:reps-cts}, we replace $\cat{A}$ by
$\cat{A}^\op$.  Thus, $\cat{A}(\dashbk, A)\from \cat{A}^\op \to \Set$
preserves limits.  A limit in $\cat{A}^\op$ is a colimit in $\cat{A}$, so
$\cat{A}(\dashbk, A)$ transforms colimits in $\cat{A}$ into limits in
$\Set$:
\begin{equation}
\label{eq:colim-out}
\cat{A}\biggl(\colt{\scat{I}} D, A\biggr)
\iso
\lt{\scat{I}^{\op}} \cat{A}(D, A).
\index{colimit!map out of}
\end{equation}
The right-hand side is a \emph{limit}, not%
\index{limit!colimit@vs.\ colimit}
a colimit!  So even though~\eqref{eq:lim-out} and~\eqref{eq:colim-out} are
dual statements, there are, in total, more limits than colimits involved.
Somehow, limits have the upper hand.

For example, let $X$, $Y$ and $A$ be objects of a category $\cat{A}$, and
suppose that the sum $X + Y$ exists.  By definition of sum, a map $X + Y
\to A$%
\index{sum!map out of}
amounts to a pair of maps $(X \to A,\, Y \to A)$.  In other words, there is
a canonical isomorphism
\[
\cat{A}(X + Y, A) 
\iso
\cat{A}(X, A) \times \cat{A}(Y, A).
\]
This is the isomorphism~\eqref{eq:colim-out} in the case where $\scat{I}$ is
the discrete category with two objects.%
\index{functor!representable!preserves limits|)}
\end{remark}

\minihead{Limits in functor categories}
\index{functor!category!limits in|(}%
\index{limit!functor category@in functor category|(}%

Earlier, we learned that it is sometimes useful to view functors as objects
in their own right, rather than as maps of categories.  For instance, when
$G$ is a group, functors $G \to \Set$ are $G$-sets
(Example~\ref{eg:functor-action}), which one would usually regard as
`things' rather than `maps'.  This point of view leads to the concept of
functor category.

We now begin an analysis of limits and colimits in functor categories
$\ftrcat{\scat{A}}{\cat{S}}$.  Here $\scat{A}$ is small and $\cat{S}$ is
locally small; these conditions together guarantee that
$\ftrcat{\scat{A}}{\cat{S}}$ is locally small.  The most important cases
for us will be $\cat{S} = \Set$ and $\cat{S} = \Set^\op$.  For that reason,
we will assume whenever necessary that $\cat{S}$ has all limits and
colimits.

We show that limits and colimits in $\ftrcat{\scat{A}}{\cat{S}}$ work in
the simplest way imaginable.  For instance, if $\cat{S}$ has binary products
then so does $\ftrcat{\scat{A}}{\cat{S}}$, and the product%
\index{functor!product of functors}
of two functors $X, Y\from \scat{A} \to \cat{S}$ is the functor $X \times
Y\from \scat{A} \to \cat{S}$ given by
\[
(X \times Y)(A) = X(A) \times Y(A)
\]
for all $A \in \scat{A}$.  

\begin{notn}
Let $\scat{A}$ and $\cat{S}$ be categories.  For each $A \in \scat{A}$, there
is a functor
\[
\begin{array}{cccc}
\ev_A\from      &\ftrcat{\scat{A}}{\cat{S}}     &
\to             &\cat{S}        \\
                &X                              &
\mapsto         &X(A),
\end{array}%
\ntn{ev}
\]
called \demph{evaluation}%
\index{evaluation}
at $A$.  We will be working with diagrams in $\ftrcat{\scat{A}}{\cat{S}}$,
and given such a diagram $D\from \scat{I} \to \ftrcat{\scat{A}}{\cat{S}}$,
we have for each $A \in \scat{A}$ a functor
\[
\begin{array}{cccc}
\ev_A \of D\from&\scat{I}       &\to            &\cat{S}        \\
                &I              &\mapsto        &D(I)(A).        
\end{array}
\]
We write $\ev_A \of D$ as $D(\dashbk)(A)$.%
\ntn{DblankA}
\end{notn}

\begin{thm}[Limits in functor categories]
\label{thm:pw}
Let $\scat{A}$ and $\scat{I}$ be small categories and $\cat{S}$ a locally
small category.  Let $D\from \scat{I} \to \ftrcat{\scat{A}}{\cat{S}}$ be a
diagram, and suppose that for each $A \in \scat{A}$, the diagram
$D(\dashbk)(A)\from \scat{I} \to \cat{S}$ has a limit.  Then there is a
cone on $D$ whose image under $\ev_A$ is a limit cone on $D(\dashbk)(A)$
for each $A \in \scat{A}$.  Moreover, any such cone on $D$ is a limit cone.
\end{thm}

Theorem~\ref{thm:pw} is often expressed as a slogan:
\begin{slogan}
Limits in a functor category are computed pointwise.%
\index{limit!computed pointwise}%
\index{pointwise}
\end{slogan}
The `points' in the word `pointwise' are the objects of $\scat{A}$.  The
slogan means, for example, that given two functors $X, Y \in
\ftrcat{\scat{A}}{\cat{S}}$, their product can be computed by first taking
the product $X(A) \times Y(A)$ in $\cat{S}$ for each `point' $A$, then
assembling them to form a functor $X \times Y$.

Of course, Theorem~\ref{thm:pw} has a dual, stating that colimits in a
functor category are also computed pointwise.

\begin{pfof}{Theorem~\ref{thm:pw}}
Take for each $A \in \scat{A}$ a limit cone 
\begin{equation}
\label{eq:pw-cone}
\Bigl( L(A) \toby{p_{I, A}} D(I)(A) \Bigr)_{I \in \scat{I}}
\end{equation}
on the diagram $D(\dashbk)(A)\from \scat{I} \to \cat{S}$.  We prove two
statements:
\begin{enumerate}[(b)]
\item 
\label{item:pw-lift} 
there is exactly one way of extending $L$ to a functor on $\scat{A}$ with
the property that $\Bigl(L \toby{p_I} D(I)\Bigr)_{I \in \scat{I}}$ is a
cone on $D$;  

\item   
\label{item:pw-lim}
this cone $\Bigl(L \toby{p_I} D(I)\Bigr)_{I \in \scat{I}}$ is a limit cone.
\end{enumerate}
The theorem will follow immediately.

For~\bref{item:pw-lift}, take a map $f\from A \to A'$ in $\scat{A}$.
Lemma~\ref{lemma:lim-functorial}\bref{lemma:lim-functorial:induced} applied
to the natural transformation
\[
\xymatrix@C+6em{
\scat{I} \rtwocell<5>^{D(\dashbk)(A)}_{D(\dashbk)(A')}%
{\hspace{2.5em}D(\dashbk)(f)}      &
\cat{S}
}
\]
implies that there is a unique map $L(f)\from L(A) \to L(A')$ such that for
all $I \in \scat{I}$, the square
\begin{equation}
\label{eq:param-lim}
\begin{array}{c}
\xymatrix{
L(A) \ar[r]^-{p_{I, A}} \ar[d]_{L(f)}   &
D(I)(A) \ar[d]^{D(I)(f)}        \\
L(A') \ar[r]_-{p_{I, A'}}       &
D(I)(A')
}
\end{array}
\end{equation}
commutes.  (This is our \emph{definition} of $L(f)$.)  We have now defined $L$
on objects and maps of $\scat{A}$.  It is easy to check that $L$ preserves
composition and identities, and is therefore a functor $L\from \scat{A}
\to \cat{S}$.  Moreover, the commutativity of diagram~\eqref{eq:param-lim}
says exactly that for each $I \in \scat{I}$, the family $\Bigl(L(A)
\toby{p_{I, A}} D(I)(A)\Bigr)_{A \in \scat{A}}$ is a natural
transformation
\[
\xymatrix@C+1em{
\scat{A} \rtwocell^{L}_{D(I)}{\hspace{.3em}p_I} &\cat{S}.
}
\]
So we have a family $\Bigl(L \toby{p_I} D(I)\Bigr)_{I \in \scat{I}}$ of
maps in $\ftrcat{\scat{A}}{\cat{S}}$, and from the fact
that~\eqref{eq:pw-cone} is a cone on $D(\dashbk)(A)$ for each $A \in
\scat{A}$, it follows immediately that $\Bigl(L \toby{p_I} D(I)\Bigr)_{I
  \in \scat{I}}$ is a cone on $D$.

For~\bref{item:pw-lim}, let $X \in \ftrcat{\scat{A}}{\cat{S}}$ and let
$\Bigl(X \toby{q_I} D(I)\Bigr)_{I \in \scat{I}}$ be a cone on $D$ in
$\ftrcat{\scat{A}}{\cat{S}}$.  For each $A \in \scat{A}$, we have a cone
\[
\Bigl(X(A) \toby{q_{I, A}} D(I)(A)\Bigr)_{I \in \scat{I}}
\]
on
$D(\dashbk)(A)$ in $\cat{S}$, so there is a unique map $\bar{q}_A\from X(A)
\to L(A)$ such that $p_{I, A} \of \bar{q}_A = q_{I, A}$ for all $I \in
\scat{I}$.  It only remains to prove that $\bar{q}_A$ is natural in $A$, and
that follows from
Lemma~\ref{lemma:lim-functorial}\bref{lemma:lim-functorial:commutes}.  
\end{pfof}

Theorem~\ref{thm:pw} has many important consequences.  We begin by recording
a cruder form of the theorem (and its dual), which we will use repeatedly.

\begin{cor}     
\label{cor:pw-main}
Let $\scat{I}$ and $\scat{A}$ be small categories, and $\cat{S}$ a locally
small category.  If $\cat{S}$ has all limits (respectively, colimits) of
shape $\scat{I}$ then so does $\ftrcat{\scat{A}}{\cat{S}}$, and for each $A
\in \scat{A}$, the evaluation functor $\ev_A\from
\ftrcat{\scat{A}}{\cat{S}} \to \cat{S}$ preserves them.
\qed
\end{cor}

\begin{warning}
If $\cat{S}$ does \emph{not} have all limits of shape $\scat{I}$ then
$\ftrcat{\scat{A}}{\cat{S}}$ may contain limits of shape $\scat{I}$
that are not%
\index{limit!non-pointwise}
computed pointwise, that is, are not preserved by all the evaluation
functors.  Examples can be constructed, as in Section~3.3 of \citeKel.
\end{warning}

Theorem~\ref{thm:pw} will also help us to prove that limits commute with
limits, in the following sense.  Take categories $\scat{I}$, $\scat{J}$ and
$\cat{S}$.  There are isomorphisms of categories
\[
\ftrcat{\scat{I}}{\ftrcat{\scat{J}}{\cat{S}}}
\iso 
\ftrcat{\scat{I} \times \scat{J}}{\cat{S}}
\iso
\ftrcat{\scat{J}}{\ftrcat{\scat{I}}{\cat{S}}}.
\]
(See Remark~\ref{rmks:global-hom}\bref{rmks:global-hom:cc} and
Exercise~\ref{ex:ftr-on-product}.)  Under these isomorphisms, a functor
$D\from \scat{I} \times \scat{J} \to \cat{S}$ corresponds to the functors
\[
\begin{array}[t]{cccc}
D^\bl\from  &\scat{I}       &\to    &\ftrcat{\scat{J}}{\cat{S}}    \\
        &I              &\mapsto&D(I, \dashbk)
\end{array}
\qquad \text{and} \qquad
\begin{array}[t]{cccc}
D_\bl\from  &\scat{J}       &\to    &\ftrcat{\scat{I}}{\cat{S}}    \\
        &J              &\mapsto&D(\dashbk, J).
\end{array}%
\ntn{D-upper-lower}
\]
Supposing that $\cat{S}$ has all limits, so do the various functor
categories, by Corollary~\ref{cor:pw-main}.  In particular, there is an
object $\lt{\scat{I}} D^\bl$ of $\ftrcat{\scat{J}}{\cat{S}}$.  This is
itself a diagram in $\cat{S}$, so we obtain in turn an object
$\lt{\scat{J}} \lt{\scat{I}} D^\bl$ of $\cat{S}$.  Alternatively, we can
take limits in the other order, producing an object $\lt{\scat{I}}
\lt{\scat{J}} D_\bl$ of $\cat{S}$.  And there is a third possibility:
taking the limit of $D$ itself, we obtain another object $\lt{\scat{I}
  \times \scat{J}} D$ of $\cat{S}$.  The next result states that these
three objects are the same.  That is, it makes no difference what order we
take limits in.

\begin{propn}[Limits commute with limits] 
\label{propn:lims-lims}
\index{limit!commutativity with limits}
Let $\scat{I}$ and $\scat{J}$ be small categories.  Let $\cat{S}$ be a
locally small category with limits of shape $\scat{I}$ and of shape
$\scat{J}$.  Then for all $D\from \scat{I} \times \scat{J} \to \cat{S}$, we
have 
\[
\lt{\scat{J}} \lt{\scat{I}} D^\bl
\iso
\lt{\scat{I} \times \scat{J}} D
\iso 
\lt{\scat{I}} \lt{\scat{J}} D_\bl,
\]
and all these limits exist.  In particular, $\cat{S}$ has limits of shape
$\scat{I} \times \scat{J}$.  
\end{propn}
This is sometimes half-jokingly called Fubini's%
\index{Fubini's theorem}
theorem, as it is something like changing the order of integration in a
double integral.  The analogy is more appealing with \emph{co}limits,%
\index{colimit!integration@and integration}
since, like integrals, colimits can be thought of as a context-sensitive
version of sums.

\begin{pf}
By symmetry, it is enough to prove the first isomorphism.  Since $\cat{S}$
has limits of shape $\scat{I}$, so does $\ftrcat{\scat{J}}{\cat{S}}$ (by
Corollary~\ref{cor:pw-main}).  So $\lt{\scat{I}} D^\bl$ exists; it is an
object of $\ftrcat{\scat{J}}{\cat{S}}$.  Since $\cat{S}$ has limits of
shape $\scat{J}$, $\lt{\scat{J}} \lt{\scat{I}} D^\bl$ exists; it is an
object of $\cat{S}$.  Then for $S \in \cat{S}$,
\begin{align*}
\cat{S}\biggl( S, \lt{\scat{J}} \lt{\scat{I}} D^\bl \biggr)     &
\iso    
\ftrcat{\scat{J}}{\cat{S}}
\biggl( \Delta S, \lt{\scat{I}} D^\bl \biggr)                   \\
        &
\iso   
\ftrcat{\scat{I}}{\ftrcat{\scat{J}}{\cat{S}}} 
(\Delta (\Delta S), D^\bl)                                      \\
        &
\iso   
\ftrcat{\scat{I} \times \scat{J}}{\cat{S}} (\Delta S, D)
\end{align*}
naturally in $S$.  The first two steps each follow from
Proposition~\ref{propn:lim-rep}.  The third uses the isomorphism
$\ftrcat{\scat{I}}{\ftrcat{\scat{J}}{\cat{S}}} \iso \ftrcat{\scat{I} \times
\scat{J}}{\cat{S}}$, under which $\Delta(\Delta S)$ corresponds to $\Delta S$
and $D^\bl$ corresponds to $D$.

Hence $\lt{\scat{J}} \lt{\scat{I}} D^\bl$ is
a representing object for the functor $\ftrcat{\scat{I} \times
\scat{J}}{\cat{S}}(\Delta\dashbk, D)$.  By Proposition~\ref{propn:lim-rep}
again, this says that $\lt{\scat{I} \times \scat{J}} D$ exists and is
isomorphic to $\lt{\scat{J}} \lt{\scat{I}} D^\bl$.
\end{pf}

\begin{example}        
\label{eg:lims-lims}
When $\scat{I} = \scat{J} = \fbox{$\bullet \hspace*{2.5em} \bullet$}\ $,
Proposition~\ref{propn:lims-lims} says that binary products commute with
binary products: if $\cat{S}$ has binary products and $S_{11}, S_{12},
S_{21}, S_{22} \in \cat{S}$ then the 4-fold product $\prod_{i, j \in \{1,
  2\}} S_{ij}$ exists and satisfies
\[
(S_{11} \times S_{21}) \times (S_{12} \times S_{22})
\iso
\prod_{i, j \in \{1, 2\}} S_{ij}
\iso
(S_{11} \times S_{12}) \times (S_{21} \times S_{22}).
\]
More generally, it makes no difference what order we write products in or
where we put the brackets: there are canonical isomorphisms
\begin{align*}
S \times T              &
\iso   
T \times S,             
\index{product!commutativity of}
\\
(S \times T) \times U   &
\iso   
S \times (T \times U)
\index{product!associativity of}%
\index{associativity}
\end{align*}
in any category with binary products.   If there is also a terminal object
$1$, there are further canonical isomorphisms
\[
S \times 1 
\iso 
S 
\iso 
1 \times S.
\]
\end{example}

\begin{warning}
The dual of Proposition~\ref{propn:lims-lims} states that colimits commute
with colimits.  For instance,
\[
(S_{11} + S_{21}) + (S_{12} + S_{22})
\iso
(S_{11} + S_{12}) + (S_{21} + S_{22})
\]
in any category $\cat{S}$ with binary sums.  But limits do \emph{not}%
\index{limit!commutativity with colimits@non-commutativity with colimits}
in general commute with colimits.  For instance, in general,
\[
(S_{11} + S_{21}) \times (S_{12} + S_{22})
\not\iso
(S_{11} \times S_{12}) + (S_{21} \times S_{22}).
\]
A counterexample is given by taking $\cat{S} = \Set$ and each $S_{ij}$ to be a
one-element set.  Then the left-hand side has $(1 + 1) \times (1 + 1) = 4$
elements, whereas the right-hand side has $(1 \times 1) + (1 \times 1) = 2$
elements.
\end{warning}

Here are two further consequences of Theorem~\ref{thm:pw}.

\begin{cor}
\index{presheaf!category of presheaves!limits in}
Let $\scat{A}$ be a small category.  Then $\pshf{\scat{A}}$ has all limits and
colimits, and for each $A \in \scat{A}$, the evaluation functor $\ev_A\from
\pshf{\scat{A}} \to \Set$ preserves them.
\end{cor}

\begin{pf}
Since $\Set$ has all limits and colimits, this is immediate from
Corollary~\ref{cor:pw-main}.
\end{pf}

\begin{cor}     
\label{cor:yoneda-cts}
\index{Yoneda embedding!preserves limits}%
\index{functor!representable!limit of representables|(}
The Yoneda embedding $\h_\bl\from \scat{A} \to \pshf{\scat{A}}$ preserves
limits, for any small category $\scat{A}$.
\end{cor}

\begin{pf}
Let $D\from \scat{I} \to \scat{A}$ be a diagram in $\scat{A}$, and let
$\biggl( \lt{\scat{I}} D \toby{p_I} D(I) \biggr)_{I \in \scat{I}}$ be a
limit cone.  For each $A \in \scat{A}$, the composite functor
\[
\scat{A} \toby{\h_\bl} \pshf{\scat{A}} \toby{\ev_A} \Set
\]
is $\h^A$, which preserves limits (Proposition~\ref{propn:reps-cts}).  So for
each $A \in \scat{A}$, 
\[
\biggl( 
\ev_A \h_\bl \biggl( \lt{\scat{I}} D \biggr) \,
\xymatrix@1@C+2em{\mbox{} \ar[r]^{\ev_A \h_\bl (p_I)}&\mbox{}}
\ev_A \h_\bl (D(I)) 
\biggr)_{I \in \scat{I}}
\]
is a limit cone.  But then, by the `moreover' part of Theorem~\ref{thm:pw}
applied to the diagram $\h_\bl \of D$ in $\pshf{\scat{A}}$, the cone
\[
\biggl( 
\h_\bl \biggl(\lt{\scat{I}} D\biggr)
\toby{\h_\bl (p_I)} 
\h_\bl (D(I)) 
\biggr)_{I \in \scat{I}}
\]
is also a limit, as required.
\end{pf}

\begin{example}
Let $\scat{A}$ be a category with binary products.
Corollary~\ref{cor:yoneda-cts} implies that for all $X, Y \in \scat{A}$,
\begin{equation}        
\label{eq:yon-pres-prods}
\h_{X \times Y} \iso \h_X \times \h_Y
\end{equation}
in $\ftrcat{\scat{A}^\op}{\Set}$.  When evaluated at a particular object
$A$, this says that
\[
\scat{A}(A, X \times Y) \iso \scat{A}(A, X) \times \scat{A}(A, Y)
\index{product!map into}
\]
(using the fact that products are computed pointwise).  This is the
isomorphism~\eqref{eq:prod-rep} that we met at the beginning of this
section.

Suppose that we view $\scat{A}$ as a subcategory of $\pshf{\scat{A}}$,
identifying $A \in \scat{A}$ with the representable $\h_A \in
\pshf{\scat{A}}$ as in Figure~\ref{fig:yoneda-embedding}.  Then the
isomorphism~\eqref{eq:yon-pres-prods} means that given two objects of
$\scat{A}$ whose product we want to form, it makes no difference whether we
think of the product as taking place in $\scat{A}$ or $\pshf{\scat{A}}$.
Similarly, if $\scat{A}$ has all limits, taking limits does not help us to
escape from $\scat{A}$ into the rest of $\pshf{\scat{A}}$: any limit of
representable presheaves is again representable.%
\index{functor!representable!limit of representables|)}
\end{example}

\begin{warning} 
\label{warning:yon-colims}
\index{Yoneda embedding!preserves colimits@does not preserve colimits}
The Yoneda embedding does \emph{not} preserve colimits.  For example, if
$\scat{A}$ has an initial object $0$ then $\h_0$ is not initial, since
$\h_0(0) = \scat{A}(0, 0)$ is a one-element set, whereas the initial object
of $\pshf{\scat{A}}$ is the presheaf with constant value $\emptyset$.  We
investigate colimits of representables next.
\index{functor!category!limits in|)}%
\index{limit!functor category@in functor category|)}%
\end{warning}

\minihead{Every presheaf is a colimit of representables}
\index{presheaf!colimit of representables@as colimit of representables|(}%
\index{functor!representable!colimit of representables|(}%

We now know that the Yoneda embedding preserves limits but not colimits.
In fact, the situation for colimits is at the opposite extreme from the
situation for limits: by taking colimits of representable presheaves, we
can obtain any presheaf we like!  This is the last main result of this
section.

Every positive integer can be expressed as a product of primes%
\index{prime numbers}
in an essentially unique way.  Somewhat similarly, every presheaf can be
expressed as a colimit of representables in a canonical (though not unique)
way.  The representables are the building blocks of presheaves.

For a different analogy, recall that any complex function holomorphic%
\index{holomorphic function}
in a neighbourhood of $0$ has a power%
\index{power!series}
series expansion, such as
\[
e^z 
=
1 + z + \frac{z^2}{2!} + \frac{z^3}{3!} + \cdots.
\]
In this sense, the power functions $z \mapsto z^n$ are the building blocks
of holomorphic functions.  We could even take the analogy further:
$\blank^n$ is like a representable $\Hom(n, \dashbk)$, and in the
categorical context, quotients and sums are types of colimit.

Before we state and prove the theorem, let us look at an easy special case.

\begin{example} 
\label{eg:discrete-dense}
Let $\scat{A}$ be the discrete category with two objects, $K$ and $L$.  A
presheaf $X$ on $\scat{A}$ is just a pair $(X(K), X(L))$ of sets, and
$\pshf{\scat{A}} \iso \Set \times \Set$.  There are two representables,
$\h_K$ and $\h_L$, given by
\[
\h_A(B) 
= 
\scat{A}(B, A)
\iso
\begin{cases}
1               &\text{if } A = B,    \\
\emptyset       &\text{if } A \neq B  
\end{cases}
\]
($A, B \in \{ K, L \}$).  Identifying $\pshf{\scat{A}}$ with $\Set \times
\Set$, we have $\h_K \iso (1, \emptyset)$ and $\h_L \iso (\emptyset, 1)$.
Every object of $\Set \times \Set$ is a sum of copies of $(1, \emptyset)$
and $(\emptyset, 1)$.  Suppose, for instance, that $X(K)$ has three elements
and $X(L)$ has two elements.  Then
\[
(X(K), X(L))
\iso
(1, \emptyset) + (1, \emptyset) + (1, \emptyset) + 
(\emptyset, 1) + (\emptyset, 1)
\]
in $\Set \times \Set$.  Equivalently,
\[
X 
\iso 
\h_K + \h_K + \h_K + \h_L + \h_L 
\]
in $\pshf{\scat{A}}$, exhibiting $X$ as a sum of representables.  
\end{example}

In this example, $X$ is expressed as a sum of five representables, that is,
a sum indexed by the set $X(K) + X(L)$ of `elements' of $X$.  A sum is a
colimit over a discrete category.  In the general case, a presheaf $X$ on a
category $\scat{A}$ is expressed as a colimit over a category whose objects
can be thought of as the `elements' of $X$.  This is made precise by the
following definition.

\begin{defn}    
\label{defn:cat-elts}
Let $\scat{A}$ be a category and $X$ a presheaf on $\scat{A}$.  The
\demph{category%
\index{element!category of elements}%
\index{category!elements@of elements}
of elements} $\elt{X}$%
\ntn{cat-elts}
of $X$ is the category in which:
\begin{itemize}
\item 
objects are pairs $(A, x)$ with $A \in \scat{A}$ and $x \in X(A)$;

\item 
maps $(A', x') \to (A, x)$ are maps $f\from A' \to A$ in $\scat{A}$ such
that $(Xf)(x) = x'$.
\end{itemize}
\end{defn}
There is a projection functor $P\from \elt{X} \to \scat{A}$ defined by $P(A, x)
= A$ and $P(f) = f$.

The following `density theorem' states that every presheaf is a colimit of
representables in a canonical way.  It is secretly dual to the Yoneda
lemma.  This becomes apparent if one expresses both in suitably lofty
categorical language (that of ends, or that of bimodules); but that is
beyond the scope of this book.

\begin{thm}[Density]
\label{thm:density}
\index{density}
Let $\scat{A}$ be a small category and $X$ a presheaf on $\scat{A}$.  Then $X$
is the colimit of the diagram
\[
\elt{X} \toby{P} \scat{A} \toby{\h_\bl} \pshf{\scat{A}}
\] 
in $\pshf{\scat{A}}$; that is, $X \iso \colt{\elt{X}} (\h_\bl \of P)$.
\end{thm}

\begin{pf}
First note that since $\scat{A}$ is small, so too is $\elt{X}$.  Hence $\h_\bl
\of P$ really is a diagram in our customary sense
(Definition~\ref{defn:diagram}). 

Now let $Y \in \pshf{\scat{A}}$.  A cocone on $\h_\bl \of P$ with vertex $Y$ is
a family
\[
\Bigl(\h_A \toby{\alpha_{A, x}} Y\Bigr)_{A \in \scat{A}, x \in X(A)}
\]
of natural transformations with the property that for all maps $A'
\toby{f} A$ in $\scat{A}$ and all $x \in X(A)$, the diagram
\[
\xymatrix@R=1ex{
\h_{A'} \ar[rd]^-{\alpha_{A', (Xf)(x)}} \ar[dd]_{\h_f}   &       \\
                                                        &Y      \\
\h_A \ar[ru]_-{\alpha_{A, x}}
}
\]
commutes.  

Equivalently (by the Yoneda lemma), a cocone on $\h_\bl \of P$ with vertex $Y$
is a family 
\[
(y_{A, x})_{A \in \scat{A}, x \in X(A)},
\]
with $y_{A, x} \in Y(A)$, such that for all maps $A' \toby{f} A$ in
$\scat{A}$ and all $x \in X(A)$,
\[
(Yf)(y_{A, x}) = y_{A', (Xf)(x)}.
\]
To see this, note that if $\alpha_{A, x} \in \pshf{\scat{A}}(\h_A, Y)$
corresponds to $y_{A, x} \in Y(A)$, then $\alpha_{A, x} \of \h_f \in
\pshf{\scat{A}}(\h_{A'}, Y)$ corresponds to $(Yf)(y_{A, x}) \in Y(A')$.

Equivalently (writing $y_{A, x}$ as $\bar{\alpha}_A(x)$), it is a family
\[
\Bigl(X(A) \toby{\bar{\alpha}_A} Y(A)\Bigr)_{A \in \scat{A}}
\] 
of functions with the property that for all maps $A' \toby{f} A$
in $\scat{A}$ and all $x \in X(A)$,
\[
(Yf)\bigl(\bar{\alpha}_A(x)\bigr) 
=
\bar{\alpha}_{A'}\bigl((Xf)(x)\bigr).
\]
But this is simply a natural transformation $\bar{\alpha}\from X \to Y$.
So we have, for each $Y \in \pshf{\scat{A}}$, a canonical bijection
\[
\ftrcat{\elt{X}}{\pshf{\scat{A}}} (\h_\bl \of P, \, \Delta Y)
\iso 
\pshf{\scat{A}} (X, Y).
\]
Hence $X$ is the colimit of $\h_\bl \of P$.
\end{pf}

\begin{example}
In Example~\ref{eg:discrete-dense}, we expressed a particular presheaf $X$
as a sum of representables.  Let us check that the way we did this is a
special case of the general construction in the density theorem.

Since $\scat{A}$ is discrete, the category of elements $\elt{X}$ is also
discrete; it is the set $X(K) + X(L)$ with five elements.  The projection
$P\from \elt{X} \to \scat{A}$ sends three of the elements to $K$ and the
other two to $L$, so the diagram $\h_\bl \of P\from \elt{X} \to
\pshf{\scat{A}}$ sends three of the elements to $\h_K$ and two to $\h_L$.
The colimit of $\h_\bl \of P$ is the sum of these five representables,
which is $X$, just as in Example~\ref{eg:discrete-dense}.
\end{example}

\begin{remarks}
\begin{enumerate}[(b)]
\item 
The term `category of elements'%
\index{element!category of elements}%
\index{category!elements@of elements}
is compatible with the generalized%
\index{element!generalized}
element terminology introduced in Definition~\ref{defn:gen-elt}.  A
generalized element of an object $X$ is just a map into $X$, say $Z \to X$;
but, as explained after that definition, we often focus on certain special
shapes $Z$.  Now suppose that we are working in a presheaf category
$\pshf{\scat{A}}$.  Among all presheaves, the representables have a special
status, so we might be especially interested in generalized elements of
representable shape.  The Yoneda lemma implies that for a presheaf $X$, the
generalized elements of $X$ of representable shape correspond to pairs $(A,
x)$ with $A \in \scat{A}$ and $x \in X(A)$.  In other words, they are the
objects of the category of elements.

\item 
In topology, a subspace $A$ of a space $B$ is called dense%
\index{density}
if every point in $B$ can be obtained as a limit of points in $A$.  This
provides some explanation for the name of Theorem~\ref{thm:density}: the
category $\scat{A}$ is `dense' in $\pshf{\scat{A}}$ because every object of
$\pshf{\scat{A}}$ can be obtained as a colimit of objects of $\scat{A}$.
\index{presheaf!colimit of representables@as colimit of representables|)}%
\index{functor!representable!colimit of representables|)}%
\end{enumerate}
\end{remarks}

\exs

\begin{question}
Fix a small category $\scat{A}$.
\begin{enumerate}[(b)]
\item
Let $\cat{S}$ be a locally small category with pullbacks.  Show that a
natural transformation
\[
\xymatrix{
\scat{A} \rtwocell^X_Y{\alpha} &\cat{S}
}
\]
is monic (as a map in $\ftrcat{\scat{A}}{\cat{S}}$) if and only if $\alpha_A$
is monic for all $A \in \scat{A}$.  (Hint: use Lemma~\ref{lemma:monic-pb}.)

\item   
\label{part:monic-epic-transf}
Describe explicitly the monics and epics in $\pshf{\scat{A}}$.%
\index{presheaf!category of presheaves!monics and epics in}%

\item 
Can you do part~\bref{part:monic-epic-transf} without relying on the fact
that limits and colimits of presheaves are computed pointwise?
\end{enumerate}
\end{question}

\begin{question}
\begin{enumerate}[(b)]
\item 
Prove that representables have the following connectedness%
\index{connectedness}
property: given a locally small category $\cat{A}$ and $A \in \cat{A}$, if
$X, Y \in \pshf{\cat{A}}$ with $\h_A \iso X + Y$, then either $X$ or $Y$ is
the constant functor $\emptyset$.

\item
Deduce that the sum%
\index{functor!representable!sum of representables}
of two representables is never representable.
\end{enumerate}
\end{question}

\begin{question}
Show how a category of elements can be described as a comma category.
\end{question}

\begin{question}
Let $X$ be a presheaf on a locally small category.  Show that $X$ is
representable if and only if its category of elements has a terminal
object.

(Since a terminal object is a limit of the empty diagram, this implies that
the concept of representability can be derived from the concept of limit.
Since a terminal object of a category $\cat{E}$ is also a right adjoint to
the unique functor $\cat{E} \to \One$, the concept of representability can
also be derived from the concept of adjoint.)
\end{question}

\begin{question}
Prove that every slice%
\index{presheaf!category of presheaves!slice of}%
\index{slice category!presheaf category@of presheaf category}
of a presheaf category is again a presheaf category.
That is, given a small category $\scat{A}$ and a presheaf $X$ on $\scat{A}$,
prove that $\pshf{\scat{A}}/X$ is equivalent to $\pshf{\scat{B}}$ for some
small category $\scat{B}$.
\end{question}

\begin{question}
\label{ex:kan}
Let $F\from \scat{A} \to \scat{B}$ be a functor between small categories.
For each object $B \in \scat{B}$, there is a comma category $\comma{F}{B}$
(defined dually to the comma category in Example~\ref{eg:comma-obj-ftr}),
and there is a projection functor $P_B \from \comma{F}{B} \to \scat{A}$.
\begin{enumerate}[(b)]
\item 
Let $X \from \scat{A} \to \cat{S}$ be a functor from $\scat{A}$ to a
category $\cat{S}$ with small colimits.  For each $B \in \scat{B}$, let
$(\Lan_F X)(B)$ be the colimit of the diagram
\[
\comma{F}{B} \toby{P_B} \scat{A} \toby{X} \cat{S}.
\]
Show that this defines a functor $\Lan_F X\from \scat{B} \to \cat{S}$, and
that for functors $Y \from \scat{B} \to \cat{S}$, there is a canonical
bijection between natural transformations $\Lan_F X \to Y$ and natural
transformations $X \to Y \of F$.

\item
\label{part:left-kan-exists}
Deduce that for any category $\cat{S}$ with small colimits, the functor
\[
\dashbk \of F \from 
\ftrcat{\scat{B}}{\cat{S}} \to
\ftrcat{\scat{A}}{\cat{S}}
\]
has a left adjoint.  (This left adjoint, $\Lan_F$, is called
\demph{left Kan%
\index{Kan extension}
extension} along $F$.)

\item
Part~\bref{part:left-kan-exists} and its dual imply that when $\cat{S}$
has small limits and colimits, the functor $\dashbk \of F$ has both left
and right adjoints.  Revisit Exercise~\ref{ex:G-set-adjns} with this
in mind, taking $F$ to be either the unique functor $\One \to G$%
\index{group!action of}%
\index{G-set@$G$-set}%
\index{representation!group or monoid@of group or monoid!linear}
or the unique functor $G \to \One$.
\end{enumerate}
\end{question}

\section{Interactions between adjoint functors and limits}
\label{sec:adj-lim}

We saw in Proposition~\ref{propn:ladj-rep} that any set-valued functor with
a left adjoint is representable, and in Proposition~\ref{propn:reps-cts}
that any representable preserves limits.  Hence, any set-valued functor
with a left adjoint preserves limits.  In fact, this conclusion holds not
only for set-valued functors, but in complete generality.

\begin{thm}     
\label{thm:adjts-cts}
\index{limit!preservation of!adjoint@by adjoint}%
\index{adjunction!limits preserved in}
Let $\hadjnli{\cat{A}}{\cat{B}}{F}{G}$ be an adjunction.  Then $F$ preserves
colimits and $G$ preserves limits.
\end{thm}

\begin{pf}
By duality, it is enough to prove that $G$ preserves limits.  Let $D\from
\scat{I} \to \cat{B}$ be a diagram for which a limit exists.  Then
\begin{align}
\cat{A} \biggl(A, G\biggl(\lt{\scat{I}} D \biggr)\biggr) &
\iso   
\cat{B}\biggl(F(A), \lt{\scat{I}} D\biggr)        
\label{eq:adj-lim-1}    \\
        &
\iso   
\lt{\scat{I}} \cat{B}(F(A), D)        
\label{eq:adj-lim-2}    \\
        &
\iso   
\lt{\scat{I}} \cat{A}(A, G \of D)      
\label{eq:adj-lim-3}    \\
        &
\iso   
\Cone(A, G \of D)
\label{eq:adj-lim-4}
\end{align}
naturally in $A \in \cat{A}$.  Here, the isomorphism~\eqref{eq:adj-lim-1}
is by adjointness, \eqref{eq:adj-lim-2} is because representables preserve
limits, \eqref{eq:adj-lim-3} is by adjointness again, and
\eqref{eq:adj-lim-4} is by Lemma~\ref{lemma:cone-rep}.  So
$G\biggl(\lt{\scat{I}} D\biggr)$ represents $\Cone(\dashbk, G \of D)$; that
is, it is a limit of $G \of D$.
\end{pf}

\begin{example}
Forgetful functors from categories of algebras to $\Set$ have left
adjoints, but hardly ever right adjoints.  Correspondingly, they preserve%
\index{functor!forgetful!preserves limits}
all limits, but rarely all colimits.
\end{example}

\begin{example}
\label{eg:arith-set}
Every set $B$ gives rise to an adjunction $(\dashbk \times B) \ladj
(\dashbk)^B$ of functors from $\Set$ to $\Set$ (Example~\ref{eg:adjn:cc}).  So
$\dashbk \times B$ preserves colimits  and
$(\dashbk)^B$ preserves limits.  In particular, $\dashbk \times B$ preserves
finite sums and $(\dashbk)^B$ preserves finite products, giving isomorphisms
\begin{align}
\label{eq:Set-dist}
0 \times B &\iso 0,
&
(A_1 + A_2) \times B &\iso (A_1 \times B) + (A_2 \times B),      \\
\label{eq:Set-codist}
1^B &\iso 1,
&
(A_1 \times A_2)^B &\iso A_1^B \times A_2^B.
\end{align}
These are the analogues of standard rules of arithmetic.%
\index{arithmetic}
(See also Example \ref{eg:lims-lims} and the `Digression on arithmetic' on
page~\pageref{p:arith}.)  Indeed, if we know~\eqref{eq:Set-dist}
and~\eqref{eq:Set-codist} for just finite sets then by taking cardinality
on both sides, we obtain exactly these standard rules.  The natural%
\index{natural numbers}
numbers are, after all, just the isomorphism classes of finite sets.
\end{example}

\begin{example}
Given a category $\cat{A}$ with all limits of shape $\scat{I}$, we have the
adjunction
$\hadjnli{\cat{A}}{\ftrcat{\scat{I}}{\cat{A}}}{\Delta}{\lt{\scat{I}}}$
(Proposition~\ref{propn:lim-const-adjn}).  Hence $\lt{\scat{I}}$ preserves
limits, or equivalently, limits of shape $\scat{I}$ commute with (all)
limits.  This gives another proof that limits commute%
\index{limit!commutativity with limits}
with limits (Proposition~\ref{propn:lims-lims}), at least in the case where
the category has all limits of one of the shapes concerned.
\end{example}

\begin{example} 
\label{eg:no-free-field}
Theorem~\ref{thm:adjts-cts} is often used to prove that a functor does
\emph{not}%
\index{adjunction!nonexistence of adjoints}
have an adjoint.  For instance, it was claimed in
Example~\ref{egs:adjns-alg}\bref{egs:adjns-alg:fields} that the forgetful
functor $U\from \Field \to \Set$%
\index{field}
does not have a left adjoint.  We can now prove this.  If $U$ had a left
adjoint $F\from \Set \to \Field$, then $F$ would preserve colimits, and in
particular, initial objects.  Hence $F(\emptyset)$ would be an initial
object of $\Field$.  But $\Field$ has no initial object, since there are no
maps between fields of different characteristic.  Further examples of
nonexistence of adjoints can be found in Exercise~\ref{ex:no-adjt}.
\end{example}

\minihead{Adjoint functor theorems}
\index{adjoint functor theorems|(}

Every functor with a left adjoint preserves limits, but limit-preservation
alone does not guarantee the existence of a left adjoint.  For example, let
$\cat{B}$ be any category.  The unique functor $\cat{B} \to \One$ always
preserves limits, but by Example~\ref{eg:init-term}, it only has a left
adjoint if $\cat{B}$ has an initial object.

On the other hand, if we have a limit-preserving functor $G\from \cat{B}
\to \cat{A}$ \emph{and $\cat{B}$ has all limits}, then there is an
excellent chance that $G$ has a left adjoint.  It is still not always true,
but counterexamples are harder to find.  For instance (taking $\cat{A} =
\One$ again), can you find a category $\cat{B}$ that has all limits but no
initial object?

The condition of having all limits is so important that it has its own word:
\begin{defn}    
\label{defn:complete}
A category is \demph{complete}%
\index{category!complete}%
\index{complete}
(or properly, \demph{small complete}) if it has all limits.
\end{defn}

There are various results called adjoint functor theorems, all of the
following form:
\begin{displaytext} \it
Let $\cat{A}$ be a category, $\cat{B}$ a complete category, and $G\from \cat{B}
\to \cat{A}$ a functor.  Suppose that $\cat{A}$, $\cat{B}$ and $G$ satisfy
certain further conditions.  Then 
\[
G \text{ has a left adjoint}
\iff
G \text{ preserves limits}.
\]
\end{displaytext}
The forwards implication is immediate from Theorem~\ref{thm:adjts-cts}.
It is the backwards implication that concerns us here.

Typically, the `further conditions' involve the distinction between small
and large collections.  But there is a special%
\index{ordered set!adjunction between|(}
case in which these complications disappear, and I will use it to explain
the main idea behind the proofs of the adjoint functor theorems.  It is the
case where the categories $\cat{A}$ and $\cat{B}$ are ordered sets.

As we saw in Section~\ref{sec:lims-basics}, limits in ordered sets are
meets.  More precisely, if $D\from \scat{I} \to \scat{B}$ is a diagram in
an ordered set $\scat{B}$, then
\[
\lt{\scat{I}} D 
= 
\Meet_{I \in \scat{I}} D(I), 
\]
with one side defined if and only if the other is.  So an ordered set is
complete if and only if every subset has a meet.  Similarly, a map $G\from
\scat{B} \to \scat{A}$ of ordered sets preserves limits if and only if
\[
G\Biggl(\Meet_{i \in I} B_i\Biggr) 
= 
\Meet_{i \in I} G(B_i)
\]
whenever $(B_i)_{i \in I}$ is a family of elements of $\scat{B}$ for which a
meet exists.

We now show that for ordered sets, there is an adjoint functor theorem of the
simplest possible kind: there are no `further conditions' at all.

\begin{propn}[Adjoint functor theorem for ordered sets]
\label{propn:OAFT}
Let $\scat{A}$ be an ordered set, $\scat{B}$ a complete ordered set, and
$G\from \scat{B} \to \scat{A}$ an order-preserving map.  Then
\[
G \text{ has a left adjoint}
\iff
G \text{ preserves meets}.
\]
\end{propn}

\begin{pf}
Suppose that $G$ preserves meets.  By Corollary~\ref{cor:pre-AFT}, it is
enough to show that for each $A \in \scat{A}$, the comma category
$\comma{A}{G}$ has an initial object.  Let $A \in \scat{A}$.  Then
$\comma{A}{G}$ is an ordered set, namely, $\{ B \in \scat{B} \such A \leq
G(B) \}$ with the order inherited from $\scat{B}$.  We have to show that
$\comma{A}{G}$ has a least element.

Since $\scat{B}$ is complete, the meet $\Meet_{B \in \scat{B}\from A \leq
  G(B)} B$ exists in $\scat{B}$.  This is the meet of all the elements of
$\comma{A}{G}$, so it suffices to show that the meet is itself an element
of $\comma{A}{G}$.  And indeed, since $G$ preserves meets, we have
\[
G \Biggl( \Meet_{B \in \scat{B}\from A \leq G(B)} B \Biggr)
=
\Meet_{B \in \scat{B}\from A \leq G(B)} G(B)
\geq
A,
\]
as required.
\end{pf}

In the general setting of Corollary~\ref{cor:pre-AFT}, the initial object
of $\comma{A}{G}$ is the pair $\Bigl(F(A),\, A \toby{\eta_A} GF(A)\Bigr)$,
where $F$ is the left adjoint and $\eta$ is the unit map.  So in
Proposition~\ref{propn:OAFT}, the left adjoint $F$ is given by
\begin{equation}        
\label{eq:ladj-recipe}
F(A) 
= 
\Meet_{B \in \scat{B}\from A \leq G(B)} B.
\end{equation}

\begin{example} 
\label{eg:oaft-least}
\index{limit!colimit@vs.\ colimit|(}
Consider Proposition~\ref{propn:OAFT} in the case $\scat{A} = \One$.  The
unique functor $G\from \scat{B} \to \One$ automatically preserves meets,
and, as observed above, a left adjoint to $G$ is an initial object of
$\scat{B}$.  So in the case $\scat{A} = \One$, the proposition states that
a complete ordered set has a least element.  This is not quite trivial,
since completeness means the existence of all meets, whereas a least
element is an empty \emph{join}.

By~\eqref{eq:ladj-recipe}, the least%
\index{least element!meet@as meet}
element of $\scat{B}$ is $\Meet_{B \in \scat{B}} B$.  Thus, a least element
is not only a colimit of the functor $\emptyset \to \scat{B}$; it is also a
limit of the identity functor $\scat{B} \to \scat{B}$.

The synonym `least upper bound' for `join' suggests a theorem: that a poset
with all meets also has all joins.  Indeed, given a poset $\scat{B}$ with
all meets, the join of a subset of $\scat{B}$ is simply the meet of its
upper bounds: quite literally, its least upper bound.%
\index{limit!colimit@vs.\ colimit|)}
\end{example}

Let us now attempt to extend Proposition~\ref{propn:OAFT} from ordered sets
to categories, starting with a limit-preserving functor $G$ from a complete
category $\cat{B}$ to a category $\cat{A}$.  In the case of ordered sets,
we had for each $A \in \cat{A}$ an inclusion map $P_A \from \comma{A}{G}
\incl \scat{B}$, and we showed that the left adjoint $F$ was given by
\begin{equation}        
\label{eq:ladj-po}
F(A) = \lt{\comma{A}{G}} P_A.  
\end{equation}
In the general case, the analogue of the inclusion functor is the
projection functor
\begin{equation}
\label{eq:pjn-aft}
\begin{array}{cccc}
P_A\from    &\comma{A}{G}                       &\to            &\cat{B}\\
            &\Bigl(B,\, A \toby{f} G(B)\Bigr)   &\mapsto        &B.
\end{array}
\end{equation}
The case of ordered sets suggests that in general,
equation~\eqref{eq:ladj-po} might define a left adjoint $F$ to $G$.  And
indeed, it can be shown that if this limit in $\cat{B}$ exists and is
preserved by $G$, then \eqref{eq:ladj-po} really does give a left adjoint
(Theorem~X.1.2 of \citeCWM).

This might seem to suggest that our adjoint functor theorem generalizes
smoothly from ordered sets to arbitrary categories, with no need for
further conditions.  But it does not, for reasons that are quite subtle.

Those reasons are more easily explained if we relax our terminology
slightly.  When we defined limits, we built in the condition that the shape
category $\scat{I}$ was small.%
\index{limit!large|(}%
\index{limit!small|(}
However, the definition of limit makes sense for an arbitrary category
$\scat{I}$.  In this discussion, we will need to refer to this more
inclusive notion of limit, so let us temporarily suspend the convention
that the shape categories $\scat{I}$ of limits are always small.

Now, in the template for adjoint functor theorems stated above (after
Definition~\ref{defn:complete}), it was only required that $\cat{B}$ has, and
$G$ preserves, \emph{small} limits.  But if $\cat{B}$ is a large category
then $\comma{A}{G}$ might also be large, since to specify an object or map in
$\comma{A}{G}$, we have to specify (among other things) an object or map in
$\cat{B}$.  So, the limit~\eqref{eq:ladj-po} defining the left adjoint is not
guaranteed to be small.  Hence there is no guarantee that this limit exists in
$\cat{B}$, nor that it is preserved by $G$.  It follows that the functor $F$
`defined' by~\eqref{eq:ladj-po} might not be defined at all, let alone a left
adjoint.

(The reader experiencing difficulty with reasoning about small and large
collections might usefully compare finite and infinite collections.  For
instance, if $\cat{B}$ is a finite category and $\cat{A}$ has finite
hom-sets then $\comma{A}{G}$ is also finite, but otherwise $\comma{A}{G}$
might be infinite.)

Proposition~\ref{propn:OAFT} still stands, since there we were dealing with
ordered \emph{sets}, which as categories are small.  We might hope to
extend it from posets to arbitrary small categories, since the problem just
described affects only large categories.  But this turns out not to be very
fruitful, since in fact, complete posets are the \emph{only}%
\index{ordered set!complete small category is}
complete small categories (Exercise~\ref{ex:small-complete}).%
\index{ordered set!adjunction between|)}

Alternatively, we could try to salvage the argument by assuming that
$\cat{B}$ has, and $G$ preserves, \emph{all} (possibly large) limits.  But
again, this is unhelpful: there are almost no such categories $\cat{B}$.

The situation therefore becomes more complicated.  Each of the best-known
adjoint functor theorems imposes further conditions implying that the large
limit $\lt{\comma{A}{G}} P_A$ can be replaced by a small limit in some
clever way.  This allows one to proceed with the argument above.%
\index{limit!large|)}%
\index{limit!small|)}

The two most famous adjoint functor theorems are the `general' and the
`special'.  Their exact statements and proofs are perhaps less significant
than their consequences.

\begin{defn}
Let $\cat{C}$ be a category.  A \demph{weakly%
\index{set!weakly initial}%
\index{weakly initial}
initial set} in $\cat{C}$ is a set $\scat{S}$ of objects with the property
that for each $C \in \cat{C}$, there exist an element $S \in \scat{S}$ and
a map $S \to C$.
\end{defn}
Note that $\scat{S}$ must be a set, that is, small.  So, the existence of a
weakly initial set is some kind of size restriction.  Such size
restrictions are comparable to finiteness conditions in algebra.

\begin{thm}[General adjoint functor theorem]     
\label{thm:gaft}
\index{adjoint functor theorems!general}%
\index{general adjoint functor theorem (GAFT)}
Let $\cat{A}$ be a category, $\cat{B}$ a complete category, and $G\from
\cat{B} \to \cat{A}$ a functor.  Suppose that $\cat{B}$ is locally small and
that for each $A \in \cat{A}$, the category $\comma{A}{G}$ has a weakly
initial set.  Then
\[
G \text{ has a left adjoint}
\iff
G \text{ preserves limits}.
\]
\end{thm}

\begin{pf}
See the appendix.
\end{pf}

\begin{example} 
\label{eg:gaft-free-alg}
The general adjoint functor theorem (GAFT) implies that for any category
$\cat{B}$ of algebras ($\Grp$, $\Vect_k$, \ldots), the forgetful%
\index{functor!forgetful!left adjoint to}
functor $U\from \cat{B} \to \Set$ has a left adjoint.  Indeed, we saw in
Example~\ref{eg:lims-alg} that $\cat{B}$ has all limits, and in
Example~\ref{eg:gp-creation} that $U$ preserves them.  Also, $\cat{B}$ is
locally small.  To apply GAFT, we now just have to check that for each $A
\in \Set$, the comma category $\comma{A}{U}$ has a weakly initial set.
This requires a little cardinal%
\index{cardinality}%
\index{arithmetic!cardinal}
arithmetic, omitted here; see Exercise~\ref{ex:free-gp-card}.

So GAFT tells us that, for instance, the free group%
\index{group!free}
functor exists.  In Examples~\ref{egs:free-functors}\bref{eg:free-group}
and~\ref{egs:adjns-alg}\bref{egs:adjns-alg:gp}, we began to see the
trickiness of explicitly constructing the free group on a generating set
$A$.  One has to define the set of `formal expressions' (such as $x^{-1} y
x^2 z y^{-3}$, with $x, y, z \in A$), then say what it means for two such
expressions to be equivalent (so that $x^{-2} x^5 y$ is equivalent to $x^3
y$), then define $F(A)$ to be the set of all equivalence classes, then
define the group structure, then check the group axioms, then prove that
the resulting group has the universal property required.  But using GAFT,
we can avoid these complications entirely.

The price to be paid is that GAFT does not give us an explicit%
\index{explicit description}
description of free groups (or left adjoints more generally).  When people
speak of knowing some object `explicitly', they usually mean knowing its
elements.  An element of an object is a map \emph{into} it, and we have no
handle on maps into $F(A)$: since $F$ is a left adjoint, it is maps
\emph{out} of $F(A)$ that we know about.  This is why explicit descriptions
of left adjoints are often hard to come by.
\end{example}

\begin{example}
More generally, GAFT guarantees that forgetful%
\index{functor!forgetful!left adjoint to}
functors between categories of algebras, such as
\[
\Ab \to \Grp, 
\quad
\Grp \to \Mon,
\quad
\Ring \to \Mon,
\quad
\Vect_\complexes \to \Vect_\reals,
\]
have left adjoints.  (Some of them are described in
Examples~\ref{egs:adjns-alg}.)  This is `more generally' because $\Set$ can
be seen as a degenerate example of a category of algebras, in the sense of
Remark~\ref{rmk:alg-thy}: a group, ring, etc., is a set equipped with some
operations satisfying some equations, and a set is a set equipped with no
operations satisfying no equations.
\end{example}

The special adjoint functor theorem (SAFT) operates under much tighter
hypotheses than GAFT, and is much less widely applicable.  Its main
advantage is that it removes the condition on weakly initial sets.  Indeed,
it removes \emph{all} further conditions on the functor $G$.

\begin{thm}[Special adjoint functor theorem]
\index{adjoint functor theorems!special}%
\index{special adjoint functor theorem}%
\index{SAFT (special adjoint functor theorem)}
Let $\cat{A}$ be a category, $\cat{B}$ a complete category, and $G\from
\cat{B} \to \cat{A}$ a functor.  Suppose that $\cat{A}$ and $\cat{B}$ are
locally small, and that $\cat{B}$ satisfies certain further conditions.  Then
\[
G \text{ has a left adjoint}
\iff
G \text{ preserves limits}.
\]
\end{thm}
A precise statement and proof can be found in Section~V.8 of \citeCWM.

\begin{example}
Here is the classic application of SAFT.  Let $\CptHff$%
\index{topological space!compact Hausdorff}
be the category of compact Hausdorff spaces, and $U\from \CptHff \to \Tp$
the forgetful functor.  SAFT tells us that $U$ has a left adjoint $F$,
turning any space into a compact Hausdorff space in a canonical way.

The existence of this left adjoint is far from obvious, and verifying the
hypotheses of SAFT (or indeed, constructing $F$ in any other way) requires
some deep theorems of topology.  Given a space $X$, the resulting compact
Hausdorff space $F(X)$ is called its \demph{Stone--\v{C}ech%
\index{Stone-Cech compactification@Stone--\v{C}ech compactification}
compactification}.  Provided that $X$ satisfies some mild separation
conditions, the unit of the adjunction at $X$ is an embedding, so that
$UF(X)$ contains $X$ as a subspace. 

Another advantage of SAFT is that one can extract from its proof a fairly
explicit formula for the left adjoint.  In this case, it tells us that
$F(X)$ is the closure of the image of the canonical map
\[
X \to [0, 1]^{\Tp(X, [0, 1])}, 
\]
where the codomain is a power of $[0, 1]$ in $\Tp$.  
\index{adjoint functor theorems|)}
\end{example}

\minihead{Cartesian closed categories}
\index{category!cartesian closed|(}%
\index{cartesian closed category|(}

We have seen that for every set $B$, there is an adjunction $(\dashbk\times B)
\ladj (\dashbk)^B$ (Example~\ref{eg:adjn:cc}), and that for every category
$\cat{B}$, there is an adjunction $(\dashbk \times \cat{B}) \ladj
\ftrcat{\cat{B}}{\dashbk}$
(Remark~\ref{rmks:global-hom}\bref{rmks:global-hom:cc}).

\begin{defn}
A category $\cat{A}$ is \demph{cartesian closed} if it has finite products
and for each $B \in \cat{A}$, the functor $\dashbk \times B\from \cat{A}
\to \cat{A}$ has a right adjoint.
\end{defn}
We write the right adjoint as $(\dashbk)^B$,%
\ntn{exp-cc}
and, for $C \in \cat{A}$, call
$C^B$ an \demph{exponential}.%
\index{exponential}
We may think of $C^B$ as the space of maps from $B$ to $C$.  Adjointness
says that for all $A, B, C \in \cat{A}$,
\[
\cat{A}(A \times B, C)
\iso
\cat{A}\bigl(A, C^B\bigr)
\]
naturally in $A$ and $C$.  In fact, the isomorphism is natural in $B$ too;
that comes for free.

\begin{example}
$\Set$ is cartesian closed; $C^B$ is the function%
\index{set!functions@of functions}%
\index{function!set of functions}
set $\Set(B, C)$.  
\end{example}

\begin{example}
\hspace*{-.5pt}%
$\CAT$ is cartesian closed; $\cat{C}^\cat{B}$ is the functor%
\index{functor!category}
category $\ftrcat{\cat{B}}{\cat{C}}$.
\end{example}

In any cartesian closed category with finite sums, the
isomorphisms~\eqref{eq:Set-dist} and~\eqref{eq:Set-codist} of
Example~\ref{eg:arith-set} hold, for the same reasons as stated there.  The
objects of a cartesian closed category therefore possess an arithmetic%
\index{arithmetic}
like that of the natural numbers.  This thought can be developed in several
interesting directions, but here we just note that these isomorphisms
provide a way of proving that a category is \emph{not} cartesian closed.

\begin{example}
$\Vect_k$%
\index{vector space!category of vector spaces!cartesian closed@is not
  cartesian closed} 
is not cartesian closed, for any field $k$.  It does have finite products,
as we saw in Example~\ref{eg:prod-vs}: binary product is direct sum
$\oplus$, and the terminal object is the trivial vector space $\{0\}$,
which is also initial.  But if $\Vect_k$ were cartesian closed then
equations~\eqref{eq:Set-dist} would hold, so that $\{0\} \oplus B \iso
\{0\}$ for all vector spaces $B$.  This is plainly false.
\end{example}

\begin{remark}
For any vector spaces $V$ and $W$, the set $\Vect_k(V, W)$ of linear maps
can itself be given the structure of a vector space, as in
Example~\ref{eg:fns-on-vs}.  Let us now call this vector space $[V, W]$.

Given that exponentials are supposed to be `spaces of maps', you might
expect $\Vect_k$ to be cartesian closed, with $[\dashbk, \dashbk]$ as its
exponential.  We have just seen that this cannot be so.  But as it turns
out, the linear maps $U \to [V, W]$ correspond to the \emph{bi}linear%
\index{map!bilinear}
maps $U \times V \to W$, or equivalently the linear maps $U \otimes V \to
W$.%
\index{tensor product}
In the jargon, $\Vect_k$ is an example of a `monoidal%
\index{monoidal closed category}%
\index{category!monoidal closed}
closed category'.  These are like cartesian closed categories, but with the
cartesian (categorical) product replaced by some other operation called
`product', in this case the tensor product of vector spaces.
\end{remark}

For any set $I$, the product category $\Set^I$ is cartesian closed, just
because $\Set$ is.  (Exponentials in $\Set^I$, as well as products, are
computed pointwise.)  Put another way, $\pshf{\scat{A}}$ is cartesian
closed whenever $\scat{A}$ is discrete.  We now show that, in fact,
$\pshf{\scat{A}}$ is cartesian closed for any small category $\scat{A}$
whatsoever.  

In preparation for proving this, let us conduct a thought%
\index{thought experiment}
experiment.  Write $\psh{\scat{A}}%
\ntn{psh}
 = \pshf{\scat{A}}$.  If $\psh{\scat{A}}$ \emph{is} cartesian
closed, what must exponentials in $\psh{\scat{A}}$ be?  In other words, given
presheaves $Y$ and $Z$, what must $Z^Y$ be in order that
\begin{equation}        
\label{eq:cc-pshf}
\psh{\scat{A}}\bigl(X, Z^Y\bigr)
\iso
\psh{\scat{A}}(X \times Y, Z)
\end{equation}
for all presheaves $X$?  If this is true for all presheaves $X$, then in
particular it is true when $X$ is representable, so
\[
Z^Y(A)
\iso
\psh{\scat{A}}\bigl(\h_A, Z^Y\bigr)
\iso
\psh{\scat{A}}(\h_A \times Y, Z) 
\]
for all $A \in \scat{A}$, the first step by Yoneda.  This tells us what
$Z^Y$ must be.  Notice that $Z^Y(A)$ is not simply $Z(A)^{Y(A)}$, as one
might at first guess: exponentials in a presheaf category are \emph{not}
generally computed pointwise.%
\index{pointwise}

\begin{thm}     
\label{thm:pshf-cc}
For any small category $\scat{A}$, the presheaf%
\index{presheaf!category of presheaves!cartesian closed@is cartesian closed}
 category $\psh{\scat{A}}$
is cartesian closed.
\end{thm}
Here is the strategy of the proof.  The argument in the thought experiment
gives us the isomorphism~\eqref{eq:cc-pshf} whenever $X$ is representable.
A general presheaf $X$ is not representable, but it is a colimit of
representables, and this allows us to bootstrap our way up.

\begin{pf}
We know that $\psh{\scat{A}}$ has all limits, and in particular, finite
products.  It remains to show that $\psh{\scat{A}}$ has exponentials.  Fix
$Y \in \psh{\scat{A}}$.

First we prove that $\dashbk \times Y\from \psh{\scat{A}} \to
\psh{\scat{A}}$ preserves colimits.  (Eventually we will prove that
$\dashbk\times Y$ has a right adjoint, from which preservation of colimits
follows, but our proof that it has a right adjoint will \emph{use}
preservation of colimits.)  Indeed, since products and colimits in
$\psh{\scat{A}}$ are computed pointwise, it is enough to prove that for any
set $S$, the functor $\dashbk \times S\from \Set \to \Set$ preserves
colimits, and this follows from the fact that $\Set$ is cartesian closed.

For each presheaf $Z$ on $\scat{A}$, let $Z^Y$ be the presheaf defined by
\[
Z^Y (A) = \psh{\scat{A}} (\h_A \times Y, Z)
\]
for all $A \in \scat{A}$.  This defines a functor $(\dashbk)^Y \from
\psh{\scat{A}} \to \psh{\scat{A}}$.  

I claim that $(\dashbk \times Y) \ladj (\dashbk)^Y$.  Let $X, Z \in
\psh{\scat{A}}$.  Write $P\from \elt{X} \to \scat{A}$ for the projection
(as in Definition~\ref{defn:cat-elts}), and write $\h_P = \h_\bl \of P$.
Then
\begin{align}
\label{eq:pcc:density}
\psh{\scat{A}}\bigl(X, Z^Y\bigr) &
\iso   
\psh{\scat{A}} \biggl( \colt{\elt{X}} \h_P, Z^Y \biggr)  \\      
\label{eq:pcc:out}
        &
\iso   
\lt{\elt{X}^{\op}} \psh{\scat{A}}\bigl(\h_P, Z^Y\bigr)   \\
\label{eq:pcc:yon}
        &
\iso   
\lt{\elt{X}^{\op}} Z^Y(P)                       \\
\label{eq:pcc:defn}
        &
\iso   
\lt{\elt{X}^{\op}} \psh{\scat{A}}(\h_P \times Y, Z)     \\
\label{eq:pcc:in}
        &
\iso   
\psh{\scat{A}}\biggl( \colt{\elt{X}} (\h_P \times Y), Z \biggr) \\
\label{eq:pcc:cocts} 
       &
\iso   
\psh{\scat{A}}\biggl(\biggl(\colt{\elt{X}} \h_P\biggr) \times Y, Z\biggr) \\
\label{eq:pcc:density-again}
        &
\iso   
\psh{\scat{A}}(X \times Y, Z)
\end{align}
naturally in $X$ and $Z$.  Here~\eqref{eq:pcc:density}
and~\eqref{eq:pcc:density-again} follow from Theorem~\ref{thm:density};
\eqref{eq:pcc:out} and~\eqref{eq:pcc:in} are because representables
preserve limits (as rephrased in Remark~\ref{rmk:rep-pres});
\eqref{eq:pcc:yon} is by Yoneda; \eqref{eq:pcc:defn} is by definition of
$Z^Y$; and~\eqref{eq:pcc:cocts} is because $\dashbk\times Y$ preserves
colimits.
\end{pf}

This result can be seen as a step along the road to topos%
\index{topos}
theory.  A topos is a category with certain special properties.  Topos
theory unifies, in an extraordinary way, important aspects of logic and
geometry.

For instance, a topos can be regarded as a `universe of sets':%
\index{set!category of sets!topos@as topos}
$\Set$ is the most basic example of a topos, and every topos shares enough
features with $\Set$ that one can reason with its objects as if they were
sets of some exotic kind.  On the other hand, a topos can be regarded as a
generalized topological space:%
\index{topological space!topos@as topos}
every space gives rise to a topos (namely, the category of sheaves%
\index{sheaf}
on it), and topological properties of the space can be reinterpreted in a
useful way as categorical properties of its associated topos.

By definition, a topos is a cartesian closed category with finite limits
and with one further property: the existence of a so-called subobject%
\index{subobject!classifier}
classifier.  For example, the two-element%
\index{set!two-element}
set $2$ is the subobject classifier of $\Set$, which means, informally,
that subsets of a set $A$ correspond one-to-one with maps $A \to 2$.
Exercises~\ref{ex:soc} and~\ref{ex:pshf-topos} give the formal definition
of subobject classifier, then guide you through the proof that $\Set$, and,
more generally, every presheaf category, is a topos.%
\index{category!cartesian closed|)}%
\index{cartesian closed category|)}

\exs

\begin{question}        
\label{ex:no-adjt}
\begin{enumerate}[(b)]
\item 
Prove that the forgetful functor $U\from \Grp \to \Set$ has no right
adjoint.

\item 
Prove that the chain of adjunctions $C \ladj D \ladj O \ladj I$%
\index{category!category of categories!adjunctions with $\Set$}
in Exercise~\ref{ex:cdoi} extends no further in either direction.

\item
Does the chain of adjunctions in Exercise~\ref{ex:pshf-adjns} extend
further in either direction?
\end{enumerate}
\end{question}

\begin{question}
Let $\cat{A}$ be a locally small category.  For functors $U\from \cat{A}
\to \Set$, consider the following three conditions: (A)~$U$ has a left
adjoint; (R)~$U$ is representable; (L)~$U$ preserves limits.  
\begin{enumerate}[(b)]
\item 
Show that (A) $\textonlyif$ (R) $\textonlyif$ (L).

\item 
Show that if $\cat{A}$ has sums then (R) $\textonlyif$ (A).%
\index{functor!representable!adjoints@and adjoints}
\end{enumerate}
(If $\cat{A}$ satisfies the hypotheses of the special adjoint functor
theorem then also (L) $\textonlyif$ (A), so the three conditions are
equivalent.)
\end{question}

\begin{question}        
\label{ex:small-complete}
\begin{enumerate}[(b)]
\item 
Prove that every preordered%
\index{ordered set!preordered set@vs.\ preordered set}
set is equivalent (as a category) to an ordered set.

\item 
Let $\cat{A}$ be a category with all small products.  Suppose that
$\cat{A}$ is not a preorder, so that there exists a parallel pair of maps
$\parpairi{A}{B}{f}{g}$ in $\cat{A}$ with $f \neq g$.  By considering the
maps $A \to B^I$ for each set $I$, prove that $\cat{A}$ is not small.

\item 
Deduce that every small category with small products is equivalent to a
complete ordered set.%
\index{ordered set!complete small category is}

\item 
Adapt the argument to prove that every finite category with finite
products is equivalent to a complete ordered set.
\end{enumerate}
\end{question}

\begin{question}        
\label{ex:free-gp-card}
Probably the most important application of the general adjoint functor
theorem is to proving that forgetful functors between categories of
algebras have left adjoints (Example~\ref{eg:gaft-free-alg}).  Verifying
the hypotheses can be done with some cardinal%
\index{cardinality}%
\index{arithmetic!cardinal}
arithmetic.  Here is a typical example.
\begin{enumerate}[(b)]
\item   
\label{part:gen-subgp-bound}
Let $A$ be a set.  Prove that for any group $G$ and family $(g_a)_{a \in
A}$ of elements of $G$, the subgroup of $G$ generated by $\{g_a \such a
\in A\}$ has cardinality at most $\max \{\crd{\nat}, \crd{A}\}$.

\item   
\label{part:small-coll-gps}
Prove that for any set $S$, the collection of isomorphism classes of groups
of cardinality at most $\crd{S}$ is small.

\item 
Let $U\from \Grp \to \Set$%
\index{group!free}
be the forgetful functor from groups to sets.  Deduce
from~\bref{part:gen-subgp-bound} and~\bref{part:small-coll-gps} that for
every set $A$, the comma category $\comma{A}{U}$ has a weakly initial set.

\item 
Use GAFT to conclude that $U$ has a left adjoint.
\end{enumerate}
\end{question}

\begin{question}
Let $\scat{A}$ be a small cartesian closed category.  Prove that the
Yoneda%
\index{Yoneda embedding!preserves exponentials}%
\index{exponential!preserved by Yoneda embedding}
embedding $\scat{A} \to \pshf{\scat{A}}$ preserves the whole cartesian
closed structure (exponentials as well as products).
\end{question}

\begin{question}        
\label{ex:soc}
Recall from Exercise~\ref{ex:subobjects} the notion of subobject.  A 
category $\cat{A}$ is \demph{well-powered}%
\index{category!well-powered}%
\index{well-powered}
if for each $A \in \cat{A}$, the class of subobjects of $A$ is small, that
is, a set.  (All of our usual examples of categories are well-powered.)
Let $\cat{A}$ be a well-powered category with pullbacks, and write
$\Sub(A)$ for the set of subobjects of an object $A \in \cat{A}$.
\begin{enumerate}[(b)]
\item
Deduce from Exercise~\ref{ex:pb-monic} that any map $A' \toby{f} A$ in
$\cat{A}$ induces a map $\Sub(f)\from \Sub(A) \to \Sub(A')$.

\item 
Show that this determines a functor $\Sub\from \cat{A}^\op \to \Set$.
(Hint: use Exercise~\ref{ex:pb-pasting}.)

\item 
For some categories $\cat{A}$, the functor $\Sub$ is representable.  A
\demph{subobject%
\index{subobject!classifier}
classifier} for $\cat{A}$ is an object $\Omega \in
\cat{A}$ such that $\Sub \iso \h_\Omega$.  Prove that $2$ is a subobject
classifier for $\Set$.
\end{enumerate}
A \demph{topos}%
\index{topos}
is a cartesian closed category with finite limits and a subobject
classifier.  You have just completed the proof that $\Set$ is a topos.
\end{question}

\begin{question}        
\label{ex:pshf-topos}
This exercise follows on from the last, culminating in the proof that every
presheaf category%
\index{presheaf!category of presheaves!topos@is topos}
is a topos.%
\index{topos}
Let $\scat{A}$ be a small category.
\begin{enumerate}[(b)]
\item
By conducting a thought%
\index{thought experiment}
experiment similar to the one before the statement of
Theorem~\ref{thm:pshf-cc}, find out what the subobject classifier $\Omega$
of $\pshf{\scat{A}}$ must be if it exists.

\item 
Prove that this $\Omega$ is indeed a subobject classifier.

\item 
Conclude that $\pshf{\scat{A}}$ is a topos.
\end{enumerate}
\end{question}


\backmatter

\addtocontents{toc}{\vspace{.2\baselineskip}} 
\oneappendix    
%
%
%

\chapter{Proof of the general adjoint functor theorem}
\index{adjoint functor theorems!general|(}%
\index{general adjoint functor theorem (GAFT)|(}

Here we prove the general adjoint functor theorem, which for convenience is
restated below.  The left-to-right implication follows immediately from
Theorem~\ref{thm:adjts-cts}; it is the right-to-left implication that we
have to prove.

\paragraph*{\upshape Theorem~\ref{thm:gaft} (General adjoint functor theorem)} 
{\itshape
Let $\cat{A}$ be a category, $\cat{B}$ a complete category, and $G\from
\cat{B} \to \cat{A}$ a functor.  Suppose that $\cat{B}$ is locally small and
that for each $A \in \cat{A}$, the category $\comma{A}{G}$ has a weakly%
\index{set!weakly initial|(}%
\index{weakly initial|(}
initial set.  Then
\[
G \text{ has a left adjoint}
\iff
G \text{ preserves limits}.
\]
}

The heart of the proof is the case $\cat{A} = \One$, where GAFT asserts
that a complete locally small category with a weakly initial set has an
initial object.  We prove this first.

The proof of this special case is illuminated by considering the even more
special case where $\cat{A} = \One$ and the category $\cat{B}$ is a poset
$\scat{B}$.  We saw in Example~\ref{eg:oaft-least} that the initial object
(least%
\index{least element}
element) of a complete poset $\scat{B}$ can be constructed as the meet of
all its elements.  Otherwise put, it is the limit of the identity functor
$1_\scat{B}\from \scat{B} \to \scat{B}$.

One might try to extend this result to arbitrary categories $\cat{B}$ by
proving that the limit of the identity%
\index{functor!identity!limit of}%
\index{object!initial!limit of identity@as limit of identity}%
\index{limit!identity@of identity}
functor $1_\cat{B}\from \cat{B} \to \cat{B}$ is (if it exists) an initial
object.  This is indeed true (Exercise~\ref{ex:lim-of-id} below).  However,
it is unhelpful: for if $\cat{B}$ is large then the limit of $1_\cat{B}$
is a large%
\index{limit!large} 
limit, but we are only given that $\cat{B}$ has small limits.

We seem to be at an impasse~-- but this is where the clever idea behind GAFT
comes in.  In order to construct the least element of a complete poset, it
is not necessary to take the meet of \emph{all} the elements.  More
economically, we could just take the meet of the elements of some weakly
initial subset (Exercise~\ref{ex:wi-poset}).  In general, for an arbitrary
complete category, the limit of any weakly initial set is an initial
object.  We prove this now.

\begin{alemma}  
\label{lemma:wi-init}
Let $\cat{C}$ be a complete locally small category with a weakly initial
set.  Then $\cat{C}$ has an initial object.
\end{alemma}

\begin{pf}
Let $\scat{S}$ be a weakly initial set in $\cat{C}$.  Regard $\scat{S}$ as
a full subcategory of $\cat{C}$; then $\scat{S}$ is small, since $\cat{C}$
is locally small.  We may therefore take a limit cone
\begin{equation}        
\label{eq:lim-wi}
\Bigl( 0 \toby{p_S} S \Bigr)_{S \in \scat{S}}
\end{equation}
of the inclusion $\scat{S} \incl \cat{C}$.  We prove that $0$ is initial.

Let $C \in \cat{C}$.  We have to show that there is exactly one map $0 \to
C$.  Certainly there is at least one, since we may choose some $S \in
\scat{S}$ and map $j \from S \to C$, and we then have the composite $j p_S
\from 0 \to C$.  To prove uniqueness, let $f, g \from 0 \to C$.  Form the
equalizer
\[
\xymatrix@1{
E \ar[r]^i      &0 \ar@<.5ex>[r]^f \ar@<-.5ex>[r]_g       &C.
}
\]
Since $\scat{S}$ is weakly initial, we may choose $S \in \scat{S}$ and
$h\from S \to E$.  We then have maps
\[
\xymatrix@1{
0 \ar[r]^{p_S} &S \ar[r]^h &E \ar[r]^i &0 \\
}
\]
with the property that for all $S' \in \scat{S}$,
\[
p_{S'}(ihp_S) = (p_{S'}ih)p_S = p_{S'} = p_{S'} 1_0
\]
(where the second equality follows from~\eqref{eq:lim-wi} being a cone).
But~\eqref{eq:lim-wi} is a \emph{limit} cone, so $ihp_S = 1_0$ by
Exercise~\ref{ex:jointly-monic}\bref{part:j-m-main}.  Hence
\[
f = fihp_S = gihp_S = g,
\]
as required.
\end{pf}

We have now proved GAFT in the special case $\cat{A} = \One$.  The rest of
the proof is comparatively routine.

\begin{alemma}  
\label{lemma:gaft-creates}
\index{limit!creation of}%
\index{creation of limits}%
\index{comma category!limits in}
Let $\cat{A}$ and $\cat{B}$ be categories.  Let $G \from \cat{B} \to
\cat{A}$ be a functor that preserves limits.  Then the projection functor
$P_A \from \comma{A}{G} \to \cat{B}$ of~\eqref{eq:pjn-aft} creates limits,
for each $A \in \cat{A}$.  In particular, if $\cat{B}$ is complete then so
is each comma category $\comma{A}{G}$.
\end{alemma}

\begin{pf}
The first statement is
Exercise~\ref{ex:gaft-creates}\bref{part:gaft-creates-creates}, and the
second follows from Lemma~\ref{lemma:creates-preserves}.
\end{pf}

We now prove GAFT.  By Corollary~\ref{cor:pre-AFT}, it is enough to show
that $\comma{A}{G}$ has an initial object for each $A \in \cat{A}$.  Let $A
\in \cat{A}$.  By Lemma~\ref{lemma:gaft-creates}, $\comma{A}{G}$ is
complete, and by hypothesis, it has a weakly initial set.  It is also
locally small, since $\cat{B}$ is.  Hence by Lemma~\ref{lemma:wi-init}, it
has an initial object, as required.

\exs

\begin{aquestion}        
\label{ex:lim-of-id}
In this exercise, we suspend the convention (made implicitly in
Definition~\ref{defn:lim}) that we only speak of the limit of a functor
$\scat{I} \to \cat{C}$ when $\scat{I}$ is small.%
\index{limit!small}%
\index{limit!large}
Let $\cat{B}$ be a
category, possibly large.  The aim is to prove that a limit%
\index{functor!identity!limit of}%
\index{object!initial!limit of identity@as limit of identity}%
\index{limit!identity@of identity}
of the identity functor on $\cat{B}$ is exactly an initial object of
$\cat{B}$.
\begin{enumerate}[(b)]
\item 
Let $0$ be an initial object of $\cat{B}$.  Show that the cone $(0 \to
B)_{B \in \cat{B}}$ on the identity functor $1_\cat{B}$ is a limit cone.

\item
Now let $\Bigl(L \toby{p_B} B\Bigr)_{B \in \cat{B}}$ be a limit cone on
$1_\cat{B}$.  Prove that $p_L$ is the identity on $L$, and deduce that $L$
is initial.  
\end{enumerate}
\end{aquestion}

\begin{aquestion}        
\label{ex:wi-poset}
Here you will prove the special case of Lemma~\ref{lemma:wi-init} in which
the category concerned is a poset.  Let $C$ be a poset and $S \sub C$.
\begin{enumerate}[(b)]
\item 
What does it mean, in purely order-theoretic terms, for $S$ to be a weakly
initial set in $C$?

\item 
Prove directly that if $S$ is weakly initial and the meet $\Meet_{s \in S}
s$ exists then $\Meet_{s \in S} s$ is a least%
\index{least element}
element of $C$.%
\index{set!weakly initial|)}%
\index{weakly initial|)}
\end{enumerate}
\end{aquestion}

\begin{aquestion}       
\label{ex:gaft-creates}
Let $G \from \cat{B} \to \cat{A}$ be a limit-preserving functor, and let $A
\in \cat{A}$.
\begin{enumerate}[(b)]
\item 
Show that for any small category $\scat{I}$, a diagram of shape $\scat{I}$
in $\comma{A}{G}$ amounts to a diagram $E$ of shape $\scat{I}$ in $\cat{B}$
together with a cone on $G \of E$ with vertex $A$.

\item
\label{part:gaft-creates-creates}
Prove that the projection functor $P_A \from \comma{A}{G} \to \cat{B}$
of~\eqref{eq:pjn-aft} creates limits.%
\index{adjoint functor theorems!general|)}%
\index{general adjoint functor theorem (GAFT)|)}
\end{enumerate}
\end{aquestion}
\endappendix

%
%
%

\chapter*{Further reading}

This book is intentionally short.  Even some topics that are included in
most introductions to category theory are omitted here.  I will indicate
some of the topics that lie beyond the scope of this book, and suggest
where you might read about them.  Since there is far more written on
category theory than anyone could read in a lifetime, these recommendations
are necessarily subjective.

The towering presence among category theory books is the classic by
one of its founders:
\begin{citedsource}
Saunders Mac~Lane,
\emph{Categories for the Working Mathematician}.\linebreak
Springer, 
1971;
second edition with two new chapters, 
1998.
\end{citedsource}
It is so well-written that more than forty years on, it is still the most
popular introduction to the subject.  It addresses a more mature readership
than this text, and covers many topics omitted here, including monads (one
formalization of the idea of algebraic theory), monoidal categories
(categories equipped with a tensor product), 2-categories (mentioned at the
end of our Chapter~\ref{ch:cfnt}), abelian categories (categories of
modules), ends (an elegant generalization of the notion of limit), and Kan
extensions (which provide the tongue-in-cheek title of the book's final
section: `All concepts are Kan extensions').

Another well-liked book, longer than the one you hold in your hands but
written for a similar readership, is:
\begin{citedsource}
Steve Awodey,
\emph{Category Theory}.
Oxford University Press, 
2010.
\end{citedsource}
Awodey's book covers less than Mac~Lane's, but is particularly strong on
connections between category theory and other parts of logic.  It has a
full chapter on cartesian closed categories, and also covers the theory of
monads.

Those who prefer lectures to books might try this library of 75 ten-minute
introductory category theory videos:
\begin{citedsource}
Eugenia Cheng and Simon Willerton,
The Catsters.
Available at\linebreak
\href{https://www.youtube.com/user/TheCatsters}{\url{https://www.youtube.com/user/TheCatsters}}, 
2007--2010.
\end{citedsource}
Other than the topics treated here, they cover monads, enriched categories,
internal groups (and other internal algebraic structures), string diagrams
(which we touched on in Remark~\ref{rmk:triangle-string}), and several more
sophisticated topics.

For inspiration as much as instruction, here are two further
recommendations.
\begin{citedsource}
Saunders Mac~Lane,
\emph{Mathematics: Form and Function}.
Springer,
1986.
\end{citedsource}
\begin{citedsource}
F. William Lawvere and Stephen H. Schanuel,
\emph{Conceptual Mathematics: A First Introduction to Categories}.
Cambridge University Press,
1997.
\end{citedsource}
\emph{Mathematics: Form and Function} is a tour through much of pure and
applied mathematics, written from a categorical perspective.  Its declared
purpose is to present the author's philosophy of mathematics, but it can
also be enjoyed for its many excellent vignettes of exposition.  (Beware of
the numerous small errors.)  \emph{Conceptual Mathematics} is a
thought-provoking text and an intriguing experiment: category theory for
high-school students, complete with classroom dialogues.

For categorical topics beyond the scope of this book, two good general
references are:
\begin{citedsource}
Francis Borceux,
\emph{Handbook of Categorical Algebra, Volumes 1--3}.
Cambridge University Press, 
1994.
\end{citedsource}
\begin{citedsource}
Various authors,
\emph{The $n$Lab}.
Available at \href{https://ncatlab.org}{\url{https://ncatlab.org}}, 2008--\linebreak
present.
\end{citedsource}
Borceux's encyclopaedic work often takes a different point of view from the
present text, but covers many, many more topics.  Apart from those just
mentioned in connection with other books, some of the more important ones
are fibrations, bimodules (also called profunctors or distributors),
Lawvere theories, Cauchy completeness, Morita equivalence, absolute
colimits, and flatness.

The $n$Lab is an ever-growing online resource for mathematics, focusing on
category theory and operating on similar principles to Wikipedia.
Individual entries can be idiosyncratic, but it has become a very useful
reference for advanced categorical topics.

Vigorous research in category theory continues to be done.  The sources
listed above provide ample onward references for anyone wishing to explore.

\minihead{Other texts cited}

\begin{citedsource}
Timothy Gowers,
\emph{Mathematics: A Very Short Introduction}.
Oxford University Press,
2002.
\end{citedsource}

\begin{citedsource}
G. M. Kelly,
\emph{Basic Concepts of Enriched Category Theory}.
Cambridge University Press,
1982. 
Also 
\emph{Reprints in Theory and Applications of Categories}
10 (2005), 1--136,
available at 
\href{http://www.tac.mta.ca/tac/reprints}{\url{http://www.tac.mta.ca/tac/reprints}}.
\end{citedsource}

\begin{citedsource}
F. William Lawvere and Robert Rosebrugh,
\emph{Sets for Mathematics}.
Cambridge University Press,
2003.
\end{citedsource}

\begin{citedsource}
Tom Leinster,
Rethinking set theory.
\emph{American Mathematical Mon\-thly} 121 (2014), no.~5, 403--415.
Also available at 
\href{https://arxiv.org/abs/1212.6543}{\url{https://arxiv.org/abs/1212.6543}}.
\end{citedsource}

%
%
%

\chapter*{Index of notation}

{%
\footnotesize%
\setlength{\parindent}{0em}%
\begin{multicols}{3}
%
%
\nuse{\nref{empty-blank}}{blank space}
\nuse{\nref{juxt}}{$gf$}
\nuse{\nref{whisker-right}}{$\alpha F$}
\nuse{\nref{whisker-left}}{$F \alpha$}
\nuse{\nref{nt-comp}}{$\alpha_A$}
\nuse{\nref{hom-set-default}}{$\cat{A}(A, B)$}
\nuse{\nref{hom-out}}{$\cat{A}(A, \dashbk)$}
\nuse{\nref{hom-in}}{$\cat{A}(\dashbk, A)$}
\nuse{\nref{hom-out-map-blank}}{$\cat{A}(f, \dashbk)$}
\nuse{\nref{hom-in-map-blank}}{$\cat{A}(\dashbk, f)$}
\nuse{\nref{AAD}}{$\cat{A}(A, D)$}
\nuse{\nref{DblankA}}{$D(\dashbk)(A)$}
\nuse{\nref{ftr-cat-power}}{$\cat{B}^\cat{A}$}
\nuse{\nref{exponential-set}, \nref{power-of-obj}, \nref{exp-cc}}{$B^A$}
\nuse{\nref{map-to-prod}}{$(f_i)_{i \in I}$}
\nuse{\nref{small-cat-face}}{$\scat{A}, \scat{B}, \ldots$ (typeface)}

%
%
\nuse{\nref{dashbk}}{$\dashbk$}
\nuse{\nref{adj-bar}, \nref{lim-bar}, \nref{colim-bar}}{$\bar{\emptybk}$}
\nuse{\pageref{eq:yoneda-fns}}{$\ynt{\emptybk}$}
\nuse{\pageref{eq:yoneda-fns}, \nref{psh}}{$\yel{\emptybk}$}
\nuse{\nref{D-upper-lower}}{$\blank^\bl$, $\blank_\bl$}
\nuse{\nref{horiz-comp}}{$*$}
\nuse{\nref{dual-vs}}{$V^*$}
\nuse{\nref{dual-map}, \nref{upper-star}}{$f^*$}
\nuse{\nref{lower-star}, \nref{lower-star-bis}}{$f_*$}
\nuse{\nref{of}, \nref{of-ftr}, \nref{of-nt}}{$\of$}
\nuse{\nref{of-blank}, \nref{of-blank-bis}}{$g \of \dashbk$}
\nuse{\nref{blank-of}}{$\dashbk \of f$}
\nuse{\nref{forall}}{$\forall$}
\nuse{\nref{ei}}{$\exists!$}
\nuse{\nref{arrow}}{$\to$}
\nuse{\nref{incl}}{$\incl$}
\nuse{\nref{iso-arrow}}{$\toiso$}
\nuse{\nref{double-arrow}, \nref{comma-cat}, \nref{comma-fix}}{$\Rightarrow$}
\nuse{\nref{ladj}}{$\ladj$}
\nuse{\nref{bot}}{$\bot$, $\top$}
\nuse{\nref{iso}, \nref{iso-cat}, \nref{iso-ftr}}{$\iso$}
\nuse{\nref{eqv}}{$\eqv$}
\nuse{\nref{leq}, \nref{card-leq}}{$\leq$}
\nuse{\nref{card-leq}}{$\crd{\emptybk}$}
\nuse{\nref{ftr-cat-bkts}}{$\ftrcat{\emptybk}{\emptybk}$}
\nuse{\nref{tensor}}{$\otimes$}
\nuse{\nref{prod-cat}, \nref{prod-set}, \nref{prod-gen}}{$\times$}
\nuse{\nref{prod-fam-set}, \nref{prod-fam-gen}}{$\prod$}
\nuse{\nref{sum-set}, \nref{sum-gen}}{$+$}
\nuse{\nref{sum-fam-set}, \nref{sum-fam-gen}}{$\sum$}
\nuse{\nref{disjt-union}}{$\amalg$}
\nuse{\nref{disjt-union-fam-gen}}{$\coprod$}
\nuse{\nref{direct-sum}}{$\oplus$}
\nuse{\nref{slice}}{$\cat{A}/A$}
\nuse{\nref{coslice}}{$A/\cat{A}$}
\nuse{\nref{qt-set}}{$\qer{A}{\sim}$}
\nuse{\nref{meet}}{$\meet$}
\nuse{\nref{Meet}}{$\Meet$}
\nuse{\nref{join}}{$\join$}
\nuse{\nref{Join}}{$\Join$}

%
%
\nuse{\nref{inverse}}{$\blank^{-1}$}
\nuse{\nref{empty-cat}}{$\emptyset$}
\nuse{\nref{initial}}{$0$}
\nuse{\nref{oneset}, \nref{id-map}, \nref{id-ftr}, \nref{id-nt},
  \nref{terminal}}{$1$} 
\nuse{\nref{terminal-cat}}{$\One$}
\nuse{\nref{two-disc-cat}, \nref{two-set}}{$2$}
\nuse{\nref{Two}}{$\Two$}

%
%
\nuse{\nref{Delta}, \nref{diag-set}, \nref{diag-gen}}{$\Delta$}
\nuse{\nref{adj-counit}}{$\epsln$}
\nuse{\nref{adj-unit}}{$\eta$}
\nuse{\nref{pi-1}}{$\pi_1$}
\nuse{\nref{char-fn}}{$\chi$}

%
%
\nuse{\nref{Ab}}{$\Ab$}
\nuse{\nref{abel}}{$\abel{\blank}$}
\nuse{\nref{Bilin}}{$\Bilin$}
\nuse{\nref{cts-ring}}{$C$}
\nuse{\nref{CAT}}{$\CAT$}
\nuse{\nref{Cat}}{$\Cat$}
\nuse{\nref{Cone}}{$\Cone$}
\nuse{\nref{CptHff}}{$\CptHff$}
\nuse{\nref{CRing}}{$\CRing$}
\nuse{\nref{discrete-space}}{$D$}
\nuse{\pageref{eq:lim-shapes}, \nref{cat-elts}}{$\scat{E}$}
\nuse{\nref{ev}}{$\ev$}
\nuse{\nref{FDVect}}{$\FDVect$}
\nuse{\nref{Field}}{$\Field$}
\nuse{\nref{FinSet}}{$\FinSet$}
\nuse{\nref{Grp}}{$\Grp$}
\nuse{\nref{hom-out}}{$\h^A$}
\nuse{\nref{hom-in}}{$\h_A$}
\nuse{\nref{hom-out-map}}{$\h^f$}
\nuse{\nref{hom-in-map}}{$\h_f$}
\nuse{\nref{hom-out-blank}}{$\h^\bl$}
\nuse{\nref{hom-in-blank}}{$\h_\bl$}
\nuse{\nref{Hom}, \nref{Hom-functor}}{$\Hom$}
\nuse{\nref{HOM}}{$\HOM$}
\nuse{\nref{indiscrete-space}}{$I$}
\nuse{\nref{lim}}{$\lt{}$}
\nuse{\nref{colim}}{$\colt{}$}
\nuse{\nref{Mon}}{$\Mon$}
\nuse{\nref{nat}}{$\nat$}
\nuse{\nref{oset}, \nref{oset-ftr}}{$\oset$}
\nuse{\nref{ob}}{$\ob$}
\nuse{\nref{op}}{$\blank^\op$}
\nuse{\pageref{eq:lim-shapes}}{$\scat{P}$}
\nuse{\nref{power-set}, \nref{power-set-precise}, \nref{power-set-ftr}}{$\pset$}
\nuse{\pageref{eq:pjn-aft}}{$P_A$}
\nuse{\nref{Ring}}{$\Ring$}
\nuse{\nref{circle}}{$S^1$}
\nuse{\nref{Set}}{$\Set$}
\nuse{\pageref{eq:lim-shapes}}{$\scat{T}$}
\nuse{\nref{Top}}{$\Tp$}
\nuse{\nref{Top-star}}{$\Tp_*$}
\nuse{\nref{Toph}}{$\Toph$}
\nuse{\nref{Toph-star}}{$\Toph_*$}
\nuse{\nref{Vect}}{$\Vect_k$}
\nuse{\nref{poly-one}}{$\integers[x]$}
\end{multicols}%
}
%
%
%

\index{bilinear|see{map, bilinear}}
\index{category!comma|see{comma category}}
\index{category!slice|see{slice category}}
\index{counit|see{unit and counit}}
\index{diagonal|see{functor, diagonal}}
\index{discrete|see{category, discrete \emph{and} topological space, discrete}}
\index{element!least|see{least element}}
\index{figure|see{element, generalized}}
\index{forgetful|see{functor, forgetful}}
\index{full|see{functor, full \emph{and} subcategory, full}}
\index{generalized element|see{element, generalized}}
\index{group!representation of|see{representation}}
\index{image!inverse|see{inverse image}}
\index{initial|see{object, initial \emph{and} set, weakly initial}}
\index{integers|see{$\integers$}}
\index{natural isomorphism|see{isomorphism, natural}}
\index{preimage|see{inverse image}}
\index{preservation|see{limit, preservation of}}
\index{representable|see{functor, representable}}
\index{terminal|see{object, terminal}}

\index{arrow|seealso{map}}
\index{category!one-object|seealso{monoid \emph{and} group}}
\index{cocone|seealso{cone}}
\index{coequalizer|seealso{equalizer}}
\index{colimit|seealso{limit}}
\index{coproduct|seealso{sum}}
\index{disjoint union|seealso{set, category of, sums in}}
\index{epic|seealso{monic}}
\index{epimorphism|seealso{epic}}
\index{exponential|seealso{set of functions}}
\index{fibred product|seealso{pullback}}
\index{G-set@$G$-set|seealso{monoid, action of}}
\index{group|seealso{monoid}}
\index{group!action of|seealso{monoid, action of}}
\index{homotopy|seealso{group, fundamental}}
\index{monoid!action of|seealso{group, action of}}
\index{monomorphism|seealso{monic}}
\index{morphism|seealso{map}}
\index{natural numbers|seealso{arithmetic}}
\index{object!terminal|seealso{object, initial}}
\index{partially ordered set|seealso{ordered set}}
\index{polynomial|seealso{ring, polynomial}}
\index{poset|seealso{ordered set}}
\index{preorder|seealso{ordered set}}
\index{pushout|seealso{pullback}}
\index{relation|seealso{equivalence relation}}
\index{sum|seealso{product}}
\index{topological space|seealso{homotopy \emph{and} group, fundamental}}
\index{vector space|seealso{bilinear map}}

\printindex             

\end{document}